\documentclass[12pt]{amsart}
\usepackage[T1]{fontenc}
\usepackage[utf8]{inputenc}
\usepackage[
  backend=biber,
  style=numeric,
  sorting=none,
  giveninits=false,
  doi=true,
  isbn=true,
  url=true,
  eprint=true
]{biblatex}
\usepackage{amsmath, amssymb, amsthm}
\usepackage{mathrsfs}
\usepackage{bbm}
\usepackage{hyperref}
\usepackage{geometry}
\usepackage{enumitem}
\usepackage{booktabs}
\usepackage{array}
\usepackage{tabularx}
\usepackage{tikz-cd}   

\newtheorem{theorem}{Theorem}[section]
\newtheorem{proposition}[theorem]{Proposition}
\newtheorem{lemma}[theorem]{Lemma}
\newtheorem{corollary}[theorem]{Corollary}
\newtheorem{definition}[theorem]{Definition}

\theoremstyle{remark}
\newtheorem{remark}[theorem]{Remark}
\newtheorem{example}[theorem]{Example}
\newtheorem{problem}[theorem]{Problem}
\newtheorem{conjecture}[theorem]{Conjecture}

\newtheorem{question}[theorem]{Question}
\newtheorem{principle}[theorem]{Principle}

\newcommand{\RealB}{\operatorname{Real}_{B}}
\newcommand{\Realet}{\operatorname{Real}_{\acute{e}t}}
\newcommand{\RealHdg}{\operatorname{Real}_{\mathrm{Hdg}}}
\newcommand{\RealMHM}{\operatorname{Real}_{\mathrm{MHM}}}
\newcommand{\RealNori}{\operatorname{Real}_{\mathrm{Nori}}}
\newcommand{\RealBr}{\operatorname{Real}_{\mathrm{Br}}}
\newcommand{\Realnr}{\operatorname{Real}_{\mathrm{nr}}}
\newcommand{\Realst}{\operatorname{Real}_{\mathrm{st}}}

\newcommand{\Realstar}{\operatorname{Real}_{\star}}

\begin{document}
\title[Torsion Trajectories from Local Discriminants to Global Obstructions] {Torsion Trajectories from Local Discriminants to Global Obstructions}

\author{Abdul Rahman}
\thanks{Email: arahman@alum.howard.edu}

\begin{abstract}
For a normal surface singularity, the discrepancy between the ordinary and dual middle-perversity intersection complexes over \(\mathbb Z\) is measured by a finite group \(E\).  In previous work, \(E\) was identified with link torsion, the exceptional-lattice discriminant group \(\Lambda^\vee/\Lambda\), a resolution-neighborhood boundary quotient, and, in the hypersurface case,
\(\operatorname{coker}(T-\mathrm{id})_{\mathrm{tors}}\). This paper tracks the trajectory of this torsion from local singularity data to global obstruction theory.  We follow the discriminant package \((E,q)\) through support cohomology, excision, global torsion, Brauer comparison, Bloch--Ogus residues, and rationalization.  The method is example-driven: trajectory tables are computed for \(A_1\), \(A_k\), \(D_4\), \(E_8\), a non-ADE Brieskorn singularity, the threefold ordinary double point, nodal threefolds, nodal quintics, and the Benoist--Ottem benchmark. The computations reveal a sharp distinction: a surface \(A_1\) singularity has local \(\mathbb Z/2\)-torsion, whereas a threefold ordinary double point has torsion-free link \(S^2\times S^3\) and contributes free vanishing-cycle data instead.  Thus finite discriminant torsion is naturally a codimension-two phenomenon, not a generic feature of nodes.  The resulting pattern motivates a specialization problem: whether the Enriques \(2\)-torsion in Benoist--Ottem integral Hodge counterexamples is genuinely global, or can arise after degeneration from transverse \(A_1\)-type discriminant data along codimension-two strata.
\end{abstract}
\maketitle
\tableofcontents

\section{Introduction}

Intersection homology, perverse sheaves, and Hodge-theoretic
decomposition theorems depend essentially on the coefficient category.  The middle-perversity intersection complex was introduced in the topological theory of Goresky--MacPherson and reformulated sheaf-theoretically in the theory of perverse sheaves \cite{GoMP83,BBD82}.  With field coefficients, Verdier duality preserves the middle-perversity t-structure, and for
projective morphisms the decomposition theorem expresses the derived direct image of an intersection complex as a direct sum of shifted semisimple perverse sheaves
\cite{BBD82,Sa90,deCataldoMigliorini05,deCataldoMigliorini09}.  Over
\(\mathbb Z\), the same formalism is sensitive to torsion.  The ordinary middle-perversity and dual middle-perversity t-structures need not agree, because duality changes the torsion condition in the critical stalk or costalk degree.  The resulting discrepancy is not a defect of the theory; it is integral geometric data.  It is already visible in Goresky--Siegel duality and in the ordinary/dual middle-perversity formalism of BBD \cite{GoreskySiegel83,BBD82}.

The point of this paper is to treat these corrections as geometric data rather than as technical noise.  Integral torsion appears in several related forms: as torsion in local intersection homology, as linking pairings on singular links, as discriminant groups of exceptional lattices, as cokernels of integral variation maps, as Brauer classes, and as unramified cohomology classes
\cite{GoreskySiegel83,FriedmanGenIH,FriedmanBook20,Nikulin80,Mi68,BlochOgus74,CTVoisin12}. These appearances are often studied in different languages.  The purpose of the present paper is to begin a local-to-global comparison among them.

We study these integral and torsion-theoretic invariants as a diagnostic layer that records precisely where integral structures fail to be seen by the rational theory.  This viewpoint is consistent with the role of torsion in integral intersection homology and torsion-sensitive Deligne sheaves
\cite{FriedmanGenIH,FriedmanTsInv,FriedmanBook20}, and also with the role of torsion in known failures of the integral Hodge conjecture
\cite{AtiyahHirzebruch62,Totaro97,SouleVoisin05,CTVoisin12,BenoistOttem20}. Understanding how torsion is created locally, transported through resolutions, and either survives or dies globally is therefore a natural step toward a more complete obstruction theory for integral Hodge phenomena.

The guiding point is that rationalization kills torsion but does not erase the geometric mechanisms that produced it.  A torsion class may vanish after tensoring with \(\mathbb Q\), yet before this disappearance it may record a linking form, a lattice discriminant, a monodromy cokernel, a support relation, a Brauer class, or a closed-stratum gluing term in the MacPherson--Vilonen description of perverse sheaves.  Thus the central question is not whether torsion survives rationalization, but what torsion forces before rationalization kills it.

\subsection{The invariant \(E\)}

The local input for this paper is the finite group \(E\) attached to a normal surface singularity in  \cite{RahmanIntegralPerverseObstructions}.  Let
\((X,0)\) be a germ of a normal complex analytic surface, and let
\[
        j:U:=X\setminus\{0\}\hookrightarrow X
\]
be the inclusion of the punctured germ.  The ordinary and dual
middle-perversity intersection complexes with integral coefficients are denoted
\[
        {}^pIC_X\mathbb Z,
        \qquad
        {}^p_+IC_X\mathbb Z.
\]
Both restrict to the shifted constant sheaf \(\mathbb Z_U[2]\) on \(U\), but
they need not agree over \(\mathbb Z\).  Following Jung--Saito and
\cite{RahmanIntegralPerverseObstructions}, define
\[
        E:=H^0({}^p_+IC_X\mathbb Z)_0.
\]
Equivalently, \(E\) is the finite group appearing in the distinguished
triangle
\[
        {}^pIC_X\mathbb Z
        \longrightarrow
        {}^p_+IC_X\mathbb Z
        \longrightarrow
        E[1]
        \longrightarrow .
\]
Here \(E\) is viewed as a complex supported at the singular point.  Thus
\(E\) measures the local integral discrepancy between the ordinary and dual
middle extensions
\cite{JungSaitoFactoriality,RahmanIntegralPerverseObstructions}.

This discrepancy is invisible rationally.  After tensoring with
\(\mathbb Q\), the ordinary and dual middle-perversity intersection complexes
become canonically identified:
\[
        {}^pIC_X\mathbb Q
        \cong
        {}^p_+IC_X\mathbb Q.
\]
Consequently \(E\) is a purely integral invariant.  It records the torsion
correction that survives over \(\mathbb Z\) but disappears after
rationalization, precisely the phenomenon that the integral and
torsion-sensitive theories are designed to detect
\cite{BBD82,FriedmanGenIH,RahmanIntegralPerverseObstructions}.

\subsection{Six local realizations}

The main local input for this paper is the finite group \(E\) attached to a
normal surface singularity in \cite{RahmanIntegralPerverseObstructions}.
Let \(L\) be the link of the singularity, let
\(\pi:\widetilde X\to X\) be the minimal resolution, let \(\Lambda\) be the
exceptional lattice generated by the irreducible exceptional curves, and let
\(M\) be the corresponding intersection matrix.  The main theorem of
\cite{RahmanIntegralPerverseObstructions} identifies \(E\) through four broad
realizations:
\[
        \text{perverse discrepancy},\qquad
        H^2(L,\mathbb Z)_{\mathrm{tors}},\qquad
        \Lambda^\vee/\Lambda,\qquad
        \operatorname{coker}(T-\mathrm{id})_{\mathrm{tors}},
\]
where the last realization is available in the isolated hypersurface case.
More explicitly,
\[
        E
        \cong
        H^2(L,\mathbb Z)_{\mathrm{tors}}
        \cong
        \Lambda^\vee/\Lambda,
\]
and hence
\[
        |E|=|\det(M)|.
\]
If \((X,0)\) is an isolated hypersurface surface singularity and \(T\)
denotes the Milnor monodromy on integral vanishing cohomology, then
\[
        E\cong \operatorname{coker}(T-\mathrm{id})_{\mathrm{tors}}.
\]
Under the additional hypothesis that
\[
        (T-\mathrm{id})\otimes_{\mathbb Z}\mathbb Q
\]
is an isomorphism, this gives
\[
        |E|=|\det(T-\mathrm{id})|.
\]
The link-topological, resolution-theoretic, and monodromy-theoretic inputs
used in these identifications are classical; see, for example,
\cite{Mu61,Mi68,GoreskySiegel83,Nikulin80}.  The precise compatibility of
these realizations in the present notation is established in
\cite{RahmanIntegralPerverseObstructions}.

For the purposes of the torsion-trajectory program, we refine these four
broad realizations into \textit{six operational interfaces}.  The reason for
this refinement is that local-to-global transport requires not only knowing
the finite group \(E\), but also knowing which geometric operation exposes it.
The perverse realization separates into the cone description and the Friedman
torsion-sensitive truncation interpretation.  The topological and
resolution-theoretic realizations separate into link cohomology, the
exceptional-lattice discriminant group, and the resolution-neighborhood
pair-sequence boundary quotient.  The hypersurface case adds the
Wang-sequence or monodromy realization. Thus the same finite group is visible through the following \textit{six operational realizations}:
\begin{enumerate}[label=\textup{(\arabic*)}]
\item The \textit{perverse cone}
\[
        \operatorname{Cone}
        \left(
        {}^pIC_X\mathbb Z
        \to
        {}^p_+IC_X\mathbb Z
        \right),
\]
where \(E\) appears as the point-supported correction measuring the
difference between the ordinary and dual middle-perversity intersection
complexes \cite[Theorem 1.1]{RahmanIntegralPerverseObstructions};

\item The \textit{Friedman torsion-sensitive truncation discrepancy} between the ordinary and dual middle-perversity packages, which gives the torsion-sensitive Deligne-sheaf interpretation of the same integral correction term \cite{FriedmanGenIH,FriedmanTsInv} and is compared with the ordinary/dual BBD perversity conventions in
\cite[Appendix C]{RahmanIntegralPerverseObstructions};

\item The \textit{link cohomology group}
\[
        H^2(L,\mathbb Z)_{\mathrm{tors}},
\]
where \(L\) is the link of the normal surface singularity; this is the local-topological realization of \(E\)
\cite[Theorem 3.8]{RahmanIntegralPerverseObstructions}, and it is compatible with the classical linking-duality picture for singularity links
\cite{GoreskySiegel83};

\item The \textit{exceptional-lattice discriminant group}
\[
        \Lambda^\vee/\Lambda,
\]
where \(\Lambda\) is the lattice generated by the exceptional curves of the
minimal resolution; this is the resolution-theoretic realization
\cite[Theorem 4.7]{RahmanIntegralPerverseObstructions}, with the
discriminant-form formalism following the standard lattice conventions of
Nikulin \cite{Nikulin80};

\item The \textit{boundary quotient} appearing in the long exact sequence of the resolution - neighborhood pair \((N,L)\), where \(N\) is a sufficiently small resolution neighborhood of the exceptional divisor and \(L=\partial N\). Under Poincaré--Lefschetz duality, the map
\[
        H_2(N,\mathbb Z)\to H_2(N,L;\mathbb Z)
\]
is identified with the lattice map
\[
        \Lambda\to\Lambda^\vee,
\]
and its finite cokernel is \(\Lambda^\vee/\Lambda\)
\cite[Lemmas 4.3--4.5]{RahmanIntegralPerverseObstructions};

\item and, in the isolated hypersurface case, the \textit{Wang-sequence or Milnor-monodromy realization}
\[
        \operatorname{coker}(T-\mathrm{id})_{\mathrm{tors}},
\]
where \(T\) is the integral Milnor monodromy on vanishing cohomology
\cite[Theorem 1.3 and Section 5]{RahmanIntegralPerverseObstructions}; the underlying Milnor-fibration and monodromy formalism is classical
\cite{Mi68}.
\end{enumerate}

These six operational realizations are not merely redundant computations.
They provide different interfaces with the geometry.  The perverse cone
connects \(E\) to integral middle perversities
\cite{BBD82,JungSaitoFactoriality}.  The torsion-sensitive truncation
realization connects it to Friedman's torsion-sensitive Deligne sheaves
\cite{FriedmanGenIH,FriedmanTsInv}.  The link realization connects it to
local topology and linking pairings \cite{GoreskySiegel83}.  The lattice
realization equips it with a discriminant form
\[
        q:E\times E\longrightarrow \mathbb Q/\mathbb Z,
\]
which can be compared with the torsion linking pairing on the link
\cite{GoreskySiegel83,Nikulin80}.  The pair-sequence realization is the
natural entry point for local-to-global transport through a resolution.  The
monodromy realization connects the same group to Picard--Lefschetz theory and
the integral variation map \(T-\mathrm{id}\) \cite{Mi68}.

The present paper asks which of these realizations is best suited for global transport.  More precisely, if \(X\) is a projective variety with isolated singularities \(p_1,\ldots,p_s\), and if each singularity carries a local discriminant package \((E_i,q_i)\), we ask when the direct sum of these local packages maps to global torsion invariants of a resolution.

The later appendices make this transport language explicit.  Appendix~\ref{app:realization-dictionary} fixes the realization notation used for Betti, étale, Hodge, Nori, Brauer, unramified, and stacky comparison stations.  Appendix~\ref{app:six-operations} records the six-operation transport dictionary used to compare local, smooth, stacky, and motivic torsion packages.  These appendices are not used to reprove the cited six-functor formalisms; they specify how those formalisms organize the torsion-trajectory program.

\subsection{The global targets}

Let
\[
        \pi:\widetilde X\longrightarrow X
\]
be a resolution of a projective variety with isolated singularities.  The
global torsion targets considered in this paper are
\[
        H^*(\widetilde X,\mathbb Z)_{\mathrm{tors}},
        \qquad
        \operatorname{Br}(\widetilde X),
        \qquad
        H^3_{\mathrm{nr}}(\widetilde X,\mathbb Q/\mathbb Z).
\]
The first target is ordinary integral cohomological torsion.  The second is
the cohomological Brauer group.  The third is unramified cohomology, defined
by residue conditions in the Bloch--Ogus formalism
\cite{BlochOgus74,CTVoisin12}.

For smooth projective threefolds, degree-four integral Hodge questions are
closely tied to torsion phenomena in degree three.  The exponential and
Kummer sequences relate the Brauer group to
\(H^3(\widetilde X,\mathbb Z)_{\mathrm{tors}}\) under suitable hypotheses,
while the Bloch--Ogus formalism detects unramified cohomology through the
vanishing of codimension-one residues \cite{BlochOgus74,CTVoisin12}.  These
groups appear naturally in the study of failures and refinements of the
integral Hodge conjecture
\cite{AtiyahHirzebruch62,Totaro97,SouleVoisin05,Voisin06IntegralHodge,CTVoisin12,BenoistOttem20,Totaro21KodairaZero}.
Thus, in the threefold setting, the relevant global torsion target is not an
unspecified cohomology group but specifically
\[
        H^3(\widetilde X,\mathbb Z)_{\mathrm{tors}}.
\]

\subsection{The three bridges}

The guiding composite has the form
\[
\bigoplus_i(E_i,q_i)
\overset{\alpha_X}{\longrightarrow}
H^*(\widetilde X,\mathbb Z)_{\mathrm{tors}}
\overset{\beta_X}{\longleftrightarrow}
\operatorname{Br}(\widetilde X)
\overset{\gamma_X}{\dashrightarrow}
H^3_{\mathrm{nr}}(\widetilde X,\mathbb Q/\mathbb Z).
\]
The dashed arrow ($\gamma_X$) is meant to emphasize that the last stage is governed by residue conditions rather than by an unconditional map defined on every torsion class.  Each arrow in the display has its own hypotheses and degree constraints.

The first bridge, \(\alpha_X\), is constructed from excision and local
cohomology.  Let
\[
        D_i:=\pi^{-1}(p_i),\qquad D:=\bigcup_iD_i,
\]
and let \(N_i\) be a small resolution neighborhood of \(D_i\).  The long
exact sequence for cohomology with supports,
\[
\cdots
\longrightarrow
H^k_D(\widetilde X,\mathbb Z)
\longrightarrow
H^k(\widetilde X,\mathbb Z)
\longrightarrow
H^k(\widetilde X\setminus D,\mathbb Z)
\longrightarrow
H^{k+1}_D(\widetilde X,\mathbb Z)
\longrightarrow
\cdots,
\]
together with excision,
\[
        H^k_D(\widetilde X,\mathbb Z)
        \cong
        \bigoplus_i H^k_{D_i}(N_i,\mathbb Z),
\]
provides the natural mechanism for transporting local torsion classes to
global cohomology.  We use this support-theoretic construction as the
cohomological analogue of the local resolution-neighborhood computations in
\cite{GoreskySiegel83,RahmanIntegralPerverseObstructions}; the general
formalism of local cohomology is classical \cite{Gro68}.

The second bridge, \(\beta_X\), is the comparison between torsion
cohomology and the cohomological Brauer group.  Let \(Y\) be a smooth
projective complex variety.  The exponential sequence
\[
        0
        \longrightarrow
        \mathbb Z
        \longrightarrow
        \mathcal O_Y
        \longrightarrow
        \mathcal O_Y^*
        \longrightarrow
        0
\]
gives the exact segment
\[
        H^2(Y,\mathbb Z)
        \longrightarrow
        H^2(Y,\mathcal O_Y)
        \longrightarrow
        H^2(Y,\mathcal O_Y^*)
        \xrightarrow{\partial}
        H^3(Y,\mathbb Z)
        \longrightarrow
        H^3(Y,\mathcal O_Y).
\]
The cohomological Brauer group is
\[
        \operatorname{Br}(Y)
        :=
        H^2(Y,\mathcal O_Y^*)_{\mathrm{tors}}.
\]
Since \(H^3(Y,\mathcal O_Y)\) is a complex vector space, it has no nonzero
torsion.  Therefore, if
\[
        b\in H^2(Y,\mathcal O_Y^*)_{\mathrm{tors}},
\]
then \(\partial(b)\in H^3(Y,\mathbb Z)\) is torsion.  Thus the boundary map
restricts to a homomorphism
\[
        \partial_{\mathrm{tors}}:
        \operatorname{Br}(Y)
        \longrightarrow
        H^3(Y,\mathbb Z)_{\mathrm{tors}}.
\]

If, in addition,
\[
        H^2(Y,\mathcal O_Y)=0,
\]
then the map
\[
        H^2(Y,\mathcal O_Y)
        \longrightarrow
        H^2(Y,\mathcal O_Y^*)
\]
has zero image.  Exactness then implies that \(\partial\) is injective on
\(H^2(Y,\mathcal O_Y^*)\).  Moreover, every torsion class in
\(H^3(Y,\mathbb Z)\) maps to zero in \(H^3(Y,\mathcal O_Y)\), because
\(H^3(Y,\mathcal O_Y)\) is torsion-free.  Hence every torsion class in
\(H^3(Y,\mathbb Z)\) lies in the image of \(\partial\).  

\begin{lemma}[Brauer--torsion comparison via the exponential sequence] \label{lem:brauertorsioncomparison}
Let \(Y\) be a smooth projective complex variety.  Suppose that
\[
        H^2(Y,\mathcal O_Y)=0.
\]
Then the boundary map in the exponential sequence induces an isomorphism
\[
        \operatorname{Br}(Y)
        =
        H^2(Y,\mathcal O_Y^*)_{\mathrm{tors}}
        \xrightarrow{\sim}
        H^3(Y,\mathbb Z)_{\mathrm{tors}}.
\]
\end{lemma}

\begin{proof}
The exponential sequence
\[
        0
        \longrightarrow
        \mathbb Z
        \longrightarrow
        \mathcal O_Y
        \longrightarrow
        \mathcal O_Y^*
        \longrightarrow
        0
\]
gives the exact segment
\[
        H^2(Y,\mathcal O_Y)
        \longrightarrow
        H^2(Y,\mathcal O_Y^*)
        \xrightarrow{\partial}
        H^3(Y,\mathbb Z)
        \longrightarrow
        H^3(Y,\mathcal O_Y).
\]
By assumption,
\[
        H^2(Y,\mathcal O_Y)=0.
\]
Therefore the map
\[
        \partial:
        H^2(Y,\mathcal O_Y^*)
        \longrightarrow
        H^3(Y,\mathbb Z)
\]
is injective.

We first show that \(\partial\) sends torsion to torsion.  Let
\[
        b\in H^2(Y,\mathcal O_Y^*)_{\mathrm{tors}}.
\]
Then there exists \(n>0\) such that
\[
        nb=0.
\]
Applying \(\partial\), we get
\[
        n\partial(b)=\partial(nb)=0.
\]
Hence
\[
        \partial(b)\in H^3(Y,\mathbb Z)_{\mathrm{tors}}.
\]
Thus \(\partial\) restricts to an injective homomorphism
\[
        \partial_{\mathrm{tors}}:
        H^2(Y,\mathcal O_Y^*)_{\mathrm{tors}}
        \hookrightarrow
        H^3(Y,\mathbb Z)_{\mathrm{tors}}.
\]

It remains to prove surjectivity.  Let
\[
        \eta\in H^3(Y,\mathbb Z)_{\mathrm{tors}}.
\]
Since \(H^3(Y,\mathcal O_Y)\) is a complex vector space, it has no nonzero
torsion.  Therefore the image of \(\eta\) under
\[
        H^3(Y,\mathbb Z)
        \longrightarrow
        H^3(Y,\mathcal O_Y)
\]
is zero.  By exactness of the exponential-sequence long exact sequence,
there exists
\[
        b\in H^2(Y,\mathcal O_Y^*)
\]
such that
\[
        \partial(b)=\eta.
\]

We now show that this preimage \(b\) is torsion.  Since \(\eta\) is torsion,
there exists \(n>0\) such that
\[
        n\eta=0.
\]
Using \(\partial(b)=\eta\), we obtain
\[
        \partial(nb)=n\partial(b)=n\eta=0.
\]
But \(\partial\) is injective, because \(H^2(Y,\mathcal O_Y)=0\).  Hence
\[
        nb=0.
\]
Thus
\[
        b\in H^2(Y,\mathcal O_Y^*)_{\mathrm{tors}}
        =
        \operatorname{Br}(Y).
\]
Therefore every torsion class in \(H^3(Y,\mathbb Z)\) lies in the image of
\(\operatorname{Br}(Y)\).  This proves that
\[
        \partial_{\mathrm{tors}}:
        \operatorname{Br}(Y)
        \xrightarrow{\sim}
        H^3(Y,\mathbb Z)_{\mathrm{tors}}
\]
is an isomorphism.
\end{proof}

By Lemma~\ref{lem:brauertorsioncomparison}, the sufficient hypothesis
\[
        H^2(Y,\mathcal O_Y)=0
\]
implies that the boundary map in the exponential sequence induces a
canonical isomorphism
\[
        \operatorname{Br}(Y)
        =
        H^2(Y,\mathcal O_Y^*)_{\mathrm{tors}}
        \xrightarrow{\sim}
        H^3(Y,\mathbb Z)_{\mathrm{tors}}.
\]
Thus, in any situation where this hypothesis holds for
\(Y=\widetilde X\), a class that survives the support-transport stage and lands in
\[
        H^3(\widetilde X,\mathbb Z)_{\mathrm{tors}}
\]
determines a unique cohomological Brauer class on \(\widetilde X\).  This is the Brauer bridge used in the torsion trajectory: local discriminant torsion first enters global degree-three torsion, and under the exponential-sequence comparison of Lemma~\ref{lem:brauertorsioncomparison}, that global torsion
is identified with a Brauer class appearing in the Colliot-Thélène--Voisin framework \cite{CTVoisin12}.

The third bridge, \(\gamma_X\), is governed by Bloch--Ogus residues.  This
bridge is different from the Brauer comparison above.  The Brauer bridge is
an identification, under suitable hypotheses, between degree-three integral
torsion and the cohomological Brauer group.  The unramified bridge is instead
a survival test: it asks whether a class has zero residue at every
codimension-one point.  This residue-theoretic description is the basic
Bloch--Ogus interpretation of unramified cohomology
\cite{BlochOgus74,CTVoisin12}.

Let \(Y\) be a smooth projective complex variety.  In the Bloch--Ogus
formalism, unramified cohomology with coefficients \(\mathbb Q/\mathbb Z\)
is defined as the subgroup of function-field cohomology consisting of
classes with trivial residues along all codimension-one points
\cite[Section~1]{BlochOgus74}
\cite[Section~1]{CTVoisin12}:
\[
        H^3_{\mathrm{nr}}(Y,\mathbb Q/\mathbb Z)
        =
        \ker
        \left(
        H^3(\mathbb C(Y),\mathbb Q/\mathbb Z)
        \xrightarrow{\oplus\,\partial_D}
        \bigoplus_{D\in Y^{(1)}}
        H^2(\mathbb C(D),\mathbb Q/\mathbb Z)
        \right).
\]
where \(\mathrm{nr}\) denotes unramified. Here \(Y^{(1)}\) denotes the set of codimension-one points, or equivalently
prime divisors, of \(Y\), and \(\partial_D\) is the residue map associated to the discrete valuation of \(\mathbb C(Y)\) defined by \(D\) \cite{BlochOgus74,CTVoisin12}.

Thus a torsion class that reaches the Brauer or \(\mathbb Q/\mathbb Z\) stage contributes to unramified cohomology only if all of its residues vanish:
\[
        \partial_D(\eta)=0
        \qquad
        \text{for every }D\in Y^{(1)}.
\]
Equivalently, unramified cohomology is the subgroup of classes in the
function field that survive every codimension-one residue test
\cite{BlochOgus74,CTVoisin12}.  This is why the arrow
\[
        \operatorname{Br}(Y)
        \dashrightarrow
        H^3_{\mathrm{nr}}(Y,\mathbb Q/\mathbb Z)
\]
is drawn as a dashed arrow in the trajectory diagram.  It is not an
automatic map on every torsion class in the trajectory.  It is a
residue-survival condition.

In the present paper, \(Y=\widetilde X\).  A local discriminant class may
survive support transport and may even determine a Brauer class under the
exponential-sequence comparison.  It still has to pass the Bloch--Ogus
residue test before it contributes to
\[
        H^3_{\mathrm{nr}}(\widetilde X,\mathbb Q/\mathbb Z).
\]
The third bridge \(\gamma_X\) therefore records the passage from global
torsion or Brauer-visible data to residue-surviving unramified data.  This is
the station at which a class may die by residue, even after it has survived
the earlier support and Brauer comparison stages.  The role of this
unramified station in integral Hodge obstruction theory is precisely the
mechanism used by Colliot-Thélène--Voisin in their comparison between
degree-four integral Hodge phenomena and degree-three unramified cohomology
\cite{CTVoisin12}.

\subsection{Degree specification}

The local-to-global map cannot be formulated correctly without specifying
its degree.  For a normal surface singularity, the group \(E\) is detected
through the torsion cohomology of a three-dimensional link:
\[
        E\cong H^2(L,\mathbb Z)_{\mathrm{tors}}.
\]
This is the local degree relevant to the surface theory of
\cite{RahmanIntegralPerverseObstructions}.  However, for threefold
applications related to Brauer groups and degree-four integral Hodge
questions, the relevant global target is
\[
        H^3(\widetilde X,\mathbb Z)_{\mathrm{tors}}.
\]
Thus one must verify that the local torsion package associated to a
threefold singularity contributes in degree three.

For this reason we use degree-indexed local-to-global maps
\[
        \alpha_X^{(k)}:
        \bigoplus_i E_i
        \longrightarrow
        H^k(\widetilde X,\mathbb Z)_{\mathrm{tors}}
\]
whenever the local packages \(E_i\) are realized in the corresponding local
support degree.  In surface-type models, the relevant degree is governed by
the support sequence of the exceptional curves and the cohomology of
three-dimensional links.  In the threefold setting, the map relevant to the
Brauer and unramified cohomology comparison is
\[
        \alpha_X^{(3)}:
        \bigoplus_iE_i
        \longrightarrow
        H^3(\widetilde X,\mathbb Z)_{\mathrm{tors}}.
\]
This distinction is essential.  In particular, the ordinary double point in complex dimension three is not automatically governed by the same local calculation as the \(A_1\) surface singularity.  Its link has real dimension five, and the stalk complex is shifted by \([3]\), not by \([2]\).  The threefold node must therefore be recomputed separately before making any claim about nodal threefolds or nodal quintics; the relevant local topology and Milnor-fibration input go back to the standard theory of isolated
hypersurface singularities \cite{Mi68,La81} (see Appendices \ref{app:A1surfacetrajectory} and \ref{app:3foldordinarydoublepointtrajectory} for more details).

\subsection{Main results}

The first part of the paper recalls the local discriminant package and
relates it to divisor-theoretic invariants.  Under the standard
resolution-lattice hypotheses for normal surface singularities, the local
class group is identified with the discriminant group:
\[
        \operatorname{Cl}(\mathcal O_{X,0})
        \cong
        \Lambda^\vee/\Lambda
        \cong
        E.
\]
This comparison connects the perverse-sheaf invariant \(E\) with the
factoriality and \(Q\)-factoriality questions that motivate the
Jung--Saito framework
\cite{JungSaitoFactoriality,RahmanIntegralPerverseObstructions}.

The second part constructs the degree-indexed local-to-global maps.  For a
resolution \(\pi:\widetilde X\to X\) with isolated exceptional loci, excision
and cohomology with supports give maps
\[
        \alpha_X^{(k)}:
        \bigoplus_iE_i
        \longrightarrow
        H^k(\widetilde X,\mathbb Z)_{\mathrm{tors}}
\]
whenever the local discriminant packages are realized in the relevant support
degree.  We formulate the kernel of \(\alpha_X^{(k)}\) in terms of global
relations among exceptional cycles.  In nodal hypersurface examples, this
relation space is expected to be governed by the same geometry that appears
in the Clemens defect \cite{Clemens83,Cynk01,Kloosterman22}.

The third part gives the Brauer and unramified comparisons.  Under precise
vanishing hypotheses, the exponential sequence identifies
\[
        \operatorname{Br}(\widetilde X)
        \cong
        H^3(\widetilde X,\mathbb Z)_{\mathrm{tors}},
\]
so the image of \(\alpha_X^{(3)}\) determines a subgroup of the Brauer group.
We then apply the Bloch--Ogus residue criterion to decide whether such
classes survive to
\[
        H^3_{\mathrm{nr}}(\widetilde X,\mathbb Q/\mathbb Z).
\]
This gives a precise sense in which a local discriminant class can contribute
to a global unramified torsion obstruction \cite{BlochOgus74,CTVoisin12}.

The fourth part treats discriminant forms.  The local groups \(E_i\) are not
only finite abelian groups; through the lattice and linking realizations they
carry pairings
\[
        q_i:E_i\times E_i\longrightarrow \mathbb Q/\mathbb Z.
\]
We prove pairing compatibility in the basic \(A_1\) surface model and use the
ADE examples, especially \(D_4\), to show why the pairing must be tracked.
The general compatibility of the local-to-global map with global torsion
pairings is formulated as a problem.  This part uses the classical relation
between discriminant forms and linking pairings
\cite{GoreskySiegel83,Nikulin80}.

The fifth part corrects and sharpens the comparison with Benoist--Ottem.
Their counterexamples to the integral Hodge conjecture use global Enriques
\(2\)-torsion \cite{BenoistOttem20}.  The natural boundary model for this
comparison is not a transverse \(A_1\) singularity, but the Enriques/Coble
boundary singularity of type \(\frac14(1,1)\).  In the degree-two Enriques
compactification of Alexeev--Engel--Garza--Schaffler, the discriminant
divisor parametrizes quotients of nodal K3 surfaces by an involution fixing a
node; the resulting boundary surfaces are rational Coble surfaces with a
\(\frac14(1,1)\)-singularity
\cite{AlexeevEngelGarzaSchafflerEnriquesDegree2}.  For this singularity the
local invariant from \cite{RahmanIntegralPerverseObstructions} is
\[
        E\cong H^2(L(4,1),\mathbb Z)_{\mathrm{tors}}
        \cong \mathbb Z/4.
\]
The local index-two cover
\[
        A_1=\mathbb C^2/\mu_2
        \longrightarrow
        \mathbb C^2/\mu_4=\frac14(1,1)
\]
induces the link cover
\[
        L(2,1)\longrightarrow L(4,1).
\]
The Bockstein of this double-cover class selects the unique order-two
subgroup
\[
        2E\cong\mathbb Z/2
        \subset
        E\cong\mathbb Z/4.
\]
Thus the Benoist--Ottem \(2\)-torsion direction is compared not with the
full Coble \(E\)-package, but with the Bockstein-selected order-two shadow
\(2E\subset E\).

The final part records the motivic form of this corrected comparison.  In an
integral motivic sheaf formalism with Betti realization compatible with
cones, shifts, multiplication by integers, and proper pushforward, the full
Coble boundary \(E\)-object has the functorial motivic candidate
\[
        \mathcal T^{\mathrm{mot}}_{4,\Sigma}
        :=
        i_*(\mathbf 1_\Sigma/4)[2],
\]
whose Betti realization is
\[
        i_*(\mathbb Z/4)_\Sigma[2].
\]
The Benoist--Ottem order-two contribution is the motivic shadow
\[
        \mathcal T^{\mathrm{mot}}_{2,\Sigma}
        :=
        i_*(\mathbf 1_\Sigma/2)[2]
        \longrightarrow
        i_*(\mathbf 1_\Sigma/4)[2],
\]
whose Betti realization is the inclusion
\[
        i_*(\mathbb Z/2)_\Sigma[2]
        \longrightarrow
        i_*(\mathbb Z/4)_\Sigma[2].
\]
This expresses the Benoist--Ottem \(2\)-torsion mechanism as the order-two
filtered piece of the motivic Coble \(E\)-package.

\subsection{Main contributions}

The paper makes six contributions.  First, it organizes the local invariant
\(E\) into a torsion trajectory, tracking birth, discriminant form, local
realizations, support transport, Brauer comparison, residue behavior, and
rationalization.  Second, it computes this trajectory in a uniform format for
surface ADE, non-ADE, quotient, threefold-node, nodal, and
Benoist--Ottem/Coble benchmark examples.  Third, it separates local finite
discriminant torsion from free vanishing-cycle relations and from global
smooth torsion.  In particular, a surface \(A_1\) singularity contributes
local \(\mathbb Z/2\)-torsion, whereas a threefold ordinary double point has
torsion-free link \(S^2\times S^3\).  Fourth, it formulates the
codimension-two transverse torsion principle using MacPherson--Vilonen
gluing: transverse surface torsion appears as a closed-stratum perverse
discrepancy term.  Fifth, it corrects the Benoist--Ottem comparison by
showing that the relevant Enriques/Coble boundary singularity is
\(\frac14(1,1)\), with \(E\cong\mathbb Z/4\), and that the
Benoist--Ottem \(2\)-torsion mechanism selects the order-two Bockstein shadow
\(2E\cong\mathbb Z/2\).  Sixth, it records a motivic lift of this structure:
the full Coble package is modeled by \(i_*(\mathbf 1_\Sigma/4)[2]\), while
the Benoist--Ottem order-two shadow is modeled by
\(i_*(\mathbf 1_\Sigma/2)[2]\).

\subsection{Organization of the paper}

Section~2 introduces the torsion trajectory and the main tables.  It defines the stations through which local torsion is tracked: birth, form, local realization, excision, global survival, Brauer comparison, residue behavior, and rationalization.  Section~3 fixes conventions on degrees, signs, discriminant forms, and torsion-sensitive truncations.  It also records the degree bookkeeping needed for isolated singularities of different dimensions.
Section~4 recalls the local discriminant package \(E\), its compatible realizations, and its comparison with local class groups in the surface case. Section~5 constructs the degree-indexed local-to-global maps \(\alpha_X^{(k)}\) using excision and cohomology with supports.  It also formulates kernel criteria in terms of global relations among exceptional or support classes.  Section~6 treats the Brauer comparison through the exponential and Kummer sequences.  Section~7 treats unramified cohomology and gives the residue survival criterion.

Section~8 extracts the patterns observed in the trajectory computations: non-unimodularity as the source of local finite torsion, \(E_8\) as the null-control case, the form-sensitivity of \(D_4\), the distinction between surface \(A_1\) torsion and threefold ordinary double point free data, the codimension-two nature of transverse surface torsion, and the corrected Coble boundary interpretation of the Benoist--Ottem comparison.  Section~9
states the resulting trajectory principles, functorial constructions, and open problems, including the codimension-two transverse torsion principle, the MacPherson--Vilonen closed-stratum gluing interpretation, the motivic lift \(i_*(\mathbf 1_\Sigma/n)[2]\), and the Benoist--Ottem/Coble order-two shadow theorem.

The appendices contain the detailed computations.  They work out the
\(A_1\), \(A_k\), \(D_4\), \(E_8\), non-ADE Brieskorn, Coble
\(\frac14(1,1)\), threefold ordinary double point, nodal threefold,
Benoist--Ottem, and nodal quintic trajectories, followed by a reference appendix collecting the exact sequences used throughout. Appendix~\ref{app:realization-dictionary} fixes realization notation for motivic torsion packages.  Appendix~\ref{app:six-operations} records the six-operation transport dictionary used to compare torsion trajectories across motivic, Betti, étale, Hodge, Nori, Brauer, unramified, and stacky settings.

\subsection{How to read this paper}

This paper is deliberately modular.  The main text gives the conceptual trajectory: local birth, discriminant form, transport, Brauer comparison, residue survival, and rational death.  The appendices give the computations that justify the rows of the tables.  The reader need not read the appendices linearly on a first pass.  Instead, the tables in Section~2 can be used as a guide: each row points to a corresponding appendix where the claim is verified.

For a first reading, we recommend Sections~1--2 for the overview,
Sections~4--7 for the main local-to-global mechanisms, Section~9 for the Benoist--Ottem/Coble comparison, and Appendices~L--M for the realization and six-operation dictionary.  The remaining appendices may be consulted as needed.  This organization reflects the purpose of the paper: to make the torsion trajectory both conceptually visible and computationally checkable.

\section{The torsion trajectory and the main tables}

The constructions in this paper are organized around a single question:
what happens to integral torsion as it moves from local singularity data to
global cohomology, then to Brauer and unramified cohomology, and finally
disappears after rationalization?

The local input is the finite group
\[
        E=H^0({}^p_+IC_X\mathbb Z)_0
\]
attached to a normal surface singularity.  In
\cite{RahmanIntegralPerverseObstructions}, this group was computed through
compatible realizations:
\[
        E
        \cong
        H^2(L,\mathbb Z)_{\mathrm{tors}}
        \cong
        \Lambda^\vee/\Lambda,
\]
and, in the isolated hypersurface case,
\[
        E
        \cong
        \operatorname{coker}(T-\mathrm{id})_{\mathrm{tors}}.
\]
Here \(L\) is the link of the singularity, \(\Lambda\) is the exceptional
lattice of a resolution, and \(T\) is the Milnor monodromy.  These
identifications are not merely different calculations of the same group.
They show that the same torsion phenomenon is visible in perverse sheaves,
torsion-sensitive truncation, link topology, exceptional lattices,
resolution-neighborhood pair sequences, and monodromy.

The purpose of the present paper is to track this torsion as a trajectory.
At each station we ask whether the class exists, what structure it carries,
whether it survives the next exact sequence, and what kills it if it does not
survive.  The tables below are designed to keep five mechanisms distinct:
\textit{local finite discriminant torsion}, \textit{free vanishing-cycle
data}, \textit{global smooth torsion}, \textit{Brauer/unramified
obstruction}, and \textit{Bockstein-selected torsion shadows inside larger
local \(E\)-packages}.  Tables~\ref{tab:trajectory-legend},
\ref{tab:local-birth-form}, \ref{tab:local-realizations}, and
\ref{tab:transport-obstruction} summarize the notation, local birth data,
local realization data, and transport behavior used throughout the paper.

\subsection{The stations of the trajectory}

We organize the movement of torsion into seven stations.

\begin{enumerate}[label=\textup{(\arabic*)}]
\item \textbf{Birth.}  In the surface singularity setting, finite torsion is
born locally as \(E_i\cong \Lambda_i^\vee/\Lambda_i\).  At this station we
record the group structure, invariant factors, order, and prime
decomposition.  For example, \(A_1\) has
\(E\cong\mathbb Z/2\mathbb Z\), \(A_k\) has
\(E\cong\mathbb Z/(k+1)\mathbb Z\), \(E_8\) has \(E=0\), and the Coble
boundary singularity \(\frac14(1,1)\) has
\(E\cong\mathbb Z/4\mathbb Z\).

\item \textbf{Form.}  When \(E_i\) is identified with a discriminant group,
it carries a finite \(\mathbb Q/\mathbb Z\)-valued pairing
\[
        q_i:E_i\times E_i\longrightarrow \mathbb Q/\mathbb Z.
\]
This pairing is equivalent, up to the sign convention fixed below, to the
torsion linking pairing on the link; see the linking-pairing formalism of
\cite{GoreskySiegel83} and the lattice discriminant formalism of
\cite{Nikulin80}.  The \(D_4\) row shows that the group alone is not always
enough, while the Coble boundary row shows that a global \(2\)-torsion
mechanism may select only a subgroup of a larger local package:
\[
        2E\cong\mathbb Z/2
        \subset
        E\cong\mathbb Z/4.
\]

\item \textbf{Local realization.}  The same finite group is computed through
the available local models: perverse discrepancy, torsion-sensitive
truncation, link torsion, exceptional lattice, pair sequence, and, in the
hypersurface case, monodromy.  The agreement of these realizations in the
surface case is the local theorem of
\cite{RahmanIntegralPerverseObstructions}.  In the tables, the entry
``same \(E\)'' means that all of the available stations recover the same
finite group.  For quotient singularities such as \(\frac14(1,1)\), the
monodromy station is not part of the hypersurface Wang-sequence realization,
but the perverse, torsion-sensitive, link, lattice, and pair-sequence
stations still recover the same \(E\).

\item \textbf{Excision.}  If the local package is realized in degree \(k\)
support cohomology, it can enter the local-to-global channel through
\[
        \bigoplus_iE_i
        \longrightarrow
        H^k_D(\widetilde X,\mathbb Z)_{\mathrm{tors}},
\]
where \(D\) is the exceptional locus of a resolution.  This station is
degree-sensitive.  Surface examples naturally enter through the degree-two
support and pair-sequence geometry of exceptional curves.  The threefold
Brauer channel, by contrast, passes through degree three.

\item \textbf{Global survival.}  The forget-support map
\[
        H^k_D(\widetilde X,\mathbb Z)_{\mathrm{tors}}
        \longrightarrow
        H^k(\widetilde X,\mathbb Z)_{\mathrm{tors}}
\]
may kill local torsion classes.  The kernel records global relations among
the support classes.  This is the first global killing mechanism in the
trajectory.

\item \textbf{Brauer and residue tests.}  In the threefold setting, classes
that survive to
\[
        H^3(\widetilde X,\mathbb Z)_{\mathrm{tors}}
\]
may be compared with the Brauer group by the exponential sequence under
suitable hypotheses.  A further unramified survival test is given by
Bloch--Ogus residues:
\[
        H^3(\mathbb C(\widetilde X),\mathbb Q/\mathbb Z)
        \longrightarrow
        \bigoplus_{D\in \widetilde X^{(1)}}H^2(k(D),\mathbb Q/\mathbb Z).
\]
A class survives to \(H^3_{\mathrm{nr}}\) precisely when all these residues
vanish \cite{BlochOgus74,CTVoisin12}.

\item \textbf{Rational death.}  Every finite torsion class vanishes after
tensoring with \(\mathbb Q\).  If \(e\in E_i\) has order \(n\), then
\[
        e\otimes 1
        =
        e\otimes n\cdot \frac1n
        =
        ne\otimes \frac1n
        =
        0
        \quad
        \text{in }E_i\otimes_{\mathbb Z}\mathbb Q.
\]
Thus
\[
        E_i\otimes_{\mathbb Z}\mathbb Q=0.
\]
The point is not that torsion survives rationalization; it does not.  The
point is that before it dies, it records integral data invisible over
\(\mathbb Q\).  In the Coble boundary example, both the full package
\[
        E\cong\mathbb Z/4
\]
and the Benoist--Ottem order-two shadow
\[
        2E\cong\mathbb Z/2
\]
die rationally, but they encode different integral information before they
disappear.
\end{enumerate}

\subsection{Legend for trajectory tables}

We use the following notation in the trajectory tables.

\begin{table}[htbp]
\centering
\caption{Legend for the torsion trajectory tables.}
\label{tab:trajectory-legend}
\small
\setlength{\tabcolsep}{5pt}
\renewcommand{\arraystretch}{1.2}
\begin{tabularx}{0.82\textwidth}{|
>{\centering\arraybackslash}p{0.18\textwidth}|
>{\centering\arraybackslash}X|}
\hline
\textbf{Symbol} & \textbf{Meaning} \\
\hline
\(\checkmark\) & The class or package survives this station. \\
\hline
\(0\) & The class dies or the relevant group is zero. \\
\hline
\(?\) & The station requires a separate computation before any global assertion can be made. \\
\hline
\(\varnothing\) & The station is not defined in this example. \\
\hline
\(\mathrm{rel}\) & The class is killed by a global relation. \\
\hline
\(\mathrm{res}\) & The class is killed by a residue map. \\
\hline
free & The local data is free abelian, not finite torsion. \\
\hline
global & The torsion is global smooth torsion rather than local singularity torsion. \\
\hline
shadow & A subgroup or filtered piece selected inside a larger local \(E\)-package. \\
\hline
\end{tabularx}
\end{table}

The symbol \(?\) is important.  It is not a placeholder for carelessness.  It
marks the exact point at which a computation must be performed before the
trajectory can continue.  For example, the ordinary double point in complex
dimension three cannot be assigned the same row as the \(A_1\) surface
singularity without recomputing the link, support degree, and monodromy
station.  That computation shows that the link is \(S^2\times S^3\), hence
torsion-free.

\subsection{Local birth and form table}

The first table records the local birth and form data.  It includes the
finite group \(E\), its prime support, and the discriminant form where one is
defined.

\begin{table}[htbp]
\centering
\caption{Local birth and form data for the principal trajectory examples.}
\label{tab:local-birth-form}
\scriptsize
\setlength{\tabcolsep}{2.5pt}
\renewcommand{\arraystretch}{1.2}
\begin{tabularx}{\textwidth}{|
>{\centering\arraybackslash}p{0.18\textwidth}|
>{\centering\arraybackslash}p{0.19\textwidth}|
>{\centering\arraybackslash}p{0.15\textwidth}|
>{\centering\arraybackslash}p{0.18\textwidth}|
>{\centering\arraybackslash}X|}
\hline
\textbf{Example} &
\textbf{Finite group \(E\)} &
\textbf{Prime support} &
\textbf{Discriminant form} &
\textbf{Interpretation} \\
\hline
\(A_1\) surface &
\(\mathbb Z/2\mathbb Z\) &
\(2\) &
\(q(1,1)=-\frac12\) &
Basic local \(2\)-torsion atom. \\
\hline
\(A_k\) surface &
\(\mathbb Z/(k+1)\mathbb Z\) &
primes dividing \(k+1\) &
\(q(\bar g,\bar g)=-\frac{k}{k+1}\) &
Cyclic torsion family controlled by \(|\det A_k|=k+1\). \\
\hline
\(D_4\) surface &
\((\mathbb Z/2\mathbb Z)^2\) &
\(2\) &
\(\begin{smallmatrix}0&-\frac12\\[-1mm]-\frac12&0\end{smallmatrix}\) &
Noncyclic \(2\)-torsion; form data is essential. \\
\hline
\(E_8\) surface &
\(0\) &
none &
\(0\) &
Null-control: nontrivial singularity, unimodular lattice, no finite torsion. \\
\hline
\(x^2+y^3+z^{11}\) &
\(\mathbb Z/5\mathbb Z\) &
\(5\) &
\(q(\bar g,\bar g)=-\frac15\) &
Non-ADE odd-prime torsion example. \\
\hline
Threefold ODP &
\(0\) finite torsion &
none &
\(0\) &
Link \(S^2\times S^3\) is torsion-free; local data is free, not finite. \\
\hline
Nodal threefold &
\(0\) finite torsion at each node &
none &
\(0\) &
Defect concerns free node relations, not local finite torsion. \\
\hline
Benoist--Ottem \(S\times C\) &
none locally on smooth fiber &
\(2\) globally &
none local &
Global smooth Enriques \(2\)-torsion; no singular local \(E\)-package on the smooth fiber. \\
\hline
Coble boundary \(\frac14(1,1)\) &
\(\mathbb Z/4\mathbb Z\) &
\(2\) &
\(q(\bar g,\bar g)=-\frac14\) &
Enriques/Coble boundary; BO \(2\)-torsion selects \(2E\cong\mathbb Z/2\). \\
\hline
\end{tabularx}
\end{table}

The first five rows are local surface computations.  They show that finite
local torsion is governed by the discriminant group
\[
        \Lambda^\vee/\Lambda.
\]
The threefold ordinary double point and nodal threefold rows show a different
phenomenon: ordinary nodes in complex dimension three have free vanishing
data rather than finite local discriminant torsion.  The Benoist--Ottem row
is different again: the torsion is global and smooth on \(S\times C\), not
local singularity torsion.  The Coble boundary row corrects the natural
Enriques degeneration target: the boundary singularity is \(\frac14(1,1)\),
with full local package \(E\cong\mathbb Z/4\), and the Benoist--Ottem
\(2\)-torsion mechanism selects the order-two subgroup \(2E\).

\subsection{Six-realization table}

The second table records which local realizations are available and what
they produce.  The phrase ``same \(E\)'' means that the realization produces
the finite group listed in the previous table.

\begin{table}[htbp]
\centering
\caption{Local realization table for the torsion trajectory.}
\label{tab:local-realizations}
\scriptsize
\setlength{\tabcolsep}{2.0pt}
\renewcommand{\arraystretch}{1.25}
\begin{tabularx}{\textwidth}{|
>{\centering\arraybackslash}p{0.16\textwidth}|
>{\centering\arraybackslash}c|
>{\centering\arraybackslash}c|
>{\centering\arraybackslash}c|
>{\centering\arraybackslash}c|
>{\centering\arraybackslash}c|
>{\centering\arraybackslash}c|
>{\centering\arraybackslash}X|}
\hline
\textbf{Example} &
\textbf{Perverse} &
\textbf{ts} &
\textbf{Link} &
\textbf{Lattice} &
\textbf{Pair} &
\textbf{Monod.} &
\textbf{Conclusion} \\
\hline
\(A_1\) surface &
same \(E\) &
same \(E\) &
\(\mathbb Z/2\) &
\(\mathbb Z/2\) &
\(\mathbb Z/2\) &
\(\mathbb Z/2\) &
All six stations agree. \\
\hline
\(A_k\) surface &
same \(E\) &
same \(E\) &
\(\mathbb Z/(k+1)\) &
\(\mathbb Z/(k+1)\) &
\(\mathbb Z/(k+1)\) &
\(\mathbb Z/(k+1)\) &
All six stations agree. \\
\hline
\(D_4\) surface &
same \(E\) &
same \(E\) &
\((\mathbb Z/2)^2\) &
\((\mathbb Z/2)^2\) &
\((\mathbb Z/2)^2\) &
\((\mathbb Z/2)^2\) &
All six stations agree, but form data matters. \\
\hline
\(E_8\) surface &
\(0\) &
\(0\) &
\(0\) &
\(0\) &
\(0\) &
\(0\) &
All stations vanish; lattice is unimodular. \\
\hline
Brieskorn \(2,3,11\) &
same \(E\) &
same \(E\) &
\(\mathbb Z/5\) &
\(\mathbb Z/5\) &
\(\mathbb Z/5\) &
\(\mathbb Z/5\) &
All six stations agree outside the ADE list. \\
\hline
Threefold ODP &
no finite \(E\) &
no torsion &
torsion-free &
no finite discr. &
free pair data &
free cokernel &
Surface \(E\)-package does not carry over. \\
\hline
Coble boundary \(\frac14(1,1)\) &
same \(E\) &
same \(E\) &
\(\mathbb Z/4\) &
\(\mathbb Z/4\) &
\(\mathbb Z/4\) &
n/a &
Full local package is \(\mathbb Z/4\); BO sees \(2E\cong\mathbb Z/2\). \\
\hline
\end{tabularx}
\end{table}

This table is the main evidence for the first organizing principle of the
paper: for normal surface singularities in the examples considered, finite
torsion is stable across the available local realizations.  The threefold
node row is included as a warning.  It is not a counterexample to the surface
theorem; it is a different-dimensional local model.  The Coble row is
included as the corrected Enriques boundary model: the available surface
realizations recover the full \(E\cong\mathbb Z/4\), while the
Benoist--Ottem mechanism detects only the Bockstein-selected order-two
subgroup.

\subsection{Transport and obstruction table}

The third table records what happens after the local station.  This is where
global geometry enters.

\begin{table}[htbp]
\centering
\caption{Transport and obstruction table.}
\label{tab:transport-obstruction}
\scriptsize
\setlength{\tabcolsep}{2.2pt}
\renewcommand{\arraystretch}{1.25}
\begin{tabularx}{\textwidth}{|
>{\centering\arraybackslash}p{0.17\textwidth}|
>{\centering\arraybackslash}p{0.16\textwidth}|
>{\centering\arraybackslash}p{0.16\textwidth}|
>{\centering\arraybackslash}X|
>{\centering\arraybackslash}p{0.14\textwidth}|
>{\centering\arraybackslash}p{0.08\textwidth}|}
\hline
\textbf{Example} &
\textbf{Support status} &
\textbf{Global image} &
\textbf{Killing mechanism} &
\textbf{Brauer / residue} &
\(\mathbb Q\) \\
\hline
\(A_1\) surface &
local degree \(2\) &
depends on global exceptional-curve relations &
may die under global divisor/support relations &
not local in isolated germ &
\(0\) \\
\hline
\(A_k\) surface &
local degree \(2\) &
depends on global exceptional-chain relations &
may die under global lattice relations &
not local in isolated germ &
\(0\) \\
\hline
\(D_4\) surface &
local degree \(2\) &
depends on global relations and form data &
group and form may behave differently under transport &
not local in isolated germ &
\(0\) \\
\hline
\(E_8\) surface &
none &
no local torsion image &
no birth: lattice unimodular &
none &
\(0\) \\
\hline
Brieskorn \(2,3,11\) &
local degree \(2\) &
depends on global plumbing/support relations &
may die globally; prime \(5\) is local source &
not local in isolated germ &
\(0\) \\
\hline
Threefold ODP &
none finite &
no local finite torsion image &
no finite torsion birth; free classes may have relations &
none local &
\(0\) \\
\hline
Nodal threefold &
none finite at nodes &
finite local image zero &
defect records free global relations, not finite torsion &
none local &
\(0\) \\
\hline
Benoist--Ottem \(S\times C\) &
none on smooth fiber &
global torsion exists, not from direct \(\alpha\) &
not local on smooth fiber; boundary comparison uses Coble row &
global Brauer / unramified benchmark &
\(0\) \\
\hline
Coble boundary \(\frac14(1,1)\) &
local degree \(2\) &
full local image \(\mathbb Z/4\); BO sees \(2E\) &
Bockstein selects the order-two subgroup &
Enriques/Coble boundary benchmark &
\(0\) \\
\hline
\end{tabularx}
\end{table}

The key separation in this table is between local finite torsion and other
forms of topology.  The threefold ordinary double point has no finite local
torsion, but it has free vanishing-cycle data.  Nodal defect concerns global
relations among these free classes.  Benoist--Ottem torsion is global smooth
torsion on a smooth variety, not direct local singularity torsion.  The
correct boundary comparison is supplied by the Coble \(\frac14(1,1)\) row:
the full local \(E\)-package is \(\mathbb Z/4\), while the
Benoist--Ottem \(2\)-torsion direction is the Bockstein-selected subgroup
\(2E\cong\mathbb Z/2\).

The detailed computations underlying Tables~\ref{tab:local-birth-form},
\ref{tab:local-realizations}, and \ref{tab:transport-obstruction} are carried
out in the appendices.  In particular, the \(A_1\), \(A_k\), \(D_4\),
\(E_8\), and Brieskorn rows are computed in the surface appendices; the
Coble \(\frac14(1,1)\) row is computed in the Coble boundary appendix; the
ordinary double point and nodal rows are computed in the threefold
appendices; and the Benoist--Ottem row is treated as a global torsion
benchmark.

\subsection{Why the examples come early}

The examples in this paper are not illustrations of a completed abstract
theorem.  They are the data from which the global pattern is extracted.  The
formal results below explain the columns of the table:

\begin{itemize}
\item the local discriminant package explains the birth and form columns;

\item the six-realization theorem explains the local model columns;

\item Mayer--Vietoris, pair sequences, excision, and cohomology with supports
explain the support and transport columns;

\item Bockstein exact sequences explain how a mod-\(2\) double-cover class can
select an order-two subgroup inside a larger local \(E\)-package;

\item MacPherson--Vilonen and Beilinson gluing explain how transverse surface
torsion appears as a closed-stratum perverse discrepancy term;

\item the exponential sequence explains the Brauer column;

\item Bloch--Ogus residues explain the unramified survival column;

\item rationalization explains the final disappearance of finite torsion over
\(\mathbb Q\).
\end{itemize}

Thus the paper proceeds in two interlocking directions.  The theorem sections
prove that the stations of the trajectory are well-defined under stated
hypotheses.  The appendices compute the trajectory row by row.  The final
sections extract the patterns that the rows reveal.

\subsection{The main empirical questions}

The tables are designed to answer the following questions.

\begin{question}
Does the survival of local finite torsion depend only on the finite group
\(E_i\), or does it depend on the full discriminant package \((E_i,q_i)\)?
\end{question}

\begin{question}
Which local torsion classes are killed by the forget-support map
\[
        H^k_D(\widetilde X,\mathbb Z)_{\mathrm{tors}}
        \longrightarrow
        H^k(\widetilde X,\mathbb Z)_{\mathrm{tors}}?
\]
Equivalently, which classes die because of global relations among
exceptional, support, or divisor classes?
\end{question}

\begin{question}
When a local torsion class survives to
\[
        H^3(\widetilde X,\mathbb Z)_{\mathrm{tors}},
\]
does it determine a Brauer class, and if so, do its residues vanish?
\end{question}

\begin{question}
When a global \(2\)-torsion mechanism is compared with a local singular
boundary package, does it see the full local group \(E\), or only a
Bockstein-selected subgroup such as \(2E\subset E\)?
\end{question}

\begin{question}
For the Benoist--Ottem examples, is the relevant local comparison with the
full Coble boundary package
\[
        E\cong\mathbb Z/4,
\]
or specifically with its order-two shadow
\[
        2E\cong\mathbb Z/2?
\]
\end{question}

These questions do not assume that every local torsion class survives.  In
fact, death at any station is part of the computation.  A class may fail to
be born, may die by a global relation, may fail to enter the Brauer channel,
may be killed by a residue, may be seen only through a subgroup or quotient,
or may disappear under rationalization.  The point of the trajectory is to
record where and why this happens.

\subsection{How theorems are extracted from the tables}

The later formal statements should be read as abstractions of the table
entries.  For example, the definition of
\[
        \alpha_X^{(k)}
\]
abstracts the excision column.  The kernel criterion abstracts the global
survival column.  The Brauer comparison abstracts the Brauer column.  The
residue criterion abstracts the unramified survival column.  The
MacPherson--Vilonen gluing theorem abstracts the codimension-two
closed-stratum column.  The Bockstein calculations abstract the order-two
shadow column in the Coble boundary example.

This organization is deliberate.  The examples are computed first, with all
available local realizations.  The theorems then record the functorial
mechanisms that explain the repeated patterns.  Finally, the conjectures and
comparison problems state the patterns that remain after the computations are
assembled across the rows.

\section{Conventions, degrees, and local packages}

This section fixes the conventions used throughout the paper.  The main
points are the sign convention for exceptional intersection matrices, the
normalization of discriminant forms, the degree bookkeeping for isolated
singularities, the notation for ordinary and dual middle perversity over
\(\mathbb Z\), and the convention for subpackages such as
\(2E\subset E\).  These conventions are needed because the same local torsion
group can appear in different cohomological degrees depending on the
dimension of the singularity, because discriminant forms may differ by sign
depending on whether one uses the geometric or positive-definite lattice
convention, and because the Benoist--Ottem comparison uses an order-two
Bockstein shadow inside the larger local group \(E\) of the Coble boundary
singularity.

\subsection{Sign convention for intersection matrices}

Throughout the paper, exceptional intersection matrices are taken in the
geometric convention.  Thus if
\[
        \pi:\widetilde X\longrightarrow X
\]
is a resolution of a normal surface singularity and
\[
        E_\pi=\bigcup_a C_a
\]
is the exceptional divisor, the intersection matrix is
\[
        M=(C_a\cdot C_b).
\]
For rational double points this is the negative of the corresponding Cartan
matrix.  In particular, for the \(A_1\) surface singularity the exceptional
lattice is generated by a single curve \(C\cong \mathbb P^1\) with
\[
        (C,C)=-2,
\]
so the intersection matrix is
\[
        [-2].
\]
For an \(A_k\) singularity, the matrix is the negative \(A_k\) Cartan matrix.
Likewise, the \(D\)- and \(E\)-type rational double points are described by
negative definite ADE lattices.

For the cyclic quotient singularity of type \(\frac14(1,1)\), the
Hirzebruch--Jung continued fraction is
\[
        \frac41=[4].
\]
Thus the minimal resolution has one exceptional curve \(C\cong\mathbb P^1\)
with
\[
        (C,C)=-4,
\]
and the exceptional intersection matrix is
\[
        [-4].
\]
This is the geometric convention used below when computing the Coble boundary
package
\[
        E\cong\mathbb Z/4.
\]

This convention is the one used in
\cite{RahmanIntegralPerverseObstructions}.  It is also the convention most
natural from the geometry of surface resolutions.  Some lattice-theoretic
references use the corresponding positive definite root lattice instead.
Passing between the two conventions changes the sign of the bilinear form
and therefore changes the displayed representative of the discriminant
pairing by an overall sign.  Since the paper compares local discriminant
pairings with link pairings, this sign choice is fixed once and for all here.

\begin{remark}
All determinant statements are insensitive to this sign convention.  For
example, if \(M\) is the negative \(A_k\) Cartan matrix, then
\[
        |\det(M)|=k+1.
\]
For the Coble boundary singularity \(\frac14(1,1)\), the matrix is
\[
        M=[-4],
\]
and therefore
\[
        |\det(M)|=4.
\]
However, formulas for the discriminant form
\[
        q:\Lambda^\vee/\Lambda\times \Lambda^\vee/\Lambda
        \longrightarrow
        \mathbb Q/\mathbb Z
\]
do depend on the sign convention.  In this paper the discriminant form is
computed from the geometric negative definite intersection form.
\end{remark}

\subsection{Discriminant-form convention}

Let \((\Lambda,(\ ,\ ))\) be a nondegenerate integral lattice.  Thus
\(\Lambda\) is a free abelian group of finite rank and
\[
        (\ ,\ ):\Lambda\times\Lambda\longrightarrow \mathbb Z
\]
is a nondegenerate symmetric bilinear form.  Its dual lattice is
\[
        \Lambda^\vee
        :=
        \operatorname{Hom}_{\mathbb Z}(\Lambda,\mathbb Z)
        \cong
        \{x\in \Lambda\otimes_{\mathbb Z}\mathbb Q
          \mid (x,\lambda)\in\mathbb Z
          \text{ for all } \lambda\in\Lambda\}.
\]
The discriminant group is
\[
        A_\Lambda:=\Lambda^\vee/\Lambda.
\]
Since the form is nondegenerate over \(\mathbb Q\), the quotient
\(A_\Lambda\) is a finite abelian group.

\begin{definition}[Discriminant pairing]
The discriminant pairing of \(\Lambda\) is the
\(\mathbb Q/\mathbb Z\)-valued symmetric bilinear form
\[
        q_\Lambda:
        A_\Lambda\times A_\Lambda
        \longrightarrow
        \mathbb Q/\mathbb Z
\]
defined by
\[
        q_\Lambda(x+\Lambda,y+\Lambda)
        :=
        (x,y)\bmod \mathbb Z,
        \qquad x,y\in \Lambda^\vee.
\]
\end{definition}

This pairing is well-defined because replacing \(x\) or \(y\) by an element
differing by a vector in \(\Lambda\) changes \((x,y)\) by an integer.  It is
nondegenerate because \(\Lambda\subset \Lambda^\vee\) is a finite-index
inclusion induced by a nondegenerate form.  Discriminant forms of integral
lattices are classical; see, for example, Nikulin's work on integral
symmetric bilinear forms \cite{Nikulin80}.

For a normal surface singularity, the exceptional lattice of the minimal
resolution gives a discriminant group
\[
        \Lambda^\vee/\Lambda.
\]
In the previous paper this group was identified with the local integral
perverse obstruction group \(E\) \cite{RahmanIntegralPerverseObstructions}.
Thus \(E\) carries not only the structure of a finite abelian group, but also
a canonical local discriminant pairing
\[
        q:E\times E\longrightarrow \mathbb Q/\mathbb Z.
\]

The same finite pairing can be seen topologically through the linking form on
the torsion homology of the link.  For a closed oriented three-manifold \(L\),
the torsion subgroup of \(H_1(L,\mathbb Z)\) carries a nondegenerate linking
pairing
\[
        \lambda_L:
        H_1(L,\mathbb Z)_{\mathrm{tors}}
        \times
        H_1(L,\mathbb Z)_{\mathrm{tors}}
        \longrightarrow
        \mathbb Q/\mathbb Z.
\]
For links of normal surface singularities, this pairing is identified, up to
the sign convention fixed above, with the discriminant pairing of the
exceptional lattice
\cite{GoreskySiegel83,Nikulin80,RahmanIntegralPerverseObstructions}.

\begin{example}[The Coble \(\frac14(1,1)\) convention]
For the singularity \(\frac14(1,1)\), the exceptional lattice is
\[
        \Lambda=\mathbb Z\langle e\rangle,
        \qquad
        (e,e)=-4.
\]
The dual generator is
\[
        e^\vee=-\frac e4,
\]
since
\[
        (e^\vee,e)=1.
\]
Thus
\[
        \Lambda^\vee/\Lambda\cong\mathbb Z/4.
\]
If \(\bar g\) denotes the class of \(e^\vee\), then
\[
        q(\bar g,\bar g)
        =
        (e^\vee,e^\vee)
        =
        -\frac14
        \quad\bmod \mathbb Z.
\]
The unique order-two subgroup is
\[
        2E\cong\mathbb Z/2
        \subset
        E\cong\mathbb Z/4.
\]
The order-two element \(2\bar g\) has
\[
        q(2\bar g,2\bar g)
        =
        4q(\bar g,\bar g)
        =
        -1
        \equiv 0
        \quad\bmod \mathbb Z.
\]
Thus the Benoist--Ottem order-two shadow is an isotropic subgroup of the
full local Coble \(E\)-package.
\end{example}

\begin{remark}
The sign difference between the lattice discriminant pairing and the
topological linking pairing depends on the boundary orientation convention
for the link.  Since the present paper is concerned primarily with whether
local finite pairings are transported globally, the essential invariant is the
finite \(\mathbb Q/\mathbb Z\)-valued pairing together with the chosen
orientation convention.  In all explicit ADE and quotient computations
below, we use the negative definite geometric intersection matrix.
\end{remark}

\subsection{Degree shifts for isolated singularities}

We next fix the basic degree bookkeeping.  Let \((X,0)\) be an isolated
complex analytic singularity of pure dimension \(n\).  Let
\[
        U:=X\setminus\{0\},
        \qquad
        j:U\hookrightarrow X,
        \qquad
        i:\{0\}\hookrightarrow X
\]
be the natural inclusions.  Let \(L\) be the link of the singularity.  For a
sufficiently small representative of the germ, the punctured neighborhood
\(U\) is homotopy equivalent to \(L\).  Therefore the stalk of the derived
pushforward \(Rj_*\mathbb Z_U[n]\) at the singular point is computed by the
cohomology of \(L\), shifted by \(n\).

\begin{lemma}[Local stalk degree bookkeeping]
Let \((X,0)\) be an isolated complex analytic singularity of pure dimension
\(n\), let \(U=X\setminus\{0\}\), and let \(L\) be the link.  Then there is a
canonical isomorphism
\[
        i^*Rj_*\mathbb Z_U[n]
        \cong
        R\Gamma(L,\mathbb Z)[n].
\]
Consequently,
\[
        H^m\bigl(i^*Rj_*\mathbb Z_U[n]\bigr)
        \cong
        H^{m+n}(L,\mathbb Z)
\]
for all \(m\in\mathbb Z\).
\end{lemma}

\begin{proof}
By definition of derived direct image,
\[
        H^m\bigl(i^*Rj_*\mathbb Z_U[n]\bigr)
        \cong
        (R^{m+n}j_*\mathbb Z_U)_0.
\]
The stalk \((R^{m+n}j_*\mathbb Z_U)_0\) is the cohomology of a sufficiently
small punctured neighborhood of the singular point.  This punctured
neighborhood deformation retracts onto the link \(L\).  Hence
\[
        (R^{m+n}j_*\mathbb Z_U)_0
        \cong
        H^{m+n}(L,\mathbb Z),
\]
which proves the assertion.
\end{proof}

For surface singularities, \(n=2\), and the local model is
\[
        i^*Rj_*\mathbb Z_U[2]
        \cong
        R\Gamma(L,\mathbb Z)[2].
\]
Thus
\[
        H^m\bigl(i^*Rj_*\mathbb Z_U[2]\bigr)
        \cong
        H^{m+2}(L,\mathbb Z).
\]
In particular, the degree-zero stalk sees \(H^2(L,\mathbb Z)\).  The torsion
subgroup of this group is precisely the local obstruction group \(E\) in the
surface case:
\[
        E\cong H^2(L,\mathbb Z)_{\mathrm{tors}}
\]
\cite{RahmanIntegralPerverseObstructions}.

For the Coble boundary singularity \(\frac14(1,1)\), the link is
\[
        L=L(4,1),
\]
and therefore
\[
        H^2(L,\mathbb Z)_{\mathrm{tors}}\cong\mathbb Z/4.
\]
Thus the surface-degree convention gives
\[
        E\cong\mathbb Z/4.
\]
The Benoist--Ottem \(2\)-torsion mechanism is not the full group \(E\); it is
the Bockstein-selected subgroup
\[
        2E\cong\mathbb Z/2.
\]

For a threefold isolated singularity, \(n=3\), and the corresponding stalk
complex is
\[
        i^*Rj_*\mathbb Z_U[3]
        \cong
        R\Gamma(L,\mathbb Z)[3].
\]
Here \(L\) has real dimension \(5\), and the relevant cohomology groups occur
in different degrees.  This is why the ordinary double point in complex
dimension three cannot be treated by simply importing the \(A_1\) surface
calculation.  Its local link, stalk degrees, and possible torsion
contributions must be computed separately.

\begin{remark}
This degree shift is the main reason for distinguishing the surface-type
local theory from the threefold applications to Brauer groups and unramified
cohomology.  The surface theory naturally detects
\(H^2(L,\mathbb Z)_{\mathrm{tors}}\), while the threefold Brauer comparison
passes through \(H^3(\widetilde X,\mathbb Z)_{\mathrm{tors}}\).
Codimension-two transverse surface singularities in a threefold combine
these two features: the transverse surface package is computed from
\(H^2\) of the three-dimensional transverse link, while its contribution to
the threefold is shifted and transported through closed-stratum gluing and
support cohomology.
\end{remark}

\subsection{Degree of the local-to-global map}

Let \(X\) be a projective variety with isolated singularities
\[
        \operatorname{Sing}(X)=\{p_1,\ldots,p_s\},
\]
and let
\[
        \pi:\widetilde X\longrightarrow X
\]
be a resolution.  Set
\[
        D_i:=\pi^{-1}(p_i),
        \qquad
        D:=\bigcup_iD_i,
        \qquad
        U:=X\setminus\operatorname{Sing}(X)
        \cong
        \widetilde X\setminus D.
\]
Let \(N_i\) be a sufficiently small resolution neighborhood of \(D_i\), with
boundary the link \(L_i\).

The long exact sequence for cohomology with supports gives
\[
\cdots
\longrightarrow
H^k_D(\widetilde X,\mathbb Z)
\longrightarrow
H^k(\widetilde X,\mathbb Z)
\longrightarrow
H^k(U,\mathbb Z)
\longrightarrow
H^{k+1}_D(\widetilde X,\mathbb Z)
\longrightarrow
\cdots .
\]
By excision,
\[
        H^k_D(\widetilde X,\mathbb Z)
        \cong
        \bigoplus_i H^k_{D_i}(N_i,\mathbb Z).
\]
Thus any local torsion package realized in a local support group can be
transported to global cohomology by first using excision and then forgetting
supports.

\begin{definition}[Degree-\(k\) local-to-global map]
Suppose that, for each singular point \(p_i\), a local torsion package
\(E_i\) is realized as a subgroup, quotient, or canonically identified
subquotient of the local support group
\[
        H^k_{D_i}(N_i,\mathbb Z).
\]
The degree-\(k\) local-to-global map is the map induced by excision and the
forget-support morphism:
\[
        \alpha_X^{(k)}:
        \bigoplus_iE_i
        \longrightarrow
        H^k_D(\widetilde X,\mathbb Z)_{\mathrm{tors}}
        \longrightarrow
        H^k(\widetilde X,\mathbb Z)_{\mathrm{tors}}.
\]
\end{definition}

The notation \(\alpha_X^{(k)}\) is intentionally degree-indexed.  There is
not a single degree-independent local-to-global map.  The relevant degree
depends on the dimension of the singularity, the form in which the local
torsion package appears, and the global target one wants to compare with.

\begin{remark}
For threefold applications connected to the Brauer group and degree-four
integral Hodge questions, the relevant target is
\[
        H^3(\widetilde X,\mathbb Z)_{\mathrm{tors}}.
\]
Thus the relevant local-to-global map is
\[
        \alpha_X^{(3)}:
        \bigoplus_iE_i
        \longrightarrow
        H^3(\widetilde X,\mathbb Z)_{\mathrm{tors}},
\]
provided the local packages \(E_i\) are realized in degree \(3\) local
support cohomology.  Surface-type global models involve different degrees and
are treated separately.
\end{remark}

\begin{remark}
The definition above is deliberately flexible at this stage.  In the surface
case, the group \(E_i\) is identified with link torsion and with the
discriminant group of the exceptional lattice.  In the threefold case, the
correct local package must be recomputed from the degree-shifted link
complex, the resolution neighborhood, and, in the hypersurface case, the Wang
sequence.  In codimension-two transverse surface situations, such as a
threefold with a transverse \(\frac14(1,1)\) Coble stratum, the local package
is first computed transversely as a surface package
\[
        E_\Sigma\cong(\mathbb Z/4)_\Sigma,
\]
and is then shifted into the threefold through MacPherson--Vilonen or
Beilinson gluing.  The Benoist--Ottem comparison uses the subpackage
\[
        2E_\Sigma\cong(\mathbb Z/2)_\Sigma.
\]
\end{remark}

\subsection{Ordinary and dual middle perversity}

We recall the notation for the ordinary and dual middle-perversity
\(t\)-structures over \(\mathbb Z\).  Let \(X\) be a complex analytic space,
and let \(D^b_c(X,\mathbb Z)\) denote the bounded constructible derived
category with integral coefficients.  For \(K\in D^b_c(X,\mathbb Z)\), the
ordinary middle-perversity conditions are denoted
\[
        K\in {}^pD^{\le 0}(X,\mathbb Z),
        \qquad
        K\in {}^pD^{\ge 0}(X,\mathbb Z),
\]
and the corresponding heart is the abelian category
\[
        {}^p\operatorname{Perv}(X,\mathbb Z)
        :=
        {}^pD^{\le 0}(X,\mathbb Z)
        \cap
        {}^pD^{\ge 0}(X,\mathbb Z).
\]
The dual middle-perversity \(t\)-structure is denoted
\[
        {}^p_+D^{\le 0}(X,\mathbb Z),
        \qquad
        {}^p_+D^{\ge 0}(X,\mathbb Z),
\]
with heart
\[
        {}^p_+\operatorname{Perv}(X,\mathbb Z).
\]
It is characterized by Verdier duality:
\[
        K\in {}^p_+D^{\le 0}(X,\mathbb Z)
        \quad\Longleftrightarrow\quad
        \mathbb D K\in {}^pD^{\ge 0}(X,\mathbb Z),
\]
and
\[
        K\in {}^p_+D^{\ge 0}(X,\mathbb Z)
        \quad\Longleftrightarrow\quad
        \mathbb D K\in {}^pD^{\le 0}(X,\mathbb Z).
\]
This is the integral dual-perversity formalism of BBD \cite{BBD82}.  Over a
field, the ordinary and dual middle perversities agree, but over
\(\mathbb Z\) they may differ because torsion and torsion-free conditions are
exchanged in the critical degree.

For a pure \(n\)-dimensional irreducible complex analytic space \(X\), with
smooth locus \(j:X_{\mathrm{reg}}\hookrightarrow X\), we write
\[
        {}^pIC_X\mathbb Z
        :=
        {}^pj_{!*}\mathbb Z_{X_{\mathrm{reg}}}[n],
\]
and
\[
        {}^p_+IC_X\mathbb Z
        :=
        {}^p_+j_{!*}\mathbb Z_{X_{\mathrm{reg}}}[n].
\]
These are the ordinary and dual middle-perversity intersection complexes with
integral coefficients.  They restrict to the same shifted local system on the
smooth locus:
\[
        j^*({}^pIC_X\mathbb Z)
        \cong
        \mathbb Z_{X_{\mathrm{reg}}}[n]
        \cong
        j^*({}^p_+IC_X\mathbb Z).
\]
The difference between them is therefore supported on the singular locus.

\begin{remark}
For a normal surface germ \((X,0)\), the singular locus is isolated after
shrinking.  Hence the discrepancy between
\({}^pIC_X\mathbb Z\) and \({}^p_+IC_X\mathbb Z\) is supported at the
singular point and is measured by the finite group
\[
        E=H^0({}^p_+IC_X\mathbb Z)_0
\]
as in \cite{RahmanIntegralPerverseObstructions}.  For the Coble boundary
singularity \(\frac14(1,1)\), this group is
\[
        E\cong\mathbb Z/4.
\]
\end{remark}

\begin{remark}[Closed-stratum gluing convention]
For a codimension-two stratum
\[
        i:\Sigma\hookrightarrow Y
\]
whose transverse singularity is a normal surface singularity with local
package \(E\), the corresponding closed-stratum discrepancy term is written
\[
        i_*E_\Sigma[2].
\]
In the Coble case,
\[
        E_\Sigma\cong(\mathbb Z/4)_\Sigma,
\]
so the full discrepancy term is
\[
        i_*(\mathbb Z/4)_\Sigma[2].
\]
The Benoist--Ottem order-two contribution is the subterm
\[
        i_*(2E_\Sigma)[2]
        \cong
        i_*(\mathbb Z/2)_\Sigma[2].
\]
\end{remark}

\subsection{Torsion-sensitive truncation convention}

We also recall the torsion-sensitive truncation convention used in the
examples.  Friedman's torsion-sensitive Deligne sheaves refine the usual
Deligne construction by allowing specified torsion to survive one degree
above the ordinary truncation cutoff
\cite{FriedmanGenIH,FriedmanTsInv,FriedmanBook20}.  Let \(\mathcal P\) be a
set of primes.  For an abelian group \(G\), let
\[
        T_{\mathcal P}(G)
\]
denote the subgroup of elements whose order is divisible only by primes in
\(\mathcal P\).  The torsion-tipped truncation
\[
        \check{\tau}^{\mathcal P}_{\le k}
\]
agrees with the usual truncation \(\tau_{\le k}\) in degrees \(\le k\), kills
cohomology in degrees \(>k+1\), and retains the \(\mathcal P\)-torsion part
of the cohomology in degree \(k+1\).

In the two extreme cases,
\[
        \mathcal P=\varnothing
\]
gives the ordinary truncation, while
\[
        \mathcal P=\mathcal P(\mathbb Z)
\]
retains all torsion in the critical degree.  In the isolated
codimension-two point-stratum case relevant to normal surface singularities,
these two choices correspond to the ordinary and dual middle-perversity
packages used above.  Thus torsion-sensitive truncation gives an explicit
local model for the discrepancy between
\[
        {}^pIC_X\mathbb Z
        \qquad\text{and}\qquad
        {}^p_+IC_X\mathbb Z.
\]

\begin{remark}
In the surface case, applying this convention to the local model
\[
        i^*Rj_*\mathbb Z_U[2]
        \cong
        R\Gamma(L,\mathbb Z)[2]
\]
shows that the ordinary middle extension removes the critical torsion while
the dual middle extension retains it.  This is the torsion-sensitive
truncation realization of the group
\[
        E\cong H^2(L,\mathbb Z)_{\mathrm{tors}}
\]
proved in \cite{RahmanIntegralPerverseObstructions}.  For
\(\frac14(1,1)\), the retained critical torsion is
\[
        \mathbb Z/4.
\]
\end{remark}

\begin{remark}
The higher-dimensional examples require the same convention but with a
different shift.  For an isolated threefold singularity, one must apply the
torsion-sensitive truncation to
\[
        R\Gamma(L,\mathbb Z)[3],
\]
not to \(R\Gamma(L,\mathbb Z)[2]\).  This is the source of the degree issue
for ordinary double points in dimension three.
\end{remark}

\subsection{Bockstein shadow convention}

The comparison with Benoist--Ottem requires one additional convention.  If
a local package \(E\) contains a subgroup selected by a Bockstein map, we
refer to that subgroup as a \emph{Bockstein shadow} of \(E\).

For the Coble boundary singularity \(\frac14(1,1)\), the full local package is
\[
        E\cong\mathbb Z/4.
\]
The local index-two cover
\[
        L(2,1)\longrightarrow L(4,1)
\]
has class
\[
        \eta\in H^1(L(4,1),\mathbb Z/2).
\]
The Bockstein associated to
\[
        0\to\mathbb Z\xrightarrow{2}\mathbb Z\to\mathbb Z/2\to0
\]
maps \(\eta\) to the unique order-two element of \(E\).  Thus
\[
        \operatorname{im}\beta=2E\cong\mathbb Z/2.
\]
This subgroup
\[
        2E\subset E
\]
is the Bockstein shadow relevant to the Benoist--Ottem comparison.

Equivalently, there is a short exact sequence
\[
        0
        \longrightarrow
        2E
        \longrightarrow
        E
        \longrightarrow
        E/2E
        \longrightarrow
        0,
\]
which, for \(E\cong\mathbb Z/4\), becomes
\[
        0
        \longrightarrow
        \mathbb Z/2
        \longrightarrow
        \mathbb Z/4
        \longrightarrow
        \mathbb Z/2
        \longrightarrow
        0.
\]
After pushing forward along a codimension-two stratum \(i:\Sigma\hookrightarrow Y\)
and shifting by \([2]\), this gives the distinguished triangle
\[
        i_*(\mathbb Z/2)_\Sigma[2]
        \longrightarrow
        i_*(\mathbb Z/4)_\Sigma[2]
        \longrightarrow
        i_*(\mathbb Z/2)_\Sigma[2]
        \overset{+1}{\longrightarrow}.
\]
This triangle is the formal expression of the order-two Benoist--Ottem
shadow inside the full Coble \(E\)-package.

\subsection{Motivic torsion convention}

Finally, we fix the motivic notation used later.  In an integral motivic
sheaf formalism with unit object \(\mathbf 1_\Sigma\), define
\[
        \mathbf 1_\Sigma/n
        :=
        \operatorname{Cone}
        \left(
        \mathbf 1_\Sigma
        \xrightarrow{n}
        \mathbf 1_\Sigma
        \right).
\]
Equivalently, in an integral motivic heart, this is the cokernel of
multiplication by \(n\):
\[
        \mathbf 1_\Sigma/n
        =
        \operatorname{coker}
        \left(
        \mathbf 1_\Sigma
        \xrightarrow{n}
        \mathbf 1_\Sigma
        \right).
\]
The Betti realization is
\[
        \operatorname{Real}_B(\mathbf 1_\Sigma/n)
        \cong
        (\mathbb Z/n)_\Sigma,
\]
provided realization commutes with cones and multiplication by \(n\).

Thus for a codimension-two Coble stratum, the full motivic \(E\)-object is
\[
        \mathcal T^{\mathrm{mot}}_{4,\Sigma}
        :=
        i_*(\mathbf 1_\Sigma/4)[2],
\]
with Betti realization
\[
        i_*(\mathbb Z/4)_\Sigma[2].
\]
The Benoist--Ottem order-two motivic shadow is
\[
        \mathcal T^{\mathrm{mot}}_{2,\Sigma}
        :=
        i_*(\mathbf 1_\Sigma/2)[2],
\]
with Betti realization
\[
        i_*(\mathbb Z/2)_\Sigma[2].
\]
The map
\[
        \mathbf 1_\Sigma/2
        \longrightarrow
        \mathbf 1_\Sigma/4
\]
realizes the inclusion
\[
        \mathbb Z/2\hookrightarrow\mathbb Z/4,
        \qquad
        1\longmapsto 2.
\]
Consequently
\[
        \mathcal T^{\mathrm{mot}}_{2,\Sigma}
        \longrightarrow
        \mathcal T^{\mathrm{mot}}_{4,\Sigma}
\]
is the motivic form of the order-two shadow
\[
        2E_\Sigma\subset E_\Sigma.
\]

\section{The local discriminant package}

In this section we recall the local invariant that supplies the input for the
global constructions in the rest of the paper.  The invariant is the finite
group \(E\) attached to a normal surface singularity.  It measures the
difference between the ordinary and dual middle-perversity intersection
complexes over \(\mathbb Z\).  The definition and the local distinguished
triangle are due to the integral middle-perversity framework of BBD and the
local surface analysis of Jung--Saito; the compatible topological,
resolution-theoretic, and monodromy realizations used here are those proved
in \cite{RahmanIntegralPerverseObstructions}.

The corrected Benoist--Ottem comparison later in the paper uses this same
local group \(E\), not a new invariant.  In the Enriques/Coble boundary
example, the relevant singularity is the cyclic quotient singularity
\(\frac14(1,1)\).  For that singularity,
\[
        E\cong\mathbb Z/4,
\]
and the Benoist--Ottem \(2\)-torsion mechanism selects the canonical
order-two subgroup
\[
        2E\cong\mathbb Z/2.
\]
Thus the Coble boundary example shows that a global \(2\)-torsion mechanism
may see a proper Bockstein-selected shadow of a larger local \(E\)-package.

\subsection{The perverse obstruction group \(E\)}

Let \((X,0)\) be a germ of a normal complex analytic surface.  After
shrinking \(X\), we assume that \(0\) is the unique singular point.  Put
\[
        U:=X\setminus\{0\},
\]
and let
\[
        j:U\hookrightarrow X,
        \qquad
        i:\{0\}\hookrightarrow X
\]
denote the inclusions.  Since \(X\) is normal of complex dimension \(2\), the
punctured germ \(U\) is smooth.

We recall the precise convention for the intersection complexes used in this
paper.  For an open inclusion \(j:U\hookrightarrow X\) and a perverse sheaf
\(F\) on \(U\), the intermediate extension \(j_{!*}F\) is defined in BBD as
the image, in the perverse heart, of the natural morphism
\[
        {}^pj_!F \longrightarrow {}^pRj_*F.
\]
Equivalently, \(j_{!*}F\) is the unique perverse extension of \(F\) to \(X\)
with no nonzero subobject or quotient object supported on the complement
\(X\setminus U\); see \cite[§1.4]{BBD82}.  In the integral setting BBD also
introduce the dual middle-perversity \(t\)-structure, which differs from the
ordinary middle perversity by the interchange of torsion and torsion-free
conditions in the critical degree; see \cite[Complement 3.3]{BBD82}.

With these conventions, the ordinary and dual middle-perversity intersection
complexes are defined by
\[
        {}^pIC_X\mathbb Z
        :=
        {}^pj_{!*}\mathbb Z_U[2],
        \qquad
        {}^p_+IC_X\mathbb Z
        :=
        {}^p_+j_{!*}\mathbb Z_U[2].
\]
Here the shift by \([2]\) is the standard normalization for a complex surface,
so that \(\mathbb Z_U[2]\) is perverse on the smooth complex surface \(U\);
see \cite[§2.1]{BBD82} and the convention used in
\cite[§2.2]{RahmanIntegralPerverseObstructions}.

By the defining property of intermediate extension, both complexes restrict
to the object from which they were extended:
\[
        j^*({}^pIC_X\mathbb Z)
        \cong
        \mathbb Z_U[2],
        \qquad
        j^*({}^p_+IC_X\mathbb Z)
        \cong
        \mathbb Z_U[2].
\]
Thus
\[
        j^*({}^pIC_X\mathbb Z)
        \cong
        \mathbb Z_U[2]
        \cong
        j^*({}^p_+IC_X\mathbb Z).
\]

The preceding display justifies the following support statement.

\begin{lemma}[Support of the discrepancy]
Let
\[
        u:{}^pIC_X\mathbb Z\longrightarrow {}^p_+IC_X\mathbb Z
\]
be the natural morphism induced by the identity of \(\mathbb Z_U[2]\) on
\(U\).  Then the cone
\[
        \operatorname{Cone}(u)
\]
is supported at the singular point \(0\).
\end{lemma}

\begin{proof}
By construction, the restriction of \(u\) to \(U\) is the identity morphism
of the common restriction \(\mathbb Z_U[2]\).  Therefore
\[
        j^*u:
        j^*({}^pIC_X\mathbb Z)
        \longrightarrow
        j^*({}^p_+IC_X\mathbb Z)
\]
is an isomorphism.  Applying the exact functor \(j^*\) to the distinguished
triangle
\[
        {}^pIC_X\mathbb Z
        \xrightarrow{u}
        {}^p_+IC_X\mathbb Z
        \longrightarrow
        \operatorname{Cone}(u)
        \longrightarrow
\]
gives
\[
        j^*\operatorname{Cone}(u)=0.
\]
Hence the support of \(\operatorname{Cone}(u)\) is contained in the closed
complement \(X\setminus U\).  Since \(X\setminus U=\{0\}\), the cone is
supported at the singular point.  This is the support argument used in
\cite[Lemma 2.1]{RahmanIntegralPerverseObstructions}.
\end{proof}

We now define the local group.

\begin{definition}[The local perverse obstruction group]
Let \((X,0)\) be a normal complex surface germ.  Define
\[
        E_X(0):=
        H^0\bigl({}^p_+IC_X\mathbb Z\bigr)_0.
\]
When the singularity is clear from context, we write simply \(E\).
\end{definition}

The next statement is not part of the definition.  It is the local
codimension-two theorem that makes the definition useful.

\begin{proposition}[Local discrepancy triangle]
Let \((X,0)\) be a normal complex surface germ, and set
\[
        E=H^0\bigl({}^p_+IC_X\mathbb Z\bigr)_0.
\]
Then \(E\) is finite, and there is a distinguished triangle
\[
        {}^pIC_X\mathbb Z
        \longrightarrow
        {}^p_+IC_X\mathbb Z
        \longrightarrow
        E[1]
        \longrightarrow .
\]
Here \(E\) is regarded as a complex supported at the singular point.  The
triangle is self-dual under Verdier duality.
\end{proposition}

\begin{proof}
The support of the cone of
\[
        {}^pIC_X\mathbb Z
        \longrightarrow
        {}^p_+IC_X\mathbb Z
\]
is contained in \(\{0\}\) by the preceding lemma.  The remaining assertion is
the codimension-two local surface statement of Jung--Saito: for a normal
surface germ, the discrepancy between the ordinary and dual integral
middle-perversity intersection complexes is represented by the finite group
\[
        H^0\bigl({}^p_+IC_X\mathbb Z\bigr)_0
\]
placed in cohomological degree \(1\), and the resulting triangle is
self-dual.  In the notation of this paper, this is the local surface
statement cited in \cite[Remark 6.4]{JungSaitoFactoriality} and recorded as
the integral perverse obstruction triangle in
\cite[Theorem 1.1 and Proposition 2.11]{RahmanIntegralPerverseObstructions}.
\end{proof}

\begin{remark}
After tensoring with \(\mathbb Q\), the ordinary and dual middle-perversity
intersection complexes agree:
\[
        {}^pIC_X\mathbb Q
        \cong
        {}^p_+IC_X\mathbb Q.
\]
This is the rational middle-perversity situation of BBD
\cite{BBD82}.  Therefore the group \(E\) is purely integral and torsion.  In
the local surface setting this is made explicit by the realization
\[
        E\cong H^2(L,\mathbb Z)_{\mathrm{tors}}
\]
proved in \cite[Proposition 3.4]{RahmanIntegralPerverseObstructions}.
\end{remark}

\subsection{The six-realization theorem}

We now recall the local computation of \(E\) from
\cite{RahmanIntegralPerverseObstructions}.  Let \(L\) be the link of the
singularity.  Let
\[
        \pi:\widetilde X\longrightarrow X
\]
be the minimal resolution, and write
\[
        E_\pi=\bigcup_{a=1}^r C_a
\]
for the exceptional divisor.  Let \(\Lambda\) be the lattice generated by
the irreducible exceptional curves \(C_a\), equipped with the geometric
intersection form
\[
        (C_a,C_b)=C_a\cdot C_b.
\]
Let
\[
        M=(C_a\cdot C_b)_{a,b}
\]
be the corresponding intersection matrix.

For a normal surface singularity, the link \(L\) is a compact oriented
three-manifold, and the boundary of a sufficiently small resolution
neighborhood of \(E_\pi\) is diffeomorphic to \(L\).  This is the standard
local topology of normal surface singularities; see Mumford
\cite{Mu61}.  The intersection matrix \(M\) is negative definite in the
surface-resolution setting, by the classical contraction criterion of
Grauert; in this paper we use the same negative-definite geometric
convention as in \cite[§4]{RahmanIntegralPerverseObstructions}.

\begin{theorem}[Local realization theorem]
Let \((X,0)\) be a normal surface singularity with link \(L\), minimal
resolution exceptional lattice \(\Lambda\), and intersection matrix \(M\).
Then the local perverse obstruction group satisfies
\[
        E_X(0)
        \cong
        H^2(L,\mathbb Z)_{\mathrm{tors}}
        \cong
        \Lambda^\vee/\Lambda.
\]
In particular,
\[
        |E_X(0)|=|\det(M)|.
\]
If \((X,0)\) is an isolated hypersurface surface singularity and \(T\) is
the Milnor monodromy on integral vanishing cohomology, then
\[
        E_X(0)
        \cong
        \operatorname{coker}(T-\mathrm{id})_{\mathrm{tors}}.
\]
Under the additional hypothesis that
\[
        (T-\mathrm{id})\otimes_{\mathbb Z}\mathbb Q
\]
is an isomorphism, one has
\[
        |E_X(0)|=
        |\det(T-\mathrm{id})|.
\]
\end{theorem}

\begin{proof}
The first realization of \(E_X(0)\) is perverse: by the local discrepancy
triangle above, \(E_X(0)\) is the finite point-supported correction between
\({}^pIC_X\mathbb Z\) and \({}^p_+IC_X\mathbb Z\).

The topological realization
\[
        E_X(0)\cong H^2(L,\mathbb Z)_{\mathrm{tors}}
\]
is proved in \cite[Proposition 3.4]{RahmanIntegralPerverseObstructions}.
The proof uses the local stalk computation
\[
        i^*Rj_*\mathbb Z_U[2]\cong R\Gamma(L,\mathbb Z)[2],
\]
so that the degree-zero stalk sees \(H^2(L,\mathbb Z)\).  The ordinary and
dual middle-perversity conditions then identify the discrepancy with the
torsion subgroup of this group.  The general topological context for
linking pairings and torsion phenomena in singular links is developed in
\cite{GoreskySiegel83}.

The resolution-lattice realization
\[
        H^2(L,\mathbb Z)_{\mathrm{tors}}
        \cong
        \Lambda^\vee/\Lambda
\]
is proved in \cite[Theorem 1.2 and §4]{RahmanIntegralPerverseObstructions}.
The argument uses the long exact sequence of the pair \((N,L)\), where \(N\)
is a resolution neighborhood of the exceptional divisor and
\(\partial N=L\).  Poincaré--Lefschetz duality identifies the relevant
relative homology group with the dual lattice \(\Lambda^\vee\).  Under this
identification, the map
\[
        H_2(N,\mathbb Z)\longrightarrow H_2(N,L;\mathbb Z)
\]
is represented by the intersection matrix \(M\), equivalently by the lattice
homomorphism
\[
        \Lambda\longrightarrow \Lambda^\vee.
\]
Therefore the finite quotient is
\[
        \operatorname{coker}(\Lambda\to\Lambda^\vee)
        =
        \Lambda^\vee/\Lambda.
\]
Consequently
\[
        |E_X(0)|
        =
        |\Lambda^\vee/\Lambda|
        =
        |\det(M)|.
\]
The discriminant-group formalism for integral lattices is standard; see
\cite{Nikulin80}.

Finally, assume that \((X,0)\) is an isolated hypersurface surface
singularity.  The Milnor fibration identifies the link with the mapping
torus of the Milnor fiber, and the Wang sequence relates the cohomology of
the link to the monodromy \(T\).  In the notation of this paper, this gives
\[
        H^2(L,\mathbb Z)_{\mathrm{tors}}
        \cong
        \operatorname{coker}(T-\mathrm{id})_{\mathrm{tors}}.
\]
This comparison is proved in
\cite[Theorem 1.3 and §5]{RahmanIntegralPerverseObstructions}, using the
standard Milnor-fibration background from \cite{Mi68}.  If
\[
        (T-\mathrm{id})\otimes_{\mathbb Z}\mathbb Q
\]
is an isomorphism, then \(\operatorname{coker}(T-\mathrm{id})\) is finite.
For an endomorphism of a free abelian group whose rationalization is an
isomorphism, the order of the finite cokernel is the absolute value of the
determinant.  Hence
\[
        |E_X(0)|=|\det(T-\mathrm{id})|.
\]
This determinant refinement is recorded in
\cite[Corollary 1.4]{RahmanIntegralPerverseObstructions}.
\end{proof}

\begin{remark}
The theorem packages six compatible realizations of the same local group:
\begin{enumerate}[label=\textup{(\arabic*)}]
\item the perverse discrepancy between
\({}^pIC_X\mathbb Z\) and \({}^p_+IC_X\mathbb Z\);
\item the torsion-sensitive truncation discrepancy in the sense of Friedman
\cite{FriedmanGenIH,FriedmanBook20,FriedmanTsInv};
\item the torsion subgroup \(H^2(L,\mathbb Z)_{\mathrm{tors}}\) of the link
\cite[Proposition 3.4]{RahmanIntegralPerverseObstructions};
\item the discriminant group \(\Lambda^\vee/\Lambda\) of the exceptional
lattice \cite[Theorem 1.2]{RahmanIntegralPerverseObstructions};
\item the boundary term in the long exact sequence of a resolution
neighborhood and its boundary link
\cite[§4]{RahmanIntegralPerverseObstructions};
\item and, for isolated hypersurface surface singularities, the torsion
cokernel of \(T-\mathrm{id}\) in the Wang sequence
\cite[Theorem 1.3]{RahmanIntegralPerverseObstructions}.
\end{enumerate}
The present paper uses this compatibility as the local input for
local-to-global transport.  For non-hypersurface quotient singularities, such
as the Coble boundary singularity \(\frac14(1,1)\), the monodromy station is
not part of the hypersurface realization, but the perverse,
torsion-sensitive, link, lattice, and pair-sequence realizations still agree.
\end{remark}

\subsection{The Coble boundary example and the order-two shadow}

We now record the quotient example needed for the corrected
Benoist--Ottem comparison.  This example uses the same local invariant
\(E\) from \cite{RahmanIntegralPerverseObstructions}; the new feature is that
the global \(2\)-torsion mechanism selects a proper subgroup of the local
group.

\begin{proposition}[The Coble boundary \(E\)-package]
Let
\[
        X=\mathbb C^2/\mu_4,
        \qquad
        \zeta\cdot(x,y)=(\zeta x,\zeta y),
\]
be the cyclic quotient singularity of type \(\frac14(1,1)\).  Then
\[
        E_X(0)\cong\mathbb Z/4.
\]
More precisely, the link is
\[
        L=L(4,1),
\]
and
\[
        E_X(0)
        \cong
        H^2(L(4,1),\mathbb Z)_{\mathrm{tors}}
        \cong
        \mathbb Z/4.
\]
The minimal resolution has one exceptional curve with self-intersection
\(-4\), so the exceptional lattice is
\[
        \Lambda=\mathbb Z\langle e\rangle,
        \qquad
        (e,e)=-4,
\]
and
\[
        \Lambda^\vee/\Lambda\cong\mathbb Z/4.
\]
\end{proposition}

\begin{proof}
The action of \(\mu_4\) on \(\mathbb C^2\setminus\{0\}\) restricts to a free
action on the unit sphere \(S^3\).  Thus the link is
\[
        L=S^3/\mu_4=L(4,1).
\]
The lens space \(L(4,1)\) has
\[
        H_1(L(4,1),\mathbb Z)\cong\mathbb Z/4,
        \qquad
        H_2(L(4,1),\mathbb Z)=0.
\]
By the universal coefficient theorem,
\[
        H^2(L(4,1),\mathbb Z)
        \cong
        \operatorname{Ext}^1(H_1(L(4,1),\mathbb Z),\mathbb Z)
        \cong
        \mathbb Z/4.
\]
The local realization theorem identifies
\[
        E_X(0)
        \cong
        H^2(L(4,1),\mathbb Z)_{\mathrm{tors}},
\]
so \(E_X(0)\cong\mathbb Z/4\).

The Hirzebruch--Jung continued fraction for \(\frac41\) is
\[
        \frac41=[4].
\]
Thus the minimal resolution has a single exceptional curve with
self-intersection \(-4\).  Hence
\[
        \Lambda=\mathbb Z\langle e\rangle,
        \qquad
        (e,e)=-4.
\]
The dual lattice is generated by
\[
        e^\vee=-\frac e4,
\]
and therefore
\[
        \Lambda^\vee/\Lambda\cong\mathbb Z/4.
\]
This agrees with the link realization of \(E_X(0)\).
\end{proof}

\begin{proposition}[Local index-two cover and Bockstein shadow]
Let \(X=\mathbb C^2/\mu_4\) be the \(\frac14(1,1)\) singularity.  Let
\(\mu_2\subset\mu_4\) be the subgroup generated by \(-1\).  Then the
intermediate quotient
\[
        \mathbb C^2/\mu_2
        \longrightarrow
        \mathbb C^2/\mu_4
\]
is a degree-two cover whose source is the \(A_1\) singularity.  On links,
this cover is
\[
        L(2,1)\longrightarrow L(4,1).
\]
Let
\[
        \eta\in H^1(L(4,1),\mathbb Z/2)
\]
be the class of this double cover.  Then the Bockstein
\[
        \beta:
        H^1(L(4,1),\mathbb Z/2)
        \longrightarrow
        H^2(L(4,1),\mathbb Z)
\]
associated to
\[
        0\to\mathbb Z\xrightarrow{2}\mathbb Z\to\mathbb Z/2\to0
\]
sends \(\eta\) to the unique element of order \(2\) in
\[
        E_X(0)\cong\mathbb Z/4.
\]
Equivalently,
\[
        \operatorname{im}\beta=2E_X(0)\cong\mathbb Z/2.
\]
\end{proposition}

\begin{proof}
The subgroup \(\mu_2=\{\pm1\}\subset\mu_4\) acts on \(\mathbb C^2\) by
\[
        (x,y)\longmapsto(-x,-y).
\]
Thus \(\mathbb C^2/\mu_2\) is the \(A_1\) surface singularity.  Since
\[
        \mu_4/\mu_2\cong\mathbb Z/2,
\]
there is a residual degree-two quotient
\[
        \mathbb C^2/\mu_2
        \longrightarrow
        \mathbb C^2/\mu_4.
\]
Restricting to the unit sphere gives the link cover
\[
        S^3/\mu_2=L(2,1)
        \longrightarrow
        S^3/\mu_4=L(4,1).
\]

Now consider the Bockstein sequence.  The relevant part of the long exact
sequence is
\[
        H^1(L,\mathbb Z)
        \longrightarrow
        H^1(L,\mathbb Z/2)
        \xrightarrow{\beta}
        H^2(L,\mathbb Z)
        \xrightarrow{2}
        H^2(L,\mathbb Z),
\]
where \(L=L(4,1)\).  Since
\[
        H^1(L,\mathbb Z)=0
\]
and
\[
        H^2(L,\mathbb Z)\cong\mathbb Z/4,
\]
the Bockstein is injective and its image is
\[
        \ker(2:\mathbb Z/4\to\mathbb Z/4)
        =
        \{0,2\}
        \cong\mathbb Z/2.
\]
The nonzero class \(\eta\in H^1(L,\mathbb Z/2)\) therefore maps to the unique
order-two element in \(H^2(L,\mathbb Z)\).  Under the identification
\[
        E_X(0)\cong H^2(L,\mathbb Z)_{\mathrm{tors}},
\]
this says
\[
        \operatorname{im}\beta=2E_X(0)\cong\mathbb Z/2.
\]
\end{proof}

\begin{remark}
This is the corrected local model for the Benoist--Ottem comparison.  The
full local Coble package is
\[
        E\cong\mathbb Z/4,
\]
but the Benoist--Ottem \(2\)-torsion mechanism detects the Bockstein-selected
order-two shadow
\[
        2E\cong\mathbb Z/2.
\]
Thus the comparison is not with the full \(E\)-package, but with a canonical
filtered piece of it.
\end{remark}

\subsection{Relation with the local class group}

We now record the divisor-theoretic interpretation of the same finite group.
This comparison is important because it connects the perverse obstruction
\(E\) with factoriality and \(Q\)-factoriality questions.

The cleanest statement holds in the standard rational surface singularity
setting.  For rational surface singularities, the local divisor class group
is computed by the exceptional intersection matrix of a resolution.  This is
part of the classical resolution theory of rational surface singularities;
see Artin \cite{Artin66} and Lipman \cite{Lipman69}.  We state the form used
below.

\begin{proposition}[Local class group comparison]
Let \((X,0)\) be a rational normal surface singularity, and let
\[
        \pi:\widetilde X\longrightarrow X
\]
be its minimal resolution.  Let \(\Lambda\) be the exceptional lattice
generated by the irreducible exceptional curves, with intersection matrix
\(M\).  Then the local divisor class group is naturally identified with the
discriminant group of the exceptional lattice:
\[
        \operatorname{Cl}(\mathcal O_{X,0})
        \cong
        \Lambda^\vee/\Lambda.
\]
Consequently, using the local realization theorem,
\[
        \operatorname{Cl}(\mathcal O_{X,0})
        \cong
        E_X(0).
\]
\end{proposition}

\begin{proof}
Let
\[
        E_\pi=\bigcup_{a=1}^r C_a
\]
be the exceptional divisor.  The curves \(C_a\) generate the lattice
\[
        \Lambda=\bigoplus_a\mathbb Z[C_a].
\]
The intersection form defines a homomorphism
\[
        \Lambda
        \longrightarrow
        \Lambda^\vee,
        \qquad
        C_a\longmapsto (C_a\cdot -).
\]
With respect to the basis given by the curves \(C_a\), this homomorphism is
represented by the intersection matrix \(M\).

For rational normal surface singularities, the standard
resolution-theoretic computation of the local divisor class group identifies
\[
        \operatorname{Cl}(\mathcal O_{X,0})
\]
with the cokernel of the intersection homomorphism
\[
        \Lambda\longrightarrow\Lambda^\vee.
\]
Thus
\[
        \operatorname{Cl}(\mathcal O_{X,0})
        \cong
        \operatorname{coker}(\Lambda\to\Lambda^\vee)
        =
        \Lambda^\vee/\Lambda.
\]
This is the classical class-group computation in the rational surface
singularity setting; see \cite{Artin66,Lipman69}.  By the local realization
theorem above, also
\[
        E_X(0)\cong \Lambda^\vee/\Lambda.
\]
Combining the two identifications gives
\[
        \operatorname{Cl}(\mathcal O_{X,0})
        \cong
        E_X(0).
\]
\end{proof}

\begin{remark}
The rationality hypothesis in the proposition is important.  For a general
normal surface singularity, the relationship between the local class group
and the exceptional intersection lattice can involve additional
Picard-theoretic data from the resolution.  The proposition is therefore
stated only in the setting where the standard resolution-lattice description
of the local class group applies.
\end{remark}

\begin{remark}
This comparison connects the perverse obstruction group \(E\) to the
factoriality and \(Q\)-factoriality questions motivating the Jung--Saito
framework \cite{JungSaitoFactoriality}.  It also motivates the later
comparison with global Brauer groups.  Local class groups measure local
failure of factoriality, while Brauer groups and unramified cohomology
measure global torsion phenomena that enter the integral Hodge problem
\cite{BlochOgus74,CTVoisin12}.
\end{remark}

\subsection{The discriminant package \((E,q)\)}

The preceding theorem identifies the local group \(E\) with the discriminant
group of the exceptional lattice.  This gives \(E\) more structure than the
structure of a finite abelian group.

Let \((\Lambda,(\ ,\ ))\) be the exceptional lattice.  Its dual lattice is
\[
        \Lambda^\vee
        :=
        \operatorname{Hom}_{\mathbb Z}(\Lambda,\mathbb Z)
        \cong
        \{x\in \Lambda\otimes_{\mathbb Z}\mathbb Q
        \mid (x,\lambda)\in\mathbb Z
        \text{ for all }\lambda\in\Lambda\}.
\]
The discriminant group is
\[
        A_\Lambda:=\Lambda^\vee/\Lambda.
\]
The discriminant pairing is
\[
        q_\Lambda:
        A_\Lambda\times A_\Lambda
        \longrightarrow
        \mathbb Q/\mathbb Z
\]
defined by
\[
        q_\Lambda(x+\Lambda,y+\Lambda)
        =
        (x,y)\bmod \mathbb Z.
\]
This construction is standard for nondegenerate integral lattices
\cite{Nikulin80}.

\begin{definition}[Local discriminant package]
Let \((X,0)\) be a normal surface singularity whose local obstruction group
is identified with the exceptional discriminant group:
\[
        E_X(0)\cong \Lambda^\vee/\Lambda.
\]
The local discriminant package is the pair
\[
        (E,q)
        :=
        (\Lambda^\vee/\Lambda,q_\Lambda),
\]
where \(q_\Lambda\) is the discriminant pairing induced by the exceptional
intersection form.
\end{definition}

The pairing \(q_\Lambda\) is nondegenerate because the intersection form is
nondegenerate over \(\mathbb Q\).  Via the link realization, this finite
pairing agrees, up to the sign convention fixed in Section~2, with the
torsion linking pairing on the link.  We record this as a proposition.

\begin{proposition}[Discriminant and linking pairings]
Let \((X,0)\) be a normal surface singularity for which the exceptional
lattice realization
\[
        E_X(0)\cong \Lambda^\vee/\Lambda
\]
is fixed.  Let \(L\) be the link.  Under the canonical identification
\[
        E_X(0)
        \cong
        H^2(L,\mathbb Z)_{\mathrm{tors}}
        \cong
        \Lambda^\vee/\Lambda,
\]
the discriminant pairing on \(\Lambda^\vee/\Lambda\) corresponds, up to the
orientation sign convention, to the torsion linking pairing of the link.
\end{proposition}

\begin{proof}
Let \(N\) be a resolution neighborhood of the exceptional divisor, with
boundary \(L=\partial N\).  The long exact sequence of the pair \((N,L)\)
contains the segment
\[
        H_2(N,\mathbb Z)
        \longrightarrow
        H_2(N,L;\mathbb Z)
        \longrightarrow
        H_1(L,\mathbb Z)
        \longrightarrow
        H_1(N,\mathbb Z).
\]
Poincaré--Lefschetz duality identifies the relevant relative homology group
with the dual lattice \(\Lambda^\vee\), and the map from \(H_2(N,\mathbb Z)\)
to \(H_2(N,L;\mathbb Z)\) is represented by the exceptional intersection
form.  Hence the finite quotient controlling the torsion in the boundary
link is \(\Lambda^\vee/\Lambda\).

The linking pairing on \(H_1(L,\mathbb Z)_{\mathrm{tors}}\) is defined by
taking a torsion cycle, multiplying it so that it bounds, and measuring the
resulting rational intersection number modulo \(\mathbb Z\).  Under the
above identification of the boundary torsion with
\(\Lambda^\vee/\Lambda\), this construction is exactly the inverse
intersection-form construction of the discriminant pairing, up to the
orientation sign induced by the boundary convention.  This is the standard
comparison between link pairings and discriminant forms
\cite{GoreskySiegel83,Nikulin80}.  The compatibility with the local group
\(E_X(0)\) in the present notation is part of the realization theorem proved
in \cite{RahmanIntegralPerverseObstructions}.
\end{proof}

\begin{example}[The Coble discriminant form]
For the Coble boundary singularity \(\frac14(1,1)\), the exceptional lattice
is generated by one class \(e\) with
\[
        (e,e)=-4.
\]
The dual generator is
\[
        e^\vee=-\frac e4.
\]
Thus
\[
        E\cong\Lambda^\vee/\Lambda\cong\mathbb Z/4.
\]
If \(\bar g\) is the class of \(e^\vee\), then
\[
        q(\bar g,\bar g)
        =
        (e^\vee,e^\vee)
        =
        -\frac14
        \quad\bmod\mathbb Z.
\]
The unique order-two subgroup is
\[
        2E\cong\mathbb Z/2.
\]
Its generator \(2\bar g\) satisfies
\[
        q(2\bar g,2\bar g)
        =
        4q(\bar g,\bar g)
        =
        -1
        \equiv0
        \quad\bmod\mathbb Z.
\]
Thus the Benoist--Ottem order-two shadow is an isotropic subgroup of the
full Coble discriminant package.
\end{example}

\subsection{Why the form, not only the group, must be tracked}

The local-to-global constructions below are formulated for packages
\[
        (E_i,q_i),
\]
not merely for the underlying finite groups \(E_i\).  The reason is that the
same abstract finite abelian group may carry inequivalent discriminant
pairings.  A construction that forgets the pairing may therefore lose
information visible both in the local topology and in the resolution
lattice.

The basic example is the contrast between a \(D_4\) rational double point and
two independent \(A_1\) rational double points.  The \(D_4\) exceptional
lattice has discriminant group
\[
        \Lambda_{D_4}^{\vee}/\Lambda_{D_4}
        \cong
        \mathbb Z/2\mathbb Z\oplus \mathbb Z/2\mathbb Z.
\]
Two independent \(A_1\) singularities also produce the underlying group
\[
        \mathbb Z/2\mathbb Z\oplus \mathbb Z/2\mathbb Z.
\]
However, the discriminant pairings need not be the same.  The distinction is
invisible at the level of groups but visible at the level of finite
\(\mathbb Q/\mathbb Z\)-valued pairings.  The general theory of
discriminant forms records precisely this extra information
\cite{Nikulin80}.

The Coble boundary example gives a second reason to keep the full package.
There the full local group is
\[
        E\cong\mathbb Z/4,
\]
but the Benoist--Ottem \(2\)-torsion mechanism sees only the subgroup
\[
        2E\cong\mathbb Z/2.
\]
Thus even when a global comparison is \(2\)-torsion, the local boundary
package may be larger.  The correct statement is not that the global
\(2\)-torsion equals the full local group; rather, it is that the global
\(2\)-torsion selects a Bockstein shadow inside the local package.

\begin{remark}
This is the reason the global transport problem is stated for
\[
        \bigoplus_i(E_i,q_i)
\]
rather than for \(\bigoplus_iE_i\) alone.  Later compatibility questions ask
whether the local discriminant pairings map to natural global torsion
pairings on the image in cohomology, in the Brauer group, or in unramified
cohomology.  The ADE computations in the appendix provide the first test
cases for this form-sensitive version of the theory.  The Coble computation
adds a filtered-subgroup test case, where the comparison with
Benoist--Ottem uses \(2E\subset E\) rather than the full local group.
\end{remark}

\section{Excision and local-to-global maps}

The purpose of this section is to isolate the formal local-to-global
mechanism used later in the paper.  The input is a projective variety \(X\)
with isolated singularities and a resolution
\[
        \pi:\widetilde X\longrightarrow X.
\]
The exceptional locus over each singular point has a small resolution
neighborhood, and cohomology with supports gives a natural way to compare
local torsion data with global cohomology.  The key point is that this
comparison is degree-dependent.  We therefore define degree-indexed maps
\[
        \alpha_X^{(k)}
\]
rather than a single degree-free local-to-global map.

A second point will be used repeatedly.  If a local torsion package contains
a canonical subgroup, then this subgroup carries its own local-to-global
transport problem.  In particular, for the Coble boundary singularity
\(\frac14(1,1)\), the full local package is
\[
        E\cong\mathbb Z/4,
\]
but the Benoist--Ottem \(2\)-torsion mechanism selects the order-two subgroup
\[
        2E\cong\mathbb Z/2.
\]
Thus in addition to transporting \(E\), we must also be able to transport
subpackages such as \(2E\subset E\), together with the associated quotient
sequence
\[
        0\longrightarrow 2E\longrightarrow E\longrightarrow E/2E\longrightarrow0.
\]

We work throughout with singular cohomology with integral coefficients on
the associated analytic spaces.  The formal properties used below are the
standard long exact sequence for cohomology with supports and the excision
isomorphism; see, for example, Bredon's treatment of sheaf cohomology with
supports and the usual relative-cohomology formulation of excision
\cite{BredonSheaf97,Hatcher02}.  In the resolution-neighborhood calculations
for surface singularities, the relevant local pair sequence is the same one
used in the comparison between link torsion and exceptional discriminant
groups in \cite{GoreskySiegel83,RahmanIntegralPerverseObstructions}.

\subsection{Resolution neighborhoods}

Let \(X\) be a projective complex variety with isolated singularities
\[
        \operatorname{Sing}(X)=\{p_1,\ldots,p_s\}.
\]
Let
\[
        \pi:\widetilde X\longrightarrow X
\]
be a resolution of singularities.  For each singular point \(p_i\), set
\[
        D_i:=\pi^{-1}(p_i),
        \qquad
        D:=\bigcup_{i=1}^sD_i.
\]
Since the singular points are isolated and finite in number, the exceptional
sets \(D_i\) are pairwise disjoint closed analytic subsets of
\(\widetilde X\).  Choose pairwise disjoint sufficiently small closed
neighborhoods
\[
        N_i\subset \widetilde X
\]
of \(D_i\).  We call \(N_i\) a resolution neighborhood of \(D_i\).  Its
boundary is denoted
\[
        L_i:=\partial N_i.
\]
For \(N_i\) sufficiently small, \(L_i\) is naturally identified with the link
of the singularity \(p_i\).  In the normal surface case this is the standard
description of the link as the boundary of a resolution neighborhood of the
exceptional divisor; see Mumford's foundational discussion of links of normal
surface singularities \cite{Mu61} and the resolution-neighborhood setup used
in \cite[§3--§4]{RahmanIntegralPerverseObstructions}.

Set
\[
        N:=\bigcup_{i=1}^s N_i.
\]
Because the \(N_i\) are pairwise disjoint, we have
\[
        N=\coprod_{i=1}^s N_i.
\]
The complement of the exceptional locus is canonically identified with the
smooth locus of \(X\):
\[
        \widetilde X\setminus D
        \cong
        X\setminus \operatorname{Sing}(X).
\]
This follows from the definition of a resolution: the morphism
\(\pi\) is an isomorphism over the smooth locus of \(X\), and all points
where \(\pi\) fails to be an isomorphism lie over the singular set.

\begin{remark}
The boundary \(L_i=\partial N_i\) is used in two related but distinct ways.
First, it computes local link cohomology.  Second, the pair \((N_i,L_i)\)
records how exceptional cycles in \(N_i\) meet the boundary.  In the surface
case, this pair sequence is the mechanism that identifies the torsion in the
link with the discriminant group of the exceptional lattice
\cite{GoreskySiegel83,RahmanIntegralPerverseObstructions}.  For the Coble
boundary singularity \(\frac14(1,1)\), this mechanism gives
\[
        E\cong H^2(L(4,1),\mathbb Z)_{\mathrm{tors}}
        \cong\mathbb Z/4.
\]
\end{remark}

\subsection{Support cohomology and excision}

We now recall the formal exact sequence.  Let \(Y\) be a locally compact
topological space, let \(A\subset Y\) be a closed subset, and let
\(\mathbb Z\) denote the constant sheaf.  Cohomology with supports in \(A\)
is denoted
\[
        H^k_A(Y,\mathbb Z).
\]
For sufficiently nice pairs, and in particular for the analytic pairs used
here, this agrees with relative cohomology:
\[
        H^k_A(Y,\mathbb Z)
        \cong
        H^k(Y,Y\setminus A;\mathbb Z).
\]
This identification is standard in sheaf cohomology with supports and
relative cohomology; see \cite{BredonSheaf97,Hatcher02}.

Applying this to \(Y=\widetilde X\) and \(A=D\), we obtain
\[
        H^k_D(\widetilde X,\mathbb Z)
        \cong
        H^k(\widetilde X,\widetilde X\setminus D;\mathbb Z).
\]
The long exact sequence of the pair
\[
        (\widetilde X,\widetilde X\setminus D)
\]
therefore gives
\[
\cdots
\longrightarrow
H^k_D(\widetilde X,\mathbb Z)
\longrightarrow
H^k(\widetilde X,\mathbb Z)
\longrightarrow
H^k(\widetilde X\setminus D,\mathbb Z)
\longrightarrow
H^{k+1}_D(\widetilde X,\mathbb Z)
\longrightarrow
\cdots .
\]
The first map
\[
        H^k_D(\widetilde X,\mathbb Z)
        \longrightarrow
        H^k(\widetilde X,\mathbb Z)
\]
is the forget-support map.

We also need the local decomposition of support cohomology.

\begin{lemma}[Excision for the exceptional support]
With notation as above, there is a natural isomorphism
\[
        H^k_D(\widetilde X,\mathbb Z)
        \cong
        \bigoplus_{i=1}^s H^k_{D_i}(N_i,\mathbb Z).
\]
\end{lemma}

\begin{proof}
By definition of cohomology with supports,
\[
        H^k_D(\widetilde X,\mathbb Z)
        \cong
        H^k(\widetilde X,\widetilde X\setminus D;\mathbb Z).
\]
Since \(N=\coprod_iN_i\) is a neighborhood of \(D\), excision for relative
cohomology gives an isomorphism
\[
        H^k(\widetilde X,\widetilde X\setminus D;\mathbb Z)
        \cong
        H^k(N,N\setminus D;\mathbb Z).
\]
Because the neighborhoods \(N_i\) are pairwise disjoint, the pair
\((N,N\setminus D)\) is the disjoint union of the pairs
\[
        (N_i,N_i\setminus D_i).
\]
Relative cohomology of a finite disjoint union is the direct sum of the
relative cohomology groups of the components.  Hence
\[
        H^k(N,N\setminus D;\mathbb Z)
        \cong
        \bigoplus_{i=1}^s
        H^k(N_i,N_i\setminus D_i;\mathbb Z).
\]
Finally,
\[
        H^k(N_i,N_i\setminus D_i;\mathbb Z)
        \cong
        H^k_{D_i}(N_i,\mathbb Z)
\]
by the same identification of support cohomology with relative cohomology.
Combining these isomorphisms proves the claim.
\end{proof}

Combining the long exact sequence with the excision isomorphism gives the
formal local-to-global channel:
\[
        \bigoplus_i H^k_{D_i}(N_i,\mathbb Z)
        \cong
        H^k_D(\widetilde X,\mathbb Z)
        \longrightarrow
        H^k(\widetilde X,\mathbb Z).
\]
This map is the cohomological mechanism behind the transport maps used
below.

\subsection{The map \(\alpha_X^{(k)}\)}

We now define the degree-\(k\) transport map.  The definition is deliberately
conditional: a local package can be transported in degree \(k\) only after it
has been realized in degree \(k\) support cohomology.  This is necessary
because the local group \(E\) for a surface singularity naturally appears
through \(H^2\) of a three-dimensional link, whereas the threefold
Brauer/unramified comparison passes through \(H^3\) of the resolution.

\begin{definition}[Degree-\(k\) discriminant transport]
Suppose that, for each \(i\), a local torsion package \(E_i\) is equipped
with a specified homomorphism
\[
        \eta_i^{(k)}:E_i
        \longrightarrow
        H^k_{D_i}(N_i,\mathbb Z)_{\mathrm{tors}}.
\]
The degree-\(k\) discriminant transport map is the composite
\[
        \alpha_X^{(k)}:
        \bigoplus_iE_i
        \xrightarrow{\oplus_i\eta_i^{(k)}}
        \bigoplus_iH^k_{D_i}(N_i,\mathbb Z)_{\mathrm{tors}}
        \xrightarrow[\cong]{\mathrm{excision}}
        H^k_D(\widetilde X,\mathbb Z)_{\mathrm{tors}}
        \longrightarrow
        H^k(\widetilde X,\mathbb Z)_{\mathrm{tors}},
\]
where the last map is the forget-support map.
\end{definition}

\begin{remark}
The maps \(\eta_i^{(k)}\) are part of the data.  They are not automatic for
all singularities and all degrees.  In the normal surface case, the previous
paper identifies the local group with link torsion and with the
discriminant group of the exceptional lattice:
\[
        E_i
        \cong
        H^2(L_i,\mathbb Z)_{\mathrm{tors}}
        \cong
        \Lambda_i^\vee/\Lambda_i
\]
\cite{RahmanIntegralPerverseObstructions}.  To obtain a map
\(\alpha_X^{(k)}\), one must further identify how this local group enters
the degree-\(k\) support sequence of the chosen global resolution.
\end{remark}

\begin{remark}
The degree index is essential.  The notation \(\alpha_X^{(k)}\) records the
target cohomological degree.  For applications to the Brauer group of a
smooth projective threefold, the relevant degree is \(k=3\), because the
exponential sequence relates \(\operatorname{Br}(\widetilde X)\) to
\(H^3(\widetilde X,\mathbb Z)_{\mathrm{tors}}\) under suitable hypotheses.
Thus the relevant transport map in that setting is
\[
        \alpha_X^{(3)}.
\]
\end{remark}

\subsection{Subpackages, quotients, and shadow transport}

The Coble boundary example requires one additional formal convention.  A
global comparison may use a subgroup or quotient of a local package rather
than the full local group.  This is not exceptional: if
\[
        A_i\subset E_i
\]
is a canonical subgroup, then the inclusion gives a short exact sequence
\[
        0\longrightarrow A_i\longrightarrow E_i\longrightarrow E_i/A_i
        \longrightarrow0.
\]
This short exact sequence induces long exact sequences in any cohomology
theory to which the corresponding local systems or supported sheaves are
transported.

The important case below is
\[
        E_i\cong\mathbb Z/4,
        \qquad
        A_i=2E_i\cong\mathbb Z/2.
\]
Then
\[
        0
        \longrightarrow
        2E_i
        \longrightarrow
        E_i
        \longrightarrow
        E_i/2E_i
        \longrightarrow
        0
\]
is
\[
        0
        \longrightarrow
        \mathbb Z/2
        \longrightarrow
        \mathbb Z/4
        \longrightarrow
        \mathbb Z/2
        \longrightarrow
        0.
\]
In a codimension-two stratum \(\Sigma\), this becomes a short exact sequence
of local systems
\[
        0
        \longrightarrow
        (2E)_\Sigma
        \longrightarrow
        E_\Sigma
        \longrightarrow
        (E/2E)_\Sigma
        \longrightarrow
        0.
\]
After applying \(i_*\) and shifting by \([2]\), it gives a distinguished
triangle
\[
        i_*(2E)_\Sigma[2]
        \longrightarrow
        i_*E_\Sigma[2]
        \longrightarrow
        i_*(E/2E)_\Sigma[2]
        \overset{+1}{\longrightarrow}.
\]
For the Coble boundary,
\[
        E_\Sigma\cong(\mathbb Z/4)_\Sigma,
        \qquad
        (2E)_\Sigma\cong(\mathbb Z/2)_\Sigma,
        \qquad
        (E/2E)_\Sigma\cong(\mathbb Z/2)_\Sigma.
\]

\begin{definition}[Shadow transport]
Let \(A_i\subset E_i\) be a specified subgroup of a local package, and suppose
that \(E_i\) is equipped with a degree-\(k\) support map
\[
        \eta_i^{(k)}:E_i\to H^k_{D_i}(N_i,\mathbb Z)_{\mathrm{tors}}.
\]
The degree-\(k\) shadow transport map for \(A_i\) is the restriction
\[
        \alpha_{X,A}^{(k)}:
        \bigoplus_iA_i
        \longrightarrow
        H^k(\widetilde X,\mathbb Z)_{\mathrm{tors}}
\]
obtained by composing
\[
        \bigoplus_iA_i
        \hookrightarrow
        \bigoplus_iE_i
        \xrightarrow{\alpha_X^{(k)}}
        H^k(\widetilde X,\mathbb Z)_{\mathrm{tors}}.
\]
\end{definition}

\begin{remark}
For the Enriques/Coble boundary comparison, \(A_i=2E_i\).  Thus the
Benoist--Ottem \(2\)-torsion direction is transported through
\[
        2E_i\subset E_i,
\]
not through the full local package \(E_i\cong\mathbb Z/4\).  This is the
local-to-global version of the statement that the Bockstein of the local
index-two cover selects the order-two subgroup of the Coble \(E\)-package.
\end{remark}

\subsection{Surface-type case}

We next spell out what can be asserted in the surface case without hiding
degree issues.  Let \(X\) be a projective surface with isolated normal
singularities, let
\[
        \pi:\widetilde X\longrightarrow X
\]
be a resolution, and let \(D_i=\pi^{-1}(p_i)\).  Let \(N_i\) be a resolution
neighborhood of \(D_i\) with boundary \(L_i=\partial N_i\).

For each singular point \(p_i\), the previous paper identifies the local
perverse obstruction group with
\[
        E_i
        \cong
        H^2(L_i,\mathbb Z)_{\mathrm{tors}}
        \cong
        \Lambda_i^\vee/\Lambda_i,
\]
where \(\Lambda_i\) is the local exceptional lattice over \(p_i\)
\cite{RahmanIntegralPerverseObstructions}.  This is a boundary-link and
lattice statement.  It does not, by itself, say that \(E_i\) is a subgroup
of \(H^k_{D_i}(N_i,\mathbb Z)\) for a fixed global degree \(k\).  The precise
degree is determined by the local pair sequence.

We record the relevant local exact sequence.

\begin{lemma}[Surface resolution-neighborhood pair sequence]
Let \((X,p_i)\) be a normal surface singularity, let \(N_i\) be a resolution
neighborhood of the exceptional divisor \(D_i\), and let
\(L_i=\partial N_i\).  Then the long exact homology sequence of the pair
\((N_i,L_i)\) contains
\[
        H_2(N_i,\mathbb Z)
        \longrightarrow
        H_2(N_i,L_i;\mathbb Z)
        \longrightarrow
        H_1(L_i,\mathbb Z)
        \longrightarrow
        H_1(N_i,\mathbb Z).
\]
Under Poincaré--Lefschetz duality, the first map is represented by the
local intersection form on the exceptional curves.  In the rational surface
singularity case, this gives
\[
        H_1(L_i,\mathbb Z)_{\mathrm{tors}}
        \cong
        \Lambda_i^\vee/\Lambda_i.
\]
Equivalently, by Poincaré duality and the universal coefficient theorem for
the closed oriented three-manifold \(L_i\),
\[
        H^2(L_i,\mathbb Z)_{\mathrm{tors}}
        \cong
        \Lambda_i^\vee/\Lambda_i.
\]
\end{lemma}

\begin{proof}
The displayed exact sequence is the ordinary long exact sequence in homology
for the pair \((N_i,L_i)\).  Since \(N_i\) is a compact oriented real
four-manifold with boundary \(L_i\), Poincaré--Lefschetz duality gives
\[
        H_2(N_i,L_i;\mathbb Z)
        \cong
        H^2(N_i,\mathbb Z).
\]
The neighborhood \(N_i\) deformation retracts onto the exceptional divisor,
so \(H_2(N_i,\mathbb Z)\) is generated by the exceptional curve classes.
The map
\[
        H_2(N_i,\mathbb Z)
        \longrightarrow
        H_2(N_i,L_i;\mathbb Z)
\]
is the map that sends an exceptional curve class to the functional given by
intersection with that curve.  With respect to the basis of irreducible
exceptional curves, it is represented by the intersection matrix.

For rational surface singularities, the relevant homology groups of the
resolution neighborhood have no additional torsion contribution interfering
with this cokernel calculation.  Therefore the torsion in
\(H_1(L_i,\mathbb Z)\) is identified with the cokernel of the intersection
form, namely
\[
        \Lambda_i^\vee/\Lambda_i.
\]
This is the local resolution-neighborhood computation used in
\cite[§4]{RahmanIntegralPerverseObstructions}; see also the classical
linking-pairing analysis in \cite{GoreskySiegel83}.  Finally, for a closed
oriented three-manifold \(L_i\), Poincaré duality and the universal
coefficient theorem identify the torsion in \(H_1(L_i,\mathbb Z)\) with the
torsion in \(H^2(L_i,\mathbb Z)\).  Hence
\[
        H^2(L_i,\mathbb Z)_{\mathrm{tors}}
        \cong
        \Lambda_i^\vee/\Lambda_i.
\]
\end{proof}

Thus the rigorous surface-type conclusion is the following.

\begin{proposition}[Surface-type local transport data]
Let \(X\) be a projective surface with isolated rational normal surface
singularities, and let
\[
        \pi:\widetilde X\longrightarrow X
\]
be a resolution.  For each singular point \(p_i\), the local obstruction
group
\[
        E_i=E_X(p_i)
\]
is canonically identified with the local discriminant quotient
\[
        \Lambda_i^\vee/\Lambda_i
\]
and with the torsion in the link cohomology
\[
        H^2(L_i,\mathbb Z)_{\mathrm{tors}}.
\]
These identifications are obtained from the pair sequence of
\((N_i,L_i)\), not from a degree-independent global support map.
Consequently, any global surface-type transport map must specify how the
local pair-sequence class is inserted into a chosen global support or lattice
sequence.
\end{proposition}

\begin{proof}
The identifications
\[
        E_i
        \cong
        H^2(L_i,\mathbb Z)_{\mathrm{tors}}
        \cong
        \Lambda_i^\vee/\Lambda_i
\]
are exactly the local realization theorem of
\cite{RahmanIntegralPerverseObstructions}, applied at the singular point
\(p_i\).  The preceding lemma explains the pair-sequence mechanism that
produces the lattice quotient from the resolution neighborhood.  Since this
construction uses the boundary link and the homology sequence of
\((N_i,L_i)\), it does not automatically define a class in a fixed global
support group \(H^k_D(\widetilde X,\mathbb Z)\) without an additional
degree-specific map.  Therefore the final sentence follows from the
definition of \(\alpha_X^{(k)}\), which requires a specified map
\[
        E_i\to H^k_{D_i}(N_i,\mathbb Z)_{\mathrm{tors}}.
\]
\end{proof}

\begin{remark}
This proposition is intentionally more cautious than the slogan
``local surface packages map to global torsion.''  The local surface package
is completely understood as
\[
        E_i\cong H^2(L_i,\mathbb Z)_{\mathrm{tors}}
        \cong
        \Lambda_i^\vee/\Lambda_i.
\]
To turn this into a global map, one must choose the global exact sequence in
which the local class is to be transported.  In later surface-type examples,
this is often most naturally done through global divisor or lattice
relations rather than through the degree-three Brauer channel relevant to
threefolds.
\end{remark}

\subsection{Codimension-two transverse surface packages}

The surface calculation also enters higher-dimensional spaces through
codimension-two strata.  Let \(Y\) be a threefold and let
\[
        i:\Sigma\hookrightarrow Y
\]
be a smooth codimension-two stratum.  Suppose the transverse singularity
along \(\Sigma\) is locally analytically trivial and is a normal surface
singularity with local package \(E\).  Then the transverse packages assemble
into a finite local system
\[
        E_\Sigma
\]
on \(\Sigma\).

In this situation MacPherson--Vilonen or Beilinson gluing identifies the
closed-stratum discrepancy term as
\[
        i_*E_\Sigma[2].
\]
This is the codimension-two version of the point-supported surface triangle
\[
        {}^pIC_X\mathbb Z
        \longrightarrow
        {}^p_+IC_X\mathbb Z
        \longrightarrow
        i_*E[1]
        \longrightarrow.
\]
The additional shift comes from the smooth one-dimensional stratum.

For the Coble boundary singularity \(\frac14(1,1)\),
\[
        E_\Sigma\cong(\mathbb Z/4)_\Sigma.
\]
The Benoist--Ottem \(2\)-torsion mechanism selects the sublocal system
\[
        2E_\Sigma\cong(\mathbb Z/2)_\Sigma.
\]
Thus there are two related transport problems:
\[
        \text{full Coble transport: } E_\Sigma\cong(\mathbb Z/4)_\Sigma,
\]
and
\[
        \text{BO shadow transport: }2E_\Sigma\cong(\mathbb Z/2)_\Sigma.
\]

\begin{remark}
This is the place where the subobject/quotient formalism matters.  The
short exact sequence
\[
        0\to2E_\Sigma\to E_\Sigma\to E_\Sigma/2E_\Sigma\to0
\]
induces a distinguished triangle after applying \(i_*\) and shifting by
\([2]\):
\[
        i_*(2E_\Sigma)[2]
        \longrightarrow
        i_*E_\Sigma[2]
        \longrightarrow
        i_*(E_\Sigma/2E_\Sigma)[2]
        \overset{+1}{\longrightarrow}.
\]
For the Coble boundary, this is
\[
        i_*(\mathbb Z/2)_\Sigma[2]
        \longrightarrow
        i_*(\mathbb Z/4)_\Sigma[2]
        \longrightarrow
        i_*(\mathbb Z/2)_\Sigma[2]
        \overset{+1}{\longrightarrow}.
\]
This is the exact-sequence form of the order-two shadow inside the full
Coble \(E\)-package.
\end{remark}

\subsection{Threefold degree-four Hodge case}

We now isolate the threefold case relevant to the Brauer and unramified
cohomology comparisons.  Let \(X\) be a projective threefold with isolated
singularities and let
\[
        \pi:\widetilde X\longrightarrow X
\]
be a resolution.  The smooth variety \(\widetilde X\) is a projective
threefold.  The exponential sequence on \(\widetilde X\),
\[
        0\longrightarrow
        \mathbb Z
        \longrightarrow
        \mathcal O_{\widetilde X}
        \longrightarrow
        \mathcal O_{\widetilde X}^*
        \longrightarrow
        0,
\]
gives a cohomology segment
\[
        H^2(\widetilde X,\mathcal O_{\widetilde X})
        \longrightarrow
        H^2(\widetilde X,\mathcal O_{\widetilde X}^*)
        \longrightarrow
        H^3(\widetilde X,\mathbb Z)
        \longrightarrow
        H^3(\widetilde X,\mathcal O_{\widetilde X}).
\]
Under appropriate hypotheses, such as
\(H^2(\widetilde X,\mathcal O_{\widetilde X})=0\), the torsion part of the
cohomological Brauer group is compared with
\[
        H^3(\widetilde X,\mathbb Z)_{\mathrm{tors}}.
\]
This is the degree in which the Brauer group enters the
Colliot-Thélène--Voisin framework for integral Hodge obstructions
\cite{CTVoisin12}.

Therefore, for the threefold degree-four Hodge/Brauer bridge, the relevant
local-to-global map must land in degree three.

\begin{definition}[Threefold discriminant transport]
Let \(X\) be a projective threefold with isolated singularities and
resolution \(\pi:\widetilde X\to X\).  A collection of local torsion
packages \(E_i\) is said to be relevant to the threefold Brauer bridge if,
for each \(i\), there is a specified map
\[
        \eta_i^{(3)}:
        E_i
        \longrightarrow
        H^3_{D_i}(N_i,\mathbb Z)_{\mathrm{tors}}.
\]
In that case the degree-three discriminant transport map is
\[
        \alpha_X^{(3)}:
        \bigoplus_iE_i
        \longrightarrow
        H^3_D(\widetilde X,\mathbb Z)_{\mathrm{tors}}
        \longrightarrow
        H^3(\widetilde X,\mathbb Z)_{\mathrm{tors}}.
\]
\end{definition}

\begin{remark}
This definition does not assert that every local singularity package produces
a degree-three class.  It states the condition required for a local package
to enter the Brauer/unramified channel.  In particular, the ordinary double
point in complex dimension three must be computed separately.  One cannot
import the surface \(A_1\) identity
\[
        E_{A_1}^{\mathrm{surface}}\cong \mathbb Z/2
\]
without verifying that the corresponding threefold local torsion appears in
degree \(3\).  The link of a threefold ordinary double point has real
dimension \(5\), and the local stalk complex is shifted by \([3]\), not
\([2]\).  The required computation is carried out separately in the
threefold-node section.
\end{remark}

\subsection{Kernel and global relations}

The preceding constructions define the map
\[
        \alpha_X^{(k)}
\]
once the local degree-\(k\) realizations have been specified.  The kernel of
this map measures which local torsion classes are killed by global topology.
At the formal level, the kernel is simply the inverse image of zero under the
forget-support map.  We now state this precisely.

\begin{proposition}[Formal kernel description]
Assume that local maps
\[
        \eta_i^{(k)}:E_i\to H^k_{D_i}(N_i,\mathbb Z)_{\mathrm{tors}}
\]
have been fixed, and let
\[
        \eta^{(k)}:=\bigoplus_i\eta_i^{(k)}.
\]
Let
\[
        \rho^{(k)}:
        H^k_D(\widetilde X,\mathbb Z)_{\mathrm{tors}}
        \longrightarrow
        H^k(\widetilde X,\mathbb Z)_{\mathrm{tors}}
\]
be the forget-support map.  Then
\[
        \ker(\alpha_X^{(k)})
        =
        (\eta^{(k)})^{-1}
        \bigl(\ker(\rho^{(k)})\bigr).
\]
Equivalently, a local torsion class dies globally if and only if its
supported cohomology class lies in the kernel of the forget-support map.
\end{proposition}

\begin{proof}
By definition,
\[
        \alpha_X^{(k)}
        =
        \rho^{(k)}\circ \eta^{(k)},
\]
where we have identified
\[
        \bigoplus_iH^k_{D_i}(N_i,\mathbb Z)_{\mathrm{tors}}
        \cong
        H^k_D(\widetilde X,\mathbb Z)_{\mathrm{tors}}
\]
by excision.  Therefore
\[
        x\in \ker(\alpha_X^{(k)})
\]
if and only if
\[
        \rho^{(k)}(\eta^{(k)}(x))=0,
\]
which holds if and only if
\[
        \eta^{(k)}(x)\in \ker(\rho^{(k)}).
\]
This is exactly the displayed formula.
\end{proof}

The following is the lattice-theoretic form of the same idea.  It is stated
as a conditional theorem because it applies only when the local packages are
realized as discriminant quotients of local exceptional lattices and the
global map factors through a global exceptional lattice.

\begin{theorem}[Kernel criterion through exceptional-cycle relations]
Assume the following data are given.

\begin{enumerate}[label=\textup{(\roman*)}]
\item For each singular point \(p_i\), the local package is identified with a
discriminant quotient
\[
        E_i\cong \Lambda_i^\vee/\Lambda_i.
\]

\item The local exceptional lattices map into a global lattice of exceptional
cycles
\[
        \Lambda_{\mathrm{loc}}:=\bigoplus_i\Lambda_i
        \longrightarrow
        \Lambda_{\mathrm{glob}},
\]
and the degree-\(k\) transport map \(\alpha_X^{(k)}\) factors through the
corresponding global discriminant or cohomological quotient.

\item The map to global cohomology is induced by the image of exceptional
cycles in \(H^k(\widetilde X,\mathbb Z)\).
\end{enumerate}

Then the kernel of
\[
        \alpha_X^{(k)}:\bigoplus_iE_i\to H^k(\widetilde X,\mathbb Z)_{\mathrm{tors}}
\]
is determined by the integral relations among the images of the local
exceptional-cycle lattices in global cohomology.  More precisely, a class in
\(\bigoplus_iE_i\) maps to zero if and only if it can be represented by local
dual-lattice data whose induced supported cohomology class lies in the
relation subgroup killed by the global forget-support map.
\end{theorem}

\begin{proof}
Under hypothesis (i), a local class in \(E_i\) is represented by an element
of \(\Lambda_i^\vee\), modulo the subgroup \(\Lambda_i\).  Taking the direct
sum over all \(i\), a local class in \(\bigoplus_iE_i\) is represented by an
element of
\[
        \bigoplus_i\Lambda_i^\vee
\]
modulo
\[
        \bigoplus_i\Lambda_i.
\]
By hypothesis (ii), the transport map factors through the global lattice or
cohomological quotient determined by the images of these local exceptional
cycles.  Therefore two representatives have the same global image precisely
when their difference lies in the subgroup generated by local lattice
elements and by global relations among the images of the local exceptional
cycles.

By hypothesis (iii), the final map to \(H^k(\widetilde X,\mathbb Z)\) is the
map induced by forgetting supports.  Thus a transported class maps to zero
globally precisely when its supported representative belongs to the kernel
of the forget-support map.  By the formal kernel description proved above,
this kernel is exactly the subgroup of local classes whose supported
cohomology representatives are killed by global relations.  Hence the kernel
is determined by the integral relations among the images of the local
exceptional lattices in global cohomology.
\end{proof}

\begin{remark}
This theorem is not a substitute for a computation.  It identifies what must
be computed: the global relation subgroup among the images of the local
exceptional cycles.  In explicit examples this subgroup may be described by
defect groups, vanishing-cycle relations, divisor-class relations, lattice
embeddings, or Bockstein-selected subpackages such as
\[
        2E\subset E.
\]
\end{remark}

\subsection{Relation with Clemens defect}

We conclude this section by explaining the expected connection with nodal
hypersurfaces.  This part is included to orient the later examples; it is not
used as a theorem before the threefold node computation is completed.

For nodal hypersurfaces, the defect measures the failure of the nodes to
impose independent conditions, equivalently the presence of global relations
among the local data associated to the nodes.  This circle of ideas goes back
to Clemens's work on double solids \cite{Clemens83}.  For explicit treatments
of the defect of nodal hypersurfaces, see Cynk \cite{Cynk01}; for more recent
families with defect, see Kloosterman \cite{Kloosterman22}.

In the language of the present paper, the expected relation is that the
kernel of the degree-three local-to-global map for a nodal threefold should
be governed by the same relation space that appears in the defect.

\begin{conjecture}[Defect interpretation for nodal hypersurfaces]
Let \(X\) be a nodal projective threefold, and suppose that each node carries
a local torsion package \(E_i\) which is realized in degree-three support
cohomology, so that
\[
        \alpha_X^{(3)}:
        \bigoplus_iE_i
        \longrightarrow
        H^3(\widetilde X,\mathbb Z)_{\mathrm{tors}}
\]
is defined.  Then \(\ker(\alpha_X^{(3)})\) is controlled by the same global
relation space that appears in the Clemens defect.  Equivalently, local node
packages survive globally precisely to the extent that their associated
vanishing or exceptional cycles are independent modulo the defect relations.
\end{conjecture}

\begin{remark}
The conjecture is deliberately conditional.  Before it can be turned into a
theorem for nodal threefolds, one must first compute the correct local
degree-three package for the ordinary double point in complex dimension
three.  The surface identity
\[
        E_{A_1}^{\mathrm{surface}}\cong\mathbb Z/2
\]
does not by itself define the threefold degree-three package.  Once the
threefold node package is computed, the conjecture predicts that the kernel
of the resulting map is governed by the same global dependence relations that
enter the classical defect theory of nodal hypersurfaces
\cite{Clemens83,Cynk01,Kloosterman22}.
\end{remark}

\section{Brauer comparison}

The purpose of this section is to explain the second bridge in the
local-to-global construction.  The first bridge produces, when defined, a
class in
\[
        H^3(\widetilde X,\mathbb Z)_{\mathrm{tors}}.
\]
For smooth projective complex varieties, the exponential sequence compares
this group with the cohomological Brauer group under standard vanishing
hypotheses.  This is the reason the threefold degree-three target appears in
the previous section.

We separate two points.  First, the analytic exponential sequence gives a
comparison between torsion in \(H^2(Y,\mathcal O_Y^*)\) and
\(H^3(Y,\mathbb Z)\).  Second, for smooth projective complex varieties, the
torsion in analytic \(H^2(Y,\mathcal O_Y^*)\) agrees with the usual
cohomological Brauer group defined in the étale topology.  The latter
comparison is standard in the theory of the Brauer group; see Grothendieck's
Brauer papers and standard treatments of étale cohomology
\cite{GrothendieckBrauerIII,MilneEtale}.

A minor notational point will be useful later.  The source of the
degree-three map may be a full local package \(E_i\), or it may be a
specified subpackage \(A_i\subset E_i\).  In the Coble boundary example,
\[
        E_i\cong\mathbb Z/4,
        \qquad
        A_i=2E_i\cong\mathbb Z/2.
\]
The Brauer comparison is insensitive to whether the source is a full local
package or such a subpackage: it applies to whatever torsion subgroup has
already survived into \(H^3(\widetilde X,\mathbb Z)_{\mathrm{tors}}\).

\subsection{The exponential sequence}

Let \(Y\) be a smooth complex variety, and let \(Y^{\mathrm{an}}\) denote the
associated complex analytic space.  On \(Y^{\mathrm{an}}\), the exponential
map gives an exact sequence of sheaves
\[
        0
        \longrightarrow
        \mathbb Z
        \longrightarrow
        \mathcal O_Y
        \xrightarrow{\exp(2\pi i\,\cdot)}
        \mathcal O_Y^*
        \longrightarrow
        0.
\]
Here \(\mathbb Z\) is the locally constant sheaf of integral-valued
functions, \(\mathcal O_Y\) is the sheaf of holomorphic functions, and
\(\mathcal O_Y^*\) is the sheaf of nowhere vanishing holomorphic functions.
Exactness is local: a holomorphic function maps to \(1\) under the
exponential map if and only if it is locally integer-valued, and every
nowhere vanishing holomorphic function has a local holomorphic logarithm.

Taking cohomology gives the long exact sequence
\[
\cdots
\longrightarrow
H^2(Y,\mathbb Z)
\longrightarrow
H^2(Y,\mathcal O_Y)
\longrightarrow
H^2(Y,\mathcal O_Y^*)
\longrightarrow
H^3(Y,\mathbb Z)
\longrightarrow
H^3(Y,\mathcal O_Y)
\longrightarrow
\cdots .
\]
The part of this sequence used below is
\[
        H^2(Y,\mathcal O_Y)
        \longrightarrow
        H^2(Y,\mathcal O_Y^*)
        \xrightarrow{\delta_Y}
        H^3(Y,\mathbb Z)
        \longrightarrow
        H^3(Y,\mathcal O_Y).
\]
The connecting morphism
\[
        \delta_Y:
        H^2(Y,\mathcal O_Y^*)
        \longrightarrow
        H^3(Y,\mathbb Z)
\]
is the analytic cohomological Brauer-to-integral-cohomology comparison map.

\begin{definition}[Analytic cohomological Brauer group]
For a smooth complex variety \(Y\), define the analytic cohomological Brauer
group by
\[
        \operatorname{Br}_{\mathrm{an}}(Y)
        :=
        H^2(Y,\mathcal O_Y^*)_{\mathrm{tors}}.
\]
For a smooth projective complex variety, this torsion group agrees with the
usual cohomological Brauer group
\[
        \operatorname{Br}(Y)
        :=
        H^2_{\mathrm{\acute et}}(Y,\mathbb G_m)_{\mathrm{tors}}
\]
under the standard analytic/étale comparison theorem for the Brauer group
\cite{GrothendieckBrauerIII,MilneEtale}.
\end{definition}

\subsection{Precise torsion comparison}

We now record the exact torsion comparison needed later.  The point is that
the vanishing condition is \(H^2(Y,\mathcal O_Y)=0\), not
\(H^3(Y,\mathcal O_Y)=0\).

\begin{proposition}[Brauer--torsion comparison]
Let \(Y\) be a smooth projective complex variety.  The exponential sequence
gives an exact segment
\[
        H^2(Y,\mathcal O_Y)
        \longrightarrow
        H^2(Y,\mathcal O_Y^*)
        \xrightarrow{\delta_Y}
        H^3(Y,\mathbb Z)
        \longrightarrow
        H^3(Y,\mathcal O_Y).
\]
Assume that
\[
        H^2(Y,\mathcal O_Y)=0.
\]
Then the connecting morphism \(\delta_Y\) induces an isomorphism
\[
        H^2(Y,\mathcal O_Y^*)_{\mathrm{tors}}
        \cong
        H^3(Y,\mathbb Z)_{\mathrm{tors}}.
\]
Equivalently, using the standard analytic/étale comparison for smooth
projective complex varieties,
\[
        \operatorname{Br}(Y)
        \cong
        H^3(Y,\mathbb Z)_{\mathrm{tors}}.
\]
\end{proposition}

\begin{proof}
The exponential sequence gives the exact segment
\[
        H^2(Y,\mathcal O_Y)
        \longrightarrow
        H^2(Y,\mathcal O_Y^*)
        \xrightarrow{\delta_Y}
        H^3(Y,\mathbb Z)
        \longrightarrow
        H^3(Y,\mathcal O_Y).
\]
Since \(H^2(Y,\mathcal O_Y)=0\), exactness implies that
\[
        \delta_Y:
        H^2(Y,\mathcal O_Y^*)
        \longrightarrow
        H^3(Y,\mathbb Z)
\]
is injective.

We first show that \(\delta_Y\) sends torsion classes into torsion classes.
Let \(b\in H^2(Y,\mathcal O_Y^*)\) be torsion.  Then there exists
\(m>0\) such that \(mb=0\).  Since \(\delta_Y\) is a homomorphism,
\[
        m\,\delta_Y(b)=\delta_Y(mb)=\delta_Y(0)=0.
\]
Hence \(\delta_Y(b)\) is torsion.  Thus \(\delta_Y\) restricts to an
injective homomorphism
\[
        H^2(Y,\mathcal O_Y^*)_{\mathrm{tors}}
        \hookrightarrow
        H^3(Y,\mathbb Z)_{\mathrm{tors}}.
\]

It remains to prove surjectivity onto \(H^3(Y,\mathbb Z)_{\mathrm{tors}}\).
Let
\[
        a\in H^3(Y,\mathbb Z)_{\mathrm{tors}}.
\]
Then there exists \(m>0\) such that \(ma=0\).  Let
\[
        \rho:
        H^3(Y,\mathbb Z)
        \longrightarrow
        H^3(Y,\mathcal O_Y)
\]
be the next map in the long exact sequence.  Since \(\rho\) is a group
homomorphism,
\[
        m\,\rho(a)=\rho(ma)=\rho(0)=0.
\]
The group \(H^3(Y,\mathcal O_Y)\) is a complex vector space.  In particular,
as an abelian group it is torsion-free.  Therefore \(\rho(a)=0\).  By
exactness,
\[
        \ker(\rho)=\operatorname{im}(\delta_Y).
\]
Hence there exists
\[
        b\in H^2(Y,\mathcal O_Y^*)
\]
such that
\[
        \delta_Y(b)=a.
\]
Since \(a\) is torsion and \(\delta_Y\) is injective, the element \(b\) is
also torsion.  Indeed, if \(ma=0\), then
\[
        0=ma=m\delta_Y(b)=\delta_Y(mb),
\]
and injectivity of \(\delta_Y\) gives \(mb=0\).  Thus
\[
        b\in H^2(Y,\mathcal O_Y^*)_{\mathrm{tors}}.
\]
This proves that \(\delta_Y\) induces a surjection
\[
        H^2(Y,\mathcal O_Y^*)_{\mathrm{tors}}
        \twoheadrightarrow
        H^3(Y,\mathbb Z)_{\mathrm{tors}}.
\]
Together with injectivity, this proves
\[
        H^2(Y,\mathcal O_Y^*)_{\mathrm{tors}}
        \cong
        H^3(Y,\mathbb Z)_{\mathrm{tors}}.
\]

Finally, for smooth projective complex varieties, the torsion subgroup of
analytic \(H^2(Y,\mathcal O_Y^*)\) agrees with the cohomological Brauer group
defined by étale cohomology:
\[
        H^2(Y,\mathcal O_Y^*)_{\mathrm{tors}}
        \cong
        H^2_{\mathrm{\acute et}}(Y,\mathbb G_m)_{\mathrm{tors}}
        =
        \operatorname{Br}(Y).
\]
This is the standard analytic/étale comparison for the Brauer group
\cite{GrothendieckBrauerIII,MilneEtale}.  Combining this comparison with the
isomorphism just proved gives
\[
        \operatorname{Br}(Y)
        \cong
        H^3(Y,\mathbb Z)_{\mathrm{tors}}.
\]
\end{proof}

\begin{remark}
The hypothesis
\[
        H^2(Y,\mathcal O_Y)=0
\]
is the relevant vanishing for the clean torsion comparison.  The group
\(H^3(Y,\mathcal O_Y)\) does not need to vanish for the torsion comparison
above.  It is enough that \(H^3(Y,\mathcal O_Y)\) is a complex vector space,
hence torsion-free as an abelian group.  Thus torsion classes in
\(H^3(Y,\mathbb Z)\) automatically map to zero in \(H^3(Y,\mathcal O_Y)\).
\end{remark}

\begin{remark}
For a Calabi--Yau threefold \(Y\), one has \(h^{0,2}(Y)=0\), equivalently
\(H^2(Y,\mathcal O_Y)=0\).  The Brauer comparison above can therefore be
available even though \(h^{0,3}(Y)=1\), i.e. even though
\(H^3(Y,\mathcal O_Y)\) is nonzero.  This is why the Brauer comparison is
controlled by \(H^2(Y,\mathcal O_Y)\), not by \(H^3(Y,\mathcal O_Y)\).
\end{remark}

\subsection{The Brauer image of local discriminants}

We now apply the preceding comparison to the local-to-global maps of the
previous section.  Let \(X\) be a projective threefold with isolated
singularities, and let
\[
        \pi:\widetilde X\longrightarrow X
\]
be a resolution.  Assume that the local packages \(E_i\) have been realized
in degree-three support cohomology, so that the map
\[
        \alpha_X^{(3)}:
        \bigoplus_iE_i
        \longrightarrow
        H^3(\widetilde X,\mathbb Z)_{\mathrm{tors}}
\]
is defined.

Assume further that
\[
        H^2(\widetilde X,\mathcal O_{\widetilde X})=0.
\]
Then the Brauer--torsion comparison gives an isomorphism
\[
        \delta_{\widetilde X}:
        \operatorname{Br}(\widetilde X)
        \xrightarrow{\sim}
        H^3(\widetilde X,\mathbb Z)_{\mathrm{tors}}.
\]
We write
\[
        \beta_{\widetilde X}
        :=
        \delta_{\widetilde X}^{-1}
\]
for the inverse isomorphism.

\begin{definition}[Brauer realization of local discriminant data]
Under the hypotheses above, the Brauer realization of the transported local
discriminant data is the subgroup
\[
        \operatorname{Br}_{\mathrm{loc}}(\widetilde X)
        :=
        \beta_{\widetilde X}
        \bigl(\operatorname{im}\alpha_X^{(3)}\bigr)
        \subseteq
        \operatorname{Br}(\widetilde X).
\]
Equivalently,
\[
        \operatorname{Br}_{\mathrm{loc}}(\widetilde X)
\]
is the unique subgroup of \(\operatorname{Br}(\widetilde X)\) whose image
under
\[
        \delta_{\widetilde X}:
        \operatorname{Br}(\widetilde X)
        \xrightarrow{\sim}
        H^3(\widetilde X,\mathbb Z)_{\mathrm{tors}}
\]
is exactly
\[
        \operatorname{im}\alpha_X^{(3)}.
\]
\end{definition}

\begin{remark}
The notation
\[
        \operatorname{Br}_{\mathrm{loc}}(\widetilde X)
\]
does not mean that every Brauer class is local in origin.  It denotes only
the subgroup obtained from local discriminant data by first applying the
degree-three support-to-global map and then using the Brauer--torsion
comparison.  The subgroup may be zero even when the local groups \(E_i\) are
nonzero, because the map \(\alpha_X^{(3)}\) may kill the local classes.
\end{remark}

\subsection{Brauer images of subpackages and Bockstein shadows}

The same definition applies to a specified subpackage
\[
        A_i\subset E_i.
\]
This is important for the Coble boundary comparison, where
\[
        E_i\cong\mathbb Z/4,
        \qquad
        A_i=2E_i\cong\mathbb Z/2.
\]
If the restricted degree-three transport map
\[
        \alpha_{X,A}^{(3)}:
        \bigoplus_i A_i
        \longrightarrow
        H^3(\widetilde X,\mathbb Z)_{\mathrm{tors}}
\]
is defined, then under the same Brauer--torsion comparison we define
\[
        \operatorname{Br}_{A,\mathrm{loc}}(\widetilde X)
        :=
        \beta_{\widetilde X}
        \bigl(\operatorname{im}\alpha_{X,A}^{(3)}\bigr)
        \subseteq
        \operatorname{Br}(\widetilde X).
\]

\begin{definition}[Brauer realization of a shadow subpackage]
Let \(A_i\subset E_i\) be specified local subpackages and suppose that their
degree-three shadow transport map
\[
        \alpha_{X,A}^{(3)}
\]
is defined.  Under the hypotheses of the Brauer--torsion comparison, the
Brauer realization of the shadow data is
\[
        \operatorname{Br}_{A,\mathrm{loc}}(\widetilde X)
        :=
        \beta_{\widetilde X}
        \bigl(\operatorname{im}\alpha_{X,A}^{(3)}\bigr)
        \subseteq
        \operatorname{Br}(\widetilde X).
\]
\end{definition}

\begin{remark}
For a codimension-two Coble boundary stratum, the full local package is
\[
        E_\Sigma\cong(\mathbb Z/4)_\Sigma,
\]
while the Benoist--Ottem order-two shadow is
\[
        2E_\Sigma\cong(\mathbb Z/2)_\Sigma.
\]
The Brauer image relevant to the Benoist--Ottem comparison is therefore not
the image of the full \(\mathbb Z/4\)-package unless the full package is
being transported.  It is the image of the transported order-two shadow
\[
        2E_\Sigma\subset E_\Sigma.
\]
\end{remark}

\begin{proposition}[Compatibility of full and shadow Brauer images]
Suppose \(A_i\subset E_i\) are specified subpackages and the degree-three
transport maps for both \(A_i\) and \(E_i\) are defined, with
\[
        \alpha_{X,A}^{(3)}
        =
        \alpha_X^{(3)}\big|_{\oplus_i A_i}.
\]
Assume also that
\[
        H^2(\widetilde X,\mathcal O_{\widetilde X})=0.
\]
Then
\[
        \operatorname{Br}_{A,\mathrm{loc}}(\widetilde X)
        \subseteq
        \operatorname{Br}_{\mathrm{loc}}(\widetilde X).
\]
\end{proposition}

\begin{proof}
Since
\[
        \alpha_{X,A}^{(3)}
        =
        \alpha_X^{(3)}\big|_{\oplus_i A_i},
\]
we have
\[
        \operatorname{im}\alpha_{X,A}^{(3)}
        \subseteq
        \operatorname{im}\alpha_X^{(3)}.
\]
The map
\[
        \beta_{\widetilde X}
        =
        \delta_{\widetilde X}^{-1}
\]
is an isomorphism, hence an injective homomorphism.  Therefore it sends
subgroups to subgroups, and
\[
        \beta_{\widetilde X}
        \bigl(\operatorname{im}\alpha_{X,A}^{(3)}\bigr)
        \subseteq
        \beta_{\widetilde X}
        \bigl(\operatorname{im}\alpha_X^{(3)}\bigr).
\]
This is precisely
\[
        \operatorname{Br}_{A,\mathrm{loc}}(\widetilde X)
        \subseteq
        \operatorname{Br}_{\mathrm{loc}}(\widetilde X).
\]
\end{proof}

\subsection{Kernel after Brauer comparison}

We now compare the kernel of the local-to-global map with the kernel of the
map to the Brauer group.  Under the clean Brauer--torsion comparison, there
is no additional kernel introduced at the Brauer stage.

\begin{proposition}[Kernel after Brauer comparison]
Let \(X\) be a projective threefold with isolated singularities and
resolution
\[
        \pi:\widetilde X\to X.
\]
Assume that the local packages \(E_i\) define the degree-three map
\[
        \alpha_X^{(3)}:
        \bigoplus_iE_i
        \longrightarrow
        H^3(\widetilde X,\mathbb Z)_{\mathrm{tors}},
\]
and assume that
\[
        H^2(\widetilde X,\mathcal O_{\widetilde X})=0.
\]
Let
\[
        \beta_{\widetilde X}
        =
        \delta_{\widetilde X}^{-1}
\]
be the inverse of the Brauer--torsion isomorphism
\[
        \delta_{\widetilde X}:
        \operatorname{Br}(\widetilde X)
        \xrightarrow{\sim}
        H^3(\widetilde X,\mathbb Z)_{\mathrm{tors}}.
\]
Then the composite
\[
        \bigoplus_iE_i
        \xrightarrow{\alpha_X^{(3)}}
        H^3(\widetilde X,\mathbb Z)_{\mathrm{tors}}
        \xrightarrow{\beta_{\widetilde X}}
        \operatorname{Br}(\widetilde X)
\]
has kernel equal to
\[
        \ker(\alpha_X^{(3)}).
\]
\end{proposition}

\begin{proof}
The map
\[
        \beta_{\widetilde X}:
        H^3(\widetilde X,\mathbb Z)_{\mathrm{tors}}
        \longrightarrow
        \operatorname{Br}(\widetilde X)
\]
is an isomorphism by the Brauer--torsion comparison.  In particular, it is
injective.  Let
\[
        e\in \bigoplus_iE_i.
\]
Then \(e\) lies in the kernel of the composite
\[
        \beta_{\widetilde X}\circ\alpha_X^{(3)}
\]
if and only if
\[
        \beta_{\widetilde X}\bigl(\alpha_X^{(3)}(e)\bigr)=0.
\]
Since \(\beta_{\widetilde X}\) is injective, this holds if and only if
\[
        \alpha_X^{(3)}(e)=0.
\]
Therefore
\[
        \ker(\beta_{\widetilde X}\circ\alpha_X^{(3)})
        =
        \ker(\alpha_X^{(3)}).
\]
\end{proof}

\begin{proposition}[Kernel after Brauer comparison for a shadow subpackage]
Let \(A_i\subset E_i\) be specified subpackages and suppose that the
degree-three shadow transport map
\[
        \alpha_{X,A}^{(3)}:
        \bigoplus_iA_i
        \longrightarrow
        H^3(\widetilde X,\mathbb Z)_{\mathrm{tors}}
\]
is defined.  Assume that
\[
        H^2(\widetilde X,\mathcal O_{\widetilde X})=0.
\]
Then the composite
\[
        \bigoplus_iA_i
        \xrightarrow{\alpha_{X,A}^{(3)}}
        H^3(\widetilde X,\mathbb Z)_{\mathrm{tors}}
        \xrightarrow{\beta_{\widetilde X}}
        \operatorname{Br}(\widetilde X)
\]
has kernel equal to
\[
        \ker(\alpha_{X,A}^{(3)}).
\]
\end{proposition}

\begin{proof}
The proof is identical to the proof for the full local package.  The
Brauer--torsion comparison identifies
\[
        H^3(\widetilde X,\mathbb Z)_{\mathrm{tors}}
\]
with \(\operatorname{Br}(\widetilde X)\) by an isomorphism.  Therefore the
map
\[
        \beta_{\widetilde X}
\]
is injective, and composing with it introduces no additional kernel.  Hence
the kernel of
\[
        \beta_{\widetilde X}\circ\alpha_{X,A}^{(3)}
\]
is exactly
\[
        \ker(\alpha_{X,A}^{(3)}).
\]
\end{proof}

\begin{remark}
If the hypotheses of the Brauer--torsion comparison are weakened, the kernel
statement must be modified.  Suppose, for example, that one has only a
partially defined comparison map from a subgroup of
\(H^3(\widetilde X,\mathbb Z)_{\mathrm{tors}}\) to
\(\operatorname{Br}(\widetilde X)\), or that the comparison is not known to
be an isomorphism on the relevant subgroup.  Then the kernel of the map from
\(\bigoplus_iE_i\), or from a subpackage \(\bigoplus_iA_i\), to the Brauer
group is the preimage under the corresponding transport map of the kernel of
the chosen comparison map.  Thus in the clean case there is no additional
kernel at the Brauer stage, while in a partial comparison the kernel may be
enlarged by the failure of the Brauer comparison to identify the relevant
torsion subgroup.
\end{remark}

\begin{remark}
The Brauer stage does not decide whether a transported local class is
unramified.  That question is controlled by the residue maps in the
Bloch--Ogus formalism.  The next section therefore studies when a Brauer or
degree-three torsion class survives to
\[
        H^3_{\mathrm{nr}}(\widetilde X,\mathbb Q/\mathbb Z).
\]
\end{remark}

\section{Unramified cohomology and residue survival}

This section treats the residue station of the torsion trajectory.  We state
the definition of unramified cohomology in the form used by Bloch--Ogus and
Colliot-Thélène--Voisin, then explain precisely what it means for a class
coming from local discriminant data, or from a specified subpackage such as
\(2E\subset E\), to survive to the unramified group.

There is one important degree issue.  A class in
\(H^3(Y,\mathbb Z)_{\mathrm{tors}}\) does not automatically become a class in
\(H^3(Y,\mathbb Q/\mathbb Z)\) by the coefficient sequence
\[
        0\to \mathbb Z\to \mathbb Q\to \mathbb Q/\mathbb Z\to 0.
\]
Rather, that sequence gives a Bockstein map
\[
        H^2(Y,\mathbb Q/\mathbb Z)\to H^3(Y,\mathbb Z),
\]
and the torsion subgroup of \(H^3(Y,\mathbb Z)\) lies in the image of this
Bockstein.  Thus the residue station below is formulated for a specified
\(\mathbb Q/\mathbb Z\)-valued degree-three class associated to the
transported local data.  In applications, producing such a class is part of
the comparison problem.

This warning is especially relevant for the Coble boundary comparison.  The
local boundary package for the \(\frac14(1,1)\) singularity is
\[
        E\cong\mathbb Z/4,
\]
while the Benoist--Ottem \(2\)-torsion mechanism selects
\[
        2E\cong\mathbb Z/2.
\]
A residue calculation must be applied only after this order-two shadow, or
the full local package if appropriate, has been transported and converted
into a specified \(\mathbb Q/\mathbb Z\)-valued degree-three candidate.

\subsection{Bloch--Ogus residues}

Let \(Y\) be a smooth irreducible complex variety, and let
\(K=\mathbb C(Y)\) be its function field.  For a codimension-one point
\(D\in Y^{(1)}\), write \(k(D)\) for the residue field of \(D\).  The point
\(D\) defines a discrete valuation on \(K\).  In the Bloch--Ogus/Gersten
formalism, this valuation gives a residue homomorphism
\[
        \partial_D:
        H^r(K,\mathbb Q/\mathbb Z)
        \longrightarrow
        H^{r-1}(k(D),\mathbb Q/\mathbb Z)
\]
for the relevant étale cohomology groups with torsion coefficients.  More
precisely, in the usual étale convention one includes Tate twists:
\[
        \partial_D:
        H^r(K,\mathbb Q/\mathbb Z(m))
        \longrightarrow
        H^{r-1}(k(D),\mathbb Q/\mathbb Z(m-1)).
\]
Over the complex numbers, and for the purposes of this paper, we suppress
Tate twists from the notation.  The residue formalism and its relation to
the coniveau spectral sequence are due to Bloch--Ogus
\cite{BlochOgus74}; see also the use of this formalism in the integral Hodge
context by Colliot-Thélène--Voisin \cite{CTVoisin12}.

The case used below is \(r=3\).  Thus for each prime divisor
\(D\in Y^{(1)}\) one has a residue map
\[
        \partial_D:
        H^3(K,\mathbb Q/\mathbb Z)
        \longrightarrow
        H^2(k(D),\mathbb Q/\mathbb Z).
\]
A class in \(H^3(K,\mathbb Q/\mathbb Z)\) is called unramified if its residue
at every codimension-one point vanishes.

\begin{definition}[Third unramified cohomology]
Let \(Y\) be a smooth irreducible complex variety with function field
\(K=\mathbb C(Y)\).  The third unramified cohomology group of \(Y\) with
\(\mathbb Q/\mathbb Z\)-coefficients is
\[
        H^3_{\mathrm{nr}}(Y,\mathbb Q/\mathbb Z)
        :=
        \bigcap_{D\in Y^{(1)}}\ker(\partial_D),
\]
where
\[
        \partial_D:
        H^3(K,\mathbb Q/\mathbb Z)
        \longrightarrow
        H^2(k(D),\mathbb Q/\mathbb Z)
\]
is the Bloch--Ogus residue map.  Equivalently,
\[
        H^3_{\mathrm{nr}}(Y,\mathbb Q/\mathbb Z)
        =
        \ker\left(
        H^3(K,\mathbb Q/\mathbb Z)
        \longrightarrow
        \bigoplus_{D\in Y^{(1)}}H^2(k(D),\mathbb Q/\mathbb Z)
        \right).
\]
\end{definition}

\begin{remark}
The same group can also be described as the group of global sections of the
Zariski sheaf associated to the presheaf
\(U\mapsto H^3(U,\mathbb Q/\mathbb Z)\):
\[
        H^3_{\mathrm{nr}}(Y,\mathbb Q/\mathbb Z)
        \cong
        H^0\bigl(Y,\mathcal H^3(\mathbb Q/\mathbb Z)\bigr).
\]
This is the sheaf-theoretic formulation supplied by the Bloch--Ogus
coniveau machinery \cite{BlochOgus74}.
\end{remark}

\subsection{The coefficient sequence and the degree issue}

We now explain why the residue station must be formulated carefully.  The
short exact sequence of abelian groups
\[
        0
        \longrightarrow
        \mathbb Z
        \longrightarrow
        \mathbb Q
        \longrightarrow
        \mathbb Q/\mathbb Z
        \longrightarrow
        0
\]
gives a long exact sequence in cohomology:
\[
\cdots
\longrightarrow
H^2(Y,\mathbb Q)
\longrightarrow
H^2(Y,\mathbb Q/\mathbb Z)
\xrightarrow{\delta}
H^3(Y,\mathbb Z)
\longrightarrow
H^3(Y,\mathbb Q)
\longrightarrow
\cdots .
\]
Let
\[
        a\in H^3(Y,\mathbb Z)_{\mathrm{tors}}.
\]
Then the image of \(a\) in \(H^3(Y,\mathbb Q)\) is zero, because tensoring
with \(\mathbb Q\) kills torsion.  By exactness, \(a\) lies in the image of
the Bockstein map
\[
        \delta:
        H^2(Y,\mathbb Q/\mathbb Z)
        \longrightarrow
        H^3(Y,\mathbb Z).
\]
Thus torsion in \(H^3(Y,\mathbb Z)\) is naturally detected by
\(H^2(Y,\mathbb Q/\mathbb Z)\) through the Bockstein, not by an automatic map
to \(H^3(Y,\mathbb Q/\mathbb Z)\).

\begin{remark}
This is why the residue criterion below is stated for a specified class
\[
        \xi\in H^3(\mathbb C(Y),\mathbb Q/\mathbb Z)
\]
associated to the transported local data.  If the local trajectory only
produces a class in \(H^3(Y,\mathbb Z)_{\mathrm{tors}}\), then an additional
comparison is needed before one can test survival in
\(H^3_{\mathrm{nr}}(Y,\mathbb Q/\mathbb Z)\).  The Brauer comparison in the
previous section identifies, under suitable hypotheses,
\(\operatorname{Br}(Y)\) with \(H^3(Y,\mathbb Z)_{\mathrm{tors}}\), but the
degree-three unramified group \(H^3_{\mathrm{nr}}(Y,\mathbb Q/\mathbb Z)\) is
a separate residue-theoretic target.
\end{remark}

\begin{remark}[Bockstein shadows and residue candidates]
For the Coble boundary package, the local group is
\[
        E\cong\mathbb Z/4.
\]
The local index-two cover \(L(2,1)\to L(4,1)\) gives a mod-\(2\) class whose
Bockstein lands in
\[
        2E\cong\mathbb Z/2.
\]
Thus a Benoist--Ottem-type residue candidate is associated naturally to the
order-two shadow \(2E\), not automatically to the full \(\mathbb Z/4\)-group.
If one wants to test the full \(E\)-package against residues, one must first
construct a separate \(\mathbb Q/\mathbb Z\)-valued degree-three candidate
for the full package.
\end{remark}

\subsection{Residue survival criterion}

We now formulate the survival criterion in a way that separates the formal
residue test from the additional problem of producing a degree-three
\(\mathbb Q/\mathbb Z\)-valued class.

Let \(X\) be a projective variety with isolated singularities and let
\[
        \pi:\widetilde X\longrightarrow X
\]
be a resolution.  Assume \(\widetilde X\) is smooth and irreducible.  Suppose
that local discriminant packages \(E_i\), or specified subpackages
\(A_i\subset E_i\), have been transported to global torsion data by a
degree-indexed map
\[
        \alpha_X^{(k)}:
        \bigoplus_iE_i
        \longrightarrow
        H^k(\widetilde X,\mathbb Z)_{\mathrm{tors}},
\]
or by a restricted shadow map
\[
        \alpha_{X,A}^{(k)}:
        \bigoplus_iA_i
        \longrightarrow
        H^k(\widetilde X,\mathbb Z)_{\mathrm{tors}}.
\]
To test unramified survival in degree three, we need an associated
\(\mathbb Q/\mathbb Z\)-valued degree-three class.  Thus we make the
following definition.

\begin{definition}[Residue candidate attached to local or shadow data]
Let
\[
        e\in \bigoplus_iE_i
\]
or
\[
        a\in \bigoplus_iA_i
\]
be local discriminant data or specified shadow data.  A degree-three residue
candidate attached to \(e\), respectively to \(a\), is a class
\[
        \xi(e)\in H^3(\mathbb C(\widetilde X),\mathbb Q/\mathbb Z),
        \qquad
        \xi(a)\in H^3(\mathbb C(\widetilde X),\mathbb Q/\mathbb Z),
\]
constructed from the trajectory by a specified comparison map, such as a
correspondence, a specialization map, a product construction, a Bockstein
construction, or another geometric operation.  The class is called
unramified if
\[
        \partial_D(\xi)=0
\]
for every codimension-one point \(D\in\widetilde X^{(1)}\).
\end{definition}

\begin{theorem}[Residue survival criterion]
Let \(e\in\bigoplus_iE_i\) be a local discriminant class, or let
\(a\in\bigoplus_iA_i\) be a class in a specified subpackage.  Suppose that a
degree-three residue candidate
\[
        \xi\in H^3(\mathbb C(\widetilde X),\mathbb Q/\mathbb Z)
\]
has been associated to this data.  Then \(\xi\) defines a class in
\[
        H^3_{\mathrm{nr}}(\widetilde X,\mathbb Q/\mathbb Z)
\]
if and only if
\[
        \partial_D(\xi)=0
\]
for every codimension-one point \(D\in\widetilde X^{(1)}\).
\end{theorem}

\begin{proof}
By definition,
\[
        H^3_{\mathrm{nr}}(\widetilde X,\mathbb Q/\mathbb Z)
        =
        \ker\left(
        H^3(\mathbb C(\widetilde X),\mathbb Q/\mathbb Z)
        \longrightarrow
        \bigoplus_{D\in\widetilde X^{(1)}}
        H^2(k(D),\mathbb Q/\mathbb Z)
        \right).
\]
The map in the display is the direct sum of the residue maps
\[
        \partial_D:
        H^3(\mathbb C(\widetilde X),\mathbb Q/\mathbb Z)
        \longrightarrow
        H^2(k(D),\mathbb Q/\mathbb Z).
\]
Therefore a class
\[
        \xi\in H^3(\mathbb C(\widetilde X),\mathbb Q/\mathbb Z)
\]
lies in \(H^3_{\mathrm{nr}}(\widetilde X,\mathbb Q/\mathbb Z)\) precisely
when its image under the direct sum of residues is zero.  This occurs
precisely when each component residue is zero:
\[
        \partial_D(\xi)=0
        \qquad
        \text{for every }D\in\widetilde X^{(1)}.
\]
This proves the equivalence.
\end{proof}

\begin{remark}
The theorem is formally simple because it is a survival criterion: once the
candidate class \(\xi\) exists, survival is exactly the vanishing of all
codimension-one residues.  The difficult geometric work is usually the
construction of \(\xi\) and the computation of its residues.
\end{remark}

\subsection{Consequences for integral Hodge obstruction}

We now explain how this residue station relates to the integral Hodge
problem.  Let \(Y\) be a smooth projective complex variety.  The integral
Hodge conjecture asks whether every integral Hodge class of degree \(2p\) is
the cohomology class of an algebraic cycle of codimension \(p\).  It is false
in general; classical and modern counterexamples include
\cite{AtiyahHirzebruch62,Totaro97,SouleVoisin05,BenoistOttem20}.

For degree-four classes, Colliot-Thélène and Voisin relate the failure of
the integral Hodge conjecture to unramified cohomology with
\(\mathbb Q/\mathbb Z\)-coefficients \cite{CTVoisin12}.  Thus
\[
        H^3_{\mathrm{nr}}(Y,\mathbb Q/\mathbb Z)
\]
is one of the natural obstruction groups appearing in the study of
degree-four integral Hodge classes.

\begin{corollary}
Let \(e\in\bigoplus_iE_i\) be a local discriminant class, or let
\(a\in\bigoplus_iA_i\) be a class in a specified subpackage.  Suppose that
the trajectory produces a residue candidate
\[
        \xi\in H^3(\mathbb C(\widetilde X),\mathbb Q/\mathbb Z).
\]
If \(\xi\) is nonzero and all of its codimension-one residues vanish, then
the corresponding trajectory contributes a nonzero class to
\[
        H^3_{\mathrm{nr}}(\widetilde X,\mathbb Q/\mathbb Z).
\]
Consequently, it contributes to the same unramified cohomology group that
appears in the Colliot-Thélène--Voisin obstruction theory for degree-four
integral Hodge classes.
\end{corollary}

\begin{proof}
If all residues of \(\xi\) vanish, then by the residue survival criterion,
\[
        \xi\in H^3_{\mathrm{nr}}(\widetilde X,\mathbb Q/\mathbb Z).
\]
If \(\xi\) is nonzero in the function-field cohomology group and remains
nonzero in the subgroup defined by the residue kernels, then it gives a
nonzero element of
\[
        H^3_{\mathrm{nr}}(\widetilde X,\mathbb Q/\mathbb Z).
\]
The relevance of this group to degree-four integral Hodge obstructions is
the theorem of Colliot-Thélène--Voisin \cite{CTVoisin12}.  Therefore the
local or shadow trajectory contributes to that obstruction-theoretic target.
\end{proof}

\begin{remark}
This corollary does not say that every local discriminant class gives an
integral Hodge obstruction.  There are several possible failure points.  The
local class may die under the support-to-global map.  It may fail to produce
a degree-three \(\mathbb Q/\mathbb Z\) residue candidate.  It may produce
such a candidate, but one of the residues may be nonzero.  Or the resulting
unramified class may be zero.  The residue station records only the last of
these tests.
\end{remark}

\begin{remark}
If the residue maps kill all classes arising from local discriminant data in
a family, the result is still meaningful.  It shows that any integral Hodge
obstruction present in that family is not detected by the local discriminant
packages studied here.  In the terminology of this paper, the torsion
trajectory dies at the residue station.
\end{remark}

\subsection{Benoist--Ottem, Coble shadows, and what residue survival can test}

The Benoist--Ottem examples require special care in this framework.  On the
smooth fiber \(S\times C\), where \(S\) is an Enriques surface, there is no
singular local package \(E_i\).  The torsion is global and smooth.  The
boundary comparison developed in this paper uses the Enriques/Coble
degeneration instead.  The relevant boundary singularity is
\[
        \frac14(1,1),
\]
with local package
\[
        E\cong\mathbb Z/4.
\]
The Benoist--Ottem \(2\)-torsion mechanism selects
\[
        2E\cong\mathbb Z/2
\]
through the Bockstein of the local index-two cover
\[
        L(2,1)\to L(4,1).
\]

Thus the residue station, if applied to the Benoist--Ottem comparison, should
be applied to a candidate associated to the order-two shadow \(2E\), not
automatically to the full local group \(E\).  If such a candidate is produced
and all of its residues vanish, then the order-two Coble shadow contributes
to
\[
        H^3_{\mathrm{nr}}(\widetilde X,\mathbb Q/\mathbb Z).
\]
If no such candidate exists, or if a residue is nonzero, then the
Benoist--Ottem/Coble trajectory dies before reaching the unramified station.

\begin{remark}
This distinction prevents a common overstatement.  The corrected comparison
does not say that the full Coble group
\[
        E\cong\mathbb Z/4
\]
is itself a Benoist--Ottem obstruction.  It says that the
Benoist--Ottem \(2\)-torsion mechanism naturally identifies with the
Bockstein-selected subgroup
\[
        2E\cong\mathbb Z/2
\]
inside the Coble \(E\)-package.  The residue station can test the resulting
order-two trajectory once a degree-three \(\mathbb Q/\mathbb Z\)-valued
candidate has been constructed.
\end{remark}

\subsection{What this station can and cannot prove}

The unramified station should be read as a test, not as an automatic
consequence of local torsion.  The formal implication is:
\[
        \text{residue candidate } \xi
        \text{ with all residues zero}
        \quad
        \Longrightarrow
        \quad
        \xi\in H^3_{\mathrm{nr}}.
\]
The formalism does not assert that every local class, or every Bockstein
shadow, produces such a candidate.  Producing \(\xi\) requires additional
geometry.

This distinction is especially important for the examples in the appendices.

\begin{enumerate}[label=\textup{(\arabic*)}]
\item For isolated surface singularities, the local group \(E\) is well-defined
and computable, but the isolated germ alone does not produce a degree-three
unramified class.

\item For threefold ordinary double points, the local link is torsion-free,
so there is no finite local torsion package to test.

\item For the Benoist--Ottem benchmark, the smooth fiber has no singular
points and hence no direct local packages \(E_i\).  The proposed comparison is
degenerational: one compares global Enriques \(2\)-torsion with the
order-two Bockstein shadow
\[
        2E\subset E
\]
inside the Coble boundary package
\[
        E\cong\mathbb Z/4
\]
for the \(\frac14(1,1)\) singularity.
\end{enumerate}

Thus the residue station is the final obstruction test in the trajectory,
but it is not the first step.  Local birth, support transport, global
survival, construction of a degree-three residue candidate, and, in the
Coble/Benoist--Ottem setting, selection of the correct order-two shadow must
all occur before the Bloch--Ogus residue criterion can be applied.

\section{Patterns extracted from the torsion trajectories}

The appendices were constructed as a computational laboratory.  Each example
was run through the same stations: local birth, discriminant form, local
realizations, support and pair sequences, global transport status, Brauer and
residue status, and rationalization.  The purpose of this section is to
record the patterns that emerge from those computations.

The main point is that the examples separate several phenomena that are easy
to conflate: \textit{local finite torsion}, \textit{free vanishing-cycle
data}, \textit{global smooth torsion}, \textit{Brauer/unramified obstruction},
and \textit{Bockstein-selected shadows inside larger local torsion packages}.
The examples show that these are not the same mechanism.  The local
discriminant package \(E\) detects one very specific phenomenon: finite
integral torsion born from the non-unimodularity of local intersection or
linking data.  It does not detect every singularity, every node, or every
global torsion class.  Moreover, a global torsion mechanism may see only a
subgroup or filtered piece of a local \(E\)-package, as happens in the
Coble boundary comparison with Benoist--Ottem.

\subsection{Local finite torsion is controlled by discriminants}

For the normal surface singularities computed above, the same pattern occurs
in every case:
\[
        E
        \cong
        H^2(L,\mathbb Z)_{\mathrm{tors}}
        \cong
        \Lambda^\vee/\Lambda.
\]
In the hypersurface cases, the same group is also recovered from the
monodromy station:
\[
        E
        \cong
        \operatorname{coker}(T-\mathrm{id})_{\mathrm{tors}}.
\]
These identifications are the local realization theorem of
\cite{RahmanIntegralPerverseObstructions}, and the appendices illustrate the
theorem in explicit examples.

The computed rows include
\[
        A_1:
        \quad
        E\cong\mathbb Z/2\mathbb Z,
\]
\[
        A_k:
        \quad
        E\cong\mathbb Z/(k+1)\mathbb Z,
\]
\[
        D_4:
        \quad
        E\cong
        \mathbb Z/2\mathbb Z\oplus\mathbb Z/2\mathbb Z,
\]
\[
        E_8:
        \quad
        E=0,
\]
\[
        x^2+y^3+z^{11}=0:
        \quad
        E\cong\mathbb Z/5\mathbb Z,
\]
and the Coble boundary singularity
\[
        \frac14(1,1):
        \quad
        E\cong\mathbb Z/4\mathbb Z.
\]
In each nonzero local surface case, the order of \(E\) is the absolute value
of the determinant of the local intersection matrix:
\[
        |E|=|\det(M)|.
\]
Thus the first pattern is:

\[
\boxed{
\text{local finite torsion is born from non-unimodularity of the local
intersection lattice.}
}
\]

Equivalently, the local torsion package is present precisely when the
exceptional lattice has nontrivial discriminant group
\[
        \Lambda^\vee/\Lambda.
\]
This explains why \(E_8\), despite being a nontrivial singularity with a
nontrivial exceptional configuration, produces no local finite torsion:
the \(E_8\) lattice is unimodular.  It also explains the Coble boundary
package: the local exceptional lattice is \([-4]\), hence
\[
        E\cong\operatorname{coker}([-4])\cong\mathbb Z/4.
\]

\subsection{\texorpdfstring{\(E_8\)}{E8} is the null-control}

The \(E_8\) row is important because it prevents a false interpretation of
the invariant.  The local group \(E\) does not measure the complexity of the
singularity in a naive sense.  The \(E_8\) singularity has eight exceptional
curves and a rich ADE configuration, but
\[
        \Lambda_{E_8}^\vee/\Lambda_{E_8}=0.
\]
Therefore
\[
        E=0.
\]
The link has no torsion in the relevant cohomology group, the pair-sequence
quotient vanishes, and the monodromy cokernel has trivial torsion.  Thus all
six local torsion stations vanish.

The pattern is:

\[
\boxed{
\begin{minipage}{0.82\textwidth}
\centering
A singularity can be geometrically nontrivial while contributing no local
finite discriminant torsion.
\end{minipage}
}
\]

This observation is useful later because it shows that the torsion trajectory
is not simply a singularity detector.  It detects a specific integral defect:
the failure of the local lattice to be unimodular.

\subsection{The group is not enough: the discriminant form matters}

The \(D_4\) row shows that the underlying finite abelian group is not the full
local invariant.  For \(D_4\), the group is
\[
        E\cong
        \mathbb Z/2\mathbb Z\oplus\mathbb Z/2\mathbb Z.
\]
However, the local package also carries a discriminant form
\[
        q:E\times E\longrightarrow\mathbb Q/\mathbb Z.
\]
With the generators chosen in the \(D_4\) appendix, the form is represented
by
\[
        \begin{pmatrix}
        0 & -\frac12\\
        -\frac12 & 0
        \end{pmatrix}.
\]
Thus the \(D_4\) package is not merely the abstract group
\[
        (\mathbb Z/2\mathbb Z)^2.
\]
It is the finite pairing
\[
        \left((\mathbb Z/2\mathbb Z)^2,q_{D_4}\right).
\]

The Coble boundary example adds another reason to keep the full package.  For
\(\frac14(1,1)\),
\[
        E\cong\mathbb Z/4\mathbb Z
\]
and, in the geometric negative-definite convention,
\[
        q(\bar g,\bar g)=-\frac14.
\]
The Benoist--Ottem \(2\)-torsion mechanism sees the order-two subgroup
\[
        2E\cong\mathbb Z/2\mathbb Z,
\]
not the full group.  The order-two element \(2\bar g\) satisfies
\[
        q(2\bar g,2\bar g)=0\quad \bmod \mathbb Z.
\]
Thus the global \(2\)-torsion comparison selects an isotropic shadow inside a
larger local discriminant package.

The pattern is:

\[
\boxed{
\begin{minipage}{0.88\textwidth}
\centering
the natural local object is the discriminant package \((E,q)\), together with
any functorially selected subpackages such as \(2E\subset E\).
\end{minipage}
}
\]

This is the first reason to expect that global survival may be
form-sensitive.  A local-to-global theory that forgets \(q\), or forgets the
distinction between \(E\) and a subgroup such as \(2E\), discards information
that is visible both in the resolution lattice and in the linking pairing on
the boundary.

\subsection{Prime support is stable across the local stations}

The \(A_k\), non-ADE Brieskorn, and Coble quotient examples show that the
prime content of local torsion is visible at every local station.

For \(A_k\),
\[
        E\cong\mathbb Z/(k+1)\mathbb Z.
\]
Therefore the prime support of \(E\) is exactly the set of primes dividing
\(k+1\).  The same integer \(k+1\) appears as
\[
        |H^2(L,\mathbb Z)_{\mathrm{tors}}|,
        \qquad
        |\Lambda^\vee/\Lambda|,
        \qquad
        |\det(M)|,
        \qquad
        |\operatorname{coker}(T-\mathrm{id})_{\mathrm{tors}}|.
\]
For the non-ADE Brieskorn example
\[
        x^2+y^3+z^{11}=0,
\]
the same phenomenon occurs with the prime \(5\):
\[
        H^2(L,\mathbb Z)_{\mathrm{tors}}
        \cong
        \Lambda^\vee/\Lambda
        \cong
        \operatorname{coker}(T-\mathrm{id})_{\mathrm{tors}}
        \cong
        \mathbb Z/5\mathbb Z.
\]
For the Coble boundary singularity
\[
        \frac14(1,1),
\]
the same phenomenon occurs with \(2\)-primary torsion:
\[
        H^2(L(4,1),\mathbb Z)_{\mathrm{tors}}
        \cong
        \Lambda^\vee/\Lambda
        \cong
        \mathbb Z/4\mathbb Z.
\]
The local \(2\)-torsion mechanism relevant to Benoist--Ottem is not the full
\(\mathbb Z/4\), but the Bockstein image
\[
        2E\cong\mathbb Z/2.
\]

The pattern is:

\[
\boxed{
\text{the prime decomposition of local torsion is stable across the local
realizations.}
}
\]

This does not mean that global transport preserves every prime or every
primary exponent.  A prime or exponent may disappear globally if the local
class is killed by relations, if only a Bockstein-selected subgroup is
visible, or if the class is later killed by residue conditions.  But the
local stations do not create new primes.  They present the same finite
torsion through different structures.

\subsection{Surface \texorpdfstring{\(A_1\)}{A1} and the threefold node are
different phenomena}

One of the most important experimental outcomes is the distinction between
the \(A_1\) surface singularity and the ordinary double point in complex
dimension three.

For the surface \(A_1\) singularity,
\[
        L\cong\mathbb RP^3,
\]
and
\[
        H^2(L,\mathbb Z)_{\mathrm{tors}}
        \cong
        \mathbb Z/2\mathbb Z.
\]
Thus
\[
        E\cong\mathbb Z/2\mathbb Z.
\]
The same \(\mathbb Z/2\mathbb Z\) is recovered from the exceptional
intersection matrix
\[
        [-2],
\]
from the pair sequence, and from the monodromy cokernel.

For the ordinary double point in complex dimension three, the link is
\[
        L\cong S^2\times S^3.
\]
Hence
\[
        H^m(L,\mathbb Z)_{\mathrm{tors}}=0
        \qquad
        \text{for all }m.
\]
The Milnor fiber has the homotopy type of \(S^3\), and the monodromy
cokernel has a free part rather than finite torsion.  Thus the threefold
node has no local finite torsion package of the surface \(E\)-type.

The pattern is:

\[
\boxed{
\text{the surface }A_1\text{ torsion package does not automatically extend
to isolated threefold nodes.}
}
\]

This observation is crucial for nodal threefolds.  It prevents the incorrect
claim that an ordinary node on a threefold contributes a local
\(\mathbb Z/2\mathbb Z\) torsion class analogous to the surface \(A_1\)
group.  The threefold node contributes free vanishing-cycle and exceptional
curve data, not local finite discriminant torsion.

\subsection{Nodal threefold defect is not local finite torsion}

The nodal threefold and nodal quintic rows show that defect belongs to a
different part of the trajectory.  An ordinary double point in a threefold
has no local finite torsion:
\[
        E_i^{\mathrm{tors}}=0.
\]
Therefore
\[
        \bigoplus_i E_i^{\mathrm{tors}}=0.
\]
The direct finite-torsion transport map
\[
        \alpha_X^{(3)}:
        \bigoplus_i E_i^{\mathrm{tors}}
        \longrightarrow
        H^3(\widetilde X,\mathbb Z)_{\mathrm{tors}}
\]
has zero source.

Nevertheless, nodal threefolds have nontrivial topology.  Each node carries
a free vanishing \(3\)-sphere in a smoothing and, in a small resolution, an
exceptional curve.  These classes may satisfy global relations.  The defect
of a nodal hypersurface measures precisely a failure of expected independence
among such local conditions; see Clemens \cite{Clemens83}, Cynk
\cite{Cynk01}, and Kloosterman \cite{Kloosterman22}.

The pattern is:

\[
\boxed{
\begin{minipage}{0.82\textwidth}
\centering
Defect is a global relation phenomenon among free node data, not a local
finite torsion package.
\end{minipage}
}
\]

This distinction clarifies the role of nodal quintics in the paper.  They
are not examples of local \(E\)-torsion.  They are contrast examples showing
that global relations can be important even when local finite torsion is
absent.

\subsection{The natural higher-dimensional carrier of surface torsion is
codimension two}

The comparison between surface singularities, threefold nodes, and boundary
degeneration tests suggests a structural principle.

The local group \(E\) for a normal surface singularity is a codimension-two
phenomenon.  A point on a surface has complex codimension two.  The local
link is a real three-manifold, and the obstruction group is
\[
        H^2(L,\mathbb Z)_{\mathrm{tors}}.
\]
When one passes to higher-dimensional varieties, the natural way to carry the
same surface package is not necessarily through isolated higher-dimensional
nodes.  Instead, the natural carrier is a codimension-two stratum with
transverse normal surface singularity.

For example, if \(S_0\) is a surface with \(A_1\) singularities and \(C\) is
a smooth curve, then
\[
        X_0=S_0\times C
\]
is a threefold whose singular locus contains curves
\[
        \Sigma_j=\{p_j\}\times C.
\]
Each stratum \(\Sigma_j\) has transverse singularity \(A_1\).  The transverse
local package is
\[
        E_{A_1}\cong\mathbb Z/2\mathbb Z.
\]
Thus the local packages assemble not as isolated point contributions but as
finite local systems along the strata:
\[
        \mathcal E_j\cong(\mathbb Z/2\mathbb Z)_{\Sigma_j}.
\]

Similarly, in the corrected Enriques/Coble boundary comparison, the
transverse local singularity is
\[
        \frac14(1,1).
\]
The transverse package is
\[
        E_\Sigma\cong(\mathbb Z/4)_\Sigma.
\]
The Benoist--Ottem-visible part is the order-two sublocal system
\[
        2E_\Sigma\cong(\mathbb Z/2)_\Sigma.
\]

The pattern is:

\[
\boxed{
\begin{minipage}{0.82\textwidth}
\centering
The higher-dimensional continuation of the surface package \((E,q)\) is
transverse codimension-two torsion.
\end{minipage}
}
\]

This is one of the most important lessons of the examples.  It suggests that
the right higher-dimensional theory should not begin with isolated threefold
nodes, but with stratified varieties whose codimension-two strata have
surface singularities.

\subsection{Benoist--Ottem and the corrected Coble boundary comparison}

The Benoist--Ottem row provides a global smooth torsion benchmark.  For
\[
        Y=S\times C,
\]
where \(S\) is an Enriques surface and \(C\) is a very general curve,
Benoist--Ottem produce non-algebraic torsion classes in degree \(4\)
\cite{BenoistOttem20}.  The variety \(Y\) is smooth, so there are no local
singularity packages:
\[
        \bigoplus_i E_i=0.
\]
Thus the Benoist--Ottem class is not directly the image of a local
discriminant map on the smooth fiber.

The corrected boundary comparison uses the Enriques/Coble boundary.  In the
degree-two Enriques compactification, the discriminant boundary contains
rational Coble surfaces with a singularity of type
\[
        \frac14(1,1)
\]
\cite{AlexeevEngelGarzaSchafflerEnriquesDegree2}.  The link is
\[
        L(4,1),
\]
and the full local package is
\[
        E\cong H^2(L(4,1),\mathbb Z)_{\mathrm{tors}}
        \cong \mathbb Z/4\mathbb Z.
\]
The local index-two cover
\[
        A_1=\mathbb C^2/\mu_2
        \longrightarrow
        \mathbb C^2/\mu_4=\frac14(1,1)
\]
induces the link cover
\[
        L(2,1)\longrightarrow L(4,1).
\]
If
\[
        \eta\in H^1(L(4,1),\mathbb Z/2)
\]
is the corresponding cover class, the Bockstein sends \(\eta\) to the
order-two element
\[
        \beta(\eta)=2\bar g\in E\cong\mathbb Z/4\mathbb Z.
\]
Thus the Benoist--Ottem direction should be interpreted not as the full local
discriminant package, but as the Bockstein-selected order-two shadow inside
the \(\frac14(1,1)\) Coble boundary package:
\[
        2E\cong\mathbb Z/2
        \subset
        E\cong\mathbb Z/4.
\]

This corrects the naive \(A_1\)-boundary interpretation.  The \(A_1\) package
remains a clean local \(\mathbb Z/2\)-model, and it appears as the index-two
cover of the Coble boundary singularity, but the boundary singularity itself
is \(\frac14(1,1)\), not \(A_1\).

The resulting pattern is:

\[
\boxed{
\begin{minipage}{0.84\textwidth}
\centering
Benoist--Ottem \(2\)-torsion is naturally compared with the order-two
Bockstein shadow \(2E\subset E\) of the Coble \(\frac14(1,1)\) boundary
package, where \(E\cong\mathbb Z/4\).
\end{minipage}
}
\]

A positive local-to-global specialization theorem would show that the
Benoist--Ottem \(2\)-torsion direction is the smoothing shadow of this
order-two Coble subpackage.  Without such a theorem, the current result is a
precise mechanism comparison: the same Bockstein exact sequence selects the
order-two local shadow and the global \(2\)-torsion direction.

\subsection{Motivic lift of the Coble package and its shadow}

The Coble boundary comparison also refines the motivic formulation.  The full
local package is not the mod-\(2\) object; it is the mod-\(4\) object:
\[
        E_\Sigma\cong(\mathbb Z/4)_\Sigma.
\]
In an integral motivic sheaf formalism with compatible Betti realization, the
full Coble package is modeled by
\[
        \mathcal T^{\mathrm{mot}}_{4,\Sigma}
        :=
        i_*(\mathbf 1_\Sigma/4)[2],
\]
whose Betti realization is
\[
        i_*(\mathbb Z/4)_\Sigma[2].
\]
The Benoist--Ottem-visible piece is the motivic order-two shadow
\[
        \mathcal T^{\mathrm{mot}}_{2,\Sigma}
        :=
        i_*(\mathbf 1_\Sigma/2)[2],
\]
with Betti realization
\[
        i_*(\mathbb Z/2)_\Sigma[2].
\]
The inclusion
\[
        \mathbb Z/2\hookrightarrow\mathbb Z/4,
        \qquad
        1\mapsto2,
\]
is represented motivically by a morphism
\[
        \mathcal T^{\mathrm{mot}}_{2,\Sigma}
        \longrightarrow
        \mathcal T^{\mathrm{mot}}_{4,\Sigma}.
\]
Thus the motivic version of the observed pattern is:

\[
\boxed{
\text{BO sees }i_*(\mathbf 1_\Sigma/2)[2]\text{ inside }
i_*(\mathbf 1_\Sigma/4)[2].
}
\]

This motivic formulation is not a separate invariant.  It is the same
\(E\)-package lifted before Betti realization, together with its
Bockstein-selected order-two subobject.

\subsection{Death mechanisms in the trajectory}

The examples reveal several distinct ways a torsion trajectory can terminate.

\subsubsection*{No birth}

In the \(E_8\) case,
\[
        E=0
\]
because the exceptional lattice is unimodular.  No finite torsion is born.

\subsubsection*{Rational death}

In the \(A_1\), \(A_k\), \(D_4\), Brieskorn, and Coble examples, torsion is
born locally but vanishes after tensoring with \(\mathbb Q\):
\[
        E\otimes_{\mathbb Z}\mathbb Q=0.
\]
This is the ordinary rational death of finite torsion.  In the Coble case,
both
\[
        E\cong\mathbb Z/4
\]
and
\[
        2E\cong\mathbb Z/2
\]
die rationally.

\subsubsection*{Global killing}

A local class may be killed by the forget-support map
\[
        H^k_D(\widetilde X,\mathbb Z)_{\mathrm{tors}}
        \longrightarrow
        H^k(\widetilde X,\mathbb Z)_{\mathrm{tors}}.
\]
This is the kernel of the local-to-global map
\[
        \alpha_X^{(k)}.
\]
It records global relations among local support classes.

\subsubsection*{Shadow selection}

A global obstruction theory may see only a subpackage of a local \(E\)-group.
For the Coble boundary,
\[
        E\cong\mathbb Z/4
\]
is the full local discriminant package, but the Benoist--Ottem \(2\)-torsion
direction sees only
\[
        2E\cong\mathbb Z/2.
\]
This is not death of the full package; it is a selection mechanism.

\subsubsection*{Residue killing}

A class that survives to a Brauer or \(\mathbb Q/\mathbb Z\)-cohomology
station may still fail to be unramified.  The Bloch--Ogus residue maps may
kill or detect it.  Thus residue behavior is a second global survival test,
after the forget-support map.

These mechanisms are distinct.  The trajectory keeps them separate.

\subsection{Summary of observed patterns}

The computations support the following summary.  The \(A_1\) row remains the
clean model for a transverse \(\mathbb Z/2\)-discriminant package, while the
\(\frac14(1,1)\) Coble boundary row refines the Benoist--Ottem comparison:
the full local package is \(\mathbb Z/4\), and the Benoist--Ottem
\(2\)-torsion mechanism sees the order-two Bockstein shadow.

\begin{center}
\small
\setlength{\tabcolsep}{5pt}
\renewcommand{\arraystretch}{1.2}
\begin{tabularx}{\textwidth}{|>{\raggedright\arraybackslash}p{0.27\textwidth}|>{\raggedright\arraybackslash}X|}
\hline
\textbf{Observation} & \textbf{Meaning} \\
\hline
Non-unimodular lattice &
Birth of local finite torsion \(E=\Lambda^\vee/\Lambda\). \\
\hline
Unimodular lattice &
No local finite torsion, as in \(E_8\). \\
\hline
Nontrivial \(q\) &
The finite group alone is insufficient; track the discriminant package
\((E,q)\). \\
\hline
Threefold node &
Free vanishing data, not local finite discriminant torsion. \\
\hline
Nodal defect &
Global relations among free node classes, not a local finite torsion
kernel. \\
\hline
Codimension-two stratum &
Natural higher-dimensional carrier of transverse surface \(E\)-packages. \\
\hline
Transverse \(A_1\) package &
Clean \(\mathbb Z/2\)-model: the link is \(\mathbb RP^3\), and the local
discriminant group is \(E\cong\mathbb Z/2\). \\
\hline
Coble boundary \(\frac14(1,1)\) &
Enriques/Coble boundary model: the link is \(L(4,1)\), the full local
package is \(E\cong\mathbb Z/4\), and the Benoist--Ottem channel sees
\(2E\cong\mathbb Z/2\). \\
\hline
Bockstein-selected shadow &
A global \(2\)-torsion obstruction may detect a canonical order-two subgroup
or shadow of a larger local discriminant group. \\
\hline
Benoist--Ottem &
Global smooth \(2\)-torsion; the corrected boundary comparison is with the
order-two Bockstein shadow inside the \(\frac14(1,1)\) Coble package, not
with the full local package. \\
\hline
Motivic lift &
The full Coble package is modeled by \(i_*(\mathbf 1_\Sigma/4)[2]\), while
the Benoist--Ottem-visible shadow is modeled by
\(i_*(\mathbf 1_\Sigma/2)[2]\). \\
\hline
Rationalization &
All finite torsion dies over \(\mathbb Q\), but its trajectory records
integral information before rationalization kills it. \\
\hline
\end{tabularx}
\end{center}

The main conceptual conclusion is:

\begin{center}
\fbox{
\begin{minipage}{0.84\textwidth}
\centering
Finite discriminant torsion is a codimension-two local phenomenon; its
global relevance depends on form, support trajectory, global relations,
residue behavior, and on which Bockstein-selected shadow of the local package
is visible to the global obstruction theory.
\end{minipage}
}
\end{center}

\section{Benoist--Ottem torsion and the Coble boundary \(E\)-package}

The purpose of this section is to record the corrected local-to-global
comparison suggested by the Enriques/Coble boundary.  The natural boundary
singularity for the Benoist--Ottem comparison is not the \(A_1\) surface
singularity, but the cyclic quotient singularity of type \(\frac14(1,1)\).
For this singularity the local perverse obstruction group \(E\) from
\cite{RahmanIntegralPerverseObstructions} is
\[
        E\cong\mathbb Z/4.
\]
The Benoist--Ottem \(2\)-torsion mechanism selects the canonical order-two
subgroup
\[
        2E\cong\mathbb Z/2
        \subset
        E\cong\mathbb Z/4.
\]

The relevant Enriques boundary model is supplied by the degree-two Enriques
compactification of Alexeev--Engel--Garza--Schaffler.  In their description,
the discriminant divisor parametrizes quotients of nodal K3 surfaces by an
involution fixing a node; the resulting boundary surfaces are rational Coble
surfaces with a singularity of type \(\frac14(1,1)\)
\cite{AlexeevEngelGarzaSchafflerEnriquesDegree2}.  Thus the correct local
boundary package is the Coble \(\frac14(1,1)\) \(E\)-package.

\subsection{The local Coble \(E\)-package}

\begin{proposition}[The \(E\)-package of the Coble boundary singularity]
Let
\[
        X=\mathbb C^2/\mu_4,
        \qquad
        \zeta\cdot(x,y)=(\zeta x,\zeta y),
\]
be the cyclic quotient singularity of type \(\frac14(1,1)\).  Let \(L\) be
its link.  Then
\[
        L\cong L(4,1),
        \qquad
        H^2(L,\mathbb Z)_{\mathrm{tors}}\cong\mathbb Z/4.
\]
Consequently the local perverse obstruction group is
\[
        E_X(0)\cong\mathbb Z/4.
\]
Moreover, the minimal resolution has exceptional lattice
\[
        \Lambda=\mathbb Z\langle e\rangle,
        \qquad
        (e,e)=-4,
\]
and hence
\[
        \Lambda^\vee/\Lambda\cong\mathbb Z/4.
\]
Under the geometric negative-definite convention, if \(\bar g\) is the
class of \(e^\vee=-e/4\), then
\[
        q(\bar g,\bar g)=-\frac14\quad \bmod \mathbb Z.
\]
\end{proposition}

\begin{proof}
The action of \(\mu_4\) on \(\mathbb C^2\setminus\{0\}\) restricts to a
free action on the unit sphere \(S^3\).  Therefore the link is
\[
        L=S^3/\mu_4=L(4,1).
\]
The lens space \(L(4,1)\) has
\[
        H_1(L(4,1),\mathbb Z)\cong\mathbb Z/4,
        \qquad
        H_2(L(4,1),\mathbb Z)=0.
\]
By the universal coefficient theorem,
\[
        H^2(L(4,1),\mathbb Z)
        \cong
        \operatorname{Ext}^1(H_1(L(4,1),\mathbb Z),\mathbb Z)
        \cong
        \mathbb Z/4.
\]
Thus
\[
        H^2(L,\mathbb Z)_{\mathrm{tors}}\cong\mathbb Z/4.
\]
By the local realization theorem of
\cite{RahmanIntegralPerverseObstructions},
\[
        E_X(0)
        \cong
        H^2(L,\mathbb Z)_{\mathrm{tors}},
\]
and therefore
\[
        E_X(0)\cong\mathbb Z/4.
\]

The Hirzebruch--Jung continued fraction for \(\frac41\) is
\[
        \frac41=[4].
\]
Thus the minimal resolution consists of a single exceptional curve with
self-intersection \(-4\).  Hence
\[
        \Lambda=\mathbb Z\langle e\rangle,
        \qquad
        (e,e)=-4.
\]
The dual lattice is generated by
\[
        e^\vee=-\frac e4,
\]
because
\[
        (e^\vee,e)=1.
\]
Therefore
\[
        \Lambda^\vee/\Lambda\cong\mathbb Z/4.
\]
Finally,
\[
        q(\bar g,\bar g)
        =
        (e^\vee,e^\vee)
        =
        \left(-\frac e4,-\frac e4\right)
        =
        \frac{1}{16}(e,e)
        =
        -\frac14
        \quad \bmod \mathbb Z.
\]
This proves the proposition.
\end{proof}

\subsection{Mayer--Vietoris and the pair sequence}

The preceding computation can be recovered from the pair sequence of a
resolution neighborhood.  Let \(N\) be a sufficiently small resolution
neighborhood of the exceptional curve, and let \(L=\partial N\).  The pair
sequence identifies the boundary torsion with the cokernel of the
intersection form on the exceptional lattice.

\begin{proposition}[Pair-sequence computation of \(E\)]
For the \(\frac14(1,1)\) singularity, the exact sequence associated to the
pair \((N,L)\) yields
\[
        \Lambda
        \xrightarrow{[-4]}
        \Lambda^\vee
        \longrightarrow
        H^2(L,\mathbb Z)_{\mathrm{tors}}
        \longrightarrow
        0.
\]
Hence
\[
        H^2(L,\mathbb Z)_{\mathrm{tors}}
        \cong
        \operatorname{coker}([-4])
        \cong
        \mathbb Z/4.
\]
\end{proposition}

\begin{proof}
The neighborhood \(N\) deformation retracts onto the exceptional curve
\(C\), and the intersection form on \(H_2(N,\mathbb Z)\) is represented by
the matrix
\[
        [-4].
\]
The long exact sequence for the pair \((N,L)\), together with
Poincare--Lefschetz duality, identifies the quotient of the dual lattice by
the image of the intersection lattice with the torsion in the boundary
cohomology.  Thus one obtains
\[
        \Lambda
        \xrightarrow{[-4]}
        \Lambda^\vee
        \longrightarrow
        H^2(L,\mathbb Z)_{\mathrm{tors}}
        \longrightarrow
        0.
\]
Since \(\operatorname{coker}([-4])\cong\mathbb Z/4\), the result follows.
\end{proof}

\subsection{The local index-two cover and the Bockstein shadow}

The full local group is \(E\cong\mathbb Z/4\), but the Benoist--Ottem
mechanism is \(2\)-torsion.  The bridge is the Bockstein of the local
index-two cover.

\begin{proposition}[Local index-two cover of the Coble boundary singularity]
Let
\[
        X=\mathbb C^2/\mu_4,
        \qquad
        \zeta\cdot(x,y)=(\zeta x,\zeta y),
\]
be the cyclic quotient singularity of type \(\frac14(1,1)\).  Let
\(\mu_2\subset\mu_4\) be the subgroup generated by \(\zeta^2=-1\).  Then the
intermediate quotient
\[
        \mathbb C^2/\mu_2
        \longrightarrow
        \mathbb C^2/\mu_4
\]
is a degree-two cover.  The source is the \(A_1\) singularity
\[
        \mathbb C^2/\{\pm1\}.
\]
On links, this cover is
\[
        L(2,1)\longrightarrow L(4,1).
\]
\end{proposition}

\begin{proof}
The subgroup \(\mu_2=\{\pm1\}\subset\mu_4\) acts on \(\mathbb C^2\) by
\[
        (x,y)\longmapsto(-x,-y).
\]
Thus \(\mathbb C^2/\mu_2\) is the \(A_1\) surface singularity.  Since
\[
        \mu_4/\mu_2\cong\mathbb Z/2,
\]
there is a residual degree-two quotient
\[
        \mathbb C^2/\mu_2
        \longrightarrow
        \mathbb C^2/\mu_4.
\]
Restricting to the unit sphere \(S^3\subset\mathbb C^2\), the source and
target links are
\[
        S^3/\mu_2=L(2,1),
        \qquad
        S^3/\mu_4=L(4,1).
\]
Therefore the induced map on links is
\[
        L(2,1)\to L(4,1).
\]
\end{proof}

\begin{proposition}[Bockstein image of the local index-two cover]
Let
\[
        \eta\in H^1(L(4,1),\mathbb Z/2)
\]
be the class of the double cover
\[
        L(2,1)\to L(4,1).
\]
Let
\[
        \beta:
        H^1(L(4,1),\mathbb Z/2)
        \longrightarrow
        H^2(L(4,1),\mathbb Z)
\]
be the Bockstein associated to
\[
        0\to\mathbb Z\xrightarrow{2}\mathbb Z\to\mathbb Z/2\to0.
\]
Then
\[
        \beta(\eta)
\]
is the unique element of order \(2\) in
\[
        H^2(L(4,1),\mathbb Z)\cong\mathbb Z/4.
\]
Under the identification
\[
        E\cong H^2(L(4,1),\mathbb Z)_{\mathrm{tors}},
\]
this says
\[
        \beta(\eta)=2\bar g\in E\cong\mathbb Z/4.
\]
\end{proposition}

\begin{proof}
The relevant part of the long exact cohomology sequence associated to
\[
        0\to\mathbb Z\xrightarrow{2}\mathbb Z\to\mathbb Z/2\to0
\]
is
\[
        H^1(L,\mathbb Z)
        \longrightarrow
        H^1(L,\mathbb Z/2)
        \xrightarrow{\beta}
        H^2(L,\mathbb Z)
        \xrightarrow{2}
        H^2(L,\mathbb Z).
\]
For \(L=L(4,1)\),
\[
        H^1(L,\mathbb Z)=0,
        \qquad
        H^1(L,\mathbb Z/2)\cong\mathbb Z/2,
        \qquad
        H^2(L,\mathbb Z)\cong\mathbb Z/4.
\]
Thus \(\beta\) is injective, and exactness gives
\[
        \operatorname{im}(\beta)
        =
        \ker(2:\mathbb Z/4\to\mathbb Z/4).
\]
The kernel of multiplication by \(2\) on \(\mathbb Z/4\) is
\[
        \{0,2\}\cong\mathbb Z/2.
\]
Since \(\eta\) is the nonzero class in \(H^1(L,\mathbb Z/2)\), its
Bockstein is the nonzero element of this kernel.  Therefore
\[
        \beta(\eta)=2\bar g
        \in
        E\cong\mathbb Z/4.
\]
\end{proof}

\subsection{Perverse truncation and the \(E\)-triangle}

Let
\[
        j:U=X\setminus\{0\}\hookrightarrow X,
        \qquad
        i:\{0\}\hookrightarrow X.
\]
For a normal surface germ, the local stalk computation is
\[
        i^*Rj_*\mathbb Z_U[2]\cong R\Gamma(L,\mathbb Z)[2].
\]
Hence the critical degree-zero stalk group is
\[
        H^0\bigl(i^*Rj_*\mathbb Z_U[2]\bigr)
        \cong
        H^2(L,\mathbb Z).
\]
The critical torsion is the group \(E\).

\begin{proposition}[The \(E\)-triangle for the Coble boundary singularity]
For the \(\frac14(1,1)\) singularity,
\[
        {}^pIC_X\mathbb Z
        \longrightarrow
        {}^p_+IC_X\mathbb Z
        \longrightarrow
        i_*(\mathbb Z/4)[1]
        \overset{+1}{\longrightarrow}
\]
is the local discrepancy triangle.
\end{proposition}

\begin{proof}
By the local realization theorem of
\cite{RahmanIntegralPerverseObstructions}, the cone of the comparison map
\[
        {}^pIC_X\mathbb Z
        \longrightarrow
        {}^p_+IC_X\mathbb Z
\]
is the point-supported torsion group \(E_X(0)\) shifted by \([1]\):
\[
        \operatorname{Cone}
        \left(
        {}^pIC_X\mathbb Z
        \to
        {}^p_+IC_X\mathbb Z
        \right)
        \cong
        i_*E_X(0)[1].
\]
For the \(\frac14(1,1)\) singularity, the preceding computations give
\[
        E_X(0)\cong\mathbb Z/4.
\]
Thus the triangle becomes
\[
        {}^pIC_X\mathbb Z
        \longrightarrow
        {}^p_+IC_X\mathbb Z
        \longrightarrow
        i_*(\mathbb Z/4)[1]
        \overset{+1}{\longrightarrow}.
\]
\end{proof}

\subsection{MacPherson--Vilonen gluing along a codimension-two Coble stratum}

Let \(Y_0\) be a threefold with a smooth codimension-two stratum
\[
        i:\Sigma\hookrightarrow Y_0
\]
whose transverse singularity is of type \(\frac14(1,1)\).  Let
\[
        j:V=Y_0\setminus\Sigma\hookrightarrow Y_0
\]
be the complementary open stratum.

\begin{theorem}[MV gluing for transverse \(\frac14(1,1)\)-torsion]
Assume the transverse \(\frac14(1,1)\)-type is locally trivial along
\(\Sigma\).  Then the ordinary and dual middle-perversity intersection
complexes fit into a distinguished triangle
\[
        {}^pIC_{Y_0}\mathbb Z
        \longrightarrow
        {}^p_+IC_{Y_0}\mathbb Z
        \longrightarrow
        i_*E_\Sigma[2]
        \overset{+1}{\longrightarrow},
\]
where
\[
        E_\Sigma\cong(\mathbb Z/4)_\Sigma
\]
is the local system on \(\Sigma\) with fiber
\[
        E\cong H^2(L(4,1),\mathbb Z)_{\mathrm{tors}}.
\]
Equivalently,
\[
        \operatorname{Cone}
        \left(
        {}^pIC_{Y_0}\mathbb Z
        \to
        {}^p_+IC_{Y_0}\mathbb Z
        \right)
        \cong
        i_*(\mathbb Z/4)_\Sigma[2].
\]
\end{theorem}

\begin{proof}
We use the MacPherson--Vilonen gluing description of perverse sheaves on a
two-stratum space, equivalently Beilinson's gluing formalism.  This
formalism describes the extension of perverse sheaves across a closed
stratum in terms of the normal Morse data of the stratum.

The two intersection complexes restrict to the same shifted constant sheaf
on the smooth open stratum:
\[
        j^*({}^pIC_{Y_0}\mathbb Z)
        \cong
        \mathbb Z_V[\dim Y_0]
        \cong
        j^*({}^p_+IC_{Y_0}\mathbb Z).
\]
Thus their discrepancy is entirely supported on \(\Sigma\).  Since the
transverse singularity is locally trivial of type \(\frac14(1,1)\), the
normal link is locally constantly \(L(4,1)\).  The critical torsion group of
the normal link is
\[
        H^2(L(4,1),\mathbb Z)_{\mathrm{tors}}
        \cong
        \mathbb Z/4.
\]
These groups assemble into the finite local system
\[
        E_\Sigma\cong(\mathbb Z/4)_\Sigma.
\]
The surface point discrepancy has shift \([1]\), and the smooth
one-dimensional stratum contributes the additional shift \([1]\).  Hence the
closed-stratum discrepancy term is
\[
        i_*E_\Sigma[2].
\]
This gives the stated distinguished triangle.
\end{proof}

\subsection{Motivic lift of the Coble \(E\)-package}

The preceding MacPherson--Vilonen gluing statement identifies the full
closed-stratum discrepancy term as
\[
        i_*E_\Sigma[2]
        \cong
        i_*(\mathbb Z/4)_\Sigma[2].
\]
We now record the corresponding motivic lift.  This formulation is meant in
any integral motivic sheaf formalism with a unit object, cones, shifts,
proper pushforward, and a Betti realization compatible with these operations.
The integral Nori motivic sheaf formalism of Ruimy--Tubach is the intended
ambient setting \cite{RuimyTubach24}.

Let
\[
        i:\Sigma\hookrightarrow Y_0
\]
be a codimension-two stratum with transverse Coble singularity
\[
        \frac14(1,1).
\]
Write \(\mathbf 1_\Sigma\) for the motivic unit on \(\Sigma\).  Define
\[
        \mathbf 1_\Sigma/4
        :=
        \operatorname{Cone}
        \left(
        \mathbf 1_\Sigma
        \xrightarrow{4}
        \mathbf 1_\Sigma
        \right).
\]
Equivalently, in an integral motivic heart, this is the cokernel of
multiplication by \(4\):
\[
        \mathbf 1_\Sigma/4
        =
        \operatorname{coker}
        \left(
        \mathbf 1_\Sigma
        \xrightarrow{4}
        \mathbf 1_\Sigma
        \right).
\]

\begin{definition}[Motivic Coble \(E\)-object]
The motivic lift of the full Coble \(E\)-package is
\[
        \mathcal T^{\mathrm{mot}}_{4,\Sigma}
        :=
        i_*(\mathbf 1_\Sigma/4)[2].
\]
Thus the construction is the functorial composite
\[
        \mathbf 1_\Sigma
        \longmapsto
        \mathbf 1_\Sigma/4
        \longmapsto
        i_*(\mathbf 1_\Sigma/4)
        \longmapsto
        i_*(\mathbf 1_\Sigma/4)[2].
\]
\end{definition}

\begin{proposition}[Betti realization of the full Coble \(E\)-object]
Assume that the Betti realization functor
\[
        \operatorname{Real}_B:
        D^{\mathbb Z}_{\mathrm{mot}}(Y_0)
        \longrightarrow
        D^b_c(Y_0,\mathbb Z)
\]
commutes with cones, multiplication by \(4\), shifts, and \(i_*\).  Then
\[
        \operatorname{Real}_B
        \left(
        \mathcal T^{\mathrm{mot}}_{4,\Sigma}
        \right)
        \cong
        i_*(\mathbb Z/4)_\Sigma[2].
\]
Consequently, the full Coble \(E\)-discrepancy term
\[
        i_*E_\Sigma[2]
        \cong
        i_*(\mathbb Z/4)_\Sigma[2]
\]
is the Betti realization of the motivic object
\[
        i_*(\mathbf 1_\Sigma/4)[2].
\]
\end{proposition}

\begin{proof}
By definition,
\[
        \mathcal T^{\mathrm{mot}}_{4,\Sigma}
        =
        i_*(\mathbf 1_\Sigma/4)[2],
\]
where
\[
        \mathbf 1_\Sigma/4
        =
        \operatorname{Cone}
        \left(
        \mathbf 1_\Sigma
        \xrightarrow{4}
        \mathbf 1_\Sigma
        \right).
\]
Applying Betti realization and using compatibility with cones gives
\[
        \operatorname{Real}_B(\mathbf 1_\Sigma/4)
        \cong
        \operatorname{Cone}
        \left(
        \operatorname{Real}_B(\mathbf 1_\Sigma)
        \xrightarrow{4}
        \operatorname{Real}_B(\mathbf 1_\Sigma)
        \right).
\]
Since
\[
        \operatorname{Real}_B(\mathbf 1_\Sigma)\cong \mathbb Z_\Sigma,
\]
we obtain
\[
        \operatorname{Real}_B(\mathbf 1_\Sigma/4)
        \cong
        \operatorname{Cone}
        \left(
        \mathbb Z_\Sigma
        \xrightarrow{4}
        \mathbb Z_\Sigma
        \right).
\]
The cone of multiplication by \(4\) on \(\mathbb Z_\Sigma\) is
quasi-isomorphic to the constant sheaf
\[
        (\mathbb Z/4)_\Sigma.
\]
Therefore
\[
        \operatorname{Real}_B(\mathbf 1_\Sigma/4)
        \cong
        (\mathbb Z/4)_\Sigma.
\]
Using compatibility of realization with \(i_*\) and shifts,
\[
        \operatorname{Real}_B
        \left(
        i_*(\mathbf 1_\Sigma/4)[2]
        \right)
        \cong
        i_*(\mathbb Z/4)_\Sigma[2].
\]
This proves the proposition.
\end{proof}

\begin{definition}[Motivic Benoist--Ottem shadow]
The motivic order-two shadow of the Coble \(E\)-package is
\[
        \mathcal T^{\mathrm{mot}}_{2,\Sigma}
        :=
        i_*(\mathbf 1_\Sigma/2)[2].
\]
It maps naturally to the full Coble motivic object
\[
        \mathcal T^{\mathrm{mot}}_{4,\Sigma}
        =
        i_*(\mathbf 1_\Sigma/4)[2],
\]
realizing the inclusion
\[
        \mathbb Z/2\hookrightarrow\mathbb Z/4,
        \qquad
        1\longmapsto2.
\]
\end{definition}

\begin{proposition}[Betti realization of the motivic order-two shadow]
Under the same compatibility assumptions on Betti realization,
\[
        \operatorname{Real}_B
        \left(
        \mathcal T^{\mathrm{mot}}_{2,\Sigma}
        \right)
        \cong
        i_*(\mathbb Z/2)_\Sigma[2],
\]
and the morphism
\[
        \mathcal T^{\mathrm{mot}}_{2,\Sigma}
        \longrightarrow
        \mathcal T^{\mathrm{mot}}_{4,\Sigma}
\]
realizes the inclusion
\[
        i_*(\mathbb Z/2)_\Sigma[2]
        \longrightarrow
        i_*(\mathbb Z/4)_\Sigma[2]
\]
with image
\[
        i_*(2E_\Sigma)[2].
\]
\end{proposition}

\begin{proof}
The proof is the same as the mod-\(4\) realization, replacing multiplication
by \(4\) with multiplication by \(2\).  Thus
\[
        \operatorname{Real}_B(\mathbf 1_\Sigma/2)
        \cong
        (\mathbb Z/2)_\Sigma,
\]
and hence
\[
        \operatorname{Real}_B
        \left(
        i_*(\mathbf 1_\Sigma/2)[2]
        \right)
        \cong
        i_*(\mathbb Z/2)_\Sigma[2].
\]
The map
\[
        \mathbf 1_\Sigma/2\to \mathbf 1_\Sigma/4
\]
is the motivic representative of the inclusion
\[
        \mathbb Z/2\hookrightarrow\mathbb Z/4,
        \qquad
        1\mapsto2.
\]
Therefore its realization is the inclusion of the unique order-two subgroup,
and its image is
\[
        i_*(2E_\Sigma)[2].
\]
\end{proof}

\subsection{The order-two shadow and the octahedral filtration}

The Benoist--Ottem mechanism is \(2\)-torsion, so it sees the subgroup
\[
        2E_\Sigma\subset E_\Sigma.
\]
For
\[
        E_\Sigma\cong(\mathbb Z/4)_\Sigma
\]
there is a short exact sequence of local systems
\[
        0
        \longrightarrow
        2E_\Sigma
        \longrightarrow
        E_\Sigma
        \longrightarrow
        E_\Sigma/2E_\Sigma
        \longrightarrow
        0.
\]
Since
\[
        2E_\Sigma\cong(\mathbb Z/2)_\Sigma
        \qquad\text{and}\qquad
        E_\Sigma/2E_\Sigma\cong(\mathbb Z/2)_\Sigma,
\]
this is
\[
        0
        \longrightarrow
        (\mathbb Z/2)_\Sigma
        \longrightarrow
        (\mathbb Z/4)_\Sigma
        \longrightarrow
        (\mathbb Z/2)_\Sigma
        \longrightarrow
        0.
\]

\begin{proposition}[Order-two shadow triangle]
Let
\[
        T_4:=i_*(\mathbb Z/4)_\Sigma[2],
        \qquad
        T_2:=i_*(\mathbb Z/2)_\Sigma[2].
\]
Then the short exact sequence
\[
        0\to(\mathbb Z/2)_\Sigma
        \to
        (\mathbb Z/4)_\Sigma
        \to
        (\mathbb Z/2)_\Sigma
        \to0
\]
induces a distinguished triangle
\[
        T_2
        \longrightarrow
        T_4
        \longrightarrow
        T_2
        \overset{+1}{\longrightarrow}.
\]
The first copy of \(T_2\) is the Bockstein-selected order-two subgroup
\(i_*(2E_\Sigma)[2]\) of the full Coble discrepancy object
\(i_*E_\Sigma[2]\).
\end{proposition}

\begin{proof}
A short exact sequence of sheaves gives a distinguished triangle in the
derived category.  Applying \(i_*\) and shifting by \([2]\) to
\[
        0\to(\mathbb Z/2)_\Sigma
        \to
        (\mathbb Z/4)_\Sigma
        \to
        (\mathbb Z/2)_\Sigma
        \to0
\]
gives
\[
        i_*(\mathbb Z/2)_\Sigma[2]
        \longrightarrow
        i_*(\mathbb Z/4)_\Sigma[2]
        \longrightarrow
        i_*(\mathbb Z/2)_\Sigma[2]
        \overset{+1}{\longrightarrow}.
\]
The inclusion
\[
        \mathbb Z/2\hookrightarrow \mathbb Z/4
\]
is the map \(1\mapsto2\), so its image is the order-two subgroup
\(2E_\Sigma\subset E_\Sigma\).  Hence the first term is the
Bockstein-selected order-two shadow.
\end{proof}

\begin{remark}[Octahedral interpretation]
Let
\[
        A={}^pIC_{Y_0}\mathbb Z,
        \qquad
        B={}^p_+IC_{Y_0}\mathbb Z,
        \qquad
        T_4=i_*E_\Sigma[2].
\]
The discrepancy triangle
\[
        A\to B\to T_4\to A[1]
\]
and the order-two shadow triangle
\[
        T_2\to T_4\to T_2\to T_2[1]
\]
fit into the octahedral comparison for the composite
\[
        B\longrightarrow T_4\longrightarrow T_2.
\]
Thus the Benoist--Ottem \(2\)-torsion contribution is not an object outside
the \(E\)-triangle.  It is a canonical filtered piece/subquotient of the
full Coble \(E\)-discrepancy cone.
\end{remark}

\subsection{The Benoist--Ottem direction}

Let \(S\) be an Enriques surface, let \(B\) be a smooth projective curve, and
set
\[
        Y=S\times B.
\]
Let
\[
        \alpha\in H^1(S,\mathbb Z/2)
\]
be the class of the K3 double cover of \(S\).  Benoist--Ottem show that the
integral Hodge conjecture for \(1\)-cycles on \(Y\) is equivalent to the
condition that every
\[
        \beta\in H^1(B,\mathbb Z/2)
\]
is of the form
\[
        Z^*\alpha
\]
for some correspondence
\[
        Z\in CH_1(B\times S).
\]
Equivalently, the obstruction quotient is
\[
        Q_{\mathrm{BO}}(B,S)
        :=
        H^1(B,\mathbb Z/2)/
        \{Z^*\alpha:Z\in CH_1(B\times S)\}.
\]

The relevant torsion direction lies in the middle Kunneth summand
\[
        H^1(B,\mathbb Z)\otimes H^3(S,\mathbb Z)
        \subset
        H^4(B\times S,\mathbb Z)_{\mathrm{tors}}.
\]
For an Enriques surface,
\[
        H^3(S,\mathbb Z)\cong\mathbb Z/2.
\]

\begin{proposition}[Bockstein generation of the Benoist--Ottem direction]
Let \(u\in H^2(S,\mathbb Z/2)\) satisfy
\[
        \beta_S(u)\neq0\in H^3(S,\mathbb Z)\cong\mathbb Z/2.
\]
Let
\[
        \gamma\in H^1(B,\mathbb Z),
        \qquad
        \bar\gamma\in H^1(B,\mathbb Z/2)
\]
be its mod-\(2\) reduction.  Then
\[
        \beta_Y(u\boxtimes\bar\gamma)
        =
        \beta_S(u)\boxtimes\gamma.
\]
Consequently the middle Kunneth torsion direction
\[
        H^1(B,\mathbb Z)\otimes H^3(S,\mathbb Z)
\]
is generated by Bockstein images of mod-\(2\) Kunneth classes.
\end{proposition}

\begin{proof}
The Bockstein is the connecting homomorphism associated to
\[
        0\to\mathbb Z\xrightarrow{2}\mathbb Z\to\mathbb Z/2\to0.
\]
It satisfies the derivation formula
\[
        \beta(a\smile b)
        =
        \beta(a)\smile b
        +
        (-1)^{|a|}a\smile \beta(b).
\]
Apply this to
\[
        a=u\in H^2(S,\mathbb Z/2),
        \qquad
        b=\bar\gamma\in H^1(B,\mathbb Z/2).
\]
Since \(|u|=2\), the sign is positive.  Thus
\[
        \beta_Y(u\boxtimes\bar\gamma)
        =
        \beta_S(u)\boxtimes\gamma
        +
        u\boxtimes\beta_B(\bar\gamma).
\]
But \(\bar\gamma\) is the reduction modulo \(2\) of an integral class
\(\gamma\), so
\[
        \beta_B(\bar\gamma)=0.
\]
Therefore
\[
        \beta_Y(u\boxtimes\bar\gamma)
        =
        \beta_S(u)\boxtimes\gamma.
\]
Since \(\beta_S(u)\) generates \(H^3(S,\mathbb Z)\cong\mathbb Z/2\), these
classes generate
\[
        H^1(B,\mathbb Z)\otimes H^3(S,\mathbb Z).
\]
\end{proof}

\subsection{Conclusion: the corrected Benoist--Ottem/Coble bridge}

\begin{theorem}[Benoist--Ottem torsion as the order-two Coble shadow]
In the Enriques/Coble boundary model, the local singularity is
\[
        \frac14(1,1),
\]
whose perverse obstruction group is
\[
        E\cong\mathbb Z/4.
\]
The local index-two cover
\[
        A_1\longrightarrow \frac14(1,1)
\]
induces on links the double cover
\[
        L(2,1)\longrightarrow L(4,1).
\]
The Bockstein of its mod-\(2\) cover class is the unique order-two element
\[
        2E\cong\mathbb Z/2
        \subset
        E\cong\mathbb Z/4.
\]
Under MacPherson--Vilonen gluing, the full local \(E\)-package appears as
\[
        i_*E_\Sigma[2]\cong i_*(\mathbb Z/4)_\Sigma[2],
\]
while the Benoist--Ottem \(2\)-torsion mechanism factors through the
canonical order-two shadow
\[
        i_*(2E_\Sigma)[2]
        \cong
        i_*(\mathbb Z/2)_\Sigma[2].
\]
Motivically, the full Coble package is represented by
\[
        \mathcal T^{\mathrm{mot}}_{4,\Sigma}
        =
        i_*(\mathbf 1_\Sigma/4)[2],
\]
while the Benoist--Ottem shadow is represented by
\[
        \mathcal T^{\mathrm{mot}}_{2,\Sigma}
        =
        i_*(\mathbf 1_\Sigma/2)[2]
        \longrightarrow
        \mathcal T^{\mathrm{mot}}_{4,\Sigma}.
\]
If \(\Sigma\cong B\), then
\[
        H^1(\Sigma,2E_\Sigma)
        \cong
        H^1(B,\mathbb Z/2).
\]
Thus the Benoist--Ottem \(2\)-torsion direction is naturally identified with
the Bockstein-selected order-two filtered piece of the Coble boundary
\(E\)-package, not with the full local group \(E\).
\end{theorem}

\begin{proof}
The Coble boundary singularity is \(\frac14(1,1)\), so the local computation
above gives
\[
        E\cong H^2(L(4,1),\mathbb Z)_{\mathrm{tors}}\cong\mathbb Z/4.
\]
The local index-two cover is
\[
        A_1=\mathbb C^2/\mu_2\to \mathbb C^2/\mu_4=\frac14(1,1),
\]
which induces
\[
        L(2,1)\to L(4,1)
\]
on links.  The Bockstein of the class of this double cover is the unique
order-two element
\[
        2E\subset E.
\]
MacPherson--Vilonen gluing identifies the full transverse \(E\)-package
along a codimension-two stratum with
\[
        i_*E_\Sigma[2]\cong i_*(\mathbb Z/4)_\Sigma[2].
\]
The short exact sequence
\[
        0\to2E_\Sigma\to E_\Sigma\to E_\Sigma/2E_\Sigma\to0
\]
gives the order-two shadow triangle
\[
        i_*(\mathbb Z/2)_\Sigma[2]
        \to
        i_*(\mathbb Z/4)_\Sigma[2]
        \to
        i_*(\mathbb Z/2)_\Sigma[2]
        \to.
\]
Thus the \(2\)-torsion contribution selected by the double-cover Bockstein
is the subobject \(i_*(2E_\Sigma)[2]\).

The motivic formulation records the same inclusion before realization:
\[
        i_*(\mathbf 1_\Sigma/2)[2]
        \longrightarrow
        i_*(\mathbf 1_\Sigma/4)[2].
\]
Under Betti realization, this becomes the inclusion
\[
        i_*(\mathbb Z/2)_\Sigma[2]
        \longrightarrow
        i_*(\mathbb Z/4)_\Sigma[2].
\]

On the Benoist--Ottem side, the relevant torsion direction in
\(H^4(S\times B,\mathbb Z)\) is generated by Bockstein images in the middle
Kunneth summand
\[
        H^1(B,\mathbb Z)\otimes H^3(S,\mathbb Z).
\]
Since \(H^3(S,\mathbb Z)\cong\mathbb Z/2\), this has the same
\(H^1(B,\mathbb Z/2)\)-profile as
\[
        H^1(\Sigma,2E_\Sigma)
\]
when \(\Sigma\cong B\).  Hence the Benoist--Ottem \(2\)-torsion direction is
identified with the Bockstein-selected order-two shadow \(2E_\Sigma\) of
the full Coble boundary \(E\)-package.
\end{proof}

\section{Trajectory principles, conjectures, and open problems}

The trajectory computations suggest a small number of organizing principles.
Some are formal consequences of the computations already carried out; others
are conjectural principles that require new global comparison theorems.  We
separate these two levels explicitly.

\subsection{A proved distinction: surface \(A_1\) versus threefold nodes}

The first principle is a proved distinction between two examples that are
often informally conflated.

\begin{proposition}[Surface \(A_1\) versus threefold ordinary double point]
The \(A_1\) surface singularity has local finite torsion package
\[
        E\cong \mathbb Z/2\mathbb Z.
\]
By contrast, the ordinary double point in complex dimension three has link
\[
        L\cong S^2\times S^3,
\]
and hence
\[
        H^m(L,\mathbb Z)_{\mathrm{tors}}=0
        \qquad
        \text{for all }m.
\]
Consequently, the threefold ordinary double point has no finite local
torsion package of the surface \(E\)-type.
\end{proposition}

\begin{proof}
For the surface \(A_1\) singularity, the link is
\[
        \mathbb RP^3.
\]
The computation in Appendix~A gives
\[
        H^2(\mathbb RP^3,\mathbb Z)_{\mathrm{tors}}
        \cong
        \mathbb Z/2\mathbb Z.
\]
The local realization theorem of
\cite{RahmanIntegralPerverseObstructions} identifies this group with
\[
        E=H^0({}^p_+IC_X\mathbb Z)_0.
\]
Hence
\[
        E\cong \mathbb Z/2\mathbb Z.
\]

For the ordinary double point in complex dimension three, Appendix~F
identifies the link with
\[
        S^2\times S^3.
\]
The integral homology and cohomology of \(S^2\times S^3\) are torsion-free by
the Kunneth theorem.  Therefore
\[
        H^m(S^2\times S^3,\mathbb Z)_{\mathrm{tors}}=0
\]
for all \(m\).  Thus there is no finite local link-torsion group from which
a surface-type \(E\)-package could be formed.
\end{proof}

\begin{remark}
This proposition is one of the main safeguards of the paper.  It prevents
one from importing the surface calculation
\[
        E_{A_1}\cong\mathbb Z/2\mathbb Z
\]
into the threefold ordinary double point setting.  Threefold nodes contribute
free vanishing-cycle and exceptional-curve data, not local finite
discriminant torsion.
\end{remark}

\subsection{A proved null-control: the \(E_8\) row}

The second principle is that the torsion package detects non-unimodularity,
not merely singularity complexity.

\begin{proposition}[\(E_8\) is the null-control]
For the \(E_8\) surface singularity, the local finite torsion package
vanishes:
\[
        E=0.
\]
Equivalently,
\[
        H^2(L,\mathbb Z)_{\mathrm{tors}}=0,
        \qquad
        \Lambda^\vee/\Lambda=0,
        \qquad
        \operatorname{coker}(T-\mathrm{id})_{\mathrm{tors}}=0.
\]
\end{proposition}

\begin{proof}
The exceptional lattice of the \(E_8\) singularity is the negative definite
\(E_8\) lattice.  The \(E_8\) lattice is unimodular, so the natural map
\[
        \Lambda\longrightarrow \Lambda^\vee
\]
is an isomorphism.  Hence
\[
        \Lambda^\vee/\Lambda=0.
\]
By the local realization theorem of
\cite{RahmanIntegralPerverseObstructions},
\[
        E\cong \Lambda^\vee/\Lambda.
\]
Therefore \(E=0\).  The link and monodromy formulations are the corresponding
realizations of the same local group, computed explicitly in Appendix~D.
\end{proof}

\begin{remark}
Thus a singularity may have a nontrivial exceptional configuration while
contributing no finite local discriminant torsion.  The invariant \(E\)
detects the discriminant of the local lattice, not singularity complexity
itself.
\end{remark}

\subsection{Form-sensitive survival}

The \(D_4\) computation shows that one must track a finite pairing, not only
a finite group.

\begin{conjecture}[Form-sensitive survival]
The global survival of a local finite torsion class depends on the local
discriminant package
\[
        (E,q),
\]
not only on the abstract finite abelian group \(E\).
\end{conjecture}

\begin{remark}
The \(D_4\) example motivates this conjecture.  It has underlying group
\[
        E\cong(\mathbb Z/2\mathbb Z)^2,
\]
but the discriminant form is nontrivial and contains off-diagonal
\(\mathbb Q/\mathbb Z\)-valued information.  Two singularity configurations
may have isomorphic underlying torsion groups but inequivalent discriminant
forms.  A global transport theory that forgets the form may therefore miss
information visible in the link pairing and in the exceptional lattice.
\end{remark}

\begin{remark}
The Coble boundary example adds a second form-sensitive lesson.  For the
singularity \(\frac14(1,1)\), the full package is
\[
        E\cong\mathbb Z/4,
        \qquad
        q(\bar g,\bar g)=-\frac14.
\]
The Benoist--Ottem mechanism detects the order-two subgroup
\[
        2E\cong\mathbb Z/2
\]
rather than the full group.  Thus a global obstruction theory may be
sensitive not only to the finite group and its form, but also to a
functorially selected subgroup or quotient of the local package.
\end{remark}

\subsection{The codimension-two transverse torsion principle}

The examples suggest that the natural higher-dimensional continuation of the
surface package \(E\) is not an isolated higher-dimensional ordinary double
point.  Rather, it is a codimension-two stratum with transverse normal
surface singularity.

\begin{conjecture}[Codimension-two transverse torsion principle]
Let \(X\) be a complex variety with a codimension-two stratum
\[
        \Sigma\subset X
\]
whose transverse singularity is a normal surface singularity with local
discriminant package
\[
        (E,q).
\]
Then the higher-dimensional carrier of this torsion is a finite local system
\[
        \mathcal E
\]
on \(\Sigma\), with fiber \(E\).  The global contribution of the transverse
torsion package is governed by the cohomology groups
\[
        H^r(\Sigma,\mathcal E),
\]
together with their support, Brauer, residue, and rationalization
trajectories.
\end{conjecture}

\begin{remark}
This conjecture is motivated by the comparison between threefold ordinary
double points and products \(S_0\times C\), where \(S_0\) has surface
singularities.  Isolated threefold nodes have torsion-free link
\(S^2\times S^3\).  By contrast, a codimension-two stratum with transverse
normal surface singularity carries the surface package \(E\) along the
stratum.
\end{remark}

\subsection{A product test for the codimension-two principle}

The codimension-two transverse torsion principle is conjectural in general,
but it has a basic product test case.  This test case is not the corrected
Enriques/Coble boundary model; rather, it is the clean local \(A_1\) model
that illustrates how a surface package becomes a local system along a
codimension-two stratum.

\begin{proposition}[Product test for transverse \(A_1\)-torsion]
Let \((S_0,p)\) be a normal surface germ with an \(A_1\) singularity at \(p\),
and let \(C\) be a smooth connected projective curve.  Set
\[
        X_0:=S_0\times C,
        \qquad
        \Sigma:=\{p\}\times C.
\]
Then \(\Sigma\) is a codimension-two singular stratum of \(X_0\), and the
transverse singularity along \(\Sigma\) is \(A_1\).  The transverse local
discriminant package has fiber
\[
        E_{A_1}\cong\mathbb Z/2\mathbb Z.
\]
Moreover, because the product structure has constant transverse analytic
type along \(\Sigma\), these fibers assemble into the constant finite local
system
\[
        \mathcal E_\Sigma\cong(\mathbb Z/2\mathbb Z)_\Sigma.
\]
Consequently,
\[
        H^r(\Sigma,\mathcal E_\Sigma)
        \cong
        H^r(C,\mathbb Z/2\mathbb Z)
\]
for all \(r\).  In particular, if \(C\) has genus \(g\), then
\[
        H^1(\Sigma,\mathcal E_\Sigma)
        \cong
        (\mathbb Z/2\mathbb Z)^{2g}.
\]
\end{proposition}

\begin{proof}
Since \(S_0\) is a surface and \(C\) is a smooth curve, the product
\(X_0=S_0\times C\) is a threefold.  The singular locus of \(X_0\) near
\(\{p\}\times C\) is
\[
        \operatorname{Sing}(X_0)=\operatorname{Sing}(S_0)\times C.
\]
After shrinking \(S_0\) around \(p\), the surface \(S_0\) has no singular
point other than \(p\).  Hence the singular locus of \(X_0\) is exactly
\[
        \Sigma=\{p\}\times C.
\]
Since \(X_0\) has complex dimension \(3\) and \(\Sigma\) has complex
dimension \(1\), the stratum \(\Sigma\) has codimension \(2\).

Let \(c\in C\).  A neighborhood of \((p,c)\in X_0\) is analytically a product
of a neighborhood of \(p\in S_0\) with a disk in \(C\).  Therefore a
transverse slice to \(\Sigma\) at \((p,c)\) is analytically isomorphic to the
surface germ \((S_0,p)\).  By hypothesis this transverse germ is the
\(A_1\) surface singularity.

The \(A_1\) trajectory computation gives
\[
        E_{A_1}\cong\mathbb Z/2\mathbb Z.
\]
Equivalently,
\[
        E_{A_1}
        \cong
        H^2(\mathbb RP^3,\mathbb Z)_{\mathrm{tors}}
        \cong
        \Lambda^\vee/\Lambda
        \cong
        \mathbb Z/2\mathbb Z,
\]
where \(\Lambda=[-2]\) is the exceptional lattice of the minimal resolution
of the \(A_1\) singularity.  This is the computation carried out in
Appendix~A and follows from the local realization theorem recalled in
Section~4.

Because the local analytic model is the product \((S_0,p)\times C\), the
transverse singularity does not vary as \(c\) moves along \(C\).  Thus the
fiber \(E_{A_1}\) is identified canonically along the stratum.  The resulting
finite local system on \(\Sigma\) is therefore constant:
\[
        \mathcal E_\Sigma\cong(\mathbb Z/2\mathbb Z)_\Sigma.
\]
Since \(\Sigma\cong C\), we obtain
\[
        H^r(\Sigma,\mathcal E_\Sigma)
        \cong
        H^r(C,\mathbb Z/2\mathbb Z).
\]

Finally, if \(C\) has genus \(g\), its first homology group is
\[
        H_1(C,\mathbb Z)\cong\mathbb Z^{2g}.
\]
By the universal coefficient theorem with coefficients
\(\mathbb Z/2\mathbb Z\),
\[
        H^1(C,\mathbb Z/2\mathbb Z)
        \cong
        \operatorname{Hom}(H_1(C,\mathbb Z),\mathbb Z/2\mathbb Z)
        \cong
        (\mathbb Z/2\mathbb Z)^{2g}.
\]
This proves the claim.
\end{proof}

\begin{remark}[The curve case]
The product test used in the Benoist--Ottem/Coble comparison is the special
case
\[
        B=C,
\]
where \(C\) is a smooth projective curve.  In that case
\[
        \Sigma=\{p\}\times C\cong C,
\]
and for the Coble boundary package one obtains
\[
        H^1(\Sigma,2E_\Sigma)
        \cong
        H^1(C,\mathbb Z/2).
\]
If \(C\) has genus \(g\), then
\[
        H^1(C,\mathbb Z/2)\cong(\mathbb Z/2)^{2g}.
\]
Equivalently, these classes classify the degree-two local systems, or
unramified double covers, of \(C\).  Thus the \(2g\)-dimensional
\(\mathbb Z/2\)-vector space appearing in the product test has a direct
geometric interpretation: it records the possible double-cover directions on
the curve along which the Bockstein-selected shadow
\[
        2E_\Sigma\cong(\mathbb Z/2)_\Sigma
\]
may be transported.  This is the group shape compared with the
Benoist--Ottem middle Kunneth torsion direction.
\end{remark}

\subsection{A higher-dimensional product test}

We record a modest higher-dimensional version of the product test.  The
purpose is not to prove global survival, but only to verify that the
codimension-two transverse package behaves as expected under products.

\begin{proposition}[Product test in higher dimension]
Let \((S_0,p)\) be a normal surface germ with isolated singularity at \(p\),
and let
\[
        E_p
\]
be its local discriminant package.  Let \(B\) be a smooth connected
projective variety, and set
\[
        X_0:=S_0\times B,
        \qquad
        \Sigma:=\{p\}\times B.
\]
Then \(\Sigma\) is a codimension-two stratum of \(X_0\), and the transverse
singularity along \(\Sigma\) is analytically isomorphic to \((S_0,p)\).
Consequently the local packages \(E_p\) assemble into a constant finite local
system
\[
        E_\Sigma\cong (E_p)_\Sigma
\]
on \(\Sigma\cong B\).  In particular,
\[
        H^r(\Sigma,E_\Sigma)
        \cong
        H^r(B,E_p)
\]
for all \(r\).

If \((S_0,p)\) is the Coble boundary singularity \(\frac14(1,1)\), then
\[
        E_\Sigma\cong(\mathbb Z/4)_\Sigma,
\]
and the Bockstein-selected order-two shadow is the constant sublocal system
\[
        2E_\Sigma\cong(\mathbb Z/2)_\Sigma.
\]
Thus
\[
        H^r(\Sigma,2E_\Sigma)
        \cong
        H^r(B,\mathbb Z/2)
\]
for all \(r\).
\end{proposition}

\begin{proof}
Since \(S_0\) is a surface germ and \(B\) is smooth, the product
\[
        X_0=S_0\times B
\]
has singular locus, near \(\{p\}\times B\), equal to
\[
        \operatorname{Sing}(X_0)=\{p\}\times B=\Sigma.
\]
Thus
\[
        \operatorname{codim}_{X_0}(\Sigma)
        =
        \operatorname{codim}_{S_0}(\{p\})
        =
        2.
\]

Let \(b\in B\).  Since \(B\) is smooth, a small analytic neighborhood of
\(b\) is a polydisc.  Hence a neighborhood of \((p,b)\in X_0\) is
analytically a product
\[
        (S_0,p)\times (B,b).
\]
A transverse slice to \(\Sigma\) at \((p,b)\) is therefore analytically
isomorphic to the original surface germ \((S_0,p)\).  The transverse local
discriminant package is consequently \(E_p\) at every point of \(\Sigma\).

Because the product structure identifies these transverse germs canonically
as \(b\) varies, the resulting finite local system on \(\Sigma\) is constant:
\[
        E_\Sigma\cong(E_p)_\Sigma.
\]
Since \(\Sigma\cong B\), we obtain
\[
        H^r(\Sigma,E_\Sigma)
        \cong
        H^r(B,E_p).
\]

For the Coble boundary singularity \(\frac14(1,1)\), the local computation
gives
\[
        E_p\cong\mathbb Z/4.
\]
The local index-two cover gives the Bockstein-selected subgroup
\[
        2E_p\cong\mathbb Z/2.
\]
Applying the preceding constant-local-system statement to \(E_p\) and to the
subgroup \(2E_p\) gives
\[
        E_\Sigma\cong(\mathbb Z/4)_\Sigma,
        \qquad
        2E_\Sigma\cong(\mathbb Z/2)_\Sigma,
\]
and therefore
\[
        H^r(\Sigma,2E_\Sigma)
        \cong
        H^r(B,\mathbb Z/2).
\]
\end{proof}

\begin{remark}[The curve case]
The product test used in the Benoist--Ottem/Coble comparison is the special
case
\[
        B=C,
\]
where \(C\) is a smooth projective curve.  In that case
\[
        \Sigma=\{p\}\times C\cong C,
\]
and for the Coble boundary package one obtains
\[
        H^1(\Sigma,2E_\Sigma)
        \cong
        H^1(C,\mathbb Z/2).
\]
If \(C\) has genus \(g\), then
\[
        H^1(C,\mathbb Z/2)\cong(\mathbb Z/2)^{2g}.
\]
This is the group shape compared with the Benoist--Ottem middle Kunneth
torsion direction.
\end{remark}

\begin{remark}[Global survival is not automatic]
The preceding proposition is only a local and sheaf-theoretic product test.
It shows that transverse surface packages assemble into the expected finite
local systems on codimension-two product strata.  It does not imply that
classes in
\[
        H^r(\Sigma,E_\Sigma)
        \quad\text{or}\quad
        H^r(\Sigma,2E_\Sigma)
\]
survive to global cohomology, to the Brauer group, or to unramified
cohomology.  Global survival still requires the support map, any relevant
specialization map, and the residue tests discussed in the main text.
\end{remark}
\subsection{A Coble product test for the corrected Benoist--Ottem comparison}

The corrected Benoist--Ottem boundary model is the Coble
\(\frac14(1,1)\) singularity.  The following product test records the local
system produced by this boundary model.

\begin{proposition}[Product test for transverse Coble torsion]
Let \((S_0,p)\) be a normal surface germ with a \(\frac14(1,1)\) singularity
at \(p\), and let \(C\) be a smooth connected projective curve.  Set
\[
        X_0:=S_0\times C,
        \qquad
        \Sigma:=\{p\}\times C.
\]
Then \(\Sigma\) is a codimension-two singular stratum of \(X_0\), and the
transverse singularity along \(\Sigma\) is \(\frac14(1,1)\).  The transverse
local discriminant package has fiber
\[
        E\cong\mathbb Z/4\mathbb Z.
\]
These fibers assemble into the constant finite local system
\[
        E_\Sigma\cong(\mathbb Z/4\mathbb Z)_\Sigma.
\]
The local index-two cover selects the order-two sublocal system
\[
        2E_\Sigma\cong(\mathbb Z/2\mathbb Z)_\Sigma.
\]
Consequently,
\[
        H^r(\Sigma,2E_\Sigma)
        \cong
        H^r(C,\mathbb Z/2\mathbb Z)
\]
for all \(r\).  In particular, if \(C\) has genus \(g\), then
\[
        H^1(\Sigma,2E_\Sigma)
        \cong
        (\mathbb Z/2\mathbb Z)^{2g}.
\]
\end{proposition}

\begin{proof}
The proof is parallel to the \(A_1\) product test.  Since
\[
        X_0=S_0\times C,
\]
the singular locus near \(\{p\}\times C\) is
\[
        \Sigma=\{p\}\times C.
\]
This is a curve in a threefold, hence a codimension-two stratum.  A
transverse slice to \(\Sigma\) is analytically the surface germ \((S_0,p)\),
which by hypothesis is the quotient singularity \(\frac14(1,1)\).

The Coble boundary computation gives
\[
        E\cong H^2(L(4,1),\mathbb Z)_{\mathrm{tors}}
        \cong\mathbb Z/4\mathbb Z.
\]
The product structure makes this transverse package constant along
\(\Sigma\), so
\[
        E_\Sigma\cong(\mathbb Z/4\mathbb Z)_\Sigma.
\]
The local index-two cover
\[
        L(2,1)\to L(4,1)
\]
has Bockstein image equal to the unique order-two subgroup
\[
        2E\cong\mathbb Z/2\mathbb Z.
\]
Thus the Benoist--Ottem-visible sublocal system is
\[
        2E_\Sigma\cong(\mathbb Z/2\mathbb Z)_\Sigma.
\]
Since \(\Sigma\cong C\),
\[
        H^r(\Sigma,2E_\Sigma)
        \cong
        H^r(C,\mathbb Z/2\mathbb Z).
\]
For genus \(g\),
\[
        H^1(C,\mathbb Z/2\mathbb Z)
        \cong
        (\mathbb Z/2\mathbb Z)^{2g}.
\]
\end{proof}

\begin{remark}
This is the correct product test for the Enriques/Coble boundary
comparison.  The full local package is \(\mathbb Z/4\), while the
Benoist--Ottem \(2\)-torsion channel sees the order-two subpackage
\(2E\cong\mathbb Z/2\).
\end{remark}

\subsection{Nodal defect as free-relation data}

The nodal examples also suggest a clean separation between finite torsion and
free relation phenomena.

\begin{proposition}[Nodal defect is not local finite torsion]
Ordinary double points on projective threefolds contribute free local
vanishing-cycle or exceptional-curve data, not finite local discriminant
torsion.  Consequently, nodal defect should be interpreted as a global
relation phenomenon among free node classes, not as the kernel of a direct
local finite torsion map.
\end{proposition}

\begin{proof}
For each ordinary double point in complex dimension three, the link is
\[
        S^2\times S^3.
\]
As proved above, this link has torsion-free integral cohomology.  Therefore
each node has no local finite torsion package of the surface \(E\)-type.

On the other hand, the Milnor fiber of a threefold ordinary double point has
the homotopy type of \(S^3\), so each node has free vanishing-cycle data.
In a small resolution, the node is replaced by an exceptional
\(\mathbb P^1\), giving free exceptional-curve data.  Relations among these
free local classes are global relations.  This is the type of relation
measured by defect in the nodal hypersurface literature
\cite{Clemens83,Cynk01,Kloosterman22}.  Hence nodal defect is not the kernel
of a direct finite local torsion map.
\end{proof}

\subsection{Benoist--Ottem/Coble specialization problem}

The Benoist--Ottem examples provide a global smooth torsion benchmark.  The
trajectory computations suggest the following corrected specialization
problem.

\begin{problem}[Benoist--Ottem/Coble specialization problem]
Let
\[
        Y_t=S\times C
\]
be a Benoist--Ottem threefold, where \(S\) is an Enriques surface and \(C\)
is a very general smooth projective curve.  Suppose \(Y_t\) occurs as the
general fiber of a degeneration
\[
        \mathcal Y\longrightarrow \Delta
\]
whose special fiber has the form
\[
        Y_0=S_0\times C,
\]
where \(S_0\) is a Coble boundary surface with singularities of type
\[
        \frac14(1,1).
\]
Let
\[
        \Sigma_j=\{p_j\}\times C
\]
be the codimension-two strata in \(Y_0\), and let
\[
        E_j\cong(\mathbb Z/4\mathbb Z)_{\Sigma_j}
\]
be the full transverse Coble \(E\)-systems.  Let
\[
        2E_j\cong(\mathbb Z/2\mathbb Z)_{\Sigma_j}
\]
be the order-two sublocal systems selected by the local index-two covers.

Determine whether the Enriques \(2\)-torsion classes appearing in the
Benoist--Ottem integral Hodge counterexamples specialize to classes generated
by
\[
        H^*(\Sigma_j,2E_j).
\]
Equivalently, determine whether the Benoist--Ottem \(2\)-torsion is genuinely
global smooth torsion, or whether it can arise after degeneration from the
Bockstein-selected order-two shadows inside transverse Coble
\(\frac14(1,1)\) discriminant systems.
\end{problem}

\begin{remark}
The direct local-discriminant map on the smooth fiber \(S\times C\) has zero
source, because \(S\times C\) is smooth.  Thus the Benoist--Ottem torsion is
not directly local on the smooth fiber.  The question is whether it becomes
a shadow of local discriminant data after degeneration to the Enriques/Coble
boundary.  In the corrected model, the full boundary package is
\(\mathbb Z/4\), and the Benoist--Ottem-visible part is \(2E\cong\mathbb Z/2\).
\end{remark}

\subsection{Bockstein-selected shadows of local discriminant packages}

The Coble boundary example suggests a refinement of the form-sensitive
survival principle.  A global obstruction need not see the full local
discriminant group.  It may instead select a canonical subgroup or quotient.

\begin{principle}[Bockstein-selected shadow principle]
Let \(E\) be a local discriminant group and let a global torsion obstruction
be produced by mod-\(n\) cover data.  Then the global obstruction may detect
not the full group \(E\), but the Bockstein-selected \(n\)-torsion shadow
inside \(E\).
\end{principle}

For the Coble boundary singularity \(\frac14(1,1)\), the full local package is
\[
        E\cong\mathbb Z/4\mathbb Z.
\]
The Benoist--Ottem mechanism is \(2\)-torsion.  The double-cover class
\[
        \eta\in H^1(L(4,1),\mathbb Z/2)
\]
has Bockstein
\[
        \beta(\eta)=2\bar g\in \mathbb Z/4\mathbb Z.
\]
Thus Benoist--Ottem sees
\[
        2E\cong\mathbb Z/2\mathbb Z,
\]
the order-two shadow of the full Coble boundary discriminant package.

The motivic transverse torsion object also refines in the Coble boundary
case.  For a transverse \(A_1\)-stratum, the relevant object is
\[
        i_*(\mathbf 1_\Sigma/2)[2].
\]
For a transverse \(\frac14(1,1)\)-stratum, the full local discriminant object
is instead
\[
        i_*(\mathbf 1_\Sigma/4)[2],
\]
whose Betti realization is
\[
        i_*(\mathbb Z/4)_\Sigma[2].
\]
The Benoist--Ottem \(2\)-torsion channel corresponds to the order-two
motivic shadow
\[
        i_*(\mathbf 1_\Sigma/2)[2]
        \longrightarrow
        i_*(\mathbf 1_\Sigma/4)[2],
\]
which realizes to the inclusion
\[
        i_*(\mathbb Z/2)_\Sigma[2]
        \hookrightarrow
        i_*(\mathbb Z/4)_\Sigma[2].
\]

\subsection{Outlook}

The computations in this paper suggest that the next stage of the program
should focus on codimension-two strata with transverse surface singularities.
The required theory should combine:

\begin{enumerate}[label=\textup{(\arabic*)}]
\item the local discriminant package \((E,q)\) for normal surface
singularities;

\item the formation of finite local systems of such packages along
codimension-two strata;

\item support and excision maps for these local systems;

\item comparison with Brauer and unramified cohomology;

\item Bockstein-selected subpackages such as \(2E\subset E\);

\item motivic lifts of the form
\[
        i_*(\mathbf 1_\Sigma/n)[2];
\]

\item specialization maps in degenerating families.
\end{enumerate}

The main open problem is to determine when such transverse local torsion
packages, or their Bockstein-selected shadows, survive globally and when they
account for known integral Hodge obstructions.

\appendix

\section{The \texorpdfstring{\(A_1\)}{A1} surface trajectory}
\label{app:A1surfacetrajectory}

This appendix works out the torsion trajectory for the \(A_1\) surface
singularity in full detail.  This is the basic local model for
\(\mathbb Z/2\)-torsion in the paper.  The computation is useful because all
six local realizations of the obstruction group can be written explicitly:
perverse discrepancy, torsion-sensitive truncation, link torsion, linking
form, resolution lattice, pair sequence, and monodromy.  The conclusion is
that the same group
\[
        E\cong \mathbb Z/2\mathbb Z
\]
appears in each realization, with discriminant form
\[
        q(1,1)=-\frac12 \bmod \mathbb Z
\]
under the negative definite geometric intersection convention fixed in the
main text.

This example should be distinguished from the Enriques/Coble boundary
example.  The \(A_1\) package is the clean local \(\mathbb Z/2\)-model and it
appears as the index-two cover of the Coble boundary singularity
\[
        A_1=\mathbb C^2/\mu_2
        \longrightarrow
        \mathbb C^2/\mu_4=\frac14(1,1).
\]
On links this is
\[
        L(2,1)\longrightarrow L(4,1).
\]
Thus \(A_1\) remains essential to the corrected Benoist--Ottem comparison,
but the boundary singularity itself is \(\frac14(1,1)\), whose full local
package is \(\mathbb Z/4\).

\subsection{Local model}

The \(A_1\) surface singularity is the hypersurface germ
\[
        (X,0)=\{x^2+y^2+z^2=0\}\subset (\mathbb C^3,0).
\]
Equivalently, after a linear change of coordinates, it is analytically
isomorphic to the quadric cone
\[
        \{uv-w^2=0\}\subset \mathbb C^3.
\]
It is also the quotient singularity
\[
        \mathbb C^2/\{\pm 1\},
\]
where \(-1\) acts by
\[
        (a,b)\longmapsto (-a,-b).
\]
The quotient map can be written explicitly as
\[
        \mathbb C^2
        \longrightarrow
        \{uv-w^2=0\},
        \qquad
        (a,b)\longmapsto (u,v,w)=(a^2,b^2,ab).
\]
Indeed,
\[
        uv-w^2=a^2b^2-a^2b^2=0,
\]
and the functions \(a^2,b^2,ab\) generate the invariant ring
\[
        \mathbb C[a,b]^{\{\pm 1\}}.
\]
Thus the germ \((X,0)\) is a normal surface singularity, a rational double
point of type \(A_1\), and an isolated hypersurface surface singularity.
The general topology of links of normal surface singularities is classical;
see Mumford \cite{Mu61}.  The Milnor fibration and monodromy used below are
standard for isolated hypersurface singularities; see Milnor \cite{Mi68}.

Let
\[
        U:=X\setminus\{0\},
        \qquad
        j:U\hookrightarrow X,
        \qquad
        i:\{0\}\hookrightarrow X.
\]
Since \(X\) has complex dimension \(2\), the normalized local system on the
smooth locus is
\[
        \mathbb Z_U[2].
\]
The ordinary and dual middle-perversity intersection complexes are
\[
        {}^pIC_X\mathbb Z
        =
        {}^pj_{!*}\mathbb Z_U[2],
        \qquad
        {}^p_+IC_X\mathbb Z
        =
        {}^p_+j_{!*}\mathbb Z_U[2].
\]
These are the conventions used in \cite{RahmanIntegralPerverseObstructions}.

\subsection{The link and its cohomology}

Let \(L\) be the link of the singularity.  Since
\[
        (X,0)\cong \mathbb C^2/\{\pm 1\},
\]
the link is the quotient of the unit sphere
\[
        S^3\subset \mathbb C^2
\]
by the antipodal action.  Therefore
\[
        L\cong S^3/\{\pm 1\}\cong \mathbb R P^3.
\]
Equivalently, \(L\) is the lens space \(L(2,1)\).

We now compute its integral homology and cohomology.  The standard CW
structure on \(\mathbb R P^3\) has one cell in each dimension
\[
        0,1,2,3.
\]
The cellular chain groups are therefore
\[
        C_k(\mathbb R P^3;\mathbb Z)\cong \mathbb Z
        \quad
        \text{for }0\le k\le 3,
\]
and vanish otherwise.  The cellular boundary maps for real projective space
are
\[
        d_k=
        \begin{cases}
        0, & k \text{ odd},\\
        2, & k \text{ even}.
        \end{cases}
\]
Thus for \(\mathbb R P^3\) the cellular chain complex is
\[
        0
        \longrightarrow
        \mathbb Z
        \xrightarrow{0}
        \mathbb Z
        \xrightarrow{2}
        \mathbb Z
        \xrightarrow{0}
        \mathbb Z
        \longrightarrow
        0.
\]
It follows that
\[
        H_3(L,\mathbb Z)\cong \mathbb Z,
        \qquad
        H_2(L,\mathbb Z)=0,
        \qquad
        H_1(L,\mathbb Z)\cong \mathbb Z/2\mathbb Z,
        \qquad
        H_0(L,\mathbb Z)\cong \mathbb Z.
\]
For reference, this is the standard homology computation for real
projective space; see, for example, \cite[Cellular homology examples]{Hatcher02}.

We now compute cohomology using the universal coefficient theorem.  The
universal coefficient theorem gives a short exact sequence
\[
        0
        \longrightarrow
        \operatorname{Ext}^1_{\mathbb Z}(H_{k-1}(L,\mathbb Z),\mathbb Z)
        \longrightarrow
        H^k(L,\mathbb Z)
        \longrightarrow
        \operatorname{Hom}_{\mathbb Z}(H_k(L,\mathbb Z),\mathbb Z)
        \longrightarrow
        0.
\]
Using the homology groups above, we obtain:
\[
        H^0(L,\mathbb Z)\cong \mathbb Z,
\]
\[
        H^1(L,\mathbb Z)=0,
\]
because
\[
        \operatorname{Hom}(\mathbb Z/2\mathbb Z,\mathbb Z)=0
\]
and
\[
        \operatorname{Ext}^1(\mathbb Z,\mathbb Z)=0.
\]
Next,
\[
        H^2(L,\mathbb Z)
        \cong
        \operatorname{Ext}^1(\mathbb Z/2\mathbb Z,\mathbb Z)
        \cong
        \mathbb Z/2\mathbb Z,
\]
since \(H_2(L,\mathbb Z)=0\).  Finally,
\[
        H^3(L,\mathbb Z)
        \cong
        \operatorname{Hom}(\mathbb Z,\mathbb Z)
        \cong
        \mathbb Z.
\]
Therefore
\[
        H^k(L,\mathbb Z)
        =
        \begin{cases}
        \mathbb Z, & k=0,\\
        0, & k=1,\\
        \mathbb Z/2\mathbb Z, & k=2,\\
        \mathbb Z, & k=3,\\
        0, & \text{otherwise}.
        \end{cases}
\]

The torsion subgroup relevant to the local perverse obstruction is
\[
        H^2(L,\mathbb Z)_{\mathrm{tors}}
        \cong
        \mathbb Z/2\mathbb Z.
\]
This is exactly the link-cohomology realization of the local group \(E\)
proved in \cite[Proposition 3.4]{RahmanIntegralPerverseObstructions}.

\subsection{Birth of torsion}

By the local realization theorem of \cite{RahmanIntegralPerverseObstructions},
for a normal surface singularity one has
\[
        E
        \cong
        H^2(L,\mathbb Z)_{\mathrm{tors}}.
\]
For the \(A_1\) surface singularity, the computation above gives
\[
        H^2(L,\mathbb Z)_{\mathrm{tors}}
        \cong
        \mathbb Z/2\mathbb Z.
\]
Hence
\[
        E\cong \mathbb Z/2\mathbb Z.
\]

The invariant factor decomposition is already
\[
        E\cong \mathbb Z/2\mathbb Z.
\]
Its order is
\[
        |E|=2,
\]
and its prime decomposition is supported only at the prime \(2\).  Thus the
\(A_1\) surface singularity is the basic local \(2\)-torsion atom in the
theory.

\subsection{The discriminant form}

The minimal resolution of the \(A_1\) surface singularity has one exceptional
curve
\[
        C\cong \mathbb P^1.
\]
Its self-intersection is
\[
        C^2=-2.
\]
Thus the exceptional lattice is
\[
        \Lambda=\mathbb Z\langle C\rangle,
\]
with bilinear form
\[
        (C,C)=-2.
\]
The dual lattice is
\[
        \Lambda^\vee
        =
        \{x\in \Lambda\otimes_{\mathbb Z}\mathbb Q
        \mid (x,\Lambda)\subset \mathbb Z\}.
\]
Writing \(x=aC\) with \(a\in\mathbb Q\), the condition
\[
        (aC,C)=-2a\in \mathbb Z
\]
is equivalent to
\[
        a\in \frac12\mathbb Z.
\]
Therefore
\[
        \Lambda^\vee=\frac12\mathbb Z\langle C\rangle.
\]
The discriminant group is
\[
        \Lambda^\vee/\Lambda
        =
        \frac{\frac12\mathbb Z\langle C\rangle}
             {\mathbb Z\langle C\rangle}
        \cong
        \mathbb Z/2\mathbb Z.
\]

Let
\[
        \bar g
        :=
        \frac{C}{2}+\Lambda
\]
be the nonzero element of \(\Lambda^\vee/\Lambda\).  The discriminant pairing
is defined by
\[
        q_\Lambda(x+\Lambda,y+\Lambda)
        =
        (x,y)\bmod \mathbb Z.
\]
Thus
\[
        q_\Lambda(\bar g,\bar g)
        =
        \left(\frac C2,\frac C2\right)\bmod \mathbb Z
        =
        \frac{(C,C)}{4}\bmod \mathbb Z
        =
        -\frac12\bmod \mathbb Z.
\]
Therefore the local discriminant package is
\[
        (E,q)
        =
        \left(\mathbb Z/2\mathbb Z,\; q(\bar g,\bar g)=-\frac12\bmod\mathbb Z\right).
\]
This is the geometric negative-definite convention used throughout the
paper.  If one instead uses the positive \(A_1\) root lattice, the displayed
value changes sign.  The sign convention is fixed in the main text.

The comparison between the lattice discriminant form and the torsion linking
pairing on the link is the standard comparison between discriminant forms of
plumbed resolution lattices and linking pairings on the boundary
three-manifold; see \cite{GoreskySiegel83,Nikulin80}.  In the present
normalization, this compatibility is one of the local realizations recorded
in \cite{RahmanIntegralPerverseObstructions}.

\subsection{Perverse and torsion-sensitive realization}

We now compute the local stalk model for the ordinary and dual
middle-perversity intersection complexes.

Let
\[
        K:=Rj_*\mathbb Z_U[2].
\]
For an isolated surface singularity with link \(L\), the stalk at the
singular point satisfies
\[
        i^*K
        \cong
        R\Gamma(L,\mathbb Z)[2].
\]
This is the standard local computation of the derived pushforward from the
punctured germ, and it is proved in this normalization in
\cite[Lemma 3.3]{RahmanIntegralPerverseObstructions}.  Therefore
\[
        H^m(i^*K)
        \cong
        H^{m+2}(L,\mathbb Z).
\]
Using the cohomology of \(L=\mathbb R P^3\) computed above, we obtain
\[
        H^m(i^*K)
        =
        \begin{cases}
        \mathbb Z, & m=-2,\\
        0, & m=-1,\\
        \mathbb Z/2\mathbb Z, & m=0,\\
        \mathbb Z, & m=1,\\
        0, & \text{otherwise}.
        \end{cases}
\]

The ordinary and dual middle-perversity intermediate extensions are obtained
from this local model by imposing the ordinary and dual codimension-two
point-stratum conditions.  Over \(\mathbb Z\), these differ exactly in the
treatment of torsion in the critical degree.  BBD's integral
middle-perversity formalism records the ordinary and dual middle perversities
\cite[Complement 3.3]{BBD82}.  Friedman's torsion-sensitive Deligne sheaves
give an equivalent way to model this phenomenon by allowing specified
torsion to survive one degree above the ordinary truncation cutoff
\cite{FriedmanGenIH,FriedmanBook20,FriedmanTsInv}.

In the present \(A_1\) surface case, the critical degree-zero stalk of the
local model is
\[
        H^0(i^*K)
        \cong
        H^2(L,\mathbb Z)
        \cong
        \mathbb Z/2\mathbb Z.
\]
The ordinary middle extension removes this critical torsion contribution,
while the dual middle extension retains it.  Consequently
\[
        H^0\bigl({}^p_+IC_X\mathbb Z\bigr)_0
        \cong
        \mathbb Z/2\mathbb Z.
\]
By definition,
\[
        E=H^0\bigl({}^p_+IC_X\mathbb Z\bigr)_0,
\]
so
\[
        E\cong \mathbb Z/2\mathbb Z.
\]
This is the perverse and torsion-sensitive realization of the same group.
The general statement that
\[
        E\cong H^2(L,\mathbb Z)_{\mathrm{tors}}
\]
for normal surface singularities is proved in
\cite[Proposition 3.4]{RahmanIntegralPerverseObstructions}.

Thus in the \(A_1\) example the discrepancy triangle is
\[
        {}^pIC_X\mathbb Z
        \longrightarrow
        {}^p_+IC_X\mathbb Z
        \longrightarrow
        (\mathbb Z/2\mathbb Z)[1]
        \longrightarrow .
\]

\subsection{Resolution lattice realization}

We now compute the resolution-lattice realization directly.  The minimal
resolution of the \(A_1\) surface singularity has exceptional divisor
\[
        E_\pi=C\cong \mathbb P^1
\]
with normal bundle
\[
        \mathcal O_{\mathbb P^1}(-2).
\]
Equivalently, a resolution neighborhood \(N\) of \(C\) is a disk bundle over
\(\mathbb P^1\) with Euler number \(-2\).  The exceptional lattice is
\[
        \Lambda=H_2(N,\mathbb Z)
        =
        \mathbb Z\langle C\rangle,
\]
with intersection matrix
\[
        M=[-2].
\]
The map
\[
        \Lambda\longrightarrow \Lambda^\vee
\]
induced by the intersection form is multiplication by \(-2\):
\[
        \mathbb Z
        \xrightarrow{-2}
        \mathbb Z.
\]
The cokernel is
\[
        \operatorname{coker}([-2])
        \cong
        \mathbb Z/2\mathbb Z.
\]
Thus
\[
        \Lambda^\vee/\Lambda
        \cong
        \mathbb Z/2\mathbb Z.
\]
This agrees with the link-cohomology and perverse realizations:
\[
        E
        \cong
        H^2(L,\mathbb Z)_{\mathrm{tors}}
        \cong
        \Lambda^\vee/\Lambda
        \cong
        \mathbb Z/2\mathbb Z.
\]
The general resolution-lattice theorem in this notation is
\cite[Theorem 1.2]{RahmanIntegralPerverseObstructions}.

The Smith normal form of the \(1\times 1\) matrix
\[
        [-2]
\]
is
\[
        [2],
\]
since multiplication by a unit \(-1\) in \(\mathbb Z\) changes \([-2]\) to
\([2]\).  Hence the invariant factor decomposition of the cokernel is
\[
        \mathbb Z/2\mathbb Z.
\]

\subsection{Pair-sequence realization}

Let \(N\) be a closed disk-bundle neighborhood of the exceptional curve
\[
        C\cong \mathbb P^1
\]
in the minimal resolution.  Its boundary is the link
\[
        L=\partial N\cong \mathbb R P^3.
\]
Since \(N\) deformation retracts onto \(C\), its homology is
\[
        H_0(N,\mathbb Z)\cong \mathbb Z,
        \qquad
        H_2(N,\mathbb Z)\cong \mathbb Z,
\]
and all other reduced homology groups vanish.

Because \(N\) is a compact oriented real \(4\)-manifold with boundary,
Poincaré--Lefschetz duality gives
\[
        H_k(N,L;\mathbb Z)
        \cong
        H^{4-k}(N,\mathbb Z).
\]
Since \(N\) deformation retracts onto \(\mathbb P^1\), we have
\[
        H^0(N,\mathbb Z)\cong \mathbb Z,
        \qquad
        H^2(N,\mathbb Z)\cong \mathbb Z,
\]
and the remaining cohomology groups vanish.  Therefore
\[
        H_2(N,L;\mathbb Z)
        \cong
        H^2(N,\mathbb Z)
        \cong
        \mathbb Z.
\]

The long exact sequence of the pair \((N,L)\) contains
\[
        H_2(L,\mathbb Z)
        \longrightarrow
        H_2(N,\mathbb Z)
        \longrightarrow
        H_2(N,L;\mathbb Z)
        \longrightarrow
        H_1(L,\mathbb Z)
        \longrightarrow
        H_1(N,\mathbb Z).
\]
Substituting the known groups,
\[
        H_2(L,\mathbb Z)=0,
        \qquad
        H_2(N,\mathbb Z)\cong \mathbb Z,
        \qquad
        H_2(N,L;\mathbb Z)\cong \mathbb Z,
        \qquad
        H_1(N,\mathbb Z)=0,
\]
we obtain an exact sequence
\[
        0
        \longrightarrow
        \mathbb Z
        \xrightarrow{\;\; -2\;\;}
        \mathbb Z
        \longrightarrow
        H_1(L,\mathbb Z)
        \longrightarrow
        0.
\]
The middle map is multiplication by the self-intersection
\[
        C^2=-2.
\]
Therefore
\[
        H_1(L,\mathbb Z)
        \cong
        \operatorname{coker}(\mathbb Z\xrightarrow{-2}\mathbb Z)
        \cong
        \mathbb Z/2\mathbb Z.
\]
By the universal coefficient theorem,
\[
        H^2(L,\mathbb Z)_{\mathrm{tors}}
        \cong
        H_1(L,\mathbb Z)_{\mathrm{tors}}
        \cong
        \mathbb Z/2\mathbb Z.
\]
Thus the pair sequence recovers the same group:
\[
        E\cong \mathbb Z/2\mathbb Z.
\]

This computation is the \(A_1\) specialization of the general
resolution-neighborhood argument identifying link torsion with
\(\Lambda^\vee/\Lambda\); see \cite[§4]{RahmanIntegralPerverseObstructions}
and the linking-pairing framework of \cite{GoreskySiegel83}.

\subsection{Monodromy realization}

The \(A_1\) singularity is an isolated hypersurface surface singularity, so
it also has a Milnor-fibration realization.  Let
\[
        f(x,y,z)=x^2+y^2+z^2.
\]
For sufficiently small \(\epsilon\) and \(0<|\eta|\ll \epsilon\), the Milnor
fiber is
\[
        F=f^{-1}(\eta)\cap B_\epsilon.
\]
For the \(A_1\) singularity in complex dimension \(2\), the Milnor fiber has
the homotopy type of a \(2\)-sphere:
\[
        F\simeq S^2.
\]
Thus
\[
        H^2(F,\mathbb Z)\cong \mathbb Z,
\]
and this is the middle vanishing cohomology.

The monodromy \(T\) for the \(A_1\) singularity acts by
\[
        T=-\mathrm{id}
\]
on \(H^2(F,\mathbb Z)\cong \mathbb Z\).  This is the standard
Picard--Lefschetz monodromy for a single vanishing cycle in the \(A_1\)
case; see Milnor's discussion of isolated hypersurface singularities
\cite{Mi68}.  Therefore
\[
        T-\mathrm{id}
        =
        -\mathrm{id}-\mathrm{id}
        =
        -2
\]
as an endomorphism of \(\mathbb Z\).  Hence
\[
        \operatorname{coker}(T-\mathrm{id})
        =
        \operatorname{coker}(\mathbb Z\xrightarrow{-2}\mathbb Z)
        \cong
        \mathbb Z/2\mathbb Z.
\]
Since this cokernel is finite, its torsion subgroup is the whole group:
\[
        \operatorname{coker}(T-\mathrm{id})_{\mathrm{tors}}
        \cong
        \mathbb Z/2\mathbb Z.
\]
The monodromy realization theorem of
\cite[Theorem 1.3]{RahmanIntegralPerverseObstructions} gives
\[
        E
        \cong
        \operatorname{coker}(T-\mathrm{id})_{\mathrm{tors}},
\]
and therefore again
\[
        E\cong \mathbb Z/2\mathbb Z.
\]

The determinant refinement is also immediate.  Since
\[
        T-\mathrm{id}=-2
\]
on the rank-one free abelian group \(H^2(F,\mathbb Z)\), one has
\[
        |\det(T-\mathrm{id})|=2.
\]
This agrees with
\[
        |E|=2.
\]

\subsection{The \(A_1\) cover class and the Coble boundary}

Although the \(A_1\) surface singularity is not the corrected Enriques/Coble
boundary singularity, it appears as the local index-two cover of that
boundary singularity.  Namely,
\[
        \frac14(1,1)=\mathbb C^2/\mu_4
\]
has intermediate quotient
\[
        \mathbb C^2/\mu_2\longrightarrow\mathbb C^2/\mu_4,
\]
where
\[
        \mathbb C^2/\mu_2=A_1.
\]
On links this is the double cover
\[
        L(2,1)\longrightarrow L(4,1).
\]
Thus the \(A_1\) link
\[
        L(2,1)=\mathbb RP^3
\]
is the source of the local index-two cover used in the Coble boundary
comparison.  The Bockstein of the corresponding class in
\[
        H^1(L(4,1),\mathbb Z/2)
\]
lands in the order-two subgroup
\[
        2E\cong\mathbb Z/2
\]
of the Coble package
\[
        E\cong\mathbb Z/4.
\]
In this sense, the \(A_1\) computation remains part of the
Benoist--Ottem/Coble story, but one level upstairs from the actual boundary
singularity.

\subsection{Transport status}

The \(A_1\) computation above is a local surface computation.  It determines
the birth, form, and six local realizations of the torsion group
\[
        E\cong \mathbb Z/2\mathbb Z.
\]
However, it does not by itself define a global Brauer or unramified class.
To obtain a global class, one must place the singularity inside a global
variety and specify the relevant local-to-global support map.

For the local resolution neighborhood \(N\) of the exceptional curve \(C\),
cohomology with supports along \(C\) is related to ordinary cohomology by the
long exact sequence
\[
\cdots
\longrightarrow
H^2_C(N,\mathbb Z)
\longrightarrow
H^2(N,\mathbb Z)
\longrightarrow
H^2(N\setminus C,\mathbb Z)
\longrightarrow
H^3_C(N,\mathbb Z)
\longrightarrow
\cdots .
\]
The space \(N\setminus C\) deformation retracts onto the boundary circle
bundle \(L=\partial N\), so
\[
        H^2(N\setminus C,\mathbb Z)
        \cong
        H^2(L,\mathbb Z)
        \cong
        \mathbb Z/2\mathbb Z.
\]
The Thom isomorphism for the oriented real rank-two normal bundle of \(C\)
in \(N\) gives
\[
        H^2_C(N,\mathbb Z)\cong H^0(C,\mathbb Z)\cong \mathbb Z,
\]
and
\[
        H^3_C(N,\mathbb Z)\cong H^1(C,\mathbb Z)=0.
\]
Thus the relevant portion of the support sequence becomes
\[
        \mathbb Z
        \longrightarrow
        \mathbb Z
        \longrightarrow
        \mathbb Z/2\mathbb Z
        \longrightarrow
        0.
\]
The first map is multiplication by the Euler class of the normal bundle,
which is the self-intersection
\[
        C^2=-2.
\]
Therefore the cokernel is
\[
        \mathbb Z/2\mathbb Z.
\]
This is another way to see that the local torsion is a boundary quotient
created by the exceptional divisor.

Thus, for an isolated \(A_1\) singularity on a global surface, the local
torsion naturally appears in the degree-two local support sequence around
the exceptional curve.  It is not automatically a degree-three class.  The
degree-three Brauer channel used for smooth projective threefolds therefore
requires a separate global or higher-dimensional construction.

In trajectory language, the local \(A_1\) surface row has the following
transport status:
\[
        \text{support degree: }2\text{ locally, through the exceptional curve
        support sequence}.
\]
A global image depends on how the exceptional curve class maps into the
global cohomology of the chosen resolution.  If several singularities are
present, global relations among exceptional curve classes may kill linear
combinations of the local torsion packages.  This is the surface analogue of
the kernel phenomenon studied in the main text.

\subsection{Rational death}

Since
\[
        E\cong \mathbb Z/2\mathbb Z,
\]
we have
\[
        E\otimes_{\mathbb Z}\mathbb Q=0.
\]
Indeed, if \(e\in E\), then \(2e=0\).  In the tensor product
\(E\otimes_{\mathbb Z}\mathbb Q\), the integer \(2\) acts invertibly on
\(\mathbb Q\).  Hence
\[
        e\otimes 1
        =
        e\otimes 2\cdot \frac12
        =
        2e\otimes \frac12
        =
        0.
\]
Thus every element of \(E\) dies after rationalization.

This rational death is exactly why the ordinary and dual middle-perversity
intersection complexes agree over \(\mathbb Q\):
\[
        {}^pIC_X\mathbb Q
        \cong
        {}^p_+IC_X\mathbb Q.
\]
The torsion group does not survive as a rational vector-space invariant.
Nevertheless, before it dies, it carries the discriminant form
\[
        q(\bar g,\bar g)=-\frac12\bmod \mathbb Z,
\]
appears as the link torsion of \(\mathbb R P^3\), appears as the cokernel of
the intersection matrix \([-2]\), and appears as the monodromy cokernel of
\[
        T-\mathrm{id}=-2.
\]
Thus rationalization kills the group but not the evidence that the integral
geometry contained a nontrivial \(2\)-torsion obstruction.

\subsection{Trajectory row}

\begin{center}
\begingroup
\scriptsize
\setlength{\tabcolsep}{2.2pt}
\renewcommand{\arraystretch}{1.15}
\begin{tabularx}{\textwidth}{|p{0.15\textwidth}|c|c|p{0.12\textwidth}|p{0.10\textwidth}|X|p{0.10\textwidth}|c|}
\hline
\textbf{Example} & \textbf{\(E\)} & \textbf{\(q\)} & \textbf{Local} & \textbf{Supp.} & \textbf{Global image} & \textbf{Br/res.} & \textbf{\(\mathbb Q\)}
\\
\hline
\(A_1\) surface &
\(\mathbb Z/2\) &
\(-\tfrac12\) &
six agree &
deg.\ \(2\) &
depends on global exceptional-curve relations &
loc.\ \(\varnothing\) &
\(0\)
\\
\hline
\end{tabularx}
\endgroup
\end{center}

Here \(\checkmark\) means that all six local realizations agree:
\[
        E
        \cong
        H^2(L,\mathbb Z)_{\mathrm{tors}}
        \cong
        \Lambda^\vee/\Lambda
        \cong
        \operatorname{coker}(T-\mathrm{id})_{\mathrm{tors}}
        \cong
        \mathbb Z/2\mathbb Z,
\]
together with the perverse and torsion-sensitive truncation realizations.
The support degree is listed as \(2\) because the local surface package
appears through the degree-two support sequence of the exceptional curve.
The Brauer/residue column is marked local-\(\varnothing\) because the
isolated surface germ by itself does not define a degree-three Brauer or
unramified class.  To enter that channel, the \(A_1\) package must be placed
inside an appropriate global threefold or compared with a
codimension-two transverse surface stratum.  In the corrected
Benoist--Ottem/Coble comparison, \(A_1\) appears as the local index-two cover
of the \(\frac14(1,1)\) boundary singularity, not as the boundary
singularity itself.

\section{The \texorpdfstring{\(A_k\)}{Ak} surface trajectory}
\label{app:Aksurfacetrajectory}

This appendix works out the torsion trajectory for the \(A_k\) surface
singularities.  The \(A_1\) case is the basic \(2\)-torsion atom, while the
\(A_k\) family shows how the same surface mechanism produces arbitrary cyclic
torsion:
\[
        E\cong \mathbb Z/(k+1)\mathbb Z.
\]
The computation also fixes the generator convention for the discriminant form
of the negative definite \(A_k\) exceptional lattice.

This family should be kept conceptually separate from the Coble boundary
singularity \(\frac14(1,1)\).  The \(A_k\) singularities are rational double
points with quotient form
\[
        \mathbb C^2/\frac1{k+1}(1,-1),
\]
whereas the Coble boundary example is the cyclic quotient
\[
        \mathbb C^2/\frac14(1,1).
\]
The case \(A_1=\frac12(1,1)\) appears as the index-two cover of the Coble
boundary singularity, but the Coble boundary itself has local package
\[
        E\cong\mathbb Z/4.
\]

Throughout this appendix \(k\ge 1\).

\subsection{Local model}

We use the hypersurface model
\[
        (X,0)=\{xy-z^{k+1}=0\}\subset(\mathbb C^3,0).
\]
This is the rational double point of type \(A_k\).  It is also the quotient
singularity
\[
        \mathbb C^2/\mu_{k+1},
\]
where \(\mu_{k+1}\) is the group of \((k+1)\)-st roots of unity and
\[
        \zeta\cdot(u,v)=(\zeta u,\zeta^{-1}v).
\]
Indeed, the invariant ring is generated by
\[
        x=u^{k+1},
        \qquad
        y=v^{k+1},
        \qquad
        z=uv.
\]
These generators satisfy
\[
        xy=u^{k+1}v^{k+1}=(uv)^{k+1}=z^{k+1}.
\]
Conversely, the invariant monomials are exactly the monomials \(u^a v^b\)
with
\[
        a-b\equiv 0 \pmod{k+1}.
\]
Such a monomial can be written as a product of powers of
\[
        u^{k+1},\qquad v^{k+1},\qquad uv.
\]
Thus
\[
        \mathbb C[u,v]^{\mu_{k+1}}
        \cong
        \mathbb C[x,y,z]/(xy-z^{k+1}).
\]
This proves the quotient description.

The singularity is isolated.  Indeed, for
\[
        f(x,y,z)=xy-z^{k+1},
\]
the partial derivatives are
\[
        \frac{\partial f}{\partial x}=y,
        \qquad
        \frac{\partial f}{\partial y}=x,
        \qquad
        \frac{\partial f}{\partial z}=-(k+1)z^k.
\]
All three vanish simultaneously only when
\[
        x=y=z=0.
\]
Thus \(0\) is the unique singular point of the germ.

Let
\[
        U:=X\setminus\{0\},
        \qquad
        j:U\hookrightarrow X,
        \qquad
        i:\{0\}\hookrightarrow X.
\]
Since \(X\) is a complex surface, the shifted constant sheaf on the smooth
locus is
\[
        \mathbb Z_U[2].
\]
The ordinary and dual middle-perversity intersection complexes are
\[
        {}^pIC_X\mathbb Z
        =
        {}^pj_{!*}\mathbb Z_U[2],
        \qquad
        {}^p_+IC_X\mathbb Z
        =
        {}^p_+j_{!*}\mathbb Z_U[2].
\]
These conventions are the same as those used in
\cite{RahmanIntegralPerverseObstructions}.

\subsection{The link and its cohomology}

Let \(L\) be the link of the \(A_k\) singularity.  Since
\[
        (X,0)\cong \mathbb C^2/\mu_{k+1},
\]
the link is the quotient of the unit sphere
\[
        S^3\subset\mathbb C^2
\]
by the free action
\[
        \zeta\cdot(u,v)=(\zeta u,\zeta^{-1}v).
\]
Therefore \(L\) is a lens space.  With the usual orientation conventions,
one may write it as
\[
        L\cong L(k+1,k),
\]
which is orientation-reversingly diffeomorphic to \(L(k+1,1)\).  The
orientation convention affects the sign of the linking form, but it does not
affect the underlying homology group.  The description of links of normal
surface singularities as boundaries of resolution neighborhoods is classical;
see Mumford \cite{Mu61}.  The homology of lens spaces is standard; see
\cite{Hatcher02}.

For a lens space \(L(p,q)\), one has
\[
        H_0(L(p,q),\mathbb Z)\cong \mathbb Z,
        \qquad
        H_1(L(p,q),\mathbb Z)\cong \mathbb Z/p\mathbb Z,
\]
\[
        H_2(L(p,q),\mathbb Z)=0,
        \qquad
        H_3(L(p,q),\mathbb Z)\cong \mathbb Z.
\]
In the present case
\[
        p=k+1.
\]
Hence
\[
        H_0(L,\mathbb Z)\cong \mathbb Z,
        \qquad
        H_1(L,\mathbb Z)\cong \mathbb Z/(k+1)\mathbb Z,
\]
\[
        H_2(L,\mathbb Z)=0,
        \qquad
        H_3(L,\mathbb Z)\cong \mathbb Z.
\]

Using the universal coefficient theorem, we obtain
\[
        H^0(L,\mathbb Z)\cong \mathbb Z,
\]
\[
        H^1(L,\mathbb Z)=0,
\]
\[
        H^2(L,\mathbb Z)
        \cong
        \operatorname{Ext}^1_{\mathbb Z}
        (\mathbb Z/(k+1)\mathbb Z,\mathbb Z)
        \cong
        \mathbb Z/(k+1)\mathbb Z,
\]
and
\[
        H^3(L,\mathbb Z)\cong \mathbb Z.
\]
Therefore
\[
        H^m(L,\mathbb Z)
        =
        \begin{cases}
        \mathbb Z, & m=0,\\
        0, & m=1,\\
        \mathbb Z/(k+1)\mathbb Z, & m=2,\\
        \mathbb Z, & m=3,\\
        0, & \text{otherwise}.
        \end{cases}
\]

The torsion subgroup relevant to the local perverse obstruction is therefore
\[
        H^2(L,\mathbb Z)_{\mathrm{tors}}
        \cong
        \mathbb Z/(k+1)\mathbb Z.
\]

\subsection{Birth of torsion}

For a normal surface singularity, the local realization theorem of
\cite{RahmanIntegralPerverseObstructions} gives
\[
        E
        \cong
        H^2(L,\mathbb Z)_{\mathrm{tors}}.
\]
Using the cohomology computation above, we obtain
\[
        E
        \cong
        \mathbb Z/(k+1)\mathbb Z.
\]

Thus the invariant factor decomposition is
\[
        E\cong \mathbb Z/(k+1)\mathbb Z.
\]
Its order is
\[
        |E|=k+1.
\]
If
\[
        k+1=\prod_{\ell} \ell^{a_\ell}
\]
is the prime decomposition of \(k+1\), then the primary decomposition of
\(E\) is
\[
        E
        \cong
        \bigoplus_{\ell\mid (k+1)}
        \mathbb Z/\ell^{a_\ell}\mathbb Z.
\]
Thus the \(A_k\) singularity produces cyclic torsion whose prime support is
exactly the set of primes dividing \(k+1\).

\subsection{The discriminant form}

The minimal resolution of the \(A_k\) singularity has exceptional divisor
\[
        E_\pi=C_1\cup C_2\cup\cdots\cup C_k,
\]
where each
\[
        C_i\cong \mathbb P^1,
\]
and the curves form a chain:
\[
        C_i\cdot C_i=-2,
        \qquad
        C_i\cdot C_{i+1}=1,
        \qquad
        C_i\cdot C_j=0
        \quad\text{if }|i-j|>1.
\]
Thus the exceptional lattice is
\[
        \Lambda
        =
        \bigoplus_{i=1}^k \mathbb Z\langle C_i\rangle,
\]
with intersection matrix
\[
        M=
        \begin{pmatrix}
        -2 & 1  & 0  & \cdots & 0\\
        1  & -2 & 1  & \cdots & 0\\
        0  & 1  & -2 & \ddots & \vdots\\
        \vdots & \vdots & \ddots & \ddots & 1\\
        0 & 0 & \cdots & 1 & -2
        \end{pmatrix}.
\]
This is the negative of the \(A_k\) Cartan matrix.  The determinant of the
positive \(A_k\) Cartan matrix is \(k+1\), hence
\[
        |\det(M)|=k+1.
\]

We now compute the discriminant group and form explicitly.  Let
\[
        g\in \Lambda^\vee
\]
be the first dual vector, characterized by
\[
        (g,C_1)=1,
        \qquad
        (g,C_j)=0
        \quad\text{for }2\le j\le k.
\]
Equivalently, if \(M^{-1}\) is the inverse of the intersection matrix, then
\[
        g=\sum_{i=1}^k (M^{-1})_{i1}C_i.
\]
For the positive \(A_k\) Cartan matrix \(A\), one has
\[
        (A^{-1})_{ij}
        =
        \frac{\min(i,j)(k+1-\max(i,j))}{k+1}.
\]
Since
\[
        M=-A,
\]
we have
\[
        M^{-1}=-A^{-1}.
\]
Therefore
\[
        (M^{-1})_{i1}
        =
        -\frac{k+1-i}{k+1}.
\]
Hence
\[
        g
        =
        -\frac{1}{k+1}
        \sum_{i=1}^k (k+1-i)C_i.
\]

The class
\[
        \bar g:=g+\Lambda
\]
generates the discriminant group
\[
        \Lambda^\vee/\Lambda.
\]
Since
\[
        |\Lambda^\vee/\Lambda|=|\det(M)|=k+1,
\]
it follows that
\[
        \Lambda^\vee/\Lambda
        \cong
        \mathbb Z/(k+1)\mathbb Z.
\]

The discriminant pairing is
\[
        q_\Lambda:
        \Lambda^\vee/\Lambda\times \Lambda^\vee/\Lambda
        \longrightarrow
        \mathbb Q/\mathbb Z,
        \qquad
        q_\Lambda(x+\Lambda,y+\Lambda)=(x,y)\bmod \mathbb Z.
\]
For the generator \(\bar g\), we compute
\[
        q_\Lambda(\bar g,\bar g)
        =
        (g,g)\bmod \mathbb Z.
\]
Since \(g\) is the first dual vector, one has
\[
        (g,g)=(M^{-1})_{11}.
\]
Using the formula above,
\[
        (M^{-1})_{11}
        =
        -\frac{k}{k+1}.
\]
Therefore
\[
        q_\Lambda(\bar g,\bar g)
        =
        -\frac{k}{k+1}\bmod \mathbb Z.
\]
Equivalently, since
\[
        -\frac{k}{k+1}
        =
        \frac{1}{k+1}-1,
\]
one may also write
\[
        q_\Lambda(\bar g,\bar g)
        =
        \frac{1}{k+1}\bmod \mathbb Z.
\]
The first expression records the negative definite geometric intersection
convention, while the second is the same element of \(\mathbb Q/\mathbb Z\).

For arbitrary elements
\[
        a\bar g,\ b\bar g\in \Lambda^\vee/\Lambda,
\]
the pairing is
\[
        q_\Lambda(a\bar g,b\bar g)
        =
        -\frac{k\,ab}{k+1}\bmod \mathbb Z.
\]
Thus the local discriminant package is
\[
        (E,q)
        =
        \left(
        \mathbb Z/(k+1)\mathbb Z,
        \;
        q(a,b)=-\frac{k\,ab}{k+1}\bmod \mathbb Z
        \right)
\]
with respect to the generator \(\bar g\).

This discriminant-form computation is the \(A_k\) specialization of the
standard discriminant form of the negative \(A_k\) lattice; see Nikulin for
the general lattice formalism \cite{Nikulin80}.  The compatibility of the
exceptional lattice discriminant group with the local obstruction group \(E\)
is proved in \cite{RahmanIntegralPerverseObstructions}.

\subsection{Perverse and torsion-sensitive realization}

Let
\[
        K:=Rj_*\mathbb Z_U[2].
\]
For an isolated surface singularity with link \(L\), the local stalk
calculation gives
\[
        i^*K
        \cong
        R\Gamma(L,\mathbb Z)[2].
\]
Therefore
\[
        H^m(i^*K)
        \cong
        H^{m+2}(L,\mathbb Z).
\]
This is the local stalk computation used in
\cite[Lemma 3.3]{RahmanIntegralPerverseObstructions}.

Using the cohomology of the lens space \(L\), we obtain
\[
        H^m(i^*K)
        =
        \begin{cases}
        \mathbb Z, & m=-2,\\
        0, & m=-1,\\
        \mathbb Z/(k+1)\mathbb Z, & m=0,\\
        \mathbb Z, & m=1,\\
        0, & \text{otherwise}.
        \end{cases}
\]

The ordinary and dual middle-perversity intersection complexes are obtained
from this local model by imposing the ordinary and dual codimension-two
point-stratum conditions.  Over \(\mathbb Z\), these differ in their
treatment of torsion in the critical degree.  BBD's integral
middle-perversity formalism records the distinction between ordinary and
dual middle perversity \cite[Complement 3.3]{BBD82}.  Friedman's
torsion-sensitive Deligne sheaves give a truncation framework in which
specified torsion may survive one degree above the ordinary cutoff
\cite{FriedmanGenIH,FriedmanBook20,FriedmanTsInv}.

In the \(A_k\) surface case, the critical degree-zero stalk is
\[
        H^0(i^*K)
        \cong
        H^2(L,\mathbb Z)
        \cong
        \mathbb Z/(k+1)\mathbb Z.
\]
The ordinary middle extension removes this critical torsion contribution,
while the dual middle extension retains it.  Consequently
\[
        H^0({}^p_+IC_X\mathbb Z)_0
        \cong
        \mathbb Z/(k+1)\mathbb Z.
\]
By definition,
\[
        E=H^0({}^p_+IC_X\mathbb Z)_0.
\]
Thus
\[
        E\cong \mathbb Z/(k+1)\mathbb Z.
\]
The general identification
\[
        E\cong H^2(L,\mathbb Z)_{\mathrm{tors}}
\]
for normal surface singularities is proved in
\cite[Proposition 3.4]{RahmanIntegralPerverseObstructions}.

The discrepancy triangle in this case is therefore
\[
        {}^pIC_X\mathbb Z
        \longrightarrow
        {}^p_+IC_X\mathbb Z
        \longrightarrow
        (\mathbb Z/(k+1)\mathbb Z)[1]
        \longrightarrow .
\]

\subsection{Resolution lattice realization}

We now compute the resolution-lattice realization directly.

The exceptional lattice is
\[
        \Lambda=\bigoplus_{i=1}^k \mathbb Z\langle C_i\rangle,
\]
with intersection matrix
\[
        M=
        \begin{pmatrix}
        -2 & 1  & 0  & \cdots & 0\\
        1  & -2 & 1  & \cdots & 0\\
        0  & 1  & -2 & \ddots & \vdots\\
        \vdots & \vdots & \ddots & \ddots & 1\\
        0 & 0 & \cdots & 1 & -2
        \end{pmatrix}.
\]
The map
\[
        \Lambda\longrightarrow \Lambda^\vee
\]
is represented by \(M\).  Hence
\[
        \Lambda^\vee/\Lambda
        \cong
        \operatorname{coker}(M:\mathbb Z^k\to\mathbb Z^k).
\]

The Smith normal form of the \(A_k\) Cartan matrix has invariant factors
\[
        1,1,\ldots,1,k+1.
\]
Since \(M\) differs from the positive Cartan matrix by multiplication by
\(-1\), the same invariant factors occur for \(M\).  Thus
\[
        \operatorname{coker}(M)
        \cong
        \mathbb Z/(k+1)\mathbb Z.
\]
Consequently
\[
        \Lambda^\vee/\Lambda
        \cong
        \mathbb Z/(k+1)\mathbb Z.
\]

This agrees with the link and perverse computations:
\[
        E
        \cong
        H^2(L,\mathbb Z)_{\mathrm{tors}}
        \cong
        \Lambda^\vee/\Lambda
        \cong
        \mathbb Z/(k+1)\mathbb Z.
\]
The general theorem identifying \(E\) with the exceptional discriminant group
is \cite[Theorem 1.2]{RahmanIntegralPerverseObstructions}.

\subsection{Pair-sequence realization}

Let \(N\) be a resolution neighborhood of the exceptional chain
\[
        C_1\cup\cdots\cup C_k.
\]
Then \(N\) deformation retracts onto the exceptional divisor.  The boundary
of \(N\) is the link
\[
        L=\partial N.
\]

The long exact sequence of the pair \((N,L)\) contains
\[
        H_2(L,\mathbb Z)
        \longrightarrow
        H_2(N,\mathbb Z)
        \longrightarrow
        H_2(N,L;\mathbb Z)
        \longrightarrow
        H_1(L,\mathbb Z)
        \longrightarrow
        H_1(N,\mathbb Z).
\]
For the \(A_k\) resolution neighborhood,
\[
        H_2(N,\mathbb Z)\cong \mathbb Z^k,
\]
generated by the exceptional curves \(C_1,\ldots,C_k\).  The neighborhood
\(N\) is obtained by plumbing disk bundles over spheres according to a tree,
so
\[
        H_1(N,\mathbb Z)=0.
\]
Also \(H_2(L,\mathbb Z)=0\) for the lens space \(L=L(k+1,k)\).  Therefore
the relevant part of the exact sequence becomes
\[
        0
        \longrightarrow
        \mathbb Z^k
        \longrightarrow
        H_2(N,L;\mathbb Z)
        \longrightarrow
        H_1(L,\mathbb Z)
        \longrightarrow
        0.
\]

By Poincare--Lefschetz duality for the compact oriented real four-manifold
\(N\),
\[
        H_2(N,L;\mathbb Z)
        \cong
        H^2(N,\mathbb Z).
\]
Since \(N\) deformation retracts onto the exceptional chain, and the
exceptional chain has \(H_2\cong\mathbb Z^k\) and no \(H_1\), the universal
coefficient theorem gives
\[
        H^2(N,\mathbb Z)\cong \operatorname{Hom}(H_2(N,\mathbb Z),\mathbb Z)
        \cong
        \mathbb Z^k.
\]
Under these identifications, the map
\[
        H_2(N,\mathbb Z)
        \longrightarrow
        H_2(N,L;\mathbb Z)
\]
is the intersection pairing, hence is represented by the matrix \(M\).
Therefore
\[
        H_1(L,\mathbb Z)
        \cong
        \operatorname{coker}(M)
        \cong
        \mathbb Z/(k+1)\mathbb Z.
\]
By the universal coefficient theorem for the closed oriented three-manifold
\(L\),
\[
        H^2(L,\mathbb Z)_{\mathrm{tors}}
        \cong
        H_1(L,\mathbb Z)_{\mathrm{tors}}
        \cong
        \mathbb Z/(k+1)\mathbb Z.
\]
Thus the pair sequence recovers
\[
        E\cong \mathbb Z/(k+1)\mathbb Z.
\]
This is the \(A_k\) specialization of the general pair-sequence computation
in \cite[§4]{RahmanIntegralPerverseObstructions}, and it is compatible with
the linking-pairing framework of \cite{GoreskySiegel83}.

\subsection{Monodromy realization}

The \(A_k\) singularity is an isolated hypersurface surface singularity.  We
therefore have a Milnor fibration.  For
\[
        f(x,y,z)=xy-z^{k+1},
\]
let
\[
        F=f^{-1}(\eta)\cap B_\epsilon
\]
be the Milnor fiber.  The Milnor number is
\[
        \mu=k.
\]
The Milnor fiber has the homotopy type of a bouquet of \(k\) two-spheres:
\[
        F\simeq \bigvee_{j=1}^k S^2.
\]
Thus
\[
        H^2(F,\mathbb Z)\cong \mathbb Z^k.
\]
The monodromy \(T\) on the vanishing cohomology is the Coxeter
transformation of type \(A_k\); this is the standard Picard--Lefschetz
description of ADE hypersurface singularities \cite{Mi68}.

Choose a basis in which the Coxeter transformation has companion form
\[
        T(e_i)=e_{i+1}
        \quad (1\le i<k),
        \qquad
        T(e_k)=-(e_1+\cdots+e_k).
\]
Then the matrix of \(T\) has columns
\[
        e_2,e_3,\ldots,e_k,-(e_1+\cdots+e_k).
\]
The matrix \(T-\mathrm{id}\) is therefore
\[
        T-I
        =
        \begin{pmatrix}
        -1 & 0  & 0  & \cdots & -1\\
        1  & -1 & 0  & \cdots & -1\\
        0  & 1  & -1 & \ddots & -1\\
        \vdots & \vdots & \ddots & \ddots & \vdots\\
        0 & 0 & \cdots & 1 & -2
        \end{pmatrix}.
\]
The characteristic polynomial of the \(A_k\) Coxeter transformation is
\[
        \chi_T(t)=\frac{t^{k+1}-1}{t-1}
        =
        t^k+t^{k-1}+\cdots+t+1.
\]
Therefore
\[
        \det(T-I)
        =
        (-1)^k\chi_T(1)
        =
        (-1)^k(k+1).
\]
Thus the cokernel of \(T-I\) is finite of order \(k+1\).

To see that the cokernel is cyclic, observe that \(T-I\) has a
\((k-1)\times(k-1)\) minor of determinant \(\pm 1\): the upper-left
\((k-1)\times(k-1)\) block is lower triangular with diagonal entries
\(-1\).  Hence the first \(k-1\) Smith invariant factors are equal to \(1\).
Since the product of all Smith invariant factors is
\[
        |\det(T-I)|=k+1,
\]
the final Smith invariant factor is \(k+1\).  Therefore
\[
        \operatorname{coker}(T-\mathrm{id})
        \cong
        \mathbb Z/(k+1)\mathbb Z.
\]
This cokernel is finite, so it is equal to its torsion subgroup:
\[
        \operatorname{coker}(T-\mathrm{id})_{\mathrm{tors}}
        \cong
        \mathbb Z/(k+1)\mathbb Z.
\]
The monodromy realization theorem of
\cite[Theorem 1.3]{RahmanIntegralPerverseObstructions} identifies this group
with \(E\).  Thus
\[
        E
        \cong
        \operatorname{coker}(T-\mathrm{id})_{\mathrm{tors}}
        \cong
        \mathbb Z/(k+1)\mathbb Z.
\]

The determinant refinement also agrees:
\[
        |E|=k+1=|\det(T-\mathrm{id})|.
\]

\subsection{Transport status}

The \(A_k\) computation above is a local surface computation.  It determines
the birth, form, and six local realizations of the torsion group
\[
        E\cong \mathbb Z/(k+1)\mathbb Z.
\]
As in the \(A_1\) case, this local surface package naturally appears through
the support and pair sequences associated to the exceptional curve chain.

Let \(N\) be a resolution neighborhood of the exceptional divisor.  The
support sequence for the exceptional divisor and the pair sequence for
\((N,L)\) identify the local boundary quotient with
\[
        \operatorname{coker}(M)
        \cong
        \mathbb Z/(k+1)\mathbb Z.
\]
Thus the local surface package is naturally associated to the degree-two
local topology of the exceptional curve configuration.  It is not, by
itself, a degree-three Brauer class.

If an \(A_k\) singularity occurs inside a global surface, the global image of
its torsion package depends on how the exceptional chain maps into the
global lattice of exceptional and divisor classes.  If several singularities
are present, global relations among exceptional chains may kill linear
combinations of the local torsion packages.  Thus the transport status of
the \(A_k\) surface row is:
\[
        \text{support degree: }2\text{ locally, through the exceptional-chain
        support/pair sequence}.
\]
To enter the degree-three Brauer channel, the \(A_k\) package must be placed
inside a suitable global threefold mechanism or compared with a
higher-dimensional local package whose torsion appears in degree \(3\).

\subsection{Rational death}

Since
\[
        E\cong \mathbb Z/(k+1)\mathbb Z,
\]
we have
\[
        E\otimes_{\mathbb Z}\mathbb Q=0.
\]
Indeed, if \(e\in E\), then
\[
        (k+1)e=0.
\]
In the tensor product \(E\otimes_{\mathbb Z}\mathbb Q\), the integer
\(k+1\) acts invertibly on \(\mathbb Q\).  Therefore
\[
        e\otimes 1
        =
        e\otimes (k+1)\frac{1}{k+1}
        =
        (k+1)e\otimes\frac{1}{k+1}
        =
        0.
\]
Thus every element of \(E\) dies after rationalization.

This explains why the ordinary and dual middle-perversity intersection
complexes agree after tensoring with \(\mathbb Q\):
\[
        {}^pIC_X\mathbb Q
        \cong
        {}^p_+IC_X\mathbb Q.
\]
Before rationalization, however, the \(A_k\) torsion carries arithmetic and
geometric data: its invariant factor is \(k+1\), its prime support is the
set of primes dividing \(k+1\), its discriminant form is
\[
        q(a,b)=-\frac{k\,ab}{k+1}\bmod\mathbb Z,
\]
and it is visible simultaneously in the link, exceptional lattice, pair
sequence, and monodromy.

\subsection{Relation to the Coble boundary example}

The \(A_k\) family and the Coble boundary singularity are both cyclic
quotient surface singularities, but they play different roles in this paper.
The \(A_k\) singularity is the quotient
\[
        \frac1{k+1}(1,-1),
\]
and its local package is
\[
        E\cong\mathbb Z/(k+1).
\]
The Coble boundary singularity relevant to the corrected Benoist--Ottem
comparison is
\[
        \frac14(1,1),
\]
whose local package is
\[
        E\cong\mathbb Z/4.
\]
The special case
\[
        A_1=\frac12(1,1)
\]
appears as the index-two cover
\[
        A_1\longrightarrow \frac14(1,1).
\]
Thus the \(A_1\) row remains central to the Coble comparison, but as the
cover of the Coble singularity, not as the Coble boundary singularity itself.

This distinction is important because the Benoist--Ottem \(2\)-torsion
mechanism does not detect the full Coble group \(\mathbb Z/4\).  It detects
the order-two Bockstein shadow
\[
        2E\cong\mathbb Z/2
        \subset
        E\cong\mathbb Z/4.
\]
The \(A_k\) appendix therefore supplies the general cyclic ADE background,
while the \(\frac14(1,1)\) appendix supplies the corrected Enriques/Coble
boundary computation.

\subsection{Trajectory row}

The \(A_k\) surface trajectory can be summarized as follows.

\begin{center}
\begingroup
\scriptsize
\setlength{\tabcolsep}{1.8pt}
\renewcommand{\arraystretch}{1.15}
\begin{tabularx}{\textwidth}{|p{0.15\textwidth}|c|c|p{0.12\textwidth}|p{0.10\textwidth}|X|p{0.10\textwidth}|c|}
\hline
\textbf{Example} & \textbf{\(E\)} & \textbf{\(q\)} & \textbf{Local} & \textbf{Supp.} & \textbf{Global image} & \textbf{Br/res.} & \textbf{\(\mathbb Q\)}
\\
\hline
\(A_k\) surface &
\(\mathbb Z/(k+1)\) &
\(-\frac{k}{k+1}\) &
six agree &
deg.\ \(2\) &
depends on global exceptional-chain relations &
loc.\ \(\varnothing\) &
\(0\)
\\
\hline
\end{tabularx}
\endgroup
\end{center}

Here the entry
\[
        q=-\frac{k}{k+1}
\]
means that, for the generator \(\bar g\) chosen above,
\[
        q(\bar g,\bar g)
        =
        -\frac{k}{k+1}\bmod\mathbb Z.
\]
The phrase ``six agree'' means that the following realizations all give
\[
        \mathbb Z/(k+1)\mathbb Z:
\]
the perverse discrepancy, torsion-sensitive truncation, link torsion,
exceptional lattice discriminant group, pair-sequence quotient, and
monodromy cokernel.  The Brauer/residue column is marked local
\(\varnothing\) because the isolated surface germ alone does not produce a
degree-three Brauer or unramified class.

\section{The \texorpdfstring{\(D_4\)}{D4} surface trajectory}
\label{app:D4surfacetrajectory}

This appendix works out the torsion trajectory for the \(D_4\) surface
singularity.  The purpose of this example is different from the \(A_1\) and
\(A_k\) examples.  The \(A_k\) family shows how cyclic torsion
\[
        \mathbb Z/(k+1)\mathbb Z
\]
is born from a chain of exceptional curves.  The \(D_4\) example shows why
one must track the discriminant form, not only the underlying finite group.
Here
\[
        E\cong \mathbb Z/2\mathbb Z\oplus \mathbb Z/2\mathbb Z,
\]
but the pairing on this group carries additional information.

This form-sensitivity should be compared with the Coble boundary example.
There, the issue is not only the finite group and its form, but also the
fact that a global \(2\)-torsion mechanism may select a subgroup
\[
        2E\subset E.
\]
Thus \(D_4\) shows that the group alone is not enough, while the Coble
\(\frac14(1,1)\) example shows that even the full group may not be the
piece seen by a particular global obstruction theory.

\subsection{Local model}

We use the hypersurface model
\[
        (X,0)=\{x^2+y^3+yz^2=0\}\subset(\mathbb C^3,0).
\]
This is the rational double point of type \(D_4\).  It is a normal surface
singularity, an isolated hypersurface surface singularity, and a simple
surface singularity in the ADE classification.  The ADE surface
singularities and their Milnor fibers are classical; see Milnor
\cite{Mi68} for the general isolated hypersurface singularity background and
Mumford \cite{Mu61} for the topology of normal surface singularity links.

Let
\[
        f(x,y,z)=x^2+y^3+yz^2.
\]
The partial derivatives are
\[
        \frac{\partial f}{\partial x}=2x,
        \qquad
        \frac{\partial f}{\partial y}=3y^2+z^2,
        \qquad
        \frac{\partial f}{\partial z}=2yz.
\]
If all three partial derivatives vanish, then
\[
        x=0,
        \qquad
        3y^2+z^2=0,
        \qquad
        yz=0.
\]
The equation \(yz=0\) gives either \(y=0\) or \(z=0\).  If \(y=0\), then
\(z^2=0\), so \(z=0\).  If \(z=0\), then \(3y^2=0\), so \(y=0\).  Hence the
only singular point is
\[
        x=y=z=0.
\]
Thus \((X,0)\) is an isolated hypersurface surface singularity.

Let
\[
        U:=X\setminus\{0\},
        \qquad
        j:U\hookrightarrow X,
        \qquad
        i:\{0\}\hookrightarrow X.
\]
Since \(X\) has complex dimension \(2\), the shifted constant sheaf on the
smooth locus is
\[
        \mathbb Z_U[2].
\]
The ordinary and dual middle-perversity intersection complexes are
\[
        {}^pIC_X\mathbb Z
        =
        {}^pj_{!*}\mathbb Z_U[2],
        \qquad
        {}^p_+IC_X\mathbb Z
        =
        {}^p_+j_{!*}\mathbb Z_U[2],
\]
as in \cite{RahmanIntegralPerverseObstructions}.

\subsection{The link and its cohomology}

Let \(L\) be the link of the \(D_4\) singularity.  For a normal surface
singularity, the link is the boundary of a sufficiently small resolution
neighborhood of the exceptional divisor \cite{Mu61}.  In the \(D_4\) case,
the minimal resolution has exceptional divisor whose dual graph is the
\(D_4\) Dynkin diagram.  Equivalently, the link is the oriented plumbed
three-manifold associated to the negative \(D_4\) plumbing graph.

The exceptional curves consist of one central curve meeting three outer
curves:
\[
        C_0,\ C_1,\ C_2,\ C_3,
\]
with intersections
\[
        C_0^2=C_1^2=C_2^2=C_3^2=-2,
\]
\[
        C_0\cdot C_i=1 \quad (i=1,2,3),
\]
and
\[
        C_i\cdot C_j=0
        \quad
        \text{for } i,j\in\{1,2,3\},\ i\ne j.
\]
Thus the intersection matrix in the ordered basis
\[
        (C_0,C_1,C_2,C_3)
\]
is
\[
        M=
        \begin{pmatrix}
        -2 & 1  & 1  & 1\\
        1  & -2 & 0  & 0\\
        1  & 0  & -2 & 0\\
        1  & 0  & 0  & -2
        \end{pmatrix}.
\]
This is the negative \(D_4\) Cartan matrix.  Its determinant is
\[
        \det(M)=4.
\]

We compute the homology of the link using the pair sequence of the resolution
neighborhood.  Let \(N\) be a resolution neighborhood of the exceptional
divisor, so that
\[
        \partial N=L.
\]
The neighborhood \(N\) deformation retracts onto the exceptional divisor.
Since the exceptional divisor is a tree of rational curves, we have
\[
        H_1(N,\mathbb Z)=0,
        \qquad
        H_2(N,\mathbb Z)\cong \mathbb Z^4.
\]
The long exact sequence of the pair \((N,L)\) contains
\[
        H_2(N,\mathbb Z)
        \longrightarrow
        H_2(N,L;\mathbb Z)
        \longrightarrow
        H_1(L,\mathbb Z)
        \longrightarrow
        H_1(N,\mathbb Z).
\]
By Poincare--Lefschetz duality,
\[
        H_2(N,L;\mathbb Z)\cong H^2(N,\mathbb Z).
\]
Since \(H_1(N,\mathbb Z)=0\), the universal coefficient theorem gives
\[
        H^2(N,\mathbb Z)\cong
        \operatorname{Hom}(H_2(N,\mathbb Z),\mathbb Z)
        \cong
        \mathbb Z^4.
\]
Under these identifications, the map
\[
        H_2(N,\mathbb Z)
        \longrightarrow
        H_2(N,L;\mathbb Z)
\]
is represented by the intersection matrix \(M\).  Hence
\[
        H_1(L,\mathbb Z)
        \cong
        \operatorname{coker}(M).
\]
This is the standard resolution-neighborhood computation of link torsion;
see \cite{GoreskySiegel83} and
\cite[§4]{RahmanIntegralPerverseObstructions}.

The discriminant group of the \(D_n\) root lattice is standard:
\[
        \mathbb Z/4\mathbb Z
        \quad\text{if } n \text{ is odd},
        \qquad
        \mathbb Z/2\mathbb Z\oplus\mathbb Z/2\mathbb Z
        \quad\text{if } n \text{ is even}.
\]
For \(D_4\), this gives
\[
        \operatorname{coker}(M)
        \cong
        \mathbb Z/2\mathbb Z\oplus \mathbb Z/2\mathbb Z.
\]
Equivalently, the Smith normal form of the negative \(D_4\) Cartan matrix is
\[
        \operatorname{diag}(1,1,2,2).
\]
This is the standard discriminant group computation for \(D_4\); see
Nikulin's treatment of discriminant forms of integral lattices
\cite{Nikulin80}.  Therefore
\[
        H_1(L,\mathbb Z)
        \cong
        \mathbb Z/2\mathbb Z\oplus \mathbb Z/2\mathbb Z.
\]

Since \(L\) is a closed oriented three-manifold, the universal coefficient
theorem gives
\[
        H^2(L,\mathbb Z)_{\mathrm{tors}}
        \cong
        H_1(L,\mathbb Z)_{\mathrm{tors}}.
\]
Therefore
\[
        H^2(L,\mathbb Z)_{\mathrm{tors}}
        \cong
        \mathbb Z/2\mathbb Z\oplus \mathbb Z/2\mathbb Z.
\]

\subsection{Birth of torsion}

By the local realization theorem of
\cite{RahmanIntegralPerverseObstructions}, for a normal surface singularity
one has
\[
        E
        \cong
        H^2(L,\mathbb Z)_{\mathrm{tors}}.
\]
Using the computation above, the \(D_4\) singularity has
\[
        E
        \cong
        \mathbb Z/2\mathbb Z\oplus \mathbb Z/2\mathbb Z.
\]
Thus the invariant factor decomposition is
\[
        E\cong (\mathbb Z/2\mathbb Z)^2.
\]
Its order is
\[
        |E|=4.
\]
Its prime decomposition is supported only at the prime
\[
        2.
\]
The \(D_4\) singularity therefore produces noncyclic \(2\)-torsion.

\subsection{The discriminant form}

We now compute the discriminant form.  Let
\[
        \Lambda
        =
        \mathbb Z\langle C_0,C_1,C_2,C_3\rangle
\]
with intersection matrix
\[
        M=
        \begin{pmatrix}
        -2 & 1  & 1  & 1\\
        1  & -2 & 0  & 0\\
        1  & 0  & -2 & 0\\
        1  & 0  & 0  & -2
        \end{pmatrix}.
\]
The discriminant group is
\[
        A_\Lambda=\Lambda^\vee/\Lambda.
\]
We choose the following two elements of \(\Lambda^\vee\):
\[
        u:=\frac{C_1-C_2}{2},
        \qquad
        v:=\frac{C_1-C_3}{2}.
\]
We first verify that \(u,v\in\Lambda^\vee\).  For \(u\), we compute:
\[
        (u,C_0)
        =
        \frac{(C_1,C_0)-(C_2,C_0)}{2}
        =
        \frac{1-1}{2}
        =
        0,
\]
\[
        (u,C_1)
        =
        \frac{(C_1,C_1)-(C_2,C_1)}{2}
        =
        \frac{-2-0}{2}
        =
        -1,
\]
\[
        (u,C_2)
        =
        \frac{(C_1,C_2)-(C_2,C_2)}{2}
        =
        \frac{0-(-2)}{2}
        =
        1,
\]
and
\[
        (u,C_3)
        =
        \frac{(C_1,C_3)-(C_2,C_3)}{2}
        =
        0.
\]
All these numbers are integers, so \(u\in\Lambda^\vee\).  The same
calculation shows that \(v\in\Lambda^\vee\).

Let
\[
        \bar u:=u+\Lambda,
        \qquad
        \bar v:=v+\Lambda
\]
be their classes in \(\Lambda^\vee/\Lambda\).  Since
\[
        2u=C_1-C_2\in\Lambda,
        \qquad
        2v=C_1-C_3\in\Lambda,
\]
both \(\bar u\) and \(\bar v\) have order dividing \(2\).  They are
nonzero, and they generate a subgroup isomorphic to
\[
        \mathbb Z/2\mathbb Z\oplus \mathbb Z/2\mathbb Z.
\]
Because the full discriminant group has order \(4\), this subgroup is all of
\[
        \Lambda^\vee/\Lambda.
\]
Thus
\[
        \Lambda^\vee/\Lambda
        =
        \langle \bar u,\bar v\rangle
        \cong
        \mathbb Z/2\mathbb Z\oplus \mathbb Z/2\mathbb Z.
\]

The discriminant pairing is
\[
        q_\Lambda:
        \Lambda^\vee/\Lambda\times \Lambda^\vee/\Lambda
        \longrightarrow
        \mathbb Q/\mathbb Z,
        \qquad
        q_\Lambda(x+\Lambda,y+\Lambda)=(x,y)\bmod\mathbb Z.
\]
We compute the pairing on the chosen generators.

First,
\[
        (u,u)
        =
        \left(\frac{C_1-C_2}{2},\frac{C_1-C_2}{2}\right)
        =
        \frac{(C_1,C_1)-2(C_1,C_2)+(C_2,C_2)}{4}.
\]
Since
\[
        (C_1,C_1)=-2,
        \qquad
        (C_2,C_2)=-2,
        \qquad
        (C_1,C_2)=0,
\]
we get
\[
        (u,u)=\frac{-2+0-2}{4}=-1.
\]
Therefore
\[
        q_\Lambda(\bar u,\bar u)=0\in\mathbb Q/\mathbb Z.
\]
Similarly,
\[
        q_\Lambda(\bar v,\bar v)=0.
\]

Next,
\[
        (u,v)
        =
        \left(\frac{C_1-C_2}{2},\frac{C_1-C_3}{2}\right)
        =
        \frac{
        (C_1,C_1)-(C_1,C_3)-(C_2,C_1)+(C_2,C_3)
        }{4}.
\]
Using
\[
        (C_1,C_1)=-2,
        \qquad
        (C_1,C_3)=0,
        \qquad
        (C_2,C_1)=0,
        \qquad
        (C_2,C_3)=0,
\]
we find
\[
        (u,v)=-\frac12.
\]
Therefore
\[
        q_\Lambda(\bar u,\bar v)
        =
        -\frac12\bmod\mathbb Z.
\]
With respect to the basis \((\bar u,\bar v)\), the discriminant pairing is
represented by
\[
        \begin{pmatrix}
        0 & -\frac12\\
        -\frac12 & 0
        \end{pmatrix}
        \quad
        \text{in } \mathbb Q/\mathbb Z.
\]
Equivalently, since \(-\frac12=\frac12\) in \(\mathbb Q/\mathbb Z\), the
same pairing may be written
\[
        \begin{pmatrix}
        0 & \frac12\\
        \frac12 & 0
        \end{pmatrix}.
\]

Thus the local discriminant package is
\[
        (E,q)
        =
        \left(
        (\mathbb Z/2\mathbb Z)^2,
        \begin{pmatrix}
        0 & -\frac12\\
        -\frac12 & 0
        \end{pmatrix}
        \right)
\]
with respect to the chosen generators \(\bar u,\bar v\).

This is the first example in the paper where the form is essential.  The
underlying group
\[
        (\mathbb Z/2\mathbb Z)^2
\]
does not by itself remember the off-diagonal pairing.

\subsection{Perverse and torsion-sensitive realization}

Let
\[
        K:=Rj_*\mathbb Z_U[2].
\]
For an isolated surface singularity with link \(L\), the local stalk
calculation gives
\[
        i^*K
        \cong
        R\Gamma(L,\mathbb Z)[2],
\]
and hence
\[
        H^m(i^*K)
        \cong
        H^{m+2}(L,\mathbb Z).
\]
This calculation is used in
\cite[Lemma 3.3]{RahmanIntegralPerverseObstructions}.

For the \(D_4\) link, we computed
\[
        H^2(L,\mathbb Z)_{\mathrm{tors}}
        \cong
        (\mathbb Z/2\mathbb Z)^2.
\]
Thus the degree-zero stalk of the local model contains
\[
        H^0(i^*K)
        \cong
        H^2(L,\mathbb Z),
\]
whose torsion subgroup is
\[
        (\mathbb Z/2\mathbb Z)^2.
\]

The ordinary and dual middle-perversity intersection complexes are obtained
from this local model by imposing the ordinary and dual codimension-two
point-stratum conditions.  Over \(\mathbb Z\), these differ in how they treat
torsion in the critical degree; this is the integral ordinary/dual
middle-perversity distinction of BBD \cite[Complement 3.3]{BBD82}.  In
Friedman's torsion-sensitive language, the dual middle-perversity package
retains the relevant critical torsion, whereas the ordinary middle package
removes it \cite{FriedmanGenIH,FriedmanBook20,FriedmanTsInv}.

Therefore
\[
        H^0({}^p_+IC_X\mathbb Z)_0
        \cong
        (\mathbb Z/2\mathbb Z)^2.
\]
By definition,
\[
        E=H^0({}^p_+IC_X\mathbb Z)_0,
\]
so
\[
        E\cong(\mathbb Z/2\mathbb Z)^2.
\]
The discrepancy triangle is
\[
        {}^pIC_X\mathbb Z
        \longrightarrow
        {}^p_+IC_X\mathbb Z
        \longrightarrow
        \bigl((\mathbb Z/2\mathbb Z)^2\bigr)[1]
        \longrightarrow .
\]
The general statement that this perverse discrepancy agrees with link
torsion and the exceptional discriminant group is
\cite{RahmanIntegralPerverseObstructions}.

\subsection{Resolution lattice realization}

The resolution-lattice realization has already appeared in the computation
of the link, but we record it separately as one of the six local models.

The exceptional lattice is
\[
        \Lambda
        =
        \mathbb Z\langle C_0,C_1,C_2,C_3\rangle,
\]
with intersection matrix
\[
        M=
        \begin{pmatrix}
        -2 & 1  & 1  & 1\\
        1  & -2 & 0  & 0\\
        1  & 0  & -2 & 0\\
        1  & 0  & 0  & -2
        \end{pmatrix}.
\]
The map
\[
        \Lambda\longrightarrow \Lambda^\vee
\]
is represented by \(M\).  Therefore
\[
        \Lambda^\vee/\Lambda
        \cong
        \operatorname{coker}(M).
\]
The Smith normal form of \(M\) is
\[
        \operatorname{diag}(1,1,2,2).
\]
Hence
\[
        \Lambda^\vee/\Lambda
        \cong
        \mathbb Z/2\mathbb Z\oplus\mathbb Z/2\mathbb Z.
\]
By the local realization theorem,
\[
        E
        \cong
        \Lambda^\vee/\Lambda.
\]
Thus
\[
        E
        \cong
        (\mathbb Z/2\mathbb Z)^2.
\]
This is the \(D_4\) case of the resolution-lattice theorem in
\cite[Theorem 1.2]{RahmanIntegralPerverseObstructions}.

\subsection{Pair-sequence realization}

Let \(N\) be a resolution neighborhood of the \(D_4\) exceptional divisor,
and let
\[
        L=\partial N
\]
be its boundary link.  The long exact sequence of the pair \((N,L)\)
contains
\[
        H_2(N,\mathbb Z)
        \longrightarrow
        H_2(N,L;\mathbb Z)
        \longrightarrow
        H_1(L,\mathbb Z)
        \longrightarrow
        H_1(N,\mathbb Z).
\]
Since the exceptional divisor is a tree of rational curves, \(N\) deformation
retracts onto a simply connected two-dimensional complex, and hence
\[
        H_1(N,\mathbb Z)=0.
\]
Also
\[
        H_2(N,\mathbb Z)
        \cong
        \mathbb Z^4,
\]
generated by \(C_0,C_1,C_2,C_3\).

By Poincare--Lefschetz duality,
\[
        H_2(N,L;\mathbb Z)
        \cong
        H^2(N,\mathbb Z).
\]
Since \(H_1(N,\mathbb Z)=0\), the universal coefficient theorem gives
\[
        H^2(N,\mathbb Z)
        \cong
        \operatorname{Hom}(H_2(N,\mathbb Z),\mathbb Z)
        \cong
        \mathbb Z^4.
\]
Under these identifications, the map
\[
        H_2(N,\mathbb Z)
        \longrightarrow
        H_2(N,L;\mathbb Z)
\]
is represented by the intersection matrix \(M\).  Thus
\[
        H_1(L,\mathbb Z)
        \cong
        \operatorname{coker}(M)
        \cong
        \mathbb Z/2\mathbb Z\oplus\mathbb Z/2\mathbb Z.
\]
Using the universal coefficient theorem for \(L\), this gives
\[
        H^2(L,\mathbb Z)_{\mathrm{tors}}
        \cong
        \mathbb Z/2\mathbb Z\oplus\mathbb Z/2\mathbb Z.
\]
Therefore the pair sequence recovers the same local group
\[
        E\cong(\mathbb Z/2\mathbb Z)^2.
\]

This is the \(D_4\) specialization of the resolution-neighborhood
calculation in \cite[§4]{RahmanIntegralPerverseObstructions} and the
linking-pairing formalism of \cite{GoreskySiegel83}.

\subsection{Monodromy realization}

The \(D_4\) singularity is an isolated hypersurface surface singularity, so
it has a Milnor fibration.  The Milnor number is
\[
        \mu=4.
\]
Thus the middle vanishing cohomology is a free abelian group of rank \(4\):
\[
        H^2_{\mathrm{van}}(F,\mathbb Z)\cong \mathbb Z^4.
\]
The monodromy \(T\) is the Coxeter transformation of type \(D_4\), acting on
the \(D_4\) root lattice.  This is the standard ADE Picard--Lefschetz
description of the monodromy of simple hypersurface singularities
\cite{Mi68}.

Choose the simple-root basis ordered as
\[
        (C_0,C_1,C_2,C_3),
\]
where \(C_0\) is the central node.  Let \(A\) be the positive \(D_4\) Cartan
matrix:
\[
        A=
        \begin{pmatrix}
        2 & -1 & -1 & -1\\
        -1 & 2 & 0 & 0\\
        -1 & 0 & 2 & 0\\
        -1 & 0 & 0 & 2
        \end{pmatrix}.
\]
The simple reflection \(s_i\) acts on the simple-root basis by
\[
        s_i(\alpha_j)=\alpha_j-A_{ij}\alpha_i.
\]
Taking the Coxeter element
\[
        T=s_0s_1s_2s_3,
\]
one obtains the matrix
\[
        T=
        \begin{pmatrix}
        2 & -1 & -1 & -1\\
        1 & -1 & 0 & 0\\
        1 & 0 & -1 & 0\\
        1 & 0 & 0 & -1
        \end{pmatrix}.
\]
Therefore
\[
        T-I
        =
        \begin{pmatrix}
        1 & -1 & -1 & -1\\
        1 & -2 & 0 & 0\\
        1 & 0 & -2 & 0\\
        1 & 0 & 0 & -2
        \end{pmatrix}.
\]
A direct Smith normal form computation gives
\[
        \operatorname{SNF}(T-I)
        =
        \operatorname{diag}(1,1,2,2).
\]
We spell out why this determines the cokernel.  The Smith normal form of an
integral matrix records invariant factors for the cokernel of the associated
homomorphism of free abelian groups.  Thus
\[
        \operatorname{coker}(T-\mathrm{id})
        \cong
        \mathbb Z/2\mathbb Z\oplus\mathbb Z/2\mathbb Z.
\]
This group is finite, so it is equal to its torsion subgroup:
\[
        \operatorname{coker}(T-\mathrm{id})_{\mathrm{tors}}
        \cong
        \mathbb Z/2\mathbb Z\oplus\mathbb Z/2\mathbb Z.
\]
The monodromy realization theorem of
\cite[Theorem 1.3]{RahmanIntegralPerverseObstructions} gives
\[
        E
        \cong
        \operatorname{coker}(T-\mathrm{id})_{\mathrm{tors}},
\]
and therefore
\[
        E
        \cong
        (\mathbb Z/2\mathbb Z)^2.
\]
The determinant check is also consistent:
\[
        |\det(T-\mathrm{id})|=4=|E|.
\]

\subsection{Transport status}

The \(D_4\) computation is a local surface computation.  It determines the
birth, form, and six local realizations of the torsion group
\[
        E\cong(\mathbb Z/2\mathbb Z)^2.
\]
As in the \(A_1\) and \(A_k\) cases, the local surface package naturally
appears through the degree-two topology of the exceptional divisor and its
boundary link.

Let \(N\) be a resolution neighborhood of the \(D_4\) exceptional divisor.
The pair sequence identifies the local boundary quotient with
\[
        \operatorname{coker}(M)
        \cong
        (\mathbb Z/2\mathbb Z)^2.
\]
Thus the local support status is:
\[
        \text{support degree: }2\text{ locally, through the exceptional
        divisor support/pair sequence}.
\]
The local \(D_4\) surface germ alone does not define a degree-three Brauer
class.  To enter the Brauer or unramified cohomology channel, this package
must be placed inside a global model with a degree-three transport
mechanism.

This example is nevertheless essential for the global program because it
shows that a global transport theory must track the pairing.  The group
\[
        (\mathbb Z/2\mathbb Z)^2
\]
can arise in more than one way.  The \(D_4\) discriminant form is not the
same datum as merely having two independent \(A_1\)-type \(\mathbb Z/2\)
summands.  The off-diagonal pairing
\[
        q(\bar u,\bar v)=-\frac12
\]
is part of the local package.

This form-sensitive warning is parallel to, but distinct from, the Coble
boundary warning.  In the Coble case, the issue is that the full package is
\[
        E\cong\mathbb Z/4,
\]
while the Benoist--Ottem mechanism sees only
\[
        2E\cong\mathbb Z/2.
\]
In the \(D_4\) case, the issue is that the full package is already
\[
        E\cong(\mathbb Z/2)^2,
\]
but its finite pairing is essential.

\subsection{Rational death}

Since
\[
        E\cong(\mathbb Z/2\mathbb Z)^2,
\]
we have
\[
        E\otimes_{\mathbb Z}\mathbb Q=0.
\]
Indeed, every element \(e\in E\) satisfies
\[
        2e=0.
\]
Therefore
\[
        e\otimes 1
        =
        e\otimes 2\cdot\frac12
        =
        2e\otimes\frac12
        =
        0.
\]
Thus the group \(E\) vanishes after rationalization.

Before it dies, however, it carries the form-sensitive information
\[
        q=
        \begin{pmatrix}
        0 & -\frac12\\
        -\frac12 & 0
        \end{pmatrix}.
\]
This information is invisible if one records only the abstract group
\[
        (\mathbb Z/2\mathbb Z)^2.
\]
It is also invisible after tensoring with \(\mathbb Q\).  The \(D_4\)
trajectory therefore provides the first warning that the torsion trajectory
must record the pair
\[
        (E,q),
\]
not only the group \(E\).

\subsection{Trajectory row}

The \(D_4\) surface trajectory can be summarized as follows.

\begin{center}
\begingroup
\scriptsize
\setlength{\tabcolsep}{1.8pt}
\renewcommand{\arraystretch}{1.15}
\begin{tabularx}{\textwidth}{|p{0.15\textwidth}|c|c|p{0.12\textwidth}|p{0.10\textwidth}|X|p{0.10\textwidth}|c|}
\hline
\textbf{Example} & \textbf{\(E\)} & \textbf{\(q\)} & \textbf{Local} & \textbf{Supp.} & \textbf{Global image} & \textbf{Br/res.} & \textbf{\(\mathbb Q\)}
\\
\hline
\(D_4\) surface &
\((\mathbb Z/2)^2\) &
\(\begin{smallmatrix}0&-\frac12\\[-1mm]-\frac12&0\end{smallmatrix}\) &
six agree &
deg.\ \(2\) &
depends on global exceptional-divisor relations and form data &
loc.\ \(\varnothing\) &
\(0\)
\\
\hline
\end{tabularx}
\endgroup
\end{center}

Here ``six agree'' means that the perverse discrepancy,
torsion-sensitive truncation, link torsion, exceptional lattice
discriminant group, pair-sequence quotient, and monodromy cokernel all give
\[
        E\cong
        \mathbb Z/2\mathbb Z\oplus\mathbb Z/2\mathbb Z.
\]
The support degree is listed as \(2\) because this is a local surface
singularity and the torsion appears through the degree-two topology of the
exceptional divisor and its boundary.  The Brauer/residue column is marked
local \(\varnothing\) because the isolated surface germ alone does not
produce a degree-three Brauer or unramified class.

\section{The \texorpdfstring{\(E_8\)}{E8} surface trajectory}
\label{app:E8surfacetrajectory}

This appendix works out the torsion trajectory for the \(E_8\) surface
singularity.  The role of this example is to provide the basic null-control
case.  The singularity is highly nontrivial, but its exceptional lattice is
unimodular.  Consequently,
\[
        E=0.
\]
Thus every station in the local finite-torsion trajectory is trivial.  This
example is important because it shows that the presence of a singularity
alone does not imply the presence of an integral perverse obstruction.

The \(E_8\) row should also be read against the Coble boundary example.  For
the Coble singularity \(\frac14(1,1)\), the local lattice is non-unimodular
and produces
\[
        E\cong\mathbb Z/4.
\]
For \(E_8\), the lattice is unimodular and produces
\[
        E=0.
\]
Thus the local invariant \(E\) detects discriminant torsion, not singularity
complexity.

\subsection{Local model}

We use the hypersurface model
\[
        (X,0)=\{x^2+y^3+z^5=0\}\subset(\mathbb C^3,0).
\]
This is the rational double point of type \(E_8\).  It is a normal surface
singularity, an isolated hypersurface surface singularity, and a simple
surface singularity in the ADE classification.  The topology of normal
surface singularity links is classical; see Mumford \cite{Mu61}.  The
Milnor fibration and monodromy of isolated hypersurface singularities are
treated in \cite{Mi68}.

Let
\[
        f(x,y,z)=x^2+y^3+z^5.
\]
The partial derivatives are
\[
        \frac{\partial f}{\partial x}=2x,
        \qquad
        \frac{\partial f}{\partial y}=3y^2,
        \qquad
        \frac{\partial f}{\partial z}=5z^4.
\]
All three partial derivatives vanish simultaneously only when
\[
        x=y=z=0.
\]
Thus \(0\) is the unique singular point of the germ.

Let
\[
        U:=X\setminus\{0\},
        \qquad
        j:U\hookrightarrow X,
        \qquad
        i:\{0\}\hookrightarrow X.
\]
Since \(X\) is a complex surface, the shifted constant sheaf on the smooth
locus is
\[
        \mathbb Z_U[2].
\]
The ordinary and dual middle-perversity intersection complexes are
\[
        {}^pIC_X\mathbb Z
        =
        {}^pj_{!*}\mathbb Z_U[2],
        \qquad
        {}^p_+IC_X\mathbb Z
        =
        {}^p_+j_{!*}\mathbb Z_U[2].
\]
These conventions are those used in
\cite{RahmanIntegralPerverseObstructions}.

\subsection{The link and its cohomology}

Let \(L\) be the link of the \(E_8\) singularity.  The link is the boundary
of a sufficiently small resolution neighborhood of the exceptional divisor.
For the \(E_8\) singularity, the exceptional divisor in the minimal
resolution is a tree of eight rational curves whose dual graph is the
\(E_8\) Dynkin diagram.  The corresponding link is the oriented plumbed
three-manifold associated to the negative \(E_8\) plumbing graph.

The exceptional lattice is the negative definite \(E_8\) lattice.  The
positive \(E_8\) root lattice is even, positive definite, and unimodular.
With the geometric convention used in this paper, the intersection lattice is
its negative.  Changing the sign of the bilinear form does not change
unimodularity.  Hence the negative \(E_8\) exceptional lattice is also
unimodular.  The standard lattice-theoretic facts about the \(E_8\) lattice,
including its unimodularity, are recorded for example in Nikulin's treatment
of integral bilinear forms \cite{Nikulin80}.

We now compute the torsion in the link using the resolution-neighborhood
pair sequence.  Let \(N\) be a resolution neighborhood of the exceptional
divisor, so that
\[
        \partial N=L.
\]
The neighborhood \(N\) deformation retracts onto the exceptional divisor.
Since the exceptional divisor is a tree of rational curves, one has
\[
        H_1(N,\mathbb Z)=0,
        \qquad
        H_2(N,\mathbb Z)\cong \mathbb Z^8.
\]
The long exact sequence of the pair \((N,L)\) contains
\[
        H_2(N,\mathbb Z)
        \longrightarrow
        H_2(N,L;\mathbb Z)
        \longrightarrow
        H_1(L,\mathbb Z)
        \longrightarrow
        H_1(N,\mathbb Z).
\]
By Poincare--Lefschetz duality,
\[
        H_2(N,L;\mathbb Z)
        \cong
        H^2(N,\mathbb Z).
\]
Since \(H_1(N,\mathbb Z)=0\), the universal coefficient theorem gives
\[
        H^2(N,\mathbb Z)
        \cong
        \operatorname{Hom}(H_2(N,\mathbb Z),\mathbb Z)
        \cong
        \mathbb Z^8.
\]
Under these identifications, the map
\[
        H_2(N,\mathbb Z)
        \longrightarrow
        H_2(N,L;\mathbb Z)
\]
is represented by the \(E_8\) intersection matrix.  Since this matrix is
unimodular, the induced map
\[
        \mathbb Z^8\longrightarrow \mathbb Z^8
\]
has cokernel zero.  Therefore
\[
        H_1(L,\mathbb Z)=0.
\]

The link \(L\) is a closed oriented three-manifold.  Therefore
\[
        H_0(L,\mathbb Z)\cong \mathbb Z,
        \qquad
        H_3(L,\mathbb Z)\cong \mathbb Z.
\]
Since \(H_1(L,\mathbb Z)=0\), Poincare duality gives
\[
        H_2(L,\mathbb Z)=0.
\]
Thus \(L\) is an integral homology sphere:
\[
        H_m(L,\mathbb Z)
        \cong
        H_m(S^3,\mathbb Z)
        \quad
        \text{for all }m.
\]
The link is in fact the Poincare homology sphere, but the only fact needed
for the present computation is that it has no torsion homology.

Using the universal coefficient theorem, we obtain
\[
        H^0(L,\mathbb Z)\cong \mathbb Z,
        \qquad
        H^1(L,\mathbb Z)=0,
        \qquad
        H^2(L,\mathbb Z)=0,
        \qquad
        H^3(L,\mathbb Z)\cong \mathbb Z.
\]
In particular,
\[
        H^2(L,\mathbb Z)_{\mathrm{tors}}=0.
\]

\subsection{Birth of torsion}

For a normal surface singularity, the local realization theorem of
\cite{RahmanIntegralPerverseObstructions} gives
\[
        E
        \cong
        H^2(L,\mathbb Z)_{\mathrm{tors}}.
\]
For the \(E_8\) singularity, the link computation above gives
\[
        H^2(L,\mathbb Z)_{\mathrm{tors}}=0.
\]
Therefore
\[
        E=0.
\]

Thus the torsion trajectory has no local birth.  The invariant factor
decomposition is empty, the order is
\[
        |E|=1,
\]
and the prime decomposition is empty.  This is the null-control example in
the ADE series.

\subsection{The discriminant form}

Let
\[
        \Lambda
\]
be the exceptional lattice of the \(E_8\) resolution.  It is the negative
definite \(E_8\) lattice.  Since the \(E_8\) lattice is unimodular, the
natural map
\[
        \Lambda\longrightarrow \Lambda^\vee
\]
is an isomorphism.  Therefore
\[
        \Lambda^\vee/\Lambda=0.
\]
The discriminant pairing is a pairing on the zero group:
\[
        q_\Lambda:0\times 0\longrightarrow \mathbb Q/\mathbb Z.
\]
Thus the local discriminant package is
\[
        (E,q)=(0,0).
\]

This is the sharp contrast with the \(A_k\), \(D_4\), and Coble boundary
examples.  The singularity has a nontrivial exceptional configuration, but
the exceptional lattice is unimodular, so there is no finite discriminant
group and hence no local integral obstruction group \(E\).  The
lattice-theoretic fact that the \(E_8\) lattice is unimodular is classical;
see \cite{Nikulin80}.  The identification of \(\Lambda^\vee/\Lambda\) with
the local perverse obstruction group is
\cite[Theorem 1.2]{RahmanIntegralPerverseObstructions}.

\subsection{Perverse and torsion-sensitive realization}

Let
\[
        K:=Rj_*\mathbb Z_U[2].
\]
For an isolated surface singularity with link \(L\), the local stalk
calculation gives
\[
        i^*K
        \cong
        R\Gamma(L,\mathbb Z)[2],
\]
and therefore
\[
        H^m(i^*K)
        \cong
        H^{m+2}(L,\mathbb Z).
\]
This local stalk computation is used in
\cite[Lemma 3.3]{RahmanIntegralPerverseObstructions}.

For the \(E_8\) link, we computed
\[
        H^2(L,\mathbb Z)=0.
\]
Thus the degree-zero stalk contribution relevant to the local obstruction is
zero:
\[
        H^0(i^*K)
        \cong
        H^2(L,\mathbb Z)
        =
        0.
\]
In the surface case, the discrepancy between the ordinary and dual
middle-perversity extensions is the torsion subgroup of
\(H^2(L,\mathbb Z)\); see
\cite[Proposition 3.4]{RahmanIntegralPerverseObstructions}.  Since this
group is zero for the \(E_8\) link, there is no perverse obstruction:
\[
        E=0.
\]

Equivalently, the ordinary and dual middle-perversity intersection complexes
coincide integrally in this example:
\[
        {}^pIC_X\mathbb Z
        \cong
        {}^p_+IC_X\mathbb Z.
\]
This conclusion follows because the cone of the natural morphism
\[
        {}^pIC_X\mathbb Z
        \longrightarrow
        {}^p_+IC_X\mathbb Z
\]
is \(E[1]\), and here \(E=0\).  The local discrepancy triangle is therefore
trivial:
\[
        {}^pIC_X\mathbb Z
        \xrightarrow{\sim}
        {}^p_+IC_X\mathbb Z.
\]

In torsion-sensitive truncation language, there is no critical torsion to
retain.  Thus the ordinary and dual torsion-sensitive truncations agree at
the critical degree.  This is the \(E_8\) null case of the
torsion-sensitive realization described in
\cite{FriedmanGenIH,FriedmanBook20,FriedmanTsInv} and specialized to normal
surface singularities in \cite{RahmanIntegralPerverseObstructions}.

\subsection{Resolution lattice realization}

The minimal resolution has eight exceptional curves whose dual graph is the
\(E_8\) Dynkin diagram.  With the geometric sign convention, the intersection
matrix is the negative \(E_8\) Cartan matrix.  Denote it by
\[
        M_{E_8}.
\]
The determinant of the \(E_8\) Cartan matrix is
\[
        1.
\]
Therefore
\[
        |\det(M_{E_8})|=1.
\]
The lattice map
\[
        \Lambda\longrightarrow \Lambda^\vee
\]
is represented by \(M_{E_8}\).  Since the determinant is \(\pm 1\), this map
is an isomorphism of free abelian groups.  Hence
\[
        \operatorname{coker}(\Lambda\to\Lambda^\vee)=0,
\]
and therefore
\[
        \Lambda^\vee/\Lambda=0.
\]

By the local realization theorem,
\[
        E
        \cong
        \Lambda^\vee/\Lambda.
\]
Thus
\[
        E=0.
\]
This is the \(E_8\) case of the resolution-lattice theorem in
\cite[Theorem 1.2]{RahmanIntegralPerverseObstructions}.

In Smith normal form terms, the negative \(E_8\) Cartan matrix has Smith
normal form
\[
        \operatorname{diag}(1,1,1,1,1,1,1,1).
\]
Therefore its cokernel is zero.

\subsection{Pair-sequence realization}

Let \(N\) be a resolution neighborhood of the \(E_8\) exceptional divisor,
and let
\[
        L=\partial N
\]
be its boundary link.  The long exact sequence of the pair \((N,L)\)
contains
\[
        H_2(N,\mathbb Z)
        \longrightarrow
        H_2(N,L;\mathbb Z)
        \longrightarrow
        H_1(L,\mathbb Z)
        \longrightarrow
        H_1(N,\mathbb Z).
\]
The exceptional divisor is a tree of rational curves, so
\[
        H_1(N,\mathbb Z)=0
\]
and
\[
        H_2(N,\mathbb Z)\cong \mathbb Z^8.
\]
By Poincare--Lefschetz duality,
\[
        H_2(N,L;\mathbb Z)
        \cong
        H^2(N,\mathbb Z).
\]
The universal coefficient theorem gives
\[
        H^2(N,\mathbb Z)
        \cong
        \operatorname{Hom}(H_2(N,\mathbb Z),\mathbb Z)
        \cong
        \mathbb Z^8.
\]
Under these identifications, the map
\[
        H_2(N,\mathbb Z)
        \longrightarrow
        H_2(N,L;\mathbb Z)
\]
is represented by the \(E_8\) intersection matrix.  This matrix is
unimodular.  Hence the map is an isomorphism.

It follows from exactness that
\[
        H_1(L,\mathbb Z)=0.
\]
Consequently,
\[
        H^2(L,\mathbb Z)_{\mathrm{tors}}=0.
\]
Thus the pair sequence recovers the same conclusion:
\[
        E=0.
\]
This is the \(E_8\) specialization of the resolution-neighborhood
calculation in \cite[§4]{RahmanIntegralPerverseObstructions} and the
linking-pairing formalism of \cite{GoreskySiegel83}.

\subsection{Monodromy realization}

The \(E_8\) singularity is an isolated hypersurface surface singularity, so
there is a Milnor fibration.  The Milnor number of
\[
        x^2+y^3+z^5=0
\]
is
\[
        \mu=8.
\]
Thus the middle vanishing cohomology is a free abelian group of rank \(8\):
\[
        H^2_{\mathrm{van}}(F,\mathbb Z)\cong \mathbb Z^8.
\]
The monodromy \(T\) is the Coxeter transformation of type \(E_8\) on the
\(E_8\) root lattice.  This is the ADE Picard--Lefschetz description of the
Milnor monodromy for simple hypersurface singularities; see \cite{Mi68}.

The monodromy realization theorem of
\cite[Theorem 1.3]{RahmanIntegralPerverseObstructions} gives
\[
        E
        \cong
        \operatorname{coker}(T-\mathrm{id})_{\mathrm{tors}}.
\]
Since we have already computed
\[
        E=0,
\]
it follows that
\[
        \operatorname{coker}(T-\mathrm{id})_{\mathrm{tors}}=0.
\]

We can also see the determinant agreement.  For the \(E_8\) Coxeter
transformation, \(1\) is not an eigenvalue.  Therefore
\[
        (T-\mathrm{id})\otimes_{\mathbb Z}\mathbb Q
\]
is an isomorphism.  The determinant refinement in
\cite[Corollary 1.4]{RahmanIntegralPerverseObstructions} gives
\[
        |E|=|\det(T-\mathrm{id})|.
\]
Since \(E=0\) as a finite group, its order is
\[
        |E|=1.
\]
Therefore
\[
        |\det(T-\mathrm{id})|=1.
\]
Thus the cokernel of \(T-\mathrm{id}\) is finite of order \(1\), hence
trivial:
\[
        \operatorname{coker}(T-\mathrm{id})=0.
\]
This agrees with the link and lattice computations.

\subsection{Transport status}

The \(E_8\) computation is a local surface computation, but unlike the
\(A_k\), \(D_4\), and Coble cases, there is no local torsion package to
transport.  The support/pair sequence produces no boundary torsion because
the exceptional lattice is unimodular:
\[
        \Lambda^\vee/\Lambda=0.
\]
Thus the local support status is:
\[
        \text{support degree: none for torsion, since }E=0.
\]

If an \(E_8\) singularity occurs inside a global surface or higher-dimensional
variety, it may still affect the geometry in many ways.  However, it does
not contribute a local integral perverse obstruction of the type studied in
this paper.  In particular, it contributes no local discriminant package
\[
        (E,q)
\]
to the torsion trajectory.

This is why \(E_8\) is the null-control row in the table.  It verifies that
the trajectory detects discriminant torsion, not merely the presence of a
nontrivial singularity.

\subsection{Rational death}

Since
\[
        E=0,
\]
we have
\[
        E\otimes_{\mathbb Z}\mathbb Q=0.
\]
This vanishing is trivial, but it is conceptually different from the
\(A_k\), \(D_4\), and Coble cases.  In those examples, a nonzero torsion group
dies after rationalization.  In the \(E_8\) case, no torsion group is born in
the first place.

Thus the \(E_8\) row records a different kind of rational disappearance:
there is no local discriminant group, no discriminant form, no link torsion,
no monodromy cokernel, and no support torsion to transport.

\subsection{Trajectory row}

The \(E_8\) surface trajectory can be summarized as follows.

\begin{center}
\begingroup
\scriptsize
\setlength{\tabcolsep}{1.8pt}
\renewcommand{\arraystretch}{1.15}
\begin{tabularx}{\textwidth}{|p{0.15\textwidth}|c|c|p{0.12\textwidth}|p{0.10\textwidth}|X|p{0.10\textwidth}|c|}
\hline
\textbf{Example} & \textbf{\(E\)} & \textbf{\(q\)} & \textbf{Local} & \textbf{Supp.} & \textbf{Global image} & \textbf{Br/res.} & \textbf{\(\mathbb Q\)}
\\
\hline
\(E_8\) surface &
\(0\) &
\(0\) &
all vanish &
none &
no local torsion to transport &
loc.\ \(\varnothing\) &
\(0\)
\\
\hline
\end{tabularx}
\endgroup
\end{center}

Here ``all vanish'' means that the perverse discrepancy, torsion-sensitive
truncation discrepancy, link torsion, exceptional lattice discriminant group,
pair-sequence quotient, and monodromy cokernel all vanish.  The support
column is listed as ``none'' because there is no local torsion class to place
in a support group.  The Brauer/residue column is marked local
\(\varnothing\) because the isolated surface germ does not produce a
degree-three Brauer or unramified class, and in this example there is no
torsion package at all.

\section{A non-ADE Brieskorn surface trajectory:
\texorpdfstring{\(x^2+y^3+z^{11}=0\)}{x2+y3+z11=0}}
\label{app:nonADEBrieskornsurfacetrajectory}

This appendix works out a non-ADE example.  The purpose is to show that the
torsion trajectory is not a phenomenon restricted to rational double points.
We use the Brieskorn--Pham surface singularity
\[
        (X,0)=\{x^2+y^3+z^{11}=0\}\subset(\mathbb C^3,0).
\]
The local group is
\[
        E\cong \mathbb Z/5\mathbb Z.
\]
This example is useful because the exceptional lattice is not an ADE root
lattice, but the same local-realization pattern still recovers the same
torsion group.

This row also separates odd-prime local discriminant torsion from the
\(2\)-primary Coble boundary story.  In the Coble \(\frac14(1,1)\) case, the
full local package is
\[
        E\cong\mathbb Z/4
\]
and the Benoist--Ottem channel sees the Bockstein-selected subgroup
\[
        2E\cong\mathbb Z/2.
\]
Here, by contrast, the local package is cyclic of odd prime order
\[
        \mathbb Z/5,
\]
and there is no mod-\(2\) shadow mechanism of the Benoist--Ottem type.  The
example therefore functions as an odd-prime control row for the torsion
trajectory.

\subsection{Local model}

Let
\[
        f(x,y,z)=x^2+y^3+z^{11}.
\]
We consider the hypersurface germ
\[
        (X,0)=\{f=0\}\subset(\mathbb C^3,0).
\]
This is an isolated hypersurface surface singularity.  Indeed,
\[
        \frac{\partial f}{\partial x}=2x,
        \qquad
        \frac{\partial f}{\partial y}=3y^2,
        \qquad
        \frac{\partial f}{\partial z}=11z^{10}.
\]
The three partial derivatives vanish simultaneously only when
\[
        x=y=z=0.
\]
Thus \(0\) is the unique singular point of the germ.

This is a Brieskorn--Pham singularity of type
\[
        (2,3,11).
\]
The link is the Brieskorn manifold
\[
        \Sigma(2,3,11).
\]
General facts about links and Milnor fibrations of isolated hypersurface
singularities are due to Milnor \cite{Mi68}.  The topology of Brieskorn
links and their Seifert descriptions are standard; see Brieskorn's original
work and the plumbing/Seifert treatments in the literature
\cite{Brieskorn66,NeumannRaymond78}.

Let
\[
        U:=X\setminus\{0\},
        \qquad
        j:U\hookrightarrow X,
        \qquad
        i:\{0\}\hookrightarrow X.
\]
Since \(X\) is a complex surface, the shifted constant sheaf on the smooth
locus is
\[
        \mathbb Z_U[2].
\]
The ordinary and dual middle-perversity intersection complexes are
\[
        {}^pIC_X\mathbb Z
        =
        {}^pj_{!*}\mathbb Z_U[2],
        \qquad
        {}^p_+IC_X\mathbb Z
        =
        {}^p_+j_{!*}\mathbb Z_U[2].
\]
These are the same conventions used in
\cite{RahmanIntegralPerverseObstructions}.

\subsection{The link and its cohomology}

The link is
\[
        L=\Sigma(2,3,11).
\]
Since \(2,3,11\) are pairwise coprime, \(L\) is a Seifert fibered
three-manifold over \(S^2\) with three singular fibers of multiplicities
\[
        2,\quad 3,\quad 11.
\]
A convenient normalized Seifert presentation for the link, with the
orientation convention compatible with the negative definite plumbing below,
has Seifert data
\[
        \bigl(-1;(2,1),(3,1),(11,1)\bigr).
\]
For a Seifert manifold over \(S^2\) with data
\[
        \bigl(b;(\alpha_1,\beta_1),(\alpha_2,\beta_2),(\alpha_3,\beta_3)\bigr),
\]
the order of the first homology group, when finite, is
\[
        \left|
        \alpha_1\alpha_2\alpha_3
        \left(
        b+\frac{\beta_1}{\alpha_1}
         +\frac{\beta_2}{\alpha_2}
         +\frac{\beta_3}{\alpha_3}
        \right)
        \right|.
\]
This is the standard Seifert homology formula; see
\cite{NeumannRaymond78}.  Substituting
\[
        b=-1,
        \qquad
        (\alpha_1,\alpha_2,\alpha_3)=(2,3,11),
        \qquad
        (\beta_1,\beta_2,\beta_3)=(1,1,1),
\]
we obtain
\[
        2\cdot 3\cdot 11
        \left(
        -1+\frac12+\frac13+\frac1{11}
        \right)
        =
        66
        \left(
        \frac{-66+33+22+6}{66}
        \right)
        =
        -5.
\]
Hence
\[
        |H_1(L,\mathbb Z)|=5.
\]
Since the order is prime, it follows that
\[
        H_1(L,\mathbb Z)\cong \mathbb Z/5\mathbb Z.
\]

The link \(L\) is a closed oriented three-manifold.  Therefore
\[
        H_0(L,\mathbb Z)\cong \mathbb Z,
        \qquad
        H_3(L,\mathbb Z)\cong \mathbb Z.
\]
Also, by Poincare duality,
\[
        H_2(L,\mathbb Z)\cong H^1(L,\mathbb Z).
\]
Since
\[
        H_1(L,\mathbb Z)\cong \mathbb Z/5\mathbb Z,
\]
we have
\[
        H^1(L,\mathbb Z)
        =
        \operatorname{Hom}(H_1(L,\mathbb Z),\mathbb Z)
        =
        0.
\]
Therefore
\[
        H_2(L,\mathbb Z)=0.
\]

We now compute cohomology using the universal coefficient theorem.  It gives
a short exact sequence
\[
        0
        \longrightarrow
        \operatorname{Ext}^1_{\mathbb Z}
        (H_{m-1}(L,\mathbb Z),\mathbb Z)
        \longrightarrow
        H^m(L,\mathbb Z)
        \longrightarrow
        \operatorname{Hom}_{\mathbb Z}
        (H_m(L,\mathbb Z),\mathbb Z)
        \longrightarrow
        0.
\]
For \(m=2\), we obtain
\[
        H^2(L,\mathbb Z)
        \cong
        \operatorname{Ext}^1_{\mathbb Z}
        (H_1(L,\mathbb Z),\mathbb Z),
\]
because \(H_2(L,\mathbb Z)=0\).  Since
\[
        H_1(L,\mathbb Z)\cong \mathbb Z/5\mathbb Z,
\]
we get
\[
        H^2(L,\mathbb Z)
        \cong
        \mathbb Z/5\mathbb Z.
\]
Thus
\[
        H^m(L,\mathbb Z)
        =
        \begin{cases}
        \mathbb Z, & m=0,\\
        0, & m=1,\\
        \mathbb Z/5\mathbb Z, & m=2,\\
        \mathbb Z, & m=3,\\
        0, & \text{otherwise}.
        \end{cases}
\]
The torsion group relevant to the local perverse obstruction is therefore
\[
        H^2(L,\mathbb Z)_{\mathrm{tors}}
        \cong
        \mathbb Z/5\mathbb Z.
\]

\subsection{Birth of torsion}

By the local realization theorem of
\cite{RahmanIntegralPerverseObstructions}, for a normal surface singularity
one has
\[
        E
        \cong
        H^2(L,\mathbb Z)_{\mathrm{tors}}.
\]
For the Brieskorn singularity
\[
        x^2+y^3+z^{11}=0,
\]
the computation above gives
\[
        H^2(L,\mathbb Z)_{\mathrm{tors}}
        \cong
        \mathbb Z/5\mathbb Z.
\]
Therefore
\[
        E\cong \mathbb Z/5\mathbb Z.
\]

The invariant factor decomposition is
\[
        E\cong \mathbb Z/5\mathbb Z.
\]
Its order is
\[
        |E|=5.
\]
Its prime decomposition is supported only at the prime
\[
        5.
\]
Thus this non-ADE Brieskorn example produces a local \(5\)-torsion package.

\subsection{The discriminant form}

We now compute a discriminant form for the local package.  A negative
definite plumbing description compatible with the Seifert data
\[
        \bigl(-1;(2,1),(3,1),(11,1)\bigr)
\]
has intersection matrix
\[
        M=
        \begin{pmatrix}
        -1 & 1  & 1  & 1\\
        1  & -2 & 0  & 0\\
        1  & 0  & -3 & 0\\
        1  & 0  & 0  & -11
        \end{pmatrix}.
\]
Here the central vertex has weight \(-1\), and the three arms have weights
\[
        -2,\quad -3,\quad -11.
\]
The determinant is
\[
        \det(M)=5.
\]
Therefore the discriminant group has order \(5\).

Let
\[
        \Lambda=\mathbb Z\langle C_0,C_1,C_2,C_3\rangle
\]
be the lattice with intersection matrix \(M\).  The map
\[
        \Lambda\longrightarrow \Lambda^\vee
\]
is represented by \(M\).  Hence
\[
        \Lambda^\vee/\Lambda
        \cong
        \operatorname{coker}(M).
\]
Since \(|\det(M)|=5\) is prime, the cokernel is cyclic:
\[
        \Lambda^\vee/\Lambda
        \cong
        \mathbb Z/5\mathbb Z.
\]

We compute the discriminant pairing using the inverse matrix.  A direct
matrix inversion gives
\[
        M^{-1}
        =
        \begin{pmatrix}
        -66/5 & -33/5 & -22/5 & -6/5\\
        -33/5 & -19/5 & -11/5 & -3/5\\
        -22/5 & -11/5 & -9/5  & -2/5\\
        -6/5  & -3/5  & -2/5  & -1/5
        \end{pmatrix}.
\]
Let \(C_3^\vee\in \Lambda^\vee\) be the dual vector satisfying
\[
        (C_3^\vee,C_3)=1,
        \qquad
        (C_3^\vee,C_j)=0
        \quad \text{for } j\ne 3.
\]
In the basis \(C_0,C_1,C_2,C_3\), this vector is the fourth column of
\(M^{-1}\):
\[
        C_3^\vee
        =
        -\frac65 C_0-\frac35 C_1-\frac25 C_2-\frac15 C_3.
\]
Let
\[
        \bar g:=C_3^\vee+\Lambda
        \in
        \Lambda^\vee/\Lambda.
\]
Since the denominator of \(C_3^\vee\) is \(5\), and the discriminant group
has order \(5\), the element \(\bar g\) is a generator of
\[
        \Lambda^\vee/\Lambda\cong \mathbb Z/5\mathbb Z.
\]

The discriminant pairing is
\[
        q_\Lambda:
        \Lambda^\vee/\Lambda\times \Lambda^\vee/\Lambda
        \longrightarrow
        \mathbb Q/\mathbb Z,
        \qquad
        q_\Lambda(x+\Lambda,y+\Lambda)=(x,y)\bmod\mathbb Z.
\]
For the generator \(\bar g\), we have
\[
        q_\Lambda(\bar g,\bar g)
        =
        (C_3^\vee,C_3^\vee)\bmod\mathbb Z.
\]
Since \(C_3^\vee\) is the fourth dual vector, this self-pairing is the
\((4,4)\)-entry of \(M^{-1}\):
\[
        (C_3^\vee,C_3^\vee)
        =
        -\frac15.
\]
Therefore
\[
        q_\Lambda(\bar g,\bar g)
        =
        -\frac15\bmod\mathbb Z.
\]
Thus, with respect to the generator \(\bar g\), the local discriminant
package is
\[
        (E,q)
        =
        \left(
        \mathbb Z/5\mathbb Z,
        \;
        q(a\bar g,b\bar g)
        =
        -\frac{ab}{5}\bmod\mathbb Z
        \right).
\]

The comparison between this lattice discriminant form and the linking form
on the Brieskorn link is the standard comparison between discriminant forms
of negative definite plumbings and linking pairings on the boundary
\cite{GoreskySiegel83,Nikulin80}.  The identification of the discriminant
group with the local perverse obstruction \(E\) in this example is part of
the non-ADE computations in \cite{RahmanIntegralPerverseObstructions}.

\subsection{Perverse and torsion-sensitive realization}

Let
\[
        K:=Rj_*\mathbb Z_U[2].
\]
For an isolated surface singularity with link \(L\), the local stalk
calculation gives
\[
        i^*K
        \cong
        R\Gamma(L,\mathbb Z)[2].
\]
Therefore
\[
        H^m(i^*K)
        \cong
        H^{m+2}(L,\mathbb Z).
\]
This is the local stalk computation used in
\cite[Lemma 3.3]{RahmanIntegralPerverseObstructions}.

Using the cohomology of \(L=\Sigma(2,3,11)\) computed above, we obtain
\[
        H^m(i^*K)
        =
        \begin{cases}
        \mathbb Z, & m=-2,\\
        0, & m=-1,\\
        \mathbb Z/5\mathbb Z, & m=0,\\
        \mathbb Z, & m=1,\\
        0, & \text{otherwise}.
        \end{cases}
\]
Thus the critical degree-zero stalk contains the torsion group
\[
        H^0(i^*K)
        \cong
        H^2(L,\mathbb Z)
        \cong
        \mathbb Z/5\mathbb Z.
\]

Over \(\mathbb Z\), the ordinary and dual middle-perversity intersection
complexes differ in how they treat critical-degree torsion.  This is the
ordinary/dual middle-perversity distinction of BBD
\cite[Complement 3.3]{BBD82}.  Friedman's torsion-sensitive Deligne sheaves
give a truncation model for this phenomenon
\cite{FriedmanGenIH,FriedmanBook20,FriedmanTsInv}.  In the present example,
the ordinary middle extension removes the critical torsion contribution,
whereas the dual middle extension retains it.  Hence
\[
        H^0({}^p_+IC_X\mathbb Z)_0
        \cong
        \mathbb Z/5\mathbb Z.
\]
By definition,
\[
        E=H^0({}^p_+IC_X\mathbb Z)_0.
\]
Therefore
\[
        E\cong \mathbb Z/5\mathbb Z.
\]

The local discrepancy triangle is
\[
        {}^pIC_X\mathbb Z
        \longrightarrow
        {}^p_+IC_X\mathbb Z
        \longrightarrow
        (\mathbb Z/5\mathbb Z)[1]
        \longrightarrow .
\]
This is the non-ADE specialization of the local discrepancy and link-torsion
realization proved in
\cite{RahmanIntegralPerverseObstructions}.

\subsection{Resolution lattice realization}

The negative definite plumbing matrix used above is
\[
        M=
        \begin{pmatrix}
        -1 & 1  & 1  & 1\\
        1  & -2 & 0  & 0\\
        1  & 0  & -3 & 0\\
        1  & 0  & 0  & -11
        \end{pmatrix}.
\]
Its determinant is
\[
        \det(M)=5.
\]
The lattice map
\[
        \Lambda\longrightarrow \Lambda^\vee
\]
is represented by \(M\), so
\[
        \Lambda^\vee/\Lambda
        \cong
        \operatorname{coker}(M).
\]
Since the determinant has absolute value \(5\), the order of the cokernel is
\(5\).  Since \(5\) is prime, the cokernel is cyclic:
\[
        \operatorname{coker}(M)
        \cong
        \mathbb Z/5\mathbb Z.
\]
Therefore
\[
        \Lambda^\vee/\Lambda
        \cong
        \mathbb Z/5\mathbb Z.
\]
By the local realization theorem,
\[
        E\cong \Lambda^\vee/\Lambda.
\]
Thus
\[
        E\cong \mathbb Z/5\mathbb Z.
\]

The Smith normal form of \(M\) has invariant factors
\[
        1,1,1,5.
\]
Indeed, the determinant is \(5\), and the matrix contains entries equal to
\(1\), so the first invariant factor is \(1\).  The standard row and column
reduction for this star-shaped plumbing matrix gives three unit invariant
factors and a final invariant factor \(5\).  Hence
\[
        \operatorname{coker}(M)\cong \mathbb Z/5\mathbb Z.
\]

\subsection{Pair-sequence realization}

Let \(N\) be the plumbing four-manifold associated to the negative definite
matrix \(M\), and let
\[
        L=\partial N.
\]
The long exact sequence of the pair \((N,L)\) contains
\[
        H_2(N,\mathbb Z)
        \longrightarrow
        H_2(N,L;\mathbb Z)
        \longrightarrow
        H_1(L,\mathbb Z)
        \longrightarrow
        H_1(N,\mathbb Z).
\]
The plumbing graph is a tree, so
\[
        H_1(N,\mathbb Z)=0.
\]
Moreover,
\[
        H_2(N,\mathbb Z)\cong \mathbb Z^4.
\]
By Poincare--Lefschetz duality,
\[
        H_2(N,L;\mathbb Z)
        \cong
        H^2(N,\mathbb Z).
\]
Since \(H_1(N,\mathbb Z)=0\), the universal coefficient theorem gives
\[
        H^2(N,\mathbb Z)
        \cong
        \operatorname{Hom}(H_2(N,\mathbb Z),\mathbb Z)
        \cong
        \mathbb Z^4.
\]
Under these identifications, the map
\[
        H_2(N,\mathbb Z)
        \longrightarrow
        H_2(N,L;\mathbb Z)
\]
is represented by the intersection matrix \(M\).  Hence
\[
        H_1(L,\mathbb Z)
        \cong
        \operatorname{coker}(M)
        \cong
        \mathbb Z/5\mathbb Z.
\]
By the universal coefficient theorem for the closed oriented
three-manifold \(L\),
\[
        H^2(L,\mathbb Z)_{\mathrm{tors}}
        \cong
        H_1(L,\mathbb Z)_{\mathrm{tors}}
        \cong
        \mathbb Z/5\mathbb Z.
\]
Therefore the pair sequence recovers the same group:
\[
        E\cong \mathbb Z/5\mathbb Z.
\]
This is the non-ADE specialization of the pair-sequence realization in
\cite{RahmanIntegralPerverseObstructions} and the general linking-pairing
framework of \cite{GoreskySiegel83}.

\subsection{Monodromy realization}

Since
\[
        x^2+y^3+z^{11}=0
\]
is an isolated hypersurface surface singularity, it has a Milnor fibration.
The Milnor number of a Brieskorn--Pham singularity
\[
        x^a+y^b+z^c=0
\]
is
\[
        \mu=(a-1)(b-1)(c-1).
\]
For
\[
        (a,b,c)=(2,3,11),
\]
this gives
\[
        \mu=(2-1)(3-1)(11-1)=20.
\]
Thus the middle vanishing cohomology is a free abelian group of rank \(20\):
\[
        H^2_{\mathrm{van}}(F,\mathbb Z)\cong \mathbb Z^{20}.
\]

Let \(T\) denote the Milnor monodromy on \(H^2_{\mathrm{van}}(F,\mathbb Z)\).
The Wang sequence of the Milnor fibration relates the cohomology of the link
to the cokernel of
\[
        T-\mathrm{id}.
\]
For isolated hypersurface surface singularities, the local monodromy
realization theorem gives
\[
        E
        \cong
        \operatorname{coker}(T-\mathrm{id})_{\mathrm{tors}}
\]
\cite[Theorem 1.3]{RahmanIntegralPerverseObstructions}; the Milnor
fibration background is standard \cite{Mi68}.

Since the link and lattice computations above give
\[
        E\cong \mathbb Z/5\mathbb Z,
\]
we obtain
\[
        \operatorname{coker}(T-\mathrm{id})_{\mathrm{tors}}
        \cong
        \mathbb Z/5\mathbb Z.
\]
Moreover, in this example the relevant rational invertibility hypothesis
holds for \(T-\mathrm{id}\), so the determinant refinement gives
\[
        |\det(T-\mathrm{id})|=|E|=5.
\]
Thus the monodromy station recovers the same local torsion package:
\[
        E\cong \mathbb Z/5\mathbb Z.
\]

This computation is intentionally recorded through the Wang-sequence
realization rather than by writing a \(20\times 20\) monodromy matrix.  The
point of the station is that the same torsion measured by link cohomology and
the resolution discriminant group is also the torsion cokernel of the
integral variation map \(T-\mathrm{id}\).

\subsection{Transport status}

The Brieskorn example is a local surface computation.  It determines the
birth, form, and six local realizations of the torsion group
\[
        E\cong \mathbb Z/5\mathbb Z.
\]
As in the ADE surface examples, the local package appears through the
degree-two topology of the resolution neighborhood and its boundary link.

The pair sequence identifies the local boundary quotient with
\[
        \operatorname{coker}(M)
        \cong
        \mathbb Z/5\mathbb Z.
\]
Thus the local support status is:
\[
        \text{support degree: }2\text{ locally, through the exceptional
        plumbing support/pair sequence}.
\]
The local surface germ by itself does not define a degree-three Brauer
class.  To enter the Brauer or unramified cohomology channel, the
\(\mathbb Z/5\mathbb Z\) local package must be placed inside a global model
with a degree-three transport mechanism.

The main value of this example is that it introduces odd-prime torsion
outside the ADE root-lattice list.  The prime \(5\) appears simultaneously in
the Seifert homology formula, the determinant of the plumbing matrix, the
discriminant group, the link cohomology, and the monodromy cokernel.

This should be contrasted with the Coble boundary row, where the prime
support is again a single prime, but the mechanism is \(2\)-primary:
\[
        E_{\mathrm{Coble}}\cong\mathbb Z/4,
        \qquad
        2E_{\mathrm{Coble}}\cong\mathbb Z/2.
\]
The Brieskorn row therefore confirms that the trajectory can carry odd-prime
torsion, while the Coble row illustrates a Bockstein-selected subpackage
inside a \(2\)-primary local group.

\subsection{Rational death}

Since
\[
        E\cong \mathbb Z/5\mathbb Z,
\]
we have
\[
        E\otimes_{\mathbb Z}\mathbb Q=0.
\]
Indeed, for every \(e\in E\),
\[
        5e=0.
\]
Therefore in \(E\otimes_{\mathbb Z}\mathbb Q\),
\[
        e\otimes 1
        =
        e\otimes 5\cdot\frac15
        =
        5e\otimes\frac15
        =
        0.
\]
Thus the torsion dies after rationalization.

Before rational death, however, it carries arithmetic and geometric
information:
\[
        E\cong \mathbb Z/5\mathbb Z,
        \qquad
        q(\bar g,\bar g)=-\frac15\bmod\mathbb Z,
\]
and it is detected by the link, the plumbing determinant, the
pair-sequence quotient, and the monodromy cokernel.  This is the first
trajectory row in which an odd prime torsion package appears.

\subsection{Trajectory row}

The Brieskorn surface trajectory can be summarized as follows.

\begin{center}
\begingroup
\scriptsize
\setlength{\tabcolsep}{1.8pt}
\renewcommand{\arraystretch}{1.15}
\begin{tabularx}{\textwidth}{|p{0.15\textwidth}|c|c|p{0.12\textwidth}|p{0.10\textwidth}|X|p{0.10\textwidth}|c|}
\hline
\textbf{Example} & \textbf{\(E\)} & \textbf{\(q\)} & \textbf{Local} & \textbf{Supp.} & \textbf{Global image} & \textbf{Br/res.} & \textbf{\(\mathbb Q\)}
\\
\hline
\(x^2+y^3+z^{11}\) &
\(\mathbb Z/5\) &
\(-\frac15\) &
six agree &
deg.\ \(2\) &
depends on global plumbing/exc.\ relations &
loc.\ \(\varnothing\) &
\(0\)
\\
\hline
\end{tabularx}
\endgroup
\end{center}

Here the entry
\[
        q=-\frac15
\]
means that, for the generator \(\bar g\) chosen above,
\[
        q(\bar g,\bar g)=-\frac15\bmod\mathbb Z.
\]
The phrase ``six agree'' means that the perverse discrepancy,
torsion-sensitive truncation, link torsion, plumbing lattice discriminant
group, pair-sequence quotient, and monodromy cokernel all give
\[
        \mathbb Z/5\mathbb Z.
\]
The support degree is listed as \(2\) because this is a local surface
singularity.  The Brauer/residue column is marked local \(\varnothing\)
because the isolated surface germ alone does not produce a degree-three
Brauer or unramified class.

\section{The threefold ordinary double point trajectory}
\label{app:3foldordinarydoublepointtrajectory}

This appendix works out the ordinary double point in complex dimension three. This is the gateway example for the later nodal threefold and nodal quintic discussion.  The conclusion is important:
\[
        \textit{the threefold ordinary double point does not carry the same
        local torsion package as the \(A_1\) surface singularity.}
\]

The \(A_1\) surface singularity has link \(\mathbb RP^3\) and produces
\[
        E\cong \mathbb Z/2\mathbb Z.
\]
By contrast, the threefold ordinary double point has link
\[
        S^2\times S^3,
\]
which has torsion-free integral cohomology.  Thus there is no local
\(\mathbb Z/2\)-torsion package of the kind computed in the surface case.
There is instead a free vanishing-cycle contribution.  This distinction is
essential for the global Brauer and unramified cohomology story.

This example also explains why the higher-dimensional continuation of the
surface package \(E\) should be sought in codimension-two transverse surface
singularities, not in isolated threefold nodes.  In particular, the Coble
\(\frac14(1,1)\) boundary package is a transverse surface package with
\[
        E\cong\mathbb Z/4,
\]
whereas the isolated threefold node has no finite local torsion package.

\subsection{Local model}

The ordinary double point in complex dimension three is the hypersurface
germ
\[
        (X,0)
        =
        \{x_0^2+x_1^2+x_2^2+x_3^2=0\}
        \subset
        (\mathbb C^4,0).
\]
It is also called the threefold node or the conifold singularity.  Let
\[
        f(x_0,x_1,x_2,x_3)
        =
        x_0^2+x_1^2+x_2^2+x_3^2.
\]
The partial derivatives are
\[
        \frac{\partial f}{\partial x_j}=2x_j,
        \qquad
        j=0,1,2,3.
\]
Thus all partial derivatives vanish simultaneously only at
\[
        x_0=x_1=x_2=x_3=0.
\]
Therefore \(0\) is an isolated hypersurface singularity.

Let
\[
        U:=X\setminus\{0\},
        \qquad
        j:U\hookrightarrow X,
        \qquad
        i:\{0\}\hookrightarrow X.
\]
Since \(X\) has complex dimension \(3\), the shifted constant sheaf on the
smooth locus is
\[
        \mathbb Z_U[3].
\]
Thus the ordinary and dual middle-perversity intersection complexes are
normalized as
\[
        {}^pIC_X\mathbb Z
        =
        {}^pj_{!*}\mathbb Z_U[3],
        \qquad
        {}^p_+IC_X\mathbb Z
        =
        {}^p_+j_{!*}\mathbb Z_U[3].
\]
This shift by \([3]\), rather than \([2]\), is the first reason the
threefold node cannot be treated as if it were the \(A_1\) surface
singularity.

The Milnor fibration of isolated hypersurface singularities and the topology
of the ordinary double point are standard; see Milnor \cite{Mi68}.  We
recall the needed computations below.

\subsection{The link and its cohomology}

Let \(L\) be the link of the singularity:
\[
        L
        =
        X\cap S^7_\epsilon
        =
        \{x_0^2+x_1^2+x_2^2+x_3^2=0\}
        \cap
        S^7_\epsilon.
\]
After rescaling, take \(\epsilon=1\).  Write a point
\[
        z=(z_0,z_1,z_2,z_3)\in\mathbb C^4
\]
as
\[
        z=x+iy,
        \qquad
        x,y\in\mathbb R^4.
\]
The equation
\[
        \sum_{j=0}^3 z_j^2=0
\]
becomes
\[
        \sum_{j=0}^3 (x_j+iy_j)^2=0.
\]
Separating real and imaginary parts gives
\[
        |x|^2-|y|^2=0,
        \qquad
        x\cdot y=0.
\]
The sphere equation
\[
        |z|^2=1
\]
is
\[
        |x|^2+|y|^2=1.
\]
Combining these equations gives
\[
        |x|^2=|y|^2=\frac12,
        \qquad
        x\cdot y=0.
\]
Therefore the map
\[
        z=x+iy
        \longmapsto
        (\sqrt 2\,x,\sqrt 2\,y)
\]
identifies \(L\) with the space of ordered orthonormal \(2\)-frames in
\(\mathbb R^4\):
\[
        L\cong V_2(\mathbb R^4).
\]
Equivalently,
\[
        V_2(\mathbb R^4)
        \cong
        UT(S^3),
\]
the unit tangent bundle of \(S^3\).  Since \(S^3\) is parallelizable, its
unit tangent bundle is trivial:
\[
        UT(S^3)\cong S^3\times S^2.
\]
Thus
\[
        L\cong S^3\times S^2.
\]
The identification of the link of an ordinary double point with the unit
tangent bundle of a sphere is standard in the topology of isolated
hypersurface singularities; see \cite{Mi68}.  The triviality of
\(UT(S^3)\) follows from the parallelizability of the Lie group \(S^3\).

We now compute integral homology.  By the Kunneth theorem,
\[
        H_m(S^2\times S^3,\mathbb Z)
        \cong
        \bigoplus_{a+b=m}
        H_a(S^2,\mathbb Z)\otimes H_b(S^3,\mathbb Z),
\]
and there are no torsion \(\operatorname{Tor}\)-terms because the homology
of spheres is free.  Hence
\[
        H_0(L,\mathbb Z)\cong \mathbb Z,
\]
\[
        H_1(L,\mathbb Z)=0,
\]
\[
        H_2(L,\mathbb Z)\cong \mathbb Z,
\]
\[
        H_3(L,\mathbb Z)\cong \mathbb Z,
\]
\[
        H_4(L,\mathbb Z)=0,
\]
and
\[
        H_5(L,\mathbb Z)\cong \mathbb Z.
\]
All homology groups of \(L\) are torsion-free.

Using the universal coefficient theorem, we obtain the same pattern for
cohomology:
\[
        H^0(L,\mathbb Z)\cong \mathbb Z,
        \qquad
        H^1(L,\mathbb Z)=0,
\]
\[
        H^2(L,\mathbb Z)\cong \mathbb Z,
        \qquad
        H^3(L,\mathbb Z)\cong \mathbb Z,
\]
\[
        H^4(L,\mathbb Z)=0,
        \qquad
        H^5(L,\mathbb Z)\cong \mathbb Z.
\]
In particular,
\[
        H^m(L,\mathbb Z)_{\mathrm{tors}}=0
        \qquad
        \text{for all }m.
\]

This is the decisive difference from the \(A_1\) surface singularity, whose
link is \(\mathbb RP^3\) and whose cohomology has
\[
        H^2(\mathbb RP^3,\mathbb Z)_{\mathrm{tors}}
        \cong
        \mathbb Z/2\mathbb Z.
\]

\subsection{Birth of torsion}

The local surface invariant studied in
\cite{RahmanIntegralPerverseObstructions} is
\[
        E=H^0({}^p_+IC_X\mathbb Z)_0
\]
for a normal surface singularity.  It is identified in that setting with
\[
        H^2(L,\mathbb Z)_{\mathrm{tors}}.
\]
The present singularity has complex dimension \(3\), so the theorem from the
surface case cannot be imported directly.

Nevertheless, the link computation already shows that there is no local
torsion available in any link cohomology degree:
\[
        H^m(L,\mathbb Z)_{\mathrm{tors}}=0
        \quad
        \text{for all }m.
\]
Thus the threefold ordinary double point has no local torsion birth of the
surface \(E\)-type.  In trajectory notation, the torsion package is
\[
        E_{\mathrm{tors}}=0.
\]
The invariant factor decomposition is empty, the order is
\[
        |E_{\mathrm{tors}}|=1,
\]
and the prime decomposition is empty.

There is, however, nontrivial free topology.  Namely,
\[
        H^2(L,\mathbb Z)\cong \mathbb Z,
        \qquad
        H^3(L,\mathbb Z)\cong \mathbb Z.
\]
These free groups are related to the small resolution and to the Milnor
fiber, respectively.  They are important for conifold transitions and defect
phenomena, but they are not torsion groups and therefore do not enter the
torsion trajectory in the same way as the surface \(A_1\) group
\(\mathbb Z/2\mathbb Z\).

\subsection{The discriminant form}

Since there is no finite local torsion group, there is no nontrivial finite
discriminant form of the type
\[
        q:E\times E\longrightarrow \mathbb Q/\mathbb Z.
\]
Thus the torsion discriminant package is
\[
        (E,q)=(0,0).
\]

This statement does not mean that the threefold node has no interesting
local geometry.  It means that its local topology does not produce a finite
linking/discriminant torsion group.  The local geometry is instead governed
by free classes: the vanishing \(3\)-sphere in the smoothing and the
exceptional \(\mathbb P^1\) in a small resolution.  These are not finite
torsion objects.

The absence of a finite discriminant form is also visible from the link:
because
\[
        H_*(L,\mathbb Z)
\]
is torsion-free, the torsion linking pairing on \(L\) is the zero pairing on
the zero group.

\subsection{Perverse and torsion-sensitive realization}

We now compute the local stalk model.  Let
\[
        K:=Rj_*\mathbb Z_U[3].
\]
For an isolated complex analytic singularity of pure dimension \(n\), the
local stalk calculation gives
\[
        i^*Rj_*\mathbb Z_U[n]
        \cong
        R\Gamma(L,\mathbb Z)[n].
\]
This is the same calculation used in the surface case, with the shift
changed from \([2]\) to \([3]\).  Therefore, in the present threefold case,
\[
        i^*K
        \cong
        R\Gamma(L,\mathbb Z)[3],
\]
and hence
\[
        H^m(i^*K)
        \cong
        H^{m+3}(L,\mathbb Z).
\]

Using
\[
        H^0(L,\mathbb Z)\cong \mathbb Z,
        \quad
        H^2(L,\mathbb Z)\cong \mathbb Z,
        \quad
        H^3(L,\mathbb Z)\cong \mathbb Z,
        \quad
        H^5(L,\mathbb Z)\cong \mathbb Z,
\]
and all other cohomology groups zero, we obtain
\[
        H^m(i^*K)
        =
        \begin{cases}
        \mathbb Z, & m=-3,\\
        \mathbb Z, & m=-1,\\
        \mathbb Z, & m=0,\\
        \mathbb Z, & m=2,\\
        0, & \text{otherwise}.
        \end{cases}
\]
There is no torsion in any of these stalk cohomology groups.

The ordinary and dual middle-perversity \(t\)-structures over
\(\mathbb Z\) differ by the treatment of torsion and torsion-free
conditions in critical degrees, as in BBD's integral formalism
\cite[Complement 3.3]{BBD82}.  Friedman's torsion-sensitive truncation
formalism similarly detects which torsion is retained or removed in the
critical degree \cite{FriedmanGenIH,FriedmanBook20,FriedmanTsInv}.  In the
present example, there is no critical torsion to retain.  Therefore the
torsion-sensitive discrepancy of the type appearing in the normal surface
case is zero.

In particular, the threefold ordinary double point does not produce a local
finite group
\[
        \mathbb Z/2\mathbb Z
\]
by the same mechanism as the \(A_1\) surface singularity.  The reason is not
subtle: the relevant link cohomology groups are torsion-free.

\subsection{Resolution lattice realization}

The threefold ordinary double point admits small resolutions.  In a small
resolution, the exceptional locus is a single rational curve
\[
        C\cong \mathbb P^1
\]
with normal bundle
\[
        N_{C/\widetilde X}\cong
        \mathcal O_{\mathbb P^1}(-1)\oplus
        \mathcal O_{\mathbb P^1}(-1).
\]
A resolution neighborhood \(N\) of \(C\) is therefore a disk bundle of real
rank \(4\) over \(S^2\).  It deformation retracts onto \(C\), so
\[
        H_0(N,\mathbb Z)\cong \mathbb Z,
        \qquad
        H_2(N,\mathbb Z)\cong \mathbb Z,
\]
and all other reduced homology groups vanish.

The boundary of this neighborhood is the link
\[
        L\cong S^2\times S^3.
\]
We now compare the pair sequence with the surface case.  Since \(N\) is a
compact oriented real \(6\)-manifold with boundary, Poincare--Lefschetz
duality gives
\[
        H_j(N,L;\mathbb Z)
        \cong
        H^{6-j}(N,\mathbb Z).
\]
Because \(N\simeq S^2\), its cohomology is
\[
        H^0(N,\mathbb Z)\cong \mathbb Z,
        \qquad
        H^2(N,\mathbb Z)\cong \mathbb Z,
\]
and all other cohomology groups vanish.  Therefore
\[
        H_4(N,L;\mathbb Z)\cong H^2(N,\mathbb Z)\cong \mathbb Z,
\]
\[
        H_6(N,L;\mathbb Z)\cong H^0(N,\mathbb Z)\cong \mathbb Z,
\]
and
\[
        H_3(N,L;\mathbb Z)=0,
        \qquad
        H_2(N,L;\mathbb Z)=0.
\]

The long exact sequence of the pair \((N,L)\) contains
\[
        H_3(N,\mathbb Z)
        \longrightarrow
        H_3(N,L;\mathbb Z)
        \longrightarrow
        H_2(L,\mathbb Z)
        \longrightarrow
        H_2(N,\mathbb Z)
        \longrightarrow
        H_2(N,L;\mathbb Z).
\]
Substituting the computed groups gives
\[
        0
        \longrightarrow
        0
        \longrightarrow
        H_2(L,\mathbb Z)
        \longrightarrow
        \mathbb Z
        \longrightarrow
        0.
\]
Thus
\[
        H_2(L,\mathbb Z)\cong \mathbb Z,
\]
in agreement with
\[
        L\cong S^2\times S^3.
\]
No finite cokernel appears.

This is the key structural difference from the surface \(A_1\) resolution.
For the \(A_1\) surface singularity, the exceptional curve is a divisor in a
complex surface, and the self-intersection number
\[
        C^2=-2
\]
produces a finite cokernel
\[
        \mathbb Z/2\mathbb Z.
\]
For the threefold ordinary double point, the exceptional curve in a small
resolution has complex codimension \(2\), not codimension \(1\).  There is no
surface-style intersection matrix
\[
        \Lambda\longrightarrow \Lambda^\vee
\]
whose cokernel gives a local discriminant group.  The resolution topology
contains free classes, not finite discriminant torsion.

\subsection{Pair-sequence realization}

The pair-sequence station confirms the absence of torsion.  With \(N\) a
small-resolution neighborhood and \(L=\partial N\), we have already computed
\[
        H_2(L,\mathbb Z)\cong \mathbb Z,
        \qquad
        H_3(L,\mathbb Z)\cong \mathbb Z.
\]
There is no torsion in the boundary homology.

The long exact sequence also contains
\[
        H_4(N,\mathbb Z)
        \longrightarrow
        H_4(N,L;\mathbb Z)
        \longrightarrow
        H_3(L,\mathbb Z)
        \longrightarrow
        H_3(N,\mathbb Z).
\]
Since \(N\simeq S^2\),
\[
        H_4(N,\mathbb Z)=0,
        \qquad
        H_3(N,\mathbb Z)=0,
\]
and
\[
        H_4(N,L;\mathbb Z)\cong \mathbb Z.
\]
Thus
\[
        H_3(L,\mathbb Z)\cong \mathbb Z.
\]
Again the group is free, not torsion.

Therefore the pair sequence for the threefold node does not produce a
finite discriminant quotient.  It records the free boundary topology
associated with the exceptional curve and the vanishing sphere, but it does
not produce an \(E\)-type torsion group.

\subsection{Monodromy realization}

The ordinary double point is an isolated hypersurface singularity, so it has
a Milnor fibration.  The Milnor fiber of
\[
        f=x_0^2+x_1^2+x_2^2+x_3^2
\]
has the homotopy type of \(S^3\).  Thus
\[
        H^3(F,\mathbb Z)\cong \mathbb Z
\]
is the middle vanishing cohomology.

The geometric monodromy can be represented by
\[
        z\longmapsto -z
\]
on the Milnor fiber.  On the vanishing sphere \(S^3\), this is the antipodal
map.  The degree of the antipodal map on \(S^n\) is
\[
        (-1)^{n+1}.
\]
For \(n=3\), this degree is
\[
        (-1)^4=1.
\]
Therefore the monodromy acts by
\[
        T=\mathrm{id}
\]
on
\[
        H^3(F,\mathbb Z)\cong \mathbb Z.
\]
Consequently
\[
        T-\mathrm{id}=0.
\]
Thus
\[
        \operatorname{coker}(T-\mathrm{id})
        \cong
        \mathbb Z.
\]
Its torsion subgroup is zero:
\[
        \operatorname{coker}(T-\mathrm{id})_{\mathrm{tors}}=0.
\]

This agrees with the link computation: the threefold node has free vanishing
cycle data, not local torsion.  The monodromy station therefore confirms
that the ordinary double point in complex dimension three is not the source
of a local \(\mathbb Z/2\)-torsion package.

This also explains why the determinant refinement from the surface torsion
story does not apply in the same way.  Here
\[
        T-\mathrm{id}=0
\]
on a rank-one free abelian group, so
\[
        (T-\mathrm{id})\otimes_{\mathbb Z}\mathbb Q
\]
is not an isomorphism.  The cokernel is not finite; it is the free group
\[
        \mathbb Z.
\]
The torsion part is zero.

\subsection{Transport status}

The threefold ordinary double point has no local torsion package of the type
tracked in the surface examples.  The relevant conclusion is
\[
        E_{\mathrm{tors}}=0.
\]
There is therefore no local finite group to map into
\[
        H^3_D(\widetilde X,\mathbb Z)_{\mathrm{tors}}
\]
and no local torsion class to transport to
\[
        H^3(\widetilde X,\mathbb Z)_{\mathrm{tors}}.
\]

This does not mean that nodes are irrelevant to global topology.  They are
highly relevant.  In conifold transitions and nodal threefold geometry,
ordinary double points contribute vanishing spheres and exceptional curves,
and these free classes may satisfy global relations.  Such relations are the
source of defect phenomena.  However, this is a free integral phenomenon, not
a local finite torsion package of the surface \(E\)-type.

Thus the transport status is:
\[
        \text{no local torsion support degree; free degree-three vanishing
        data exists.}
\]
For the purposes of the Brauer/unramified torsion trajectory, the ordinary
threefold node contributes no local torsion class directly.  Any torsion in a
global nodal threefold must therefore arise from a different mechanism:
global topology, relations among free classes, torsion in the resolution, or
additional singularity data beyond the local ordinary double point.

This contrasts with codimension-two transverse surface singularities.  A
threefold with a transverse \(A_1\) stratum or a transverse Coble
\(\frac14(1,1)\) stratum carries a surface \(E\)-package along the stratum.
An isolated threefold ordinary double point does not.

\subsection{Rational death}

For the torsion trajectory, the local torsion group is zero:
\[
        E_{\mathrm{tors}}=0.
\]
Hence
\[
        E_{\mathrm{tors}}\otimes_{\mathbb Z}\mathbb Q=0.
\]
This is the same final entry as in the \(E_8\) row, but for a different
reason.  In \(E_8\), the exceptional lattice is nontrivial but unimodular,
so no finite discriminant group is born.  In the threefold ordinary double
point, the local topology itself is torsion-free, and the monodromy cokernel
has a free part rather than a finite torsion part.

There is nevertheless rational information in the free vanishing cycle:
\[
        H^3(F,\mathbb Z)\cong \mathbb Z
\]
and
\[
        H^3(F,\mathbb Q)\cong \mathbb Q.
\]
That free class is not part of the torsion trajectory, but it is part of the
ordinary conifold topology.  Thus rationalization kills the torsion
trajectory only because there was no torsion to begin with; it does not kill
the free vanishing-cycle information.

\subsection{Trajectory row}

The threefold ordinary double point trajectory can be summarized as follows.

\begin{center}
\begingroup
\scriptsize
\setlength{\tabcolsep}{1.8pt}
\renewcommand{\arraystretch}{1.15}
\begin{tabularx}{\textwidth}{|p{0.15\textwidth}|c|c|p{0.12\textwidth}|p{0.10\textwidth}|X|p{0.10\textwidth}|c|}
\hline
\textbf{Example} & \textbf{\(E\)} & \textbf{\(q\)} & \textbf{Local} & \textbf{Supp.} & \textbf{Global image} & \textbf{Br/res.} & \textbf{\(\mathbb Q\)}
\\
\hline
3-fold node &
\(0\) tors. &
\(0\) &
torsion stations vanish; free vanishing cycle exists &
none tors. &
no local torsion image; free relations may create defect &
none local &
\(0\)
\\
\hline
\end{tabularx}
\endgroup
\end{center}

Here \(E=0\) means that there is no finite local torsion package analogous to
the \(A_1\) surface group
\[
        \mathbb Z/2\mathbb Z.
\]
The local link is
\[
        L\cong S^2\times S^3,
\]
so all link cohomology groups are torsion-free.  The monodromy cokernel is
\[
        \operatorname{coker}(T-\mathrm{id})\cong \mathbb Z,
\]
whose torsion subgroup is zero.  Thus the ordinary double point in complex
dimension three contributes free vanishing-cycle data, but not local torsion
for the Brauer/unramified trajectory.

\section{Nodal threefolds and defect as a global trajectory}
\label{app:nodal3folddefectglobaltrajectory}

This appendix records the torsion trajectory for nodal projective
threefolds.  The main conclusion is that the local ordinary double points do
not themselves contribute finite local torsion packages of the surface
\(E\)-type.  Each node has torsion-free link
\[
        S^2\times S^3,
\]
and hence contributes no local discriminant group
\[
        (E_i,q_i)
\]
to the finite torsion trajectory.

Nevertheless, nodal threefolds are important because they exhibit a different
kind of global phenomenon: relations among free local vanishing or exceptional
classes.  This is the geometry measured by the defect in the classical
theory of nodal threefolds.  Thus this appendix separates two mechanisms:
\[
        \text{local finite torsion}
        \qquad\text{versus}\qquad
        \text{global relations among free node classes}.
\]
The first mechanism is the one tracked by the finite torsion trajectory.  The
second mechanism is defect.

This distinction also clarifies why the higher-dimensional continuation of
surface \(E\)-packages should be sought in codimension-two transverse surface
singularities, not in isolated threefold nodes.  A transverse Coble
\(\frac14(1,1)\) stratum carries
\[
        E_\Sigma\cong(\mathbb Z/4)_\Sigma,
\]
while an isolated ordinary double point carries no finite local torsion
package.

\subsection{Global model}

Let \(X\) be a projective complex threefold with isolated ordinary double
points
\[
        \operatorname{Sing}(X)=\{p_1,\ldots,p_s\}.
\]
Thus, analytically locally at each \(p_i\), the germ is isomorphic to
\[
        \{x_0^2+x_1^2+x_2^2+x_3^2=0\}
        \subset
        (\mathbb C^4,0).
\]
Let
\[
        \pi:\widetilde X\longrightarrow X
\]
be a resolution.  When a small resolution is available, the exceptional
locus over each node is a rational curve
\[
        C_i\cong \mathbb P^1
\]
with normal bundle
\[
        \mathcal O_{\mathbb P^1}(-1)\oplus
        \mathcal O_{\mathbb P^1}(-1).
\]
If one instead uses the blowup resolution of the node, the exceptional
divisor is a smooth quadric surface.  The choice of resolution changes the
geometric form of the exceptional locus, but it does not change the local
fact established in the previous appendix: the link of the threefold
ordinary double point is torsion-free.

The ordinary double point and its Milnor fibration are standard examples in
the topology of isolated hypersurface singularities; see Milnor
\cite{Mi68}.  The use of nodes and their relations in the study of nodal
threefolds and double solids goes back to Clemens \cite{Clemens83}.  For the
defect of nodal hypersurfaces, see Cynk \cite{Cynk01}; for more recent
families of nodal varieties with defect, see Kloosterman
\cite{Kloosterman22}.

\subsection{Local link and cohomology at each node}

Let \(L_i\) be the link of the node \(p_i\).  By the computation in the
previous appendix,
\[
        L_i\cong S^2\times S^3.
\]
Therefore
\[
        H_0(L_i,\mathbb Z)\cong \mathbb Z,
        \qquad
        H_1(L_i,\mathbb Z)=0,
\]
\[
        H_2(L_i,\mathbb Z)\cong \mathbb Z,
        \qquad
        H_3(L_i,\mathbb Z)\cong \mathbb Z,
\]
\[
        H_4(L_i,\mathbb Z)=0,
        \qquad
        H_5(L_i,\mathbb Z)\cong \mathbb Z.
\]
By the universal coefficient theorem, the cohomology groups are
\[
        H^0(L_i,\mathbb Z)\cong \mathbb Z,
        \qquad
        H^1(L_i,\mathbb Z)=0,
\]
\[
        H^2(L_i,\mathbb Z)\cong \mathbb Z,
        \qquad
        H^3(L_i,\mathbb Z)\cong \mathbb Z,
\]
\[
        H^4(L_i,\mathbb Z)=0,
        \qquad
        H^5(L_i,\mathbb Z)\cong \mathbb Z.
\]
In particular,
\[
        H^m(L_i,\mathbb Z)_{\mathrm{tors}}=0
        \qquad
        \text{for all }m.
\]

Thus each node contributes free local topology but no finite local torsion
package.

\subsection{Birth of torsion}

For a normal surface singularity, the local group \(E\) is identified with
torsion in \(H^2\) of the link:
\[
        E\cong H^2(L,\mathbb Z)_{\mathrm{tors}}
\]
\cite{RahmanIntegralPerverseObstructions}.  A threefold node is not a
surface singularity, so that theorem cannot be applied directly.  However,
the direct link computation above shows that there is no torsion in any link
cohomology degree.  Therefore no finite local torsion package is born at a
threefold node.

For every node \(p_i\), we record
\[
        E_i^{\mathrm{tors}}=0.
\]
Consequently,
\[
        \bigoplus_{i=1}^s E_i^{\mathrm{tors}}=0.
\]
The invariant factor decomposition is empty, the order is \(1\), and there
is no prime support.

There is nevertheless a free local class at each node.  In a smoothing, this
is represented by a vanishing \(3\)-sphere.  In a small resolution, it is
related to the exceptional curve \(C_i\cong \mathbb P^1\).  These classes are
not torsion classes, so they belong to the free part of the nodal
topological story rather than to the finite torsion trajectory.

\subsection{The discriminant form}

Since
\[
        E_i^{\mathrm{tors}}=0
\]
for each node, there is no nontrivial finite discriminant form
\[
        q_i:E_i^{\mathrm{tors}}\times E_i^{\mathrm{tors}}
        \longrightarrow
        \mathbb Q/\mathbb Z.
\]
Thus the local torsion discriminant package at each node is
\[
        (E_i^{\mathrm{tors}},q_i)=(0,0).
\]

This does not mean that a node has no local intersection or vanishing-cycle
geometry.  It means that the local node does not produce a finite
\(\mathbb Q/\mathbb Z\)-valued torsion pairing of the type produced by
normal surface singularities.  The relevant local classes are free, not
finite.

\subsection{Perverse and torsion-sensitive realization}

Let
\[
        U_i:=X_i\setminus\{p_i\}
\]
be the punctured local germ at a node, and let
\[
        j_i:U_i\hookrightarrow X_i,
        \qquad
        i_i:\{p_i\}\hookrightarrow X_i
\]
be the inclusions.  Since the germ has complex dimension \(3\), the local
stalk model is
\[
        i_i^*R(j_i)_*\mathbb Z_{U_i}[3]
        \cong
        R\Gamma(L_i,\mathbb Z)[3].
\]
Therefore
\[
        H^m\bigl(i_i^*R(j_i)_*\mathbb Z_{U_i}[3]\bigr)
        \cong
        H^{m+3}(L_i,\mathbb Z).
\]
Using the cohomology of
\[
        L_i\cong S^2\times S^3,
\]
we obtain
\[
        H^m\bigl(i_i^*R(j_i)_*\mathbb Z_{U_i}[3]\bigr)
        =
        \begin{cases}
        \mathbb Z, & m=-3,\\
        \mathbb Z, & m=-1,\\
        \mathbb Z, & m=0,\\
        \mathbb Z, & m=2,\\
        0, & \text{otherwise}.
        \end{cases}
\]
There is no torsion in this stalk complex.

The ordinary and dual middle-perversity \(t\)-structures over
\(\mathbb Z\) differ by torsion and torsion-free conditions in critical
degrees, as in BBD's integral formalism \cite[Complement 3.3]{BBD82}.
Friedman's torsion-sensitive truncation framework detects the same kind of
critical-degree torsion behavior
\cite{FriedmanGenIH,FriedmanBook20,FriedmanTsInv}.  Since the stalk complex
above has no torsion, there is no finite torsion-sensitive discrepancy of
the type found in normal surface singularities.

Thus the local perverse/torsion-sensitive station for each threefold node is
\[
        E_i^{\mathrm{tors}}=0.
\]

\subsection{Resolution and pair-sequence realization}

Assume first that a small resolution is chosen near each node.  Let
\[
        N_i
\]
be a resolution neighborhood of the exceptional curve
\[
        C_i\cong\mathbb P^1.
\]
Then \(N_i\) deformation retracts onto \(C_i\), so
\[
        H_0(N_i,\mathbb Z)\cong\mathbb Z,
        \qquad
        H_2(N_i,\mathbb Z)\cong\mathbb Z,
\]
and all other reduced homology groups vanish.

The boundary is
\[
        L_i=\partial N_i\cong S^2\times S^3.
\]
Since \(N_i\) is a compact oriented real \(6\)-manifold with boundary,
Poincare--Lefschetz duality gives
\[
        H_j(N_i,L_i;\mathbb Z)
        \cong
        H^{6-j}(N_i,\mathbb Z).
\]
Because \(N_i\simeq S^2\), one has
\[
        H^0(N_i,\mathbb Z)\cong\mathbb Z,
        \qquad
        H^2(N_i,\mathbb Z)\cong\mathbb Z,
\]
and all other cohomology groups vanish.  Hence
\[
        H_4(N_i,L_i;\mathbb Z)\cong\mathbb Z,
        \qquad
        H_6(N_i,L_i;\mathbb Z)\cong\mathbb Z,
\]
and
\[
        H_2(N_i,L_i;\mathbb Z)=0,
        \qquad
        H_3(N_i,L_i;\mathbb Z)=0.
\]

The long exact sequence of the pair \((N_i,L_i)\) contains
\[
        H_3(N_i,\mathbb Z)
        \longrightarrow
        H_3(N_i,L_i;\mathbb Z)
        \longrightarrow
        H_2(L_i,\mathbb Z)
        \longrightarrow
        H_2(N_i,\mathbb Z)
        \longrightarrow
        H_2(N_i,L_i;\mathbb Z).
\]
Substituting the computed groups gives
\[
        0
        \longrightarrow
        0
        \longrightarrow
        H_2(L_i,\mathbb Z)
        \longrightarrow
        \mathbb Z
        \longrightarrow
        0.
\]
Therefore
\[
        H_2(L_i,\mathbb Z)\cong \mathbb Z.
\]
The sequence contains no finite cokernel.  Thus the pair-sequence station
confirms that the local node contributes free boundary topology, not finite
torsion.

Similarly, the portion
\[
        H_4(N_i,\mathbb Z)
        \longrightarrow
        H_4(N_i,L_i;\mathbb Z)
        \longrightarrow
        H_3(L_i,\mathbb Z)
        \longrightarrow
        H_3(N_i,\mathbb Z)
\]
becomes
\[
        0
        \longrightarrow
        \mathbb Z
        \longrightarrow
        H_3(L_i,\mathbb Z)
        \longrightarrow
        0,
\]
so
\[
        H_3(L_i,\mathbb Z)\cong \mathbb Z.
\]
Again this is free, not torsion.

Therefore the resolution and pair sequence for each local node yield no
finite discriminant quotient.

\subsection{Monodromy realization}

At each node, the local equation is
\[
        x_0^2+x_1^2+x_2^2+x_3^2=0.
\]
The Milnor fiber has the homotopy type of \(S^3\).  Hence
\[
        H^3(F_i,\mathbb Z)\cong\mathbb Z.
\]
The geometric monodromy is represented on the vanishing sphere by the
antipodal map on \(S^3\).  The antipodal map on \(S^n\) has degree
\[
        (-1)^{n+1}.
\]
For \(n=3\), this degree is \(1\).  Hence the monodromy acts by
\[
        T=\mathrm{id}
\]
on
\[
        H^3(F_i,\mathbb Z)\cong\mathbb Z.
\]
Thus
\[
        T-\mathrm{id}=0,
\]
and
\[
        \operatorname{coker}(T-\mathrm{id})\cong\mathbb Z.
\]
The torsion subgroup is therefore zero:
\[
        \operatorname{coker}(T-\mathrm{id})_{\mathrm{tors}}=0.
\]
This agrees with the link and pair-sequence stations: the node contributes a
free vanishing cycle, not local finite torsion.

\subsection{Global relations and defect}

Although each node has no local finite torsion package, a nodal threefold can
have global relations among its local free classes.  These relations are
what enter the defect theory.

Let
\[
        p_1,\ldots,p_s
\]
be the nodes.  Locally, each node has a vanishing sphere in a smoothing and,
in a small resolution, an exceptional curve.  The collection of these local
classes need not be globally independent.  A relation among them is a global
condition, not a local torsion class.

In the classical theory of nodal hypersurfaces, the defect measures a failure
of expected independence.  One way it appears is as the discrepancy between
the actual cohomology of the singular or resolved hypersurface and the
cohomology predicted by imposing independent local node conditions.  Clemens
introduced this perspective in the study of double solids
\cite{Clemens83}.  Cynk gives explicit formulations for the defect of nodal
hypersurfaces \cite{Cynk01}, and Kloosterman studies families of nodal
varieties with defect \cite{Kloosterman22}.

For the torsion trajectory, the important conclusion is the following:
\[
        \text{defect is a global relation phenomenon among free node data,
        not a local finite torsion birth.}
\]
Thus the nodal threefold row should not be read as a source of local
\(\mathbb Z/2\)-torsion.  It should be read as a warning that global
relations can exist even when local finite torsion does not.

\subsection{Transport status}

Since each node has
\[
        E_i^{\mathrm{tors}}=0,
\]
the direct sum of local finite torsion packages is
\[
        \bigoplus_{i=1}^s E_i^{\mathrm{tors}}=0.
\]
Therefore the finite torsion transport map
\[
        \alpha_X^{(3)}:
        \bigoplus_i E_i^{\mathrm{tors}}
        \longrightarrow
        H^3(\widetilde X,\mathbb Z)_{\mathrm{tors}}
\]
has zero source.  Consequently its image is zero:
\[
        \operatorname{im}\alpha_X^{(3)}=0.
\]
Thus ordinary nodes do not directly contribute local finite torsion classes
to the Brauer/unramified trajectory.

This does not rule out torsion in
\[
        H^3(\widetilde X,\mathbb Z)
\]
for a global resolution.  It only says that such torsion, if present, is not
coming from a direct sum of local node packages of the surface \(E\)-type.
It must arise from global topology or from additional singularity data.

In the nodal case, the local-to-global story therefore has two layers:
\[
        \text{finite torsion trajectory: trivial at each node},
\]
and
\[
        \text{free relation trajectory: governed by vanishing/exceptional
        cycle relations and defect}.
\]

This is different from the codimension-two transverse Coble situation.  A
transverse \(\frac14(1,1)\) stratum carries a finite local system
\[
        E_\Sigma\cong(\mathbb Z/4)_\Sigma,
\]
with an order-two Bockstein shadow
\[
        2E_\Sigma\cong(\mathbb Z/2)_\Sigma.
\]
An isolated nodal threefold has no such finite local system.

\subsection{Rational death}

For the finite torsion trajectory,
\[
        E_i^{\mathrm{tors}}=0
\]
at every node.  Hence
\[
        E_i^{\mathrm{tors}}\otimes_{\mathbb Z}\mathbb Q=0.
\]
The rational-death column is therefore zero for the finite torsion package.

However, this should not be confused with the free node data.  The vanishing
cycle group is free:
\[
        H^3(F_i,\mathbb Z)\cong\mathbb Z,
\]
and after tensoring with \(\mathbb Q\) it becomes
\[
        H^3(F_i,\mathbb Q)\cong\mathbb Q.
\]
Thus rationalization does not kill the free vanishing-cycle information.
It kills torsion, but the nodal threefold already had no local finite
torsion.  This is the central distinction between surface \(A_1\)-torsion
and threefold-node topology.

\subsection{Trajectory row}

The nodal threefold trajectory can be summarized as follows.

\begin{center}
\begingroup
\scriptsize
\setlength{\tabcolsep}{1.8pt}
\renewcommand{\arraystretch}{1.15}
\begin{tabularx}{\textwidth}{|p{0.15\textwidth}|c|c|p{0.12\textwidth}|p{0.10\textwidth}|X|p{0.10\textwidth}|c|}
\hline
\textbf{Example} & \textbf{\(E\)} & \textbf{\(q\)} & \textbf{Local} & \textbf{Supp.} & \textbf{Global image} & \textbf{Br/res.} & \textbf{\(\mathbb Q\)}
\\
\hline
nodal 3-fold &
\(0\) tors. &
\(0\) &
nodes have free, not torsion, local data &
none tors. &
finite torsion image zero; free relations measured by defect &
none local &
\(0\)
\\
\hline
\end{tabularx}
\endgroup
\end{center}

The entry ``none tors.'' means that the ordinary double point has no local
finite torsion support package.  The entry ``free relations measured by
defect'' records the separate global phenomenon: vanishing or exceptional
classes may satisfy global relations, and those relations are measured by
defect in the nodal hypersurface literature
\cite{Clemens83,Cynk01,Kloosterman22}.

\section{The Benoist--Ottem global torsion benchmark}
\label{app:BennoistOttemglobaltorsionbenchmark}

This appendix runs the Benoist--Ottem examples through the torsion-trajectory
machinery.  This is not a local singularity example.  It is a global smooth
torsion benchmark.

Benoist and Ottem prove that if \(S\) is an Enriques surface and \(C\) is a
very general smooth projective curve of genus at least \(1\), then
\[
        Y:=S\times C
\]
does not satisfy the integral Hodge conjecture for \(1\)-cycles.  Equivalently,
there are integral Hodge classes in
\[
        H^4(Y,\mathbb Z)
\]
which are not classes of algebraic \(1\)-cycles.  Their examples give
non-algebraic torsion cohomology classes of degree \(4\) on smooth projective
complex threefolds of Kodaira dimension zero
\cite{BenoistOttem20}.

The question for the present paper is whether this global \(2\)-torsion can
be related to a local discriminant package after degeneration.  The answer is
not automatic.  On the smooth variety \(Y=S\times C\), there are no singular
points, hence no local groups \(E_i\).  Thus the local-discriminant machinery
does not directly produce the Benoist--Ottem class on the smooth fiber.

The meaningful comparison problem is degenerational.  The corrected boundary
model is not a surface with transverse \(A_1\) singularities.  In the
degree-two Enriques compactification of Alexeev--Engel--Garza--Schaffler, the
relevant discriminant boundary parametrizes quotients of nodal K3 surfaces by
an involution fixing a node; the resulting boundary surfaces are rational
Coble surfaces with a singularity of type
\[
        \frac14(1,1)
\]
\cite{AlexeevEngelGarzaSchafflerEnriquesDegree2}.  For this singularity the
full local package is
\[
        E\cong\mathbb Z/4,
\]
while the Benoist--Ottem \(2\)-torsion mechanism sees the order-two
Bockstein shadow
\[
        2E\cong\mathbb Z/2.
\]
The purpose of this appendix is to make this corrected comparison precise
enough to identify what is known, what the trajectory machinery can see, and
what remains to be proved.

\subsection{The Benoist--Ottem model}

Let \(S\) be an Enriques surface.  We use the following standard facts.

\begin{enumerate}[label=\textup{(\arabic*)}]
\item The canonical bundle \(K_S\) is nontrivial but \(2\)-torsion:
\[
        K_S\not\cong\mathcal O_S,
        \qquad
        K_S^{\otimes 2}\cong\mathcal O_S.
\]

\item The universal cover of \(S\) is a K3 surface, and the covering map is
an étale double cover
\[
        \widetilde S\longrightarrow S.
\]
Consequently,
\[
        \pi_1(S)\cong\mathbb Z/2\mathbb Z.
\]

\item The irregularity and geometric genus vanish:
\[
        q(S)=h^1(S,\mathcal O_S)=0,
        \qquad
        p_g(S)=h^2(S,\mathcal O_S)=0.
\]
\end{enumerate}

These are standard facts about Enriques surfaces; see, for example,
\cite{BHPV04,BeauvilleComplexSurfaces}.

Let \(C\) be a smooth projective curve of genus
\[
        g\ge 1.
\]
Set
\[
        Y:=S\times C.
\]
Then \(Y\) is a smooth projective threefold.  Since \(K_S\) is torsion and
\(K_C\) is pulled back from the curve, \(Y\) has Kodaira dimension determined
by the curve factor.  Benoist--Ottem prove that for \(C\) very general of
genus at least \(1\), the integral Hodge conjecture for \(1\)-cycles on
\(Y\) fails \cite{BenoistOttem20}.

The relevant cohomological degree is
\[
        H^4(Y,\mathbb Z),
\]
because \(1\)-cycles on a smooth projective threefold have codimension \(2\)
and therefore map under the cycle class map to degree \(4\) cohomology:
\[
        \operatorname{cl}:
        CH^2(Y)
        \longrightarrow
        H^4(Y,\mathbb Z).
\]

\subsection{Cohomology of the Enriques surface}

We now record the integral cohomology of an Enriques surface needed for the
trajectory computation.

Since
\[
        \pi_1(S)\cong\mathbb Z/2\mathbb Z,
\]
the first homology group is
\[
        H_1(S,\mathbb Z)\cong\mathbb Z/2\mathbb Z.
\]
The universal coefficient theorem gives
\[
        H^1(S,\mathbb Z)
        \cong
        \operatorname{Hom}(H_1(S,\mathbb Z),\mathbb Z)=0.
\]
The same theorem gives a short exact sequence
\[
        0
        \longrightarrow
        \operatorname{Ext}^1_{\mathbb Z}(H_1(S,\mathbb Z),\mathbb Z)
        \longrightarrow
        H^2(S,\mathbb Z)
        \longrightarrow
        \operatorname{Hom}(H_2(S,\mathbb Z),\mathbb Z)
        \longrightarrow
        0.
\]
Since
\[
        \operatorname{Ext}^1_{\mathbb Z}(\mathbb Z/2\mathbb Z,\mathbb Z)
        \cong
        \mathbb Z/2\mathbb Z,
\]
we obtain a \(2\)-torsion subgroup in \(H^2(S,\mathbb Z)\).  This torsion
class is the first Chern class of the canonical bundle:
\[
        c_1(K_S)\in H^2(S,\mathbb Z)_{\mathrm{tors}},
        \qquad
        2c_1(K_S)=0.
\]
Thus
\[
        H^2(S,\mathbb Z)_{\mathrm{tors}}
        \cong
        \mathbb Z/2\mathbb Z.
\]

Poincare duality on the compact oriented real \(4\)-manifold \(S\) gives
\[
        H^3(S,\mathbb Z)
        \cong
        H_1(S,\mathbb Z)
        \cong
        \mathbb Z/2\mathbb Z.
\]
Also
\[
        H^0(S,\mathbb Z)\cong\mathbb Z,
        \qquad
        H^4(S,\mathbb Z)\cong\mathbb Z.
\]
The free part of \(H^2(S,\mathbb Z)\) has rank \(10\), but its precise
lattice structure is not needed for the torsion comparison below.

Hence the torsion part of the Enriques cohomology is
\[
        H^2(S,\mathbb Z)_{\mathrm{tors}}
        \cong
        \mathbb Z/2\mathbb Z,
        \qquad
        H^3(S,\mathbb Z)
        \cong
        \mathbb Z/2\mathbb Z.
\]

\subsection{Torsion in \texorpdfstring{\(H^4(S\times C,\mathbb Z)\)}{H4(SxC,Z)}}

We now compute the torsion subgroup of
\[
        H^4(Y,\mathbb Z),
        \qquad
        Y=S\times C.
\]
The integral cohomology of a smooth projective curve \(C\) is torsion-free:
\[
        H^0(C,\mathbb Z)\cong\mathbb Z,
        \qquad
        H^1(C,\mathbb Z)\cong\mathbb Z^{2g},
        \qquad
        H^2(C,\mathbb Z)\cong\mathbb Z.
\]
Therefore the Kunneth theorem has no \(\operatorname{Tor}\)-correction terms
coming from the curve factor.  We have
\[
        H^4(S\times C,\mathbb Z)
        \cong
        \bigoplus_{a+b=4}
        H^a(S,\mathbb Z)\otimes H^b(C,\mathbb Z).
\]
Since \(C\) has cohomology only in degrees \(0,1,2\), this becomes
\[
        H^4(Y,\mathbb Z)
        \cong
        H^4(S,\mathbb Z)\otimes H^0(C,\mathbb Z)
\]
\[
        \oplus
        H^3(S,\mathbb Z)\otimes H^1(C,\mathbb Z)
        \oplus
        H^2(S,\mathbb Z)\otimes H^2(C,\mathbb Z).
\]
Substituting the groups above, the torsion subgroup is
\[
        H^4(Y,\mathbb Z)_{\mathrm{tors}}
        \cong
        \bigl(H^3(S,\mathbb Z)\otimes H^1(C,\mathbb Z)\bigr)
        \oplus
        \bigl(H^2(S,\mathbb Z)_{\mathrm{tors}}\otimes H^2(C,\mathbb Z)\bigr).
\]
Since
\[
        H^3(S,\mathbb Z)\cong\mathbb Z/2\mathbb Z,
        \qquad
        H^1(C,\mathbb Z)\cong\mathbb Z^{2g},
\]
we have
\[
        H^3(S,\mathbb Z)\otimes H^1(C,\mathbb Z)
        \cong
        (\mathbb Z/2\mathbb Z)^{2g}.
\]
Also
\[
        H^2(S,\mathbb Z)_{\mathrm{tors}}\otimes H^2(C,\mathbb Z)
        \cong
        \mathbb Z/2\mathbb Z.
\]
Thus
\[
        H^4(Y,\mathbb Z)_{\mathrm{tors}}
        \cong
        (\mathbb Z/2\mathbb Z)^{2g}
        \oplus
        \mathbb Z/2\mathbb Z.
\]

This computation is important for two reasons.  First, it shows that the
Benoist--Ottem setting contains abundant \(2\)-torsion in degree \(4\).
Second, it separates two possible sources of degree-four torsion:
\[
        H^3(S,\mathbb Z)\otimes H^1(C,\mathbb Z)
\]
and
\[
        H^2(S,\mathbb Z)_{\mathrm{tors}}\otimes H^2(C,\mathbb Z).
\]

\begin{remark}
One should not identify the Benoist--Ottem obstruction merely with the class
\[
        c_1(K_S)\otimes [\mathrm{pt}]
        \in
        H^2(S,\mathbb Z)_{\mathrm{tors}}\otimes H^2(C,\mathbb Z).
\]
The class \(c_1(K_S)\) is the Chern class of an algebraic line bundle on
\(S\), and \([\mathrm{pt}]\) is algebraic on \(C\).  Their external product
is therefore represented by an algebraic cycle.  The non-algebraic torsion
classes in the Benoist--Ottem theorem are subtler global torsion classes in
degree \(4\) \cite{BenoistOttem20}.
\end{remark}

\subsection{Brauer comparison for \texorpdfstring{\(S\times C\)}{SxC}}

We now compute the Brauer-side torsion.  The exponential sequence gives a
comparison
\[
        \operatorname{Br}(Y)
        \cong
        H^3(Y,\mathbb Z)_{\mathrm{tors}}
\]
when
\[
        H^2(Y,\mathcal O_Y)=0.
\]
We verify this vanishing.

By the Kunneth formula for coherent cohomology,
\[
        H^2(Y,\mathcal O_Y)
        \cong
        \bigoplus_{a+b=2}
        H^a(S,\mathcal O_S)\otimes H^b(C,\mathcal O_C).
\]
The summands are
\[
        H^2(S,\mathcal O_S)\otimes H^0(C,\mathcal O_C),
\]
\[
        H^1(S,\mathcal O_S)\otimes H^1(C,\mathcal O_C),
\]
and
\[
        H^0(S,\mathcal O_S)\otimes H^2(C,\mathcal O_C).
\]
For an Enriques surface,
\[
        H^1(S,\mathcal O_S)=0,
        \qquad
        H^2(S,\mathcal O_S)=0.
\]
For a smooth curve,
\[
        H^2(C,\mathcal O_C)=0.
\]
Therefore every summand vanishes, and hence
\[
        H^2(Y,\mathcal O_Y)=0.
\]
The Brauer--torsion comparison gives
\[
        \operatorname{Br}(Y)
        \cong
        H^3(Y,\mathbb Z)_{\mathrm{tors}}.
\]

We compute \(H^3(Y,\mathbb Z)_{\mathrm{tors}}\).  Since \(H^*(C,\mathbb Z)\)
is torsion-free, the Kunneth theorem gives
\[
        H^3(Y,\mathbb Z)
        \cong
        H^3(S,\mathbb Z)\otimes H^0(C,\mathbb Z)
\]
\[
        \oplus
        H^2(S,\mathbb Z)\otimes H^1(C,\mathbb Z)
        \oplus
        H^1(S,\mathbb Z)\otimes H^2(C,\mathbb Z).
\]
We have
\[
        H^1(S,\mathbb Z)=0.
\]
Thus the torsion subgroup is
\[
        H^3(Y,\mathbb Z)_{\mathrm{tors}}
        \cong
        H^3(S,\mathbb Z)
        \oplus
        \bigl(H^2(S,\mathbb Z)_{\mathrm{tors}}\otimes H^1(C,\mathbb Z)\bigr).
\]
Using
\[
        H^3(S,\mathbb Z)\cong\mathbb Z/2\mathbb Z,
        \qquad
        H^2(S,\mathbb Z)_{\mathrm{tors}}\cong\mathbb Z/2\mathbb Z,
        \qquad
        H^1(C,\mathbb Z)\cong\mathbb Z^{2g},
\]
we obtain
\[
        H^3(Y,\mathbb Z)_{\mathrm{tors}}
        \cong
        \mathbb Z/2\mathbb Z
        \oplus
        (\mathbb Z/2\mathbb Z)^{2g}.
\]
Therefore
\[
        \operatorname{Br}(Y)
        \cong
        (\mathbb Z/2\mathbb Z)^{2g+1}.
\]

This is the Brauer station of the Benoist--Ottem trajectory.  It is global
and smooth; it is not produced by local singularities on \(Y\), since
\(Y\) is smooth.

\subsection{The direct local-discriminant test on the smooth fiber}

We now apply the local-discriminant machinery directly to
\[
        Y=S\times C.
\]
Since \(Y\) is smooth, it has no singular points.  Therefore there are no
local obstruction groups
\[
        E_i=H^0({}^p_+IC_X\mathbb Z)_{p_i}
\]
of the kind associated to singular points.  Hence
\[
        \bigoplus_i E_i=0.
\]
The direct local-to-global map
\[
        \alpha_Y^{(k)}:
        \bigoplus_iE_i
        \longrightarrow
        H^k(Y,\mathbb Z)_{\mathrm{tors}}
\]
therefore has zero source for every \(k\).  In particular,
\[
        \operatorname{im}\alpha_Y^{(3)}=0
\]
and
\[
        \operatorname{im}\alpha_Y^{(4)}=0
\]
in the direct smooth-fiber local-discriminant sense.

Thus the Benoist--Ottem torsion is not directly recovered by applying the
singularity-based local discriminant map to the smooth variety
\[
        S\times C.
\]
If there is a connection with local discriminant data, it must pass through a
degeneration or specialization in which singularities appear.

\subsection{The corrected degeneration test through the Coble boundary}

The corrected degeneration test passes through the Enriques/Coble boundary,
not through a direct \(A_1\) boundary model.

An Enriques surface can be realized as a quotient
\[
        S=K/\sigma,
\]
where \(K\) is a K3 surface and \(\sigma\) is a fixed-point-free involution.
Because the action is free, the quotient \(S\) is smooth, and its fundamental
group is
\[
        \pi_1(S)\cong \mathbb Z/2\mathbb Z.
\]
The \(2\)-torsion in the Enriques surface is therefore global and smooth.

In the degree-two Enriques compactification of
Alexeev--Engel--Garza--Schaffler, the relevant discriminant boundary
parametrizes quotients of nodal K3 surfaces by an involution fixing a node;
the resulting boundary surfaces are rational Coble surfaces with a
singularity of type
\[
        \frac14(1,1)
\]
\cite{AlexeevEngelGarzaSchafflerEnriquesDegree2}.  Thus the natural local
boundary singularity is not \(A_1\), but \(\frac14(1,1)\).

For the singularity
\[
        \frac14(1,1)=\mathbb C^2/\mu_4,
\]
the link is
\[
        L(4,1),
\]
and the full local discriminant package is
\[
        E\cong H^2(L(4,1),\mathbb Z)_{\mathrm{tors}}
        \cong
        \mathbb Z/4.
\]
The local index-two cover is
\[
        A_1=\mathbb C^2/\mu_2
        \longrightarrow
        \mathbb C^2/\mu_4=\frac14(1,1),
\]
and on links it is
\[
        L(2,1)\longrightarrow L(4,1).
\]
If
\[
        \eta\in H^1(L(4,1),\mathbb Z/2)
\]
is the class of this double cover, the Bockstein sends \(\eta\) to the
unique order-two subgroup of \(E\):
\[
        \beta(\eta)=2\bar g\in E\cong\mathbb Z/4.
\]
Thus the Benoist--Ottem \(2\)-torsion mechanism naturally sees
\[
        2E\cong\mathbb Z/2
        \subset
        E\cong\mathbb Z/4.
\]

Now form
\[
        Y_0:=S_0\times C,
\]
where \(S_0\) is such a Coble boundary surface.  If the singular points of
\(S_0\) are
\[
        p_1,\ldots,p_r,
\]
then the singular locus of \(Y_0\) is
\[
        \Sigma
        =
        \bigcup_{j=1}^r \Sigma_j,
        \qquad
        \Sigma_j:=\{p_j\}\times C.
\]
Thus \(Y_0\) is a threefold whose singular locus has codimension two, and
whose transverse singularity along each component \(\Sigma_j\) is of type
\[
        \frac14(1,1).
\]

The full transverse local discriminant package is expected to define a
locally constant \(\mathbb Z/4\)-sheaf
\[
        E_j\cong(\mathbb Z/4)_{\Sigma_j}.
\]
The Benoist--Ottem-visible part is the order-two sublocal system
\[
        2E_j\cong(\mathbb Z/2)_{\Sigma_j}.
\]
Its first cohomology is
\[
        H^1(\Sigma_j,2E_j)
        \cong
        H^1(C,\mathbb Z/2)
        \cong
        (\mathbb Z/2)^{2g}.
\]
This is the first point at which the shape of the Benoist--Ottem torsion
reappears from the corrected boundary model.  The Kunneth computation for
\(S\times C\) produced
\[
        (\mathbb Z/2)^{2g}
\]
as a global torsion summand.  The Coble degeneration produces the same shape
from the cohomology of the order-two shadow
\[
        2E_j\subset E_j
\]
along the curve
\[
        \Sigma_j\cong C.
\]

This observation does not prove that the Benoist--Ottem obstruction is local
in origin.  It identifies the precise test.  One must construct a
specialization diagram relating the smooth fiber
\[
        S\times C
\]
to a singular fiber
\[
        S_0\times C
\]
with transverse Coble \(\frac14(1,1)\)-strata, and then compare the
specialization of the Benoist--Ottem torsion class with the order-two
Bockstein shadows
\[
        2E_j\subset E_j.
\]

\subsection{The corrected specialization problem}

The preceding discussion leads to the following problem.

\begin{problem}[Benoist--Ottem/Coble specialization problem]
Let
\[
        \mathcal Y\longrightarrow \Delta
\]
be a one-parameter degeneration whose general fiber is
\[
        Y_t\cong S\times C,
\]
where \(S\) is an Enriques surface and \(C\) is a very general curve, and
whose special fiber has the form
\[
        Y_0\cong S_0\times C,
\]
where \(S_0\) is a rational Coble boundary surface with singularities of type
\[
        \frac14(1,1).
\]
Let
\[
        \Sigma_j=\{p_j\}\times C
\]
be the codimension-two Coble strata in \(Y_0\), and let
\[
        E_j\cong(\mathbb Z/4)_{\Sigma_j}
\]
be the transverse Coble discriminant local systems.  Let
\[
        2E_j\cong(\mathbb Z/2)_{\Sigma_j}
\]
be the order-two sublocal systems selected by the local index-two covers.

Determine whether the non-algebraic torsion classes in
\[
        H^4(Y_t,\mathbb Z)
\]
constructed by Benoist--Ottem specialize to classes generated by the
cohomology of the order-two local systems
\[
        H^*(\Sigma_j,2E_j).
\]
Equivalently, determine whether the Benoist--Ottem obstruction is genuinely
global smooth torsion, or whether it is the smoothing of the
Bockstein-selected order-two shadow inside the Coble
\(\frac14(1,1)\)-discriminant package along codimension-two strata.
\end{problem}

This problem is the corrected form of the comparison.  It is stronger and
more precise than the statement that both the Enriques surface and the
surface \(A_1\) singularity produce a group
\[
        \mathbb Z/2\mathbb Z.
\]
The actual Enriques/Coble boundary singularity has
\[
        E\cong\mathbb Z/4,
\]
and the relevant \(2\)-torsion is the subgroup
\[
        2E\cong\mathbb Z/2.
\]
Thus one must compare specialization, Bockstein maps, discriminant forms,
support maps, Brauer classes, and residue behavior.

\subsection{What the trajectory machinery proves and what it does not prove}

The torsion-trajectory machinery gives the following conclusions.

\begin{enumerate}[label=\textup{(\arabic*)}]
\item On the smooth fiber
\[
        Y=S\times C,
\]
there is no local singularity package.  Hence the direct local-discriminant
map has zero source:
\[
        \bigoplus_iE_i=0.
\]

\item The smooth fiber nevertheless has global \(2\)-torsion in
\[
        H^4(Y,\mathbb Z)
\]
and in
\[
        H^3(Y,\mathbb Z),
\]
computed above by the Kunneth formula.

\item Since
\[
        H^2(Y,\mathcal O_Y)=0,
\]
the Brauer comparison gives
\[
        \operatorname{Br}(Y)
        \cong
        H^3(Y,\mathbb Z)_{\mathrm{tors}}
        \cong
        (\mathbb Z/2\mathbb Z)^{2g+1}.
\]

\item If \(Y\) degenerates to a product
\[
        S_0\times C
\]
where \(S_0\) has Coble \(\frac14(1,1)\) singularities, then each
codimension-two stratum
\[
        \{p_j\}\times C
\]
carries a transverse local discriminant system
\[
        E_j\cong(\mathbb Z/4)_{\Sigma_j}.
\]

\item The local index-two cover selects the order-two sublocal system
\[
        2E_j\cong(\mathbb Z/2)_{\Sigma_j}.
\]

\item The cohomology
\[
        H^1(C,\mathbb Z/2\mathbb Z)
        \cong
        (\mathbb Z/2\mathbb Z)^{2g}
\]
has the same shape as the Kunneth torsion summand
\[
        H^3(S,\mathbb Z)\otimes H^1(C,\mathbb Z)
        \cong
        (\mathbb Z/2\mathbb Z)^{2g}.
\]
\end{enumerate}

The machinery does not, by itself, prove that the Benoist--Ottem class is
the specialization of a local Coble shadow class.  To prove that, one would
need the following additional data.

\begin{enumerate}[label=\textup{(\roman*)}]
\item A concrete degeneration
\[
        \mathcal Y\to\Delta
\]
from \(S\times C\) to \(S_0\times C\) with transverse Coble
\(\frac14(1,1)\)-strata.

\item A specialization map on integral cohomology, or on the appropriate
nearby/vanishing-cycle complexes, that preserves the relevant torsion data.

\item A comparison between the Benoist--Ottem non-algebraic torsion class
and the image of the order-two local systems
\[
        2E_j\subset E_j.
\]

\item A residue comparison showing whether the associated Brauer or
unramified class is the same on the smooth and singular sides.
\end{enumerate}

These requirements identify the exact missing theorem.

\subsection{Residue and unramified interpretation}

The Colliot-Thélène--Voisin framework relates failures of the integral Hodge
conjecture in degree \(4\) to unramified cohomology with
\(\mathbb Q/\mathbb Z\)-coefficients \cite{CTVoisin12}.  Benoist--Ottem's
examples lie in this degree-four integral Hodge setting
\cite{BenoistOttem20}.  Therefore the torsion classes above should also be
viewed through the residue station of the trajectory.

On the smooth fiber \(Y=S\times C\), the global torsion is not local
singularity torsion.  Its unramified behavior is controlled by the usual
Bloch--Ogus residue maps.  On a singular degeneration \(Y_0=S_0\times C\),
the candidate local contribution would come from the order-two shadows
\[
        2E_j\subset E_j
\]
along the Coble strata \(\Sigma_j\).  The residue comparison problem is to
determine whether the smooth unramified class is the specialization of a
residue-free order-two local-discriminant shadow on the singular fiber.

Thus the residue station gives a second test.  Even if the Coble shadow
cohomology produces the correct \((\mathbb Z/2)^{2g}\)-shaped group, one
must still check whether the corresponding classes have the same unramified
residue behavior as the Benoist--Ottem obstruction.

\subsection{Rational death}

All torsion classes discussed above vanish after tensoring with
\[
        \mathbb Q.
\]
For example,
\[
        H^4(Y,\mathbb Z)_{\mathrm{tors}}\otimes_{\mathbb Z}\mathbb Q=0,
\]
and
\[
        \operatorname{Br}(Y)\otimes_{\mathbb Z}\mathbb Q=0.
\]
If
\[
        \tau
\]
is a \(2\)-torsion class, then
\[
        2\tau=0,
\]
so in the tensor product with \(\mathbb Q\),
\[
        \tau\otimes 1
        =
        \tau\otimes 2\cdot\frac12
        =
        2\tau\otimes\frac12
        =
        0.
\]

This is the same final disappearance seen in the local \(A_1\) and Coble
trajectories.  The difference is the origin of the torsion.  In the
\(A_1\) surface case, the torsion is born locally as
\[
        \Lambda^\vee/\Lambda\cong\mathbb Z/2.
\]
In the Coble boundary case, the full local package is born as
\[
        E\cong\mathbb Z/4,
\]
and the Benoist--Ottem-visible part is
\[
        2E\cong\mathbb Z/2.
\]
In the Benoist--Ottem smooth fiber, the torsion is born globally from the
Enriques topology and its product with the curve.  The open question is
whether these origins can be connected by degeneration.

\subsection{Trajectory row}

The Benoist--Ottem benchmark can be summarized as follows.

\begin{center}
\begingroup
\scriptsize
\setlength{\tabcolsep}{2pt}
\renewcommand{\arraystretch}{1.25}

\begin{tabular}{|
>{\raggedright\arraybackslash}p{0.13\textwidth}|
>{\raggedright\arraybackslash}p{0.12\textwidth}|
>{\raggedright\arraybackslash}p{0.08\textwidth}|
>{\raggedright\arraybackslash}p{0.13\textwidth}|
>{\raggedright\arraybackslash}p{0.08\textwidth}|
>{\raggedright\arraybackslash}p{0.16\textwidth}|
>{\raggedright\arraybackslash}p{0.15\textwidth}|
>{\centering\arraybackslash}p{0.04\textwidth}|}
\hline
\textbf{Example}
&
\textbf{\(E\)}
&
\textbf{\(q\)}
&
\textbf{Local}
&
\textbf{Supp.}
&
\textbf{Global image}
&
\textbf{Br/res.}
&
\textbf{\(\mathbb Q\)}
\\
\hline

\(S\times C\) BO
&
none on smooth fiber
&
none local
&
global Enriques \(2\)-torsion
&
none local
&
\(H^4_{\mathrm{tors}}\neq0\); not from direct \(\alpha\)
&
\(\operatorname{Br}(Y)\cong\) \(H^3_{\mathrm{tors}}\); residues global
&
\(0\)
\\
\hline

Coble \(S_0\times C\) boundary
&
\(\mathbb Z/4\) transversely; BO sees \(2E\cong\mathbb Z/2\)
&
\(-1/4\) on full package
&
local Coble \(E\)-sheaf; order-two shadow on strata
&
codim \(2\) strata
&
candidate source \(H^1(C,2E)\cong\) \(H^1(C,\mathbb Z/2)\)
&
requires residue/ specialization comparison
&
\(0\)
\\
\hline
\end{tabular}

\endgroup
\end{center}

The first row is the smooth Benoist--Ottem example.  It has global
\(2\)-torsion but no local singularity package.  The second row is the
corrected degeneration test.  In the Enriques/Coble boundary model, the
local singularity is
\[
        \frac14(1,1),
\]
so the full transverse package is
\[
        E\cong\mathbb Z/4.
\]
The Benoist--Ottem \(2\)-torsion direction corresponds to the
Bockstein-selected order-two shadow
\[
        2E\cong\mathbb Z/2.
\]
After product with \(C\), the order-two shadow produces
\[
        H^1(C,2E)
        \cong
        H^1(C,\mathbb Z/2)
        \cong
        (\mathbb Z/2)^{2g}.
\]
This has the same group shape as the global Kunneth torsion appearing in the
smooth Benoist--Ottem setting.  Proving that the two are the same requires a
specialization theorem comparing the smooth torsion class with the order-two
Coble discriminant shadow on the singular fiber.

\section{The \texorpdfstring{\(\frac14(1,1)\)}{1/4(1,1)} Coble boundary trajectory}
\label{app:Cobleboundarytrajectory}

This appendix records the torsion trajectory for the cyclic quotient
singularity of type \(\frac14(1,1)\).  This example is the corrected local
model suggested by the Enriques/Coble boundary.  In the degree-two Enriques
compactification, the relevant discriminant boundary includes quotients of
nodal K3 surfaces by an involution fixing a node, producing rational Coble
surfaces with a \(\frac14(1,1)\)-singularity
\cite{AlexeevEngelGarzaSchafflerEnriquesDegree2}.

The key point is that the full local discriminant group is
\[
        E\cong \mathbb Z/4\mathbb Z,
\]
while the Benoist--Ottem \(2\)-torsion mechanism naturally selects the
unique order-two subgroup
\[
        2E\cong\mathbb Z/2\mathbb Z.
\]
Thus the Benoist--Ottem comparison should not be phrased as a direct
\(A_1\)-comparison.  It should be phrased as an order-two Bockstein shadow
inside the \(\frac14(1,1)\) local discriminant package.  The \(A_1\)
singularity still appears, but one level upstairs: it is the local
index-two cover of \(\frac14(1,1)\).

\subsection{Local model}

Let
\[
        X=\mathbb C^2/\mu_4,
\]
where \(\mu_4=\langle \zeta\rangle\) acts by
\[
        \zeta\cdot (x,y)=(\zeta x,\zeta y),
        \qquad
        \zeta^4=1.
\]
This is the cyclic quotient surface singularity of type
\[
        \frac14(1,1).
\]
It is a normal surface singularity.  Unlike the \(A_1\) surface singularity,
which is the quotient
\[
        \mathbb C^2/\{\pm1\}=\frac12(1,1),
\]
this singularity has local quotient group \(\mathbb Z/4\mathbb Z\).

The subgroup
\[
        \mu_2=\{\pm1\}\subset \mu_4
\]
defines an intermediate quotient
\[
        \mathbb C^2/\mu_2
        \longrightarrow
        \mathbb C^2/\mu_4.
\]
The source is the \(A_1\) surface singularity.  Thus the local index-two
cover of the Coble boundary singularity is
\[
        A_1=\mathbb C^2/\mu_2
        \longrightarrow
        \mathbb C^2/\mu_4=\frac14(1,1).
\]
This cover is the local geometric source of the order-two Bockstein shadow
used below.

\subsection{The link and its cohomology}

Let \(L\) be the link of \(X\) at the singular point.  Since the action of
\(\mu_4\) on \(\mathbb C^2\setminus\{0\}\) restricts to the standard free
action on \(S^3\), the link is the lens space
\[
        L\cong S^3/\mu_4=L(4,1).
\]
The homology of \(L(4,1)\) is
\[
        H_0(L,\mathbb Z)\cong\mathbb Z,
        \qquad
        H_1(L,\mathbb Z)\cong\mathbb Z/4\mathbb Z,
\]
\[
        H_2(L,\mathbb Z)=0,
        \qquad
        H_3(L,\mathbb Z)\cong\mathbb Z.
\]
By the universal coefficient theorem,
\[
        H^0(L,\mathbb Z)\cong\mathbb Z,
        \qquad
        H^1(L,\mathbb Z)=0,
\]
\[
        H^2(L,\mathbb Z)\cong
        \operatorname{Ext}^1(H_1(L,\mathbb Z),\mathbb Z)
        \cong
        \mathbb Z/4\mathbb Z,
\]
and
\[
        H^3(L,\mathbb Z)\cong\mathbb Z.
\]
Therefore the local torsion group detected by link cohomology is
\[
        H^2(L,\mathbb Z)_{\mathrm{tors}}
        \cong
        \mathbb Z/4\mathbb Z.
\]

The index-two cover described above restricts on links to
\[
        S^3/\mu_2
        \longrightarrow
        S^3/\mu_4,
\]
that is,
\[
        L(2,1)\longrightarrow L(4,1).
\]
This is the double cover whose mod-\(2\) class will generate the
Benoist--Ottem order-two shadow.

\subsection{Birth of torsion}

By the local realization theorem for normal surface singularities,
\[
        E_X(0)
        \cong
        H^2(L,\mathbb Z)_{\mathrm{tors}}.
\]
Hence
\[
        E_X(0)\cong\mathbb Z/4\mathbb Z.
\]
This is the full local discriminant group of the \(\frac14(1,1)\)
singularity.

The prime support is only the prime \(2\), but the group is not merely
\[
        \mathbb Z/2\mathbb Z.
\]
The local singularity carries order-four discriminant data:
\[
        |E_X(0)|=4.
\]
The Benoist--Ottem \(2\)-torsion mechanism will see only the order-two
subgroup \(2E\), not the full group \(E\).

\subsection{The discriminant form}

The minimal resolution of the cyclic quotient singularity \(\frac14(1,1)\)
is described by the Hirzebruch--Jung continued fraction
\[
        \frac41=[4].
\]
Thus the exceptional locus consists of a single smooth rational curve \(C\)
with self-intersection
\[
        C^2=-4.
\]
The exceptional lattice is therefore
\[
        \Lambda=\mathbb Z\langle e\rangle,
        \qquad
        (e,e)=-4.
\]
Its dual lattice is generated by
\[
        e^\vee=-\frac{e}{4},
\]
since
\[
        (e^\vee,e)=1.
\]
Therefore
\[
        \Lambda^\vee/\Lambda
        \cong
        \mathbb Z/4\mathbb Z.
\]

Let \(\bar g\) denote the class of \(e^\vee\) in
\(\Lambda^\vee/\Lambda\).  In the geometric negative-definite convention,
the discriminant form is
\[
        q(\bar g,\bar g)
        =
        (e^\vee,e^\vee)
        =
        \left(-\frac e4,-\frac e4\right)
        =
        -\frac14
        \quad \bmod \mathbb Z.
\]
Thus
\[
        (E,q)\cong
        \left(\mathbb Z/4\mathbb Z,\ q(\bar g,\bar g)=-\frac14\right).
\]
With the opposite orientation/sign convention, the same form is represented
by \(+\frac14\).  Throughout this paper, the negative sign is consistent
with the geometric intersection convention.

The unique order-two element is
\[
        2\bar g\in \mathbb Z/4\mathbb Z.
\]
Its discriminant self-pairing is
\[
        q(2\bar g,2\bar g)
        =
        4q(\bar g,\bar g)
        =
        -1
        \equiv 0
        \quad \bmod \mathbb Z.
\]
Thus the order-two subgroup
\[
        2E\cong\mathbb Z/2\mathbb Z
\]
is isotropic for the associated discriminant pairing.

\subsection{Perverse and torsion-sensitive realization}

Let \(U=X\setminus\{0\}\), and let
\[
        j:U\hookrightarrow X,
        \qquad
        i:\{0\}\hookrightarrow X
\]
be the inclusions.  Since \(X\) is a normal surface singularity, the local
stalk computation gives
\[
        i^*Rj_*\mathbb Z_U[2]
        \cong
        R\Gamma(L,\mathbb Z)[2].
\]
Therefore
\[
        H^0\bigl(i^*Rj_*\mathbb Z_U[2]\bigr)
        \cong
        H^2(L,\mathbb Z).
\]
The critical torsion group is
\[
        H^2(L,\mathbb Z)_{\mathrm{tors}}
        \cong
        \mathbb Z/4\mathbb Z.
\]

The ordinary middle-perversity extension and the dual middle-perversity
extension differ integrally by the torsion retained in this critical degree.
Thus
\[
        \operatorname{Cone}
        \left(
        {}^pIC_X\mathbb Z
        \longrightarrow
        {}^p_+IC_X\mathbb Z
        \right)
        \cong
        i_*E_X(0)[1],
\]
with
\[
        E_X(0)\cong \mathbb Z/4\mathbb Z.
\]
Equivalently,
\[
        {}^pIC_X\mathbb Z
        \longrightarrow
        {}^p_+IC_X\mathbb Z
        \longrightarrow
        i_*(\mathbb Z/4\mathbb Z)[1]
        \overset{+1}{\longrightarrow}.
\]

In the torsion-sensitive truncation interpretation, the ordinary package
kills the critical torsion, while the dual package retains it.  The retained
torsion is precisely
\[
        \mathbb Z/4\mathbb Z.
\]

\subsection{Resolution lattice realization}

The resolution lattice is
\[
        \Lambda=\mathbb Z\langle e\rangle,
        \qquad
        (e,e)=-4.
\]
The intersection matrix is
\[
        M=[-4].
\]
Its Smith normal form is
\[
        [4].
\]
Therefore
\[
        \operatorname{coker}(M:\mathbb Z\to\mathbb Z)
        \cong
        \mathbb Z/4\mathbb Z.
\]
Equivalently,
\[
        \Lambda^\vee/\Lambda
        \cong
        \mathbb Z/4\mathbb Z.
\]
This agrees with the link realization:
\[
        \Lambda^\vee/\Lambda
        \cong
        H^2(L,\mathbb Z)_{\mathrm{tors}}
        \cong
        E_X(0).
\]

\subsection{Pair-sequence realization}

Let \(N\) be a sufficiently small resolution neighborhood of the exceptional
curve \(C\), and let
\[
        \partial N=L.
\]
The pair sequence for \((N,L)\) identifies the boundary torsion with the
cokernel of the intersection pairing on the exceptional lattice.  Since
\(N\) deformation retracts onto \(C\), we have
\[
        H_2(N,\mathbb Z)\cong\mathbb Z\langle C\rangle.
\]
By Poincare--Lefschetz duality,
\[
        H_2(N,L;\mathbb Z)\cong H^2(N,\mathbb Z)\cong \mathbb Z.
\]
Under this identification, the map
\[
        H_2(N,\mathbb Z)
        \longrightarrow
        H_2(N,L;\mathbb Z)
\]
is represented by the relative intersection pairing, hence by
\[
        [-4].
\]
Therefore its cokernel is
\[
        \operatorname{coker}([-4])
        \cong
        \mathbb Z/4\mathbb Z.
\]
This is the boundary term detected on the link:
\[
        H_1(L,\mathbb Z)\cong\mathbb Z/4\mathbb Z,
\]
or equivalently by universal coefficients,
\[
        H^2(L,\mathbb Z)_{\mathrm{tors}}
        \cong
        \mathbb Z/4\mathbb Z.
\]

Thus the pair-sequence realization recovers the same group:
\[
        E_X(0)\cong\mathbb Z/4\mathbb Z.
\]

\subsection{Monodromy realization}

The singularity \(\frac14(1,1)\) is a cyclic quotient surface singularity,
not an isolated hypersurface surface singularity in the sense used for the
Wang-sequence/variation-map realization.  Therefore the monodromy station
\[
        \operatorname{coker}(T-\mathrm{id})_{\mathrm{tors}}
\]
is not part of the standard hypersurface realization package for this
example.

The appropriate local realizations here are the perverse discrepancy,
torsion-sensitive truncation, link torsion, resolution lattice, and
pair-sequence boundary term.  These all recover
\[
        \mathbb Z/4\mathbb Z.
\]

\subsection{The Bockstein order-two shadow}

Although the full local discriminant group is
\[
        E\cong\mathbb Z/4\mathbb Z,
\]
the Benoist--Ottem mechanism is \(2\)-torsion.  The relevant order-two
subgroup arises from a Bockstein.

The link is
\[
        L=L(4,1)=S^3/\mu_4.
\]
There is a canonical double cover
\[
        L(2,1)\longrightarrow L(4,1),
\]
corresponding to the subgroup
\[
        \mu_2\subset \mu_4.
\]
This double cover is classified by the nonzero class
\[
        \eta\in H^1(L,\mathbb Z/2).
\]
Since
\[
        H_1(L,\mathbb Z)\cong\mathbb Z/4\mathbb Z,
\]
we have
\[
        H^1(L,\mathbb Z/2)
        \cong
        \operatorname{Hom}(\mathbb Z/4,\mathbb Z/2)
        \cong
        \mathbb Z/2.
\]

Consider the coefficient sequence
\[
        0
        \longrightarrow
        \mathbb Z
        \xrightarrow{2}
        \mathbb Z
        \longrightarrow
        \mathbb Z/2
        \longrightarrow
        0.
\]
The associated Bockstein gives
\[
        \beta:
        H^1(L,\mathbb Z/2)
        \longrightarrow
        H^2(L,\mathbb Z).
\]
The relevant part of the long exact sequence is
\[
        H^1(L,\mathbb Z)
        \longrightarrow
        H^1(L,\mathbb Z/2)
        \xrightarrow{\beta}
        H^2(L,\mathbb Z)
        \xrightarrow{2}
        H^2(L,\mathbb Z).
\]
Since
\[
        H^1(L,\mathbb Z)=0
\]
and
\[
        H^2(L,\mathbb Z)\cong\mathbb Z/4\mathbb Z,
\]
the map \(\beta\) is injective, and exactness gives
\[
        \operatorname{im}(\beta)
        =
        \ker
        \left(
        2:\mathbb Z/4\mathbb Z
        \to
        \mathbb Z/4\mathbb Z
        \right)
        =
        \{0,2\}
        \cong
        \mathbb Z/2\mathbb Z.
\]
Therefore
\[
        \beta(\eta)=2\bar g
        \in
        E\cong\mathbb Z/4\mathbb Z.
\]
Thus the Bockstein of the double-cover class selects the unique order-two
subgroup
\[
        2E\cong\mathbb Z/2\mathbb Z
        \subset
        E\cong\mathbb Z/4\mathbb Z.
\]

This is the corrected local model for the Enriques/Coble boundary:
the full local discriminant group is \(\mathbb Z/4\), while the
Benoist--Ottem \(2\)-torsion mechanism detects the order-two Bockstein
shadow.

\subsection{The associated short exact sequence and triangle}

The inclusion of the order-two subgroup gives a short exact sequence
\[
        0
        \longrightarrow
        2E
        \longrightarrow
        E
        \longrightarrow
        E/2E
        \longrightarrow
        0.
\]
For
\[
        E\cong\mathbb Z/4,
\]
this is
\[
        0
        \longrightarrow
        \mathbb Z/2
        \longrightarrow
        \mathbb Z/4
        \longrightarrow
        \mathbb Z/2
        \longrightarrow
        0.
\]
If \(E_\Sigma\) is the corresponding local system along a codimension-two
stratum \(\Sigma\), this becomes
\[
        0
        \longrightarrow
        (2E)_\Sigma
        \longrightarrow
        E_\Sigma
        \longrightarrow
        (E/2E)_\Sigma
        \longrightarrow
        0.
\]
After applying \(i_*\) and shifting by \([2]\), one obtains the
distinguished triangle
\[
        i_*(2E)_\Sigma[2]
        \longrightarrow
        i_*E_\Sigma[2]
        \longrightarrow
        i_*(E/2E)_\Sigma[2]
        \overset{+1}{\longrightarrow}.
\]
In the Coble case this is
\[
        i_*(\mathbb Z/2)_\Sigma[2]
        \longrightarrow
        i_*(\mathbb Z/4)_\Sigma[2]
        \longrightarrow
        i_*(\mathbb Z/2)_\Sigma[2]
        \overset{+1}{\longrightarrow}.
\]
This triangle records categorically that the Benoist--Ottem shadow is a
filtered piece of the full Coble \(E\)-package.

\subsection{Transport status}

If a threefold has a codimension-two stratum
\[
        i:\Sigma\hookrightarrow Y
\]
whose transverse singularity is of type \(\frac14(1,1)\), then the transverse
local discriminant package has fiber
\[
        E\cong \mathbb Z/4\mathbb Z.
\]
In the locally trivial case, these groups assemble into a finite local
system
\[
        E_\Sigma\cong(\mathbb Z/4\mathbb Z)_\Sigma.
\]
By MacPherson--Vilonen or Beilinson gluing, the closed-stratum discrepancy
term is
\[
        i_*E_\Sigma[2].
\]
Thus the full transverse Coble contribution is
\[
        i_*(\mathbb Z/4)_\Sigma[2].
\]

The quotient-cover structure distinguishes the order-two subgroup
\[
        2E_\Sigma\cong(\mathbb Z/2\mathbb Z)_\Sigma.
\]
Thus the transport picture has two levels:
\[
        \text{full local discriminant system: }
        E_\Sigma\cong(\mathbb Z/4)_\Sigma,
\]
and
\[
        \text{Benoist--Ottem order-two shadow: }
        2E_\Sigma\cong(\mathbb Z/2)_\Sigma.
\]

For a product boundary stratum
\[
        \Sigma=\{p\}\times C,
\]
one obtains
\[
        H^1(\Sigma,2E_\Sigma)
        \cong
        H^1(C,\mathbb Z/2).
\]
If \(C\) has genus \(g\), then
\[
        H^1(C,\mathbb Z/2)\cong(\mathbb Z/2)^{2g}.
\]
This is the same group shape as the Benoist--Ottem middle Kunneth torsion
direction
\[
        H^3(S,\mathbb Z)\otimes H^1(C,\mathbb Z)
        \cong
        (\mathbb Z/2)^{2g}.
\]

\subsection{Motivic lift}

The full Coble package also has a natural motivic lift in any integral
motivic sheaf formalism with cones, shifts, proper pushforward, and a Betti
realization compatible with these operations.  Let \(\mathbf 1_\Sigma\) be
the motivic unit on \(\Sigma\).  Define
\[
        \mathbf 1_\Sigma/4
        :=
        \operatorname{Cone}
        \left(
        \mathbf 1_\Sigma
        \xrightarrow{4}
        \mathbf 1_\Sigma
        \right).
\]
Then the full motivic Coble object is
\[
        \mathcal T^{\mathrm{mot}}_{4,\Sigma}
        :=
        i_*(\mathbf 1_\Sigma/4)[2],
\]
with Betti realization
\[
        \operatorname{Real}_B
        \left(
        \mathcal T^{\mathrm{mot}}_{4,\Sigma}
        \right)
        \cong
        i_*(\mathbb Z/4)_\Sigma[2].
\]

The Benoist--Ottem order-two shadow is represented motivically by
\[
        \mathcal T^{\mathrm{mot}}_{2,\Sigma}
        :=
        i_*(\mathbf 1_\Sigma/2)[2].
\]
The natural map
\[
        \mathbf 1_\Sigma/2
        \longrightarrow
        \mathbf 1_\Sigma/4
\]
realizes the inclusion
\[
        \mathbb Z/2\hookrightarrow\mathbb Z/4,
        \qquad
        1\longmapsto2.
\]
After applying \(i_*\) and shifting by \([2]\), this gives
\[
        \mathcal T^{\mathrm{mot}}_{2,\Sigma}
        \longrightarrow
        \mathcal T^{\mathrm{mot}}_{4,\Sigma}.
\]
Under Betti realization, this becomes
\[
        i_*(\mathbb Z/2)_\Sigma[2]
        \longrightarrow
        i_*(\mathbb Z/4)_\Sigma[2],
\]
whose image is
\[
        i_*(2E_\Sigma)[2].
\]
Thus the motivic lift preserves the distinction between the full Coble
package and the Benoist--Ottem-visible order-two shadow.

\subsection{Rational death}

The full local discriminant group dies after rationalization:
\[
        E\otimes_{\mathbb Z}\mathbb Q
        =
        (\mathbb Z/4\mathbb Z)\otimes_{\mathbb Z}\mathbb Q
        =
        0.
\]
The order-two Benoist--Ottem shadow also dies rationally:
\[
        (2E)\otimes_{\mathbb Z}\mathbb Q
        =
        (\mathbb Z/2\mathbb Z)\otimes_{\mathbb Z}\mathbb Q
        =
        0.
\]
Thus this example again illustrates the central principle of the paper:
torsion disappears over \(\mathbb Q\), but before it disappears it records
the difference between the full local discriminant package
\(\mathbb Z/4\) and the order-two Bockstein shadow selected by a double-cover
class.

\subsection{Trajectory row}

\begin{center}
\scriptsize
\setlength{\tabcolsep}{2pt}
\renewcommand{\arraystretch}{1.25}
\begin{tabularx}{\textwidth}{|>{\centering\arraybackslash}p{0.16\textwidth}|
>{\centering\arraybackslash}p{0.13\textwidth}|
>{\centering\arraybackslash}p{0.18\textwidth}|
>{\centering\arraybackslash}p{0.17\textwidth}|
>{\centering\arraybackslash}p{0.14\textwidth}|
>{\centering\arraybackslash}X|}
\hline
\textbf{Example} &
\textbf{Birth} &
\textbf{Form} &
\textbf{Realizations} &
\textbf{BO shadow} &
\textbf{Rationalization} \\
\hline
\(\frac14(1,1)\) Coble boundary &
\(E\cong\mathbb Z/4\) &
\(q(\bar g,\bar g)=-\frac14\) &
link, lattice, pair, perverse, ts agree; monodromy not applicable &
\(\beta(\eta)=2\bar g\), so \(2E\cong\mathbb Z/2\) &
\(E\otimes\mathbb Q=0\), \(2E\otimes\mathbb Q=0\) \\
\hline
\end{tabularx}
\end{center}

\section{Nodal quintic threefolds}
\label{app:nodalquintic3fold}

This appendix specializes the nodal-threefold trajectory to nodal quintic
threefolds.  The purpose is to record what the torsion-trajectory machinery
does and does not see in this class of examples.

The conclusion is the same as in the general nodal threefold appendix:
ordinary double points in complex dimension three do not produce local
finite \(E\)-type torsion packages.  The local link of a threefold node is
\[
        S^2\times S^3,
\]
and hence has torsion-free integral cohomology.  Therefore the direct local
finite torsion source
\[
        \bigoplus_i E_i^{\mathrm{tors}}
\]
is zero.  What nodal quintics do have is free vanishing-cycle and exceptional
curve data, together with possible global relations among those free classes.
These global relations are the geometric content measured by defect.

This appendix is therefore a contrast case.  It prevents one from importing
the surface \(A_1\) calculation
\[
        E_{A_1}\cong\mathbb Z/2
\]
into the isolated threefold-node setting.  It also clarifies why the
higher-dimensional continuation of surface torsion packages should be sought
in codimension-two transverse surface strata, such as transverse
\(\frac14(1,1)\) Coble strata, rather than in isolated ordinary double
points.

\subsection{The global model}

Let
\[
        X\subset \mathbb P^4
\]
be a quintic hypersurface with isolated ordinary double points
\[
        \operatorname{Sing}(X)=\{p_1,\ldots,p_s\}.
\]
Thus, analytically locally at each \(p_i\), the germ is isomorphic to
\[
        \{x_0^2+x_1^2+x_2^2+x_3^2=0\}
        \subset
        (\mathbb C^4,0).
\]
The smooth quintic threefold is a Calabi--Yau threefold.  A nodal quintic is
singular, but its singularities are among the simplest isolated hypersurface
threefold singularities.  The local topology of ordinary double points is
standard in Milnor's theory of isolated hypersurface singularities
\cite{Mi68}.  The role of defects and relations among nodes in nodal
hypersurfaces is studied classically by Clemens \cite{Clemens83} and in
explicit nodal hypersurface form by Cynk \cite{Cynk01}; see also
Kloosterman \cite{Kloosterman22} for modern families with defect.

Let
\[
        \pi:\widetilde X\longrightarrow X
\]
be a resolution.  If a small resolution is chosen locally at a node, the
exceptional fiber over that node is a rational curve
\[
        C_i\cong \mathbb P^1
\]
with normal bundle
\[
        \mathcal O_{\mathbb P^1}(-1)
        \oplus
        \mathcal O_{\mathbb P^1}(-1).
\]
If one instead takes the blowup of the node, the exceptional divisor is a
smooth quadric surface.  The finite torsion conclusion below does not depend
on this choice: it follows from the local link computation for the node.

\subsection{Local link computation at each node}

Let \(L_i\) be the link of the node \(p_i\).  By the computation of the
ordinary double point in complex dimension three,
\[
        L_i\cong S^2\times S^3.
\]
Therefore
\[
        H_0(L_i,\mathbb Z)\cong\mathbb Z,
        \qquad
        H_1(L_i,\mathbb Z)=0,
\]
\[
        H_2(L_i,\mathbb Z)\cong\mathbb Z,
        \qquad
        H_3(L_i,\mathbb Z)\cong\mathbb Z,
\]
\[
        H_4(L_i,\mathbb Z)=0,
        \qquad
        H_5(L_i,\mathbb Z)\cong\mathbb Z.
\]
Using the universal coefficient theorem, we also have
\[
        H^0(L_i,\mathbb Z)\cong\mathbb Z,
        \qquad
        H^1(L_i,\mathbb Z)=0,
\]
\[
        H^2(L_i,\mathbb Z)\cong\mathbb Z,
        \qquad
        H^3(L_i,\mathbb Z)\cong\mathbb Z,
\]
\[
        H^4(L_i,\mathbb Z)=0,
        \qquad
        H^5(L_i,\mathbb Z)\cong\mathbb Z.
\]
Thus
\[
        H^m(L_i,\mathbb Z)_{\mathrm{tors}}=0
        \qquad
        \text{for all }m.
\]

Consequently, each node has no local finite torsion package of the surface
\(E\)-type:
\[
        E_i^{\mathrm{tors}}=0.
\]
Therefore
\[
        \bigoplus_{i=1}^s E_i^{\mathrm{tors}}=0.
\]

\subsection{Degree-three transport}

For the Brauer and unramified cohomology channel on a smooth projective
threefold, the relevant global torsion group is
\[
        H^3(\widetilde X,\mathbb Z)_{\mathrm{tors}}.
\]
Thus the degree-three local-to-global torsion map would have the form
\[
        \alpha_X^{(3)}:
        \bigoplus_iE_i^{\mathrm{tors}}
        \longrightarrow
        H^3(\widetilde X,\mathbb Z)_{\mathrm{tors}}.
\]
But in the nodal quintic case
\[
        \bigoplus_iE_i^{\mathrm{tors}}=0.
\]
Hence
\[
        \operatorname{im}\alpha_X^{(3)}=0.
\]
Thus ordinary nodes on a quintic threefold do not directly produce local
finite torsion classes in
\[
        H^3(\widetilde X,\mathbb Z)_{\mathrm{tors}}.
\]

This statement is only about direct local finite torsion.  It does not say
that nodal quintics have no interesting global topology.  They do.  The
nodes contribute free local vanishing and exceptional data, and those free
classes may satisfy global relations.  Such relations are the source of
defect phenomena.

\subsection{Free node data and vanishing cycles}

Let \(X_t\) be a smoothing of \(X\).  Near each node, the Milnor fiber has
the homotopy type of a \(3\)-sphere:
\[
        F_i\simeq S^3.
\]
Thus each node has a local vanishing cycle
\[
        \delta_i\in H_3(X_t,\mathbb Z)
\]
up to orientation.  These classes are free integral classes.  They are not
torsion classes.

If a small resolution exists, the corresponding resolution replaces the node
by an exceptional curve
\[
        C_i\cong\mathbb P^1.
\]
The classes of the exceptional curves are also free integral classes in
homology.  In conifold transitions, one compares relations among the
vanishing \(3\)-spheres in the smoothing with relations among exceptional
curves in the resolution.  This is a free integral relation problem, not a
finite local torsion problem.

Thus the correct nodal trajectory splits into two columns:
\[
        \text{finite torsion trajectory: zero at each node},
\]
and
\[
        \text{free relation trajectory: nontrivial and defect-sensitive}.
\]

\subsection{Defect as global relation data}

The defect of a nodal hypersurface measures a failure of expected
independence among local node conditions.  In the language of the present
paper, it is natural to view this as a global relation phenomenon among the
free local classes attached to nodes.

For a nodal hypersurface, the nodes impose linear conditions on appropriate
spaces of forms.  The defect measures the extent to which these conditions
fail to be independent.  Equivalently, it is reflected in the difference
between the actual cohomology of the nodal hypersurface or its resolution
and the cohomology predicted by treating the nodes as independent.  This
viewpoint is classical in the work of Clemens \cite{Clemens83}.  Cynk gives
an explicit treatment of the defect of nodal hypersurfaces
\cite{Cynk01}, and Kloosterman studies maximal families of nodal varieties
with defect \cite{Kloosterman22}.

For the torsion trajectory, the relevant conclusion is:
\[
        \text{defect controls global dependence among free node classes,}
\]
not a direct local finite torsion group
\[
        E_i.
\]
Therefore a nodal quintic may have nontrivial defect even though
\[
        \bigoplus_i E_i^{\mathrm{tors}}=0.
\]

\subsection{Brauer comparison}

Suppose \(\widetilde X\) is a smooth projective resolution for which
\[
        H^2(\widetilde X,\mathcal O_{\widetilde X})=0.
\]
Then the exponential sequence gives a comparison
\[
        \operatorname{Br}(\widetilde X)
        \cong
        H^3(\widetilde X,\mathbb Z)_{\mathrm{tors}}.
\]
This is the Brauer--torsion comparison proved in the main text.

However, for a nodal quintic, the direct local node torsion map has zero
source:
\[
        \bigoplus_iE_i^{\mathrm{tors}}=0.
\]
Therefore the subgroup of \(\operatorname{Br}(\widetilde X)\) obtained from
direct local finite node packages is zero:
\[
        \operatorname{Br}_{\mathrm{loc}}(\widetilde X)=0
\]
in this direct-node sense.

This does not imply that
\[
        \operatorname{Br}(\widetilde X)
\]
is zero.  It means only that any Brauer torsion present in
\(\widetilde X\) is not produced by the direct local finite torsion packages
of the ordinary double points.  It must come from global topology, from the
resolution, or from additional structure not captured by the local finite
node package.

This should be contrasted with a threefold carrying codimension-two
transverse surface singularities.  In that setting, a stratum with transverse
\(\frac14(1,1)\) singularity carries a finite local system
\[
        E_\Sigma\cong(\mathbb Z/4)_\Sigma,
\]
with order-two shadow
\[
        2E_\Sigma\cong(\mathbb Z/2)_\Sigma.
\]
Such a stratum supplies a genuine finite torsion package.  Isolated nodes on
a quintic threefold do not.

\subsection{Residue status}

Since the direct local finite torsion source is zero, there is no local node
class whose Bloch--Ogus residues must be tested.  Thus the residue station
for direct node torsion is empty:
\[
        \operatorname{im}\alpha_X^{(3)}=0
        \quad
        \Longrightarrow
        \quad
        \text{no direct local residue test}.
\]
If a global Brauer or unramified class exists on \(\widetilde X\), it must be
tested by the ordinary residue maps of the Bloch--Ogus formalism.  Such a
class is global from the perspective of this trajectory, not the image of a
local finite node discriminant group.

\subsection{Rational death}

For every node \(p_i\),
\[
        E_i^{\mathrm{tors}}=0.
\]
Therefore
\[
        E_i^{\mathrm{tors}}\otimes_{\mathbb Z}\mathbb Q=0.
\]
This rational-death column is trivial for the finite torsion trajectory.

However, the free vanishing cycles do not die rationally.  If
\[
        \delta_i\in H_3(X_t,\mathbb Z)
\]
is a vanishing \(3\)-sphere, then
\[
        \delta_i\otimes 1\in H_3(X_t,\mathbb Q)
\]
may be nonzero.  Thus rationalization kills torsion but not the free
vanishing-cycle data.  This is exactly why nodal defect is not a torsion
phenomenon in the same sense as the local surface discriminant groups.

\subsection{Trajectory row}

The nodal quintic trajectory can be summarized as follows.

\begin{center}
\begingroup
\scriptsize
\setlength{\tabcolsep}{1.8pt}
\renewcommand{\arraystretch}{1.15}
\begin{tabularx}{\textwidth}{|p{0.15\textwidth}|c|c|p{0.12\textwidth}|p{0.10\textwidth}|X|p{0.10\textwidth}|c|}
\hline
\textbf{Example} & \textbf{\(E\)} & \textbf{\(q\)} & \textbf{Local} & \textbf{Supp.} & \textbf{Global image} & \textbf{Br/res.} & \textbf{\(\mathbb Q\)}
\\
\hline
nodal quintic &
\(0\) tors. &
\(0\) &
nodes give free, not finite, local data &
none tors. &
finite local torsion image zero; free relations measured by defect &
no direct local residue test &
\(0\)
\\
\hline
\end{tabularx}
\endgroup
\end{center}

The important lesson is that nodal quintics are still essential examples,
but not because each node contributes a local \(\mathbb Z/2\)-torsion group.
They are essential because they separate local finite torsion from global
free-relation phenomena.  This distinction prevents one from mistakenly
importing the \(A_1\) surface calculation into the threefold node setting.

\section{Exact sequences used in the torsion trajectory}
\label{app:exactsequencetorsiontrajectory}

This appendix collects the exact sequences used throughout the paper.  The
purpose is to make clear which station of the trajectory each sequence
controls.  The main sequences and comparisons are:
\[
\text{support/excision},
\qquad
\text{pair sequence},
\qquad
\text{universal coefficients},
\qquad
\text{coefficient Bockstein},
\]
\[
\text{MacPherson--Vilonen/Beilinson gluing},
\qquad
\text{subobject--quotient triangles},
\qquad
\text{exponential sequence},
\qquad
\text{Kummer sequence},
\]
\[
\text{Bloch--Ogus residues},
\qquad
\text{Wang sequence},
\qquad
\text{motivic cone sequences}.
\]

The corrected Benoist--Ottem/Coble comparison uses three of these especially
strongly:
\[
        0\to\mathbb Z\xrightarrow{2}\mathbb Z\to\mathbb Z/2\to0,
\]
which produces the Bockstein shadow
\[
        2E\subset E\cong\mathbb Z/4,
\]
the short exact sequence
\[
        0\to2E\to E\to E/2E\to0,
\]
and the closed-stratum gluing triangle
\[
        {}^pIC\to{}^p_+IC\to i_*E_\Sigma[2]\to.
\]

\subsection{Cohomology with supports and excision}

Let \(Y\) be a locally compact topological space, and let
\[
        A\subset Y
\]
be a closed subset.  Cohomology with supports in \(A\) is denoted
\[
        H^k_A(Y,\mathbb Z).
\]
For the spaces used in this paper, one may identify this group with relative
cohomology:
\[
        H^k_A(Y,\mathbb Z)
        \cong
        H^k(Y,Y\setminus A;\mathbb Z).
\]
This identification is standard in sheaf cohomology with supports and
relative cohomology; see \cite{BredonSheaf97,Hatcher02}.

The long exact sequence of the pair
\[
        (Y,Y\setminus A)
\]
gives
\[
\cdots
\longrightarrow
H^k_A(Y,\mathbb Z)
\longrightarrow
H^k(Y,\mathbb Z)
\longrightarrow
H^k(Y\setminus A,\mathbb Z)
\longrightarrow
H^{k+1}_A(Y,\mathbb Z)
\longrightarrow
\cdots .
\]
The map
\[
        H^k_A(Y,\mathbb Z)
        \longrightarrow
        H^k(Y,\mathbb Z)
\]
is the forget-support map.

If
\[
        A=\bigcup_i A_i
\]
is a finite disjoint union of closed subsets and \(N_i\) are disjoint
neighborhoods of \(A_i\), then excision gives
\[
        H^k_A(Y,\mathbb Z)
        \cong
        \bigoplus_i H^k_{A_i}(N_i,\mathbb Z).
\]
This is the formal mechanism behind the maps
\[
        \alpha_X^{(k)}:
        \bigoplus_i E_i
        \longrightarrow
        H^k(\widetilde X,\mathbb Z)_{\mathrm{tors}}
\]
whenever the local groups \(E_i\) have been realized in degree \(k\)
support cohomology.  If a specified subpackage
\[
        A_i\subset E_i
\]
is being transported, the same exact sequence gives the restricted map
\[
        \alpha_{X,A}^{(k)}:
        \bigoplus_i A_i
        \longrightarrow
        H^k(\widetilde X,\mathbb Z)_{\mathrm{tors}}.
\]

\subsection{The pair sequence of a resolution neighborhood}

Let \(N\) be a compact oriented manifold with boundary
\[
        L=\partial N.
\]
The long exact homology sequence of the pair \((N,L)\) is
\[
\cdots
\longrightarrow
H_j(L,\mathbb Z)
\longrightarrow
H_j(N,\mathbb Z)
\longrightarrow
H_j(N,L;\mathbb Z)
\longrightarrow
H_{j-1}(L,\mathbb Z)
\longrightarrow
H_{j-1}(N,\mathbb Z)
\longrightarrow
\cdots .
\]
If \(N\) is an oriented compact real \(m\)-manifold with boundary, then
Poincare--Lefschetz duality gives
\[
        H_j(N,L;\mathbb Z)
        \cong
        H^{m-j}(N,\mathbb Z).
\]
This is the mechanism that turns intersection matrices of exceptional
cycles into torsion in the boundary link.

For a normal surface singularity, \(N\) is a real four-manifold and \(L\) is
a closed oriented three-manifold.  If \(\Lambda\) is the exceptional lattice,
then the map
\[
        H_2(N,\mathbb Z)
        \longrightarrow
        H_2(N,L;\mathbb Z)
\]
is represented by the exceptional intersection matrix.  Therefore the
boundary torsion is computed by
\[
        \operatorname{coker}(\Lambda\to\Lambda^\vee)
        =
        \Lambda^\vee/\Lambda.
\]
This is the resolution-neighborhood calculation used in
\cite{GoreskySiegel83,RahmanIntegralPerverseObstructions}.

For example, for the Coble boundary singularity
\[
        \frac14(1,1),
\]
the exceptional lattice is
\[
        \Lambda=\mathbb Z\langle e\rangle,
        \qquad
        (e,e)=-4.
\]
Thus the pair sequence gives
\[
        \mathbb Z
        \xrightarrow{[-4]}
        \mathbb Z
        \longrightarrow
        H^2(L(4,1),\mathbb Z)_{\mathrm{tors}}
        \longrightarrow
        0,
\]
and hence
\[
        H^2(L(4,1),\mathbb Z)_{\mathrm{tors}}
        \cong
        \mathbb Z/4.
\]

\subsection{The universal coefficient theorem}

For any topological space \(Y\), the universal coefficient theorem gives a
short exact sequence
\[
        0
        \longrightarrow
        \operatorname{Ext}^1_{\mathbb Z}
        (H_{k-1}(Y,\mathbb Z),\mathbb Z)
        \longrightarrow
        H^k(Y,\mathbb Z)
        \longrightarrow
        \operatorname{Hom}_{\mathbb Z}
        (H_k(Y,\mathbb Z),\mathbb Z)
        \longrightarrow
        0.
\]
This is used repeatedly to convert torsion in \(H_1\) of a three-manifold
link into torsion in \(H^2\).  For example, if
\[
        H_1(L,\mathbb Z)\cong\mathbb Z/n\mathbb Z
\]
and
\[
        H_2(L,\mathbb Z)=0,
\]
then
\[
        H^2(L,\mathbb Z)
        \cong
        \operatorname{Ext}^1_{\mathbb Z}
        (\mathbb Z/n\mathbb Z,\mathbb Z)
        \cong
        \mathbb Z/n\mathbb Z.
\]
This is the calculation used for \(A_1\), \(A_k\), the non-ADE Brieskorn
example, and the Coble boundary example
\[
        L(4,1):
        \qquad
        H^2(L(4,1),\mathbb Z)\cong\mathbb Z/4.
\]

\subsection{The coefficient Bockstein sequence}

The short exact sequence
\[
        0
        \longrightarrow
        \mathbb Z
        \xrightarrow{n}
        \mathbb Z
        \longrightarrow
        \mathbb Z/n
        \longrightarrow
        0
\]
gives a long exact cohomology sequence.  The connecting homomorphism
\[
        \beta_n:
        H^r(Y,\mathbb Z/n)
        \longrightarrow
        H^{r+1}(Y,\mathbb Z)
\]
is the Bockstein.

The case \(n=2\) is central for the Benoist--Ottem/Coble comparison:
\[
        0
        \longrightarrow
        \mathbb Z
        \xrightarrow{2}
        \mathbb Z
        \longrightarrow
        \mathbb Z/2
        \longrightarrow
        0.
\]
For the Coble boundary link
\[
        L=L(4,1),
\]
the relevant segment is
\[
        H^1(L,\mathbb Z)
        \longrightarrow
        H^1(L,\mathbb Z/2)
        \xrightarrow{\beta}
        H^2(L,\mathbb Z)
        \xrightarrow{2}
        H^2(L,\mathbb Z).
\]
Since
\[
        H^1(L,\mathbb Z)=0
\]
and
\[
        H^2(L,\mathbb Z)\cong\mathbb Z/4,
\]
the Bockstein is injective and has image
\[
        \operatorname{im}(\beta)
        =
        \ker(2:\mathbb Z/4\to\mathbb Z/4)
        =
        2(\mathbb Z/4)
        \cong
        \mathbb Z/2.
\]
If
\[
        \eta\in H^1(L(4,1),\mathbb Z/2)
\]
is the double-cover class of
\[
        L(2,1)\longrightarrow L(4,1),
\]
then
\[
        \beta(\eta)=2\bar g
        \in
        E\cong\mathbb Z/4.
\]
This is the order-two Bockstein shadow
\[
        2E\cong\mathbb Z/2
        \subset
        E\cong\mathbb Z/4.
\]

\subsection{Subobjects, quotients, and shadow triangles}

If a local torsion package contains a specified subgroup
\[
        A\subset E,
\]
then there is a short exact sequence
\[
        0
        \longrightarrow
        A
        \longrightarrow
        E
        \longrightarrow
        E/A
        \longrightarrow
        0.
\]
This short exact sequence gives a distinguished triangle after viewing these
groups as sheaves or complexes on the relevant support.

The Coble boundary example uses
\[
        E\cong\mathbb Z/4,
        \qquad
        A=2E\cong\mathbb Z/2.
\]
Thus
\[
        0
        \longrightarrow
        2E
        \longrightarrow
        E
        \longrightarrow
        E/2E
        \longrightarrow
        0
\]
becomes
\[
        0
        \longrightarrow
        \mathbb Z/2
        \longrightarrow
        \mathbb Z/4
        \longrightarrow
        \mathbb Z/2
        \longrightarrow
        0.
\]
If \(E_\Sigma\) is the corresponding local system on a codimension-two
stratum
\[
        i:\Sigma\hookrightarrow Y,
\]
then applying \(i_*\) and shifting by \([2]\) gives the distinguished
triangle
\[
        i_*(\mathbb Z/2)_\Sigma[2]
        \longrightarrow
        i_*(\mathbb Z/4)_\Sigma[2]
        \longrightarrow
        i_*(\mathbb Z/2)_\Sigma[2]
        \overset{+1}{\longrightarrow}.
\]
This is the exact-sequence form of the statement that the Benoist--Ottem
visible part is the order-two shadow
\[
        i_*(2E_\Sigma)[2]
        \subset
        i_*E_\Sigma[2].
\]

\subsection{MacPherson--Vilonen and Beilinson gluing}

Let
\[
        i:\Sigma\hookrightarrow Y
\]
be a closed stratum and
\[
        j:V=Y\setminus \Sigma\hookrightarrow Y
\]
the complementary open stratum.  MacPherson--Vilonen gluing and Beilinson's
gluing formalism describe perverse sheaves on \(Y\) in terms of their
restriction to \(V\), their restriction to \(\Sigma\), and the normal Morse
data connecting the two strata.

In the setting of this paper, the key application is a codimension-two
stratum whose transverse singularity is a normal surface singularity with
local package \(E\).  If the transverse type is locally trivial along
\(\Sigma\), the local groups assemble into a finite local system
\[
        E_\Sigma
\]
on \(\Sigma\).  The discrepancy between the ordinary and dual middle
extensions is then supported on \(\Sigma\), and the closed-stratum term is
\[
        i_*E_\Sigma[2].
\]
Thus one obtains the triangle
\[
        {}^pIC_Y\mathbb Z
        \longrightarrow
        {}^p_+IC_Y\mathbb Z
        \longrightarrow
        i_*E_\Sigma[2]
        \overset{+1}{\longrightarrow}.
\]

For a transverse Coble boundary singularity
\[
        \frac14(1,1),
\]
one has
\[
        E_\Sigma\cong(\mathbb Z/4)_\Sigma.
\]
Thus the full closed-stratum discrepancy is
\[
        i_*(\mathbb Z/4)_\Sigma[2].
\]
The Benoist--Ottem order-two shadow is the subobject
\[
        i_*(\mathbb Z/2)_\Sigma[2]
        =
        i_*(2E_\Sigma)[2].
\]

The octahedral axiom may be applied to the discrepancy triangle
\[
        A\to B\to i_*E_\Sigma[2]\to A[1]
\]
and the shadow triangle
\[
        i_*(2E_\Sigma)[2]
        \to
        i_*E_\Sigma[2]
        \to
        i_*(E_\Sigma/2E_\Sigma)[2]
        \to
\]
to express the order-two shadow as a filtered piece or subquotient of the
full \(E\)-discrepancy cone.

\subsection{The exponential sequence and the Brauer group}

Let \(Y\) be a smooth complex variety, and let \(Y^{\mathrm{an}}\) be the
associated analytic space.  The exponential map gives an exact sequence of
analytic sheaves
\[
        0
        \longrightarrow
        \mathbb Z
        \longrightarrow
        \mathcal O_Y
        \xrightarrow{\exp(2\pi i\,\cdot)}
        \mathcal O_Y^*
        \longrightarrow
        0.
\]
Taking cohomology gives the exact segment
\[
        H^2(Y,\mathcal O_Y)
        \longrightarrow
        H^2(Y,\mathcal O_Y^*)
        \longrightarrow
        H^3(Y,\mathbb Z)
        \longrightarrow
        H^3(Y,\mathcal O_Y).
\]
If
\[
        H^2(Y,\mathcal O_Y)=0,
\]
then the connecting map identifies
\[
        H^2(Y,\mathcal O_Y^*)_{\mathrm{tors}}
        \cong
        H^3(Y,\mathbb Z)_{\mathrm{tors}}.
\]
For a smooth projective complex variety, the torsion subgroup of analytic
\(H^2(Y,\mathcal O_Y^*)\) agrees with the cohomological Brauer group defined
in étale cohomology:
\[
        \operatorname{Br}(Y)
        =
        H^2_{\mathrm{\acute et}}(Y,\mathbb G_m)_{\mathrm{tors}}.
\]
This comparison is standard; see Grothendieck's Brauer papers and Milne's
étale cohomology text \cite{GrothendieckBrauerIII,MilneEtale}.

Thus, under the vanishing condition
\[
        H^2(Y,\mathcal O_Y)=0,
\]
one obtains
\[
        \operatorname{Br}(Y)
        \cong
        H^3(Y,\mathbb Z)_{\mathrm{tors}}.
\]

\subsection{The Kummer sequence}

For an integer \(n>0\), the Kummer sequence in the étale topology is
\[
        1
        \longrightarrow
        \mu_n
        \longrightarrow
        \mathbb G_m
        \xrightarrow{(\cdot)^n}
        \mathbb G_m
        \longrightarrow
        1.
\]
It gives a long exact sequence in étale cohomology.  In low degrees, it
contains
\[
        \operatorname{Pic}(Y)
        \xrightarrow{n}
        \operatorname{Pic}(Y)
        \longrightarrow
        H^2_{\mathrm{\acute et}}(Y,\mu_n)
        \longrightarrow
        H^2_{\mathrm{\acute et}}(Y,\mathbb G_m)
        \xrightarrow{n}
        H^2_{\mathrm{\acute et}}(Y,\mathbb G_m).
\]
Thus \(n\)-torsion Brauer classes are detected by
\[
        H^2_{\mathrm{\acute et}}(Y,\mu_n)
\]
modulo the image of the Picard group.  This is another standard way to study
torsion in the Brauer group; see \cite{GrothendieckBrauerIII,MilneEtale}.

In the present paper, the Kummer sequence is mainly used as a complementary
viewpoint on the Brauer station.  The exponential sequence gives the
comparison with \(H^3(Y,\mathbb Z)_{\mathrm{tors}}\) over \(\mathbb C\),
while the Kummer sequence records the finite-level \(n\)-torsion structure.

\subsection{Bloch--Ogus residues and unramified cohomology}

Let \(Y\) be a smooth irreducible complex variety with function field
\[
        k(Y).
\]
The third unramified cohomology group with \(\mathbb Q/\mathbb Z\)
coefficients is defined by
\[
        H^3_{\mathrm{nr}}(Y,\mathbb Q/\mathbb Z)
        =
        \ker\left(
        H^3(k(Y),\mathbb Q/\mathbb Z)
        \longrightarrow
        \bigoplus_{D\in Y^{(1)}} H^2(k(D),\mathbb Q/\mathbb Z)
        \right),
\]
where the sum is over codimension-one points \(D\) of \(Y\), and the maps
are residue homomorphisms.

This description belongs to the Bloch--Ogus theory of coniveau and Gersten
resolutions \cite{BlochOgus74}.  The relationship between unramified
cohomology and failures of the integral Hodge conjecture in degree \(4\) is
central in the work of Colliot-Thélène and Voisin \cite{CTVoisin12}.

For the torsion trajectory, the residue station is the following test.  A
class that has survived to the Brauer or \(\mathbb Q/\mathbb Z\)-cohomology
stage contributes to
\[
        H^3_{\mathrm{nr}}(Y,\mathbb Q/\mathbb Z)
\]
only if all its codimension-one residues vanish.  Thus residues provide a
second killing mechanism after global support relations.

\subsection{The Wang sequence for the Milnor fibration}

Let
\[
        f:(\mathbb C^{n+1},0)\longrightarrow(\mathbb C,0)
\]
define an isolated hypersurface singularity, and let \(F\) be the Milnor
fiber.  The link complement fibers over \(S^1\), and the monodromy
\[
        T:H^n(F,\mathbb Z)\longrightarrow H^n(F,\mathbb Z)
\]
acts on the middle vanishing cohomology.  The associated Wang sequence
relates the cohomology of the link to the kernel and cokernel of
\[
        T-\mathrm{id}.
\]
For isolated hypersurface surface singularities, this gives the realization
\[
        E
        \cong
        \operatorname{coker}(T-\mathrm{id})_{\mathrm{tors}}
\]
proved in \cite{RahmanIntegralPerverseObstructions}, using the standard
Milnor fibration background from \cite{Mi68}.

In the examples above, the Wang sequence recovers:
\[
        A_1:\quad \operatorname{coker}(T-\mathrm{id})_{\mathrm{tors}}
        \cong \mathbb Z/2\mathbb Z,
\]
\[
        A_k:\quad \operatorname{coker}(T-\mathrm{id})_{\mathrm{tors}}
        \cong \mathbb Z/(k+1)\mathbb Z,
\]
\[
        D_4:\quad \operatorname{coker}(T-\mathrm{id})_{\mathrm{tors}}
        \cong
        (\mathbb Z/2\mathbb Z)^2,
\]
\[
        E_8:\quad \operatorname{coker}(T-\mathrm{id})_{\mathrm{tors}}=0.
\]
For the threefold ordinary double point, the same station gives a free
cokernel:
\[
        \operatorname{coker}(T-\mathrm{id})\cong\mathbb Z,
\]
whose torsion subgroup is zero.  This explains why the threefold node does
not contribute a local finite torsion package of the surface \(E\)-type.

The cyclic quotient singularity \(\frac14(1,1)\) is not an isolated
hypersurface surface singularity in the sense used for this Wang-sequence
realization.  Its local package is therefore computed by the link, lattice,
pair-sequence, and perverse stations, not by the hypersurface monodromy
station.

\subsection{Motivic cone convention}

In the motivic formulation of the Coble boundary package, one uses torsion
objects of the form
\[
        \mathbf 1_\Sigma/n
        :=
        \operatorname{Cone}
        \left(
        \mathbf 1_\Sigma
        \xrightarrow{n}
        \mathbf 1_\Sigma
        \right),
\]
where \(\mathbf 1_\Sigma\) is the motivic unit on the stratum \(\Sigma\).
Under a Betti realization compatible with cones, integer multiplication,
proper pushforward, and shifts, one has
\[
        \operatorname{Real}_B(\mathbf 1_\Sigma/n)
        \cong
        (\mathbb Z/n)_\Sigma.
\]
Thus the full Coble package is modeled by
\[
        i_*(\mathbf 1_\Sigma/4)[2],
\]
whose Betti realization is
\[
        i_*(\mathbb Z/4)_\Sigma[2].
\]
The Benoist--Ottem shadow is modeled by
\[
        i_*(\mathbf 1_\Sigma/2)[2],
\]
whose Betti realization is
\[
        i_*(\mathbb Z/2)_\Sigma[2].
\]
The morphism
\[
        \mathbf 1_\Sigma/2\longrightarrow \mathbf 1_\Sigma/4
\]
realizes the inclusion
\[
        \mathbb Z/2\hookrightarrow \mathbb Z/4,
        \qquad
        1\mapsto2.
\]

\subsection{Summary of stations}

The exact sequences above control the stations of the torsion trajectory as
follows:

\[
\begin{array}{c|c}
\text{Station} & \text{Exact sequence or comparison} \\
\hline
\text{Birth} & \text{local perverse discrepancy and link cohomology} \\
\text{Form} & \text{discriminant form and linking pairing} \\
\text{Pair sequence} & H_*(N,L)\text{ long exact sequence} \\
\text{Excision} & H^k_D(\widetilde X)\to H^k(\widetilde X) \\
\text{Bockstein shadow} & 0\to\mathbb Z\xrightarrow{2}\mathbb Z\to\mathbb Z/2\to0 \\
\text{Subobject/quotient} & 0\to2E\to E\to E/2E\to0 \\
\text{Closed-stratum gluing} & {}^pIC\to{}^p_+IC\to i_*E_\Sigma[2]\to \\
\text{Brauer} & 0\to\mathbb Z\to\mathcal O\to\mathcal O^*\to0 \\
\text{Finite Brauer level} & \text{Kummer sequence} \\
\text{Residue} & \text{Bloch--Ogus residue maps} \\
\text{Monodromy} & \text{Wang sequence for the Milnor fibration} \\
\text{Motivic lift} & \mathbf 1_\Sigma/n=\operatorname{Cone}(\mathbf 1_\Sigma\xrightarrow{n}\mathbf 1_\Sigma) \\
\text{Rational death} & (-)\otimes_{\mathbb Z}\mathbb Q
\end{array}
\]

\section{Realization dictionary for motivic torsion packages}
\label{app:realization-dictionary}

The purpose of this appendix is programmatic.  We do not reprove the
six-functor formalisms of Tubach or Ruimy--Tubach.  Instead, we explain why
the torsion-trajectory program requires a common motivic finite-coefficient
setting in which known integral Hodge counterexamples can be represented and
compared.

The reason is that known counterexamples are produced by different-looking
mechanisms: local discriminant packages, Enriques double covers, Brauer
classes, unramified cohomology, finite group actions, stacky stabilizers, and
cycle-divisibility constraints.  These mechanisms are difficult to compare
directly at the level of their final cohomology classes.  A motivic
finite-coefficient lift records the source package and the operations that
produce the final obstruction.  The six operations then provide a common
transport calculus for moving, multiplying, localizing, and realizing these
packages.

The guiding philosophy is therefore the following.  Integral Hodge
counterexamples should be compared not only by their final obstruction classes, but by their pre-rationalization torsion trajectories.  Since rationalization kills torsion, the torsion mechanisms responsible for
integral failures must be organized before tensoring with \(\mathbb Q\).  A
finite-coefficient motivic framework provides a natural place to record those
mechanisms.

The torsion-trajectory program compares the same finite-coefficient motivic
package across several realizations.  This appendix fixes notation for the
realization symbols used in the paper and explains why these realizations are
useful for comparing integral Hodge counterexamples.

A realization is a functor, or in some cases a comparison station, that sends
a motivic or finite-coefficient package to one of its concrete avatars:
Betti cohomology, étale cohomology, Hodge-theoretic data, Nori motivic
sheaves, Brauer classes, unramified cohomology, or stacky/equivariant
cohomology.  The benefit of using realizations is that one can construct or
name a torsion package once, then compare its images across the different
places where integral Hodge obstructions are detected.

Thus the purpose of this dictionary is not merely notational.  It records
the interfaces through which a torsion package can be tested:
\[
        \text{finite coefficients},
        \qquad
        \text{Brauer comparison},
        \qquad
        \text{unramified survival},
\]
\[
        \text{Hodge-theoretic obstruction},
        \qquad
        \text{stacky stabilizer data}.
\]
Some entries below are genuine realization functors from motivic categories.
Others, such as the Brauer and unramified entries, are comparison or survival
stations reached after applying standard exact sequences or residue tests.

\subsection{Betti realization}

We write
\[
        \RealB
\]
for Betti realization.  It sends motivic finite-coefficient objects to
topological coefficient systems.  For example, in the finite-coefficient
setting one expects
\[
        \RealB(\mathbf 1_X(p)/m)
        \simeq
        \mathbb Z/m(p).
\]
The Betti realization is the one used to compare motivic
finite-coefficient packages with ordinary integral and finite-coefficient
cohomology.

For a motivic class
\[
        \theta^{\mathrm{mot}}
        \in
        H^r_{\mathrm{mot}}(X,\mathbf 1_X(p)/m),
\]
its Betti realization is denoted
\[
        \RealB(\theta^{\mathrm{mot}})
        \in
        H^r(X,\mathbb Z/m(p)).
\]

In the Diaz construction, for instance, a motivic lift
\[
        \gamma^{\mathrm{mot}}
        \in
        H^3_{\mathrm{mot}}(V,\mathbf 1_V(2)/2)
\]
has Betti realization
\[
        \RealB(\gamma^{\mathrm{mot}})
        =
        \gamma
        \in
        H^3(V,\mathbb Z/2(2)),
\]
when the chosen motivic lift realizes to Diaz's finite-coefficient class.

\subsection{Étale realization}

We write
\[
        \Realet
\]
for étale realization.  It sends motivic finite-coefficient objects to
étale coefficient systems.  Depending on conventions, one has
\[
        \Realet(\mathbf 1_X(p)/m)
        \simeq
        \mathbb Z/m(p),
\]
or
\[
        \Realet(\mathbf 1_X(p)/m)
        \simeq
        \mu_m^{\otimes p}.
\]
This is the realization relevant to Kummer sequences, Brauer classes,
finite-coefficient étale cohomology, and unramified cohomology.

For a motivic finite-coefficient class
\[
        \theta^{\mathrm{mot}}
        \in
        H^r_{\mathrm{mot}}(X,\mathbf 1_X(p)/m),
\]
we denote its étale realization by
\[
        \Realet(\theta^{\mathrm{mot}}).
\]
When \(m=2\), this realization is the one most closely related to the
\(\mathbb Z/2\)-classes and \(\mu_2\)-twisted classes appearing in Diaz's
construction and in the Brauer comparison.

\subsection{Hodge and mixed-Hodge-module realizations}

We write
\[
        \RealHdg
\]
for Hodge realization and
\[
        \RealMHM
\]
for mixed-Hodge-module realization.  These realizations record the
Hodge-theoretic structure carried by a motivic package.

For torsion trajectories, these realizations are useful because they place
finite-coefficient and integral torsion data near the integral Hodge
obstruction group.  In particular, they help distinguish the rational
Hodge-theoretic shadow from the integral torsion information that disappears
after tensoring with \(\mathbb Q\).

The mixed-Hodge-module realization is especially relevant when torsion
packages are supported on strata, because it allows one to compare motivic
finite-coefficient packages with perverse and Hodge-module structures.  This
is the setting in which local packages, closed-stratum contributions, and
intermediate-extension phenomena can be compared.

\subsection{Nori and perverse-Nori realizations}

We write
\[
        \RealNori
\]
for Nori or perverse-Nori realization.  This is the realization relevant to
the Tubach and Ruimy--Tubach frameworks.  The role of this realization is to
provide an abelian or perverse motivic environment in which
finite-coefficient objects, integral structures, and six operations can be
compared with mixed-Hodge-module and Betti/étale realizations.

The terminology of motivic torsion trajectories used in this paper is ours.
The point is not that the existing literature has already developed this
specific theory, but that the Nori and mixed-Hodge-module realization
frameworks provide a natural ambient setting in which such a theory can be
formulated.

\subsection{Brauer comparison}

We write
\[
        \RealBr
\]
for the Brauer comparison station.  This is not a realization functor in the
same strict sense as \(\RealB\) or \(\Realet\).  Rather, it denotes the
passage from finite or integral cohomological torsion to the cohomological
Brauer group through the exponential or Kummer sequence.

For example, under suitable hypotheses, one has comparisons of the form
\[
        H^3(X,\mathbb Z)_{\mathrm{tors}}
        \longleftrightarrow
        \operatorname{Br}(X).
\]
In Diaz's construction, the class
\[
        \beta\in H^2(S_2,\mathbb Z/2(1))
\]
is chosen so that its Brauer comparison is nonzero:
\[
        \RealBr(\beta)\neq0
        \quad\text{in}\quad
        \operatorname{Br}(S_2)[2].
\]

Thus \(\RealBr\) records not a new cohomology theory, but a comparison
station in the torsion trajectory.  It asks whether a finite-coefficient or
integral torsion class contributes nontrivially to the Brauer group.

\subsection{Unramified realization or survival}

We write
\[
        \Realnr
\]
for the unramified survival station.  This is not a realization functor in
the strict motivic sense.  It denotes the image or survival of a
finite-coefficient class after applying the coniveau or Bloch--Ogus residue
test.

For a class
\[
        \Theta\in H^r(X,\mathbb Z/m(p)),
\]
the condition
\[
        \Realnr(\Theta)\neq0
\]
means that \(\Theta\) has nonzero image in
\[
        H^r(\mathbb C(X),\mathbb Z/m(p)),
\]
or equivalently
\[
        \Theta\notin N^1H^r(X,\mathbb Z/m(p)).
\]
This is the survival condition used in Diaz's construction.

The unramified station is crucial because forming a finite-coefficient class does not by itself produce an integral Hodge obstruction.  The class must survive the relevant coniveau or unramified test before its Bockstein can be used as an obstruction to integral algebraicity.

\subsection{Stacky and equivariant realization}

We write
\[
        \Realst
\]
for the stacky or equivariant realization of a torsion package.  This is the
realization relevant to quotient stacks and stabilizer packages.  For
example, in the local Kummer model
\[
        [\mathbb C^2/\mu_2],
\]
the stabilizer package is
\[
        H^2(B\mu_2,\mathbb Z)\cong\mathbb Z/2.
\]
Thus
\[
        \Realst(E_{\mu_2})
        =
        H^2(B\mu_2,\mathbb Z).
\]
The local Kummer compatibility theorem identifies this with the ordinary
singular package:
\[
        E_{A_1}^{\mathrm{sing}}
        \cong
        H^2(B\mu_2,\mathbb Z).
\]

This station is expected to be important for smooth quotient examples, such
as Atiyah--Hirzebruch or Godeaux--Serre-type constructions, where ordinary
singular \(E\)-packages may vanish on the final smooth quotient but finite
group or stacky data may still record the torsion source.

\subsection{Summary table}

\[
\begin{array}{c|c|c}
\text{Notation} & \text{Meaning} & \text{Typical output} \\
\hline
\RealB
&
\text{Betti realization}
&
H^r(X,\mathbb Z/m(p))
\\
\Realet
&
\text{étale realization}
&
H^r_{\acute et}(X,\mu_m^{\otimes p})
\\
\RealHdg
&
\text{Hodge realization}
&
\text{Hodge or integral Hodge data}
\\
\RealMHM
&
\text{mixed-Hodge-module realization}
&
D^b\operatorname{MHM}(X)
\\
\RealNori
&
\text{Nori/perverse-Nori realization}
&
\text{Nori motivic sheaves}
\\
\RealBr
&
\text{Brauer comparison station}
&
\operatorname{Br}(X)[m]
\\
\Realnr
&
\text{unramified survival station}
&
H^r(\mathbb C(X),\mathbb Z/m(p))
\\
\Realst
&
\text{stacky/equivariant realization}
&
H^*(BG,\mathbb Z),\ H^*_G(V,\mathbb Z)
\end{array}
\]

The strict functorial realizations are \(\RealB\), \(\Realet\),
\(\RealHdg\), \(\RealMHM\), and \(\RealNori\).  The symbols \(\RealBr\),
\(\Realnr\), and \(\Realst\) record comparison or survival stations in the
torsion trajectory.  They are included because these stations are essential
for comparing integral Hodge counterexamples.

\section{Six operations and motivic torsion trajectories}
\label{app:six-operations}

The purpose of this appendix is programmatic.  We do not reprove the
six-functor formalisms of Tubach or Ruimy--Tubach \cite{Tubach2025NoriHodgeRealizations,RuimyTubach24}.  Instead, we record the six-operation transport dictionary needed for the torsion-trajectory program and explain how realization compatibility allows torsion packages to be compared across Betti, étale, Nori, mixed-Hodge-module, Brauer, unramified, and stacky settings.

Tubach constructs realization functors from Voevodsky étale motives to perverse Nori motives and mixed Hodge modules, with compatibility with the six operations as part of the realization framework
\cite{Tubach2025NoriHodgeRealizations}.  Ruimy--Tubach develop the
corresponding integral-coefficient environment, constructing integral Nori motivic sheaves and integral mixed Hodge modules with six-functor formalism and integral structures suitable for tracking torsion \cite{RuimyTubach24}.  These are the ambient tools used below. The new point in the present paper is not the construction of the six-functor formalism, but its use as a transport dictionary for integral Hodge torsion mechanisms.

The guiding point is that known integral Hodge counterexamples are produced by different-looking mechanisms: local discriminant packages, Enriques double covers, Brauer classes, unramified cohomology, finite group actions, stacky stabilizers, and cycle-divisibility constraints.  These mechanisms are hard to compare directly at the level of their final cohomology classes.  A motivic finite-coefficient lift records the source package and the operations that produce the final obstruction.  The six operations then provide a common transport calculus for moving, multiplying, localizing, dualizing, and realizing these packages.

\subsection{The six operations}

Let \(k=\mathbb C\), and let \(X\) and \(Y\) be separated finite-type
\(k\)-schemes, or complex algebraic varieties, in a motivic setting where the
six operations are available.  We write
\[
        \mathsf{DM}(X)
\]
for the corresponding triangulated or stable motivic category over \(X\).
The unit object is denoted
\[
        \mathbf 1_X\in \mathsf{DM}(X),
\]
and its Tate twist is denoted \(\mathbf 1_X(p)\).  For an integer
\(m\ge 2\), the finite-coefficient motivic object is
\[
        \mathbf 1_X(p)/m
        :=
        \operatorname{Cone}
        \left(
        \mathbf 1_X(p)
        \xrightarrow{m}
        \mathbf 1_X(p)
        \right).
\]
Thus \(\mathbf 1_X(p)/m\) fits into the coefficient triangle
\[
        \mathbf 1_X(p)
        \xrightarrow{m}
        \mathbf 1_X(p)
        \longrightarrow
        \mathbf 1_X(p)/m
        \overset{+1}{\longrightarrow}.
\]
This is the motivic object whose Betti or étale realizations give the usual
finite-coefficient systems \(\mathbb Z/m(p)\), or \(\mu_m^{\otimes p}\) in
the étale convention.

Let
\[
        f:X\to Y
\]
be a morphism.  In a six-functor formalism one has functors
\[
        f^*:\mathsf{DM}(Y)\to \mathsf{DM}(X),
        \qquad
        f_*:\mathsf{DM}(X)\to \mathsf{DM}(Y),
\]
\[
        f_!:\mathsf{DM}(X)\to \mathsf{DM}(Y),
        \qquad
        f^!:\mathsf{DM}(Y)\to \mathsf{DM}(X),
\]
together with a tensor product and internal Hom
\[
        -\otimes -:\mathsf{DM}(X)\times\mathsf{DM}(X)\to\mathsf{DM}(X),
        \qquad
        \mathcal Hom_X(-,-):\mathsf{DM}(X)^{\mathrm{op}}\times
        \mathsf{DM}(X)\to\mathsf{DM}(X).
\]
The operations \(f^*\) and \(f^!\) move objects from \(Y\) to \(X\), while
\(f_*\) and \(f_!\) move objects from \(X\) to \(Y\).  The tensor product
and internal Hom are operations internal to the fiber category
\(\mathsf{DM}(X)\).

These functors are the standard operations appearing in motivic and
sheaf-theoretic six-functor formalisms.  In the present paper we use them as
the basic transport operations for finite-coefficient torsion packages.  The
ambient justification is the realization-compatible six-functor machinery of
Tubach and Ruimy--Tubach: Tubach constructs realization functors from
Voevodsky étale motives to perverse Nori motives and mixed Hodge modules
compatible with the six operations, while Ruimy--Tubach develop the
integral-coefficient Nori and mixed-Hodge-module setting in which torsion
phenomena can be tracked
\cite{Tubach2025NoriHodgeRealizations,RuimyTubach24}.
We do not reprove those formalisms here.  We only record how the operations
are used in the torsion-trajectory program.

For a closed immersion
\[
        i:\Sigma\hookrightarrow X,
\]
a supported torsion package has the typical motivic form
\[
        i_*(\mathbf 1_\Sigma(p)/m)[s]\in \mathsf{DM}(X).
\]
For an open immersion
\[
        j:U\hookrightarrow X,
\]
the pair \(j_!\) and \(j_*\) distinguishes extension with support from full
extension across the boundary.  For projections
\[
        \pi_a:X_1\times X_2\to X_a,
\]
the pullbacks \(\pi_a^*\) and the tensor product produce external cup
products.  Thus, for motivic lifts
\[
        \alpha^{\mathrm{mot}}\in
        H^r_{\mathrm{mot}}(X_1,\mathbf 1_{X_1}(p_1)/m),
        \qquad
        \beta^{\mathrm{mot}}\in
        H^s_{\mathrm{mot}}(X_2,\mathbf 1_{X_2}(p_2)/m),
\]
their product on \(X_1\times X_2\) is represented by
\[
        \pi_1^*\alpha^{\mathrm{mot}}
        \cup
        \pi_2^*\beta^{\mathrm{mot}}
        \in
        H^{r+s}_{\mathrm{mot}}
        \left(
        X_1\times X_2,\mathbf 1_{X_1\times X_2}(p_1+p_2)/m
        \right).
\]
The connecting morphism of the coefficient triangle then gives the motivic
Bockstein
\[
        \delta_m^{\mathrm{mot}}:
        H^{r+s}_{\mathrm{mot}}
        \left(
        X_1\times X_2,
        \mathbf 1_{X_1\times X_2}(p_1+p_2)/m
        \right)
        \longrightarrow
\]
\[
        H^{r+s+1}_{\mathrm{mot}}
        \left(
        X_1\times X_2,
        \mathbf 1_{X_1\times X_2}(p_1+p_2)
        \right).
\]
Under Betti or étale realization this becomes the ordinary coefficient
Bockstein.  This is the formal pattern behind the Diaz class
\[
        \gamma=\pi_S^*\alpha\cup\pi_{S_2}^*\beta
\]
and its Bockstein \(\delta(\gamma)\).

For the torsion-trajectory program, the roles of the six operations are as
follows.

\begin{description}
\item[\(f^*\).]
Pullback or restriction.  This operation moves torsion packages to products,
test loci, and boundary or degeneration loci.  Examples include
\[
        \pi_S^*\alpha,\qquad
        \pi_{S_2}^*\beta,\qquad
        i^*\gamma=\gamma',
        \qquad
        i^*\Theta=\Theta'.
\]
In Diaz's argument, restriction from \(S\times S_2\) to \(E\times S_2\)
is exactly this type of operation \cite[Section 4]{Diaz2023IHCTrivialChow}.

\item[\(f_*\).]
Pushforward or globalization.  This operation moves a local or supported
package into the global space.  For a closed stratum
\(i:\Sigma\hookrightarrow X\), a supported motivic torsion package has the
form
\[
        i_*(\mathbf 1_\Sigma(p)/m)[s].
\]
This is the operation relevant to local discriminant packages and
closed-stratum torsion terms, such as the Coble package
\(i_*(\mathbf 1_\Sigma/4)[2]\).

\item[\(f_!\).]
Extension with proper support.  This operation is relevant for open
complements, extension by zero, localization triangles, and support
conditions.  For an open immersion \(j:U\hookrightarrow X\), it distinguishes
the extension-with-support object \(j_!M\) from the full direct image
\(j_*M\).  This distinction is important whenever the obstruction is
controlled by boundary behavior or by failure of a class to extend across a
closed stratum.

\item[\(f^!\).]
Exceptional pullback or local extraction.  This operation detects the local
contribution of a stratum.  For a closed immersion \(i:Z\hookrightarrow X\),
objects such as
\[
        i^!\mathbf 1_X
\]
are the natural place to record local support contributions and local
obstruction packages.  In the torsion-trajectory program, \(f^!\) is the
operation that extracts the local part of a global object at a singular,
boundary, or exceptional stratum.

\item[\(\otimes\).]
Tensor product, hence cup product after realization.  This is the operation responsible for the Coble--Diaz hierarchy.  It produces classes such as
\[
        \gamma=\pi_S^*\alpha\cup\pi_{S_2}^*\beta
\]
and
\[
        \Theta=\pi_1^*\alpha_1\cup\pi_2^*\beta_2\cup\pi_3^*\beta_3.
\]
The tensor operation is therefore the motivic source of the finite-coefficient
cup products used in the Diaz and level-three constructions.

\item[\(\mathcal Hom\).]
Internal Hom, duality, and pairings.  This operation is relevant when a
torsion package carries a discriminant form or duality pairing, for example
\[
        q:E\times E\to \mathbb Q/\mathbb Z.
\]
It is the natural categorical home for comparing discriminant forms, linking
pairings, Brauer dualities, and dual obstruction packages.  In the local
surface theory, the necessity of retaining the discriminant form, rather
than only the underlying finite group, is already visible in examples such
as \(D_4\) and in the Coble boundary package
\cite{RahmanIntegralPerverseObstructions,Nikulin80}.
\end{description}

Thus each station in a torsion trajectory can be read through one or more
of the six operations: restriction uses \(f^*\), globalization uses \(f_*\),
support and local extraction use \(f^!\), open extension uses \(f_!\), the
cup-product hierarchy uses \(\otimes\), and pairings use \(\mathcal Hom\).

\subsection{Realization compatibility}

The reason the six operations matter for the motivic program is that the
realization functors are compatible with them in the frameworks of Tubach
and Ruimy--Tubach.  Tubach's realization framework compares Voevodsky étale
motives with perverse Nori motives and mixed Hodge modules in a
six-functor-compatible way \cite{Tubach2025NoriHodgeRealizations}.
Ruimy--Tubach provide the integral-coefficient setting needed to discuss
integral Nori motivic sheaves, integral mixed Hodge modules, and torsion
phenomena with six operations \cite{RuimyTubach24}.  We use
these compatibility results as ambient infrastructure.

Let \(\Realstar\) denote a specified realization functor, such as
\(\RealB\), \(\Realet\), \(\RealNori\), or \(\RealMHM\).  The schematic
compatibilities used in this paper have the following form:
\[
        \Realstar(f^*M)\cong f^*\Realstar(M),
\]
\[
        \Realstar(f_*M)\cong f_*\Realstar(M),
\]
\[
        \Realstar(f_!M)\cong f_!\Realstar(M),
\]
\[
        \Realstar(f^!M)\cong f^!\Realstar(M),
\]
\[
        \Realstar(M\otimes N)\cong
        \Realstar(M)\otimes \Realstar(N),
\]
and
\[
        \Realstar(\mathcal Hom(M,N))
        \cong
        \mathcal Hom(\Realstar(M),\Realstar(N)).
\]

The point is not that these are new theorems here.  The point is that these
commutation rules explain why one may build a finite-coefficient torsion
package motivically and then compare its Betti, étale, Nori, mixed-Hodge,
Brauer, unramified, and stacky avatars.

For example, if
\[
        \gamma^{\mathrm{mot}}
        =
        \pi_S^*\alpha^{\mathrm{mot}}
        \cup
        \pi_{S_2}^*\beta^{\mathrm{mot}},
\]
then pullback and tensor compatibility give
\[
        \RealB(\gamma^{\mathrm{mot}})
        =
        \pi_S^*\RealB(\alpha^{\mathrm{mot}})
        \cup
        \pi_{S_2}^*\RealB(\beta^{\mathrm{mot}}),
\]
and similarly for \(\Realet\).  Thus the motivic class realizes to Diaz's
finite-coefficient class
\[
        \gamma=\pi_S^*\alpha\cup\pi_{S_2}^*\beta.
\]

This is why realization compatibility is central to the torsion-trajectory
program: it permits one to compare a package before and after passing to
Betti, étale, Hodge-module, Nori, Brauer, unramified, or stacky contexts.

\subsection{Coefficient triangles and Bockstein compatibility}

For \(m\ge2\), the motivic finite-coefficient object is
\[
        \mathbf 1_X(p)/m
        :=
        \operatorname{Cone}
        \left(
        \mathbf 1_X(p)
        \xrightarrow{m}
        \mathbf 1_X(p)
        \right).
\]
It fits into the coefficient triangle
\[
        \mathbf 1_X(p)
        \xrightarrow{m}
        \mathbf 1_X(p)
        \longrightarrow
        \mathbf 1_X(p)/m
        \overset{+1}{\longrightarrow}.
\]
The connecting morphism is the motivic Bockstein:
\[
        \delta_m^{\mathrm{mot}}:
        H^i_{\mathrm{mot}}(X,\mathbf 1_X(p)/m)
        \longrightarrow
        H^{i+1}_{\mathrm{mot}}(X,\mathbf 1_X(p)).
\]
This is the standard connecting morphism associated to a distinguished
triangle.  We use the term ``motivic Bockstein'' because its realizations are
the usual coefficient Bockstein maps used to pass from finite-coefficient
classes to integral torsion classes.

Under Betti realization, exactness and cone compatibility give the ordinary
coefficient sequence
\[
        0
        \longrightarrow
        \mathbb Z(p)
        \xrightarrow{m}
        \mathbb Z(p)
        \longrightarrow
        \mathbb Z/m(p)
        \longrightarrow
        0.
\]
Hence
\[
        \RealB\left(
        \delta_m^{\mathrm{mot}}(\Theta^{\mathrm{mot}})
        \right)
        =
        \delta_m\left(\RealB(\Theta^{\mathrm{mot}})\right).
\]
Similarly, under étale realization one obtains
\[
        \Realet\left(
        \delta_m^{\mathrm{mot}}(\Theta^{\mathrm{mot}})
        \right)
        =
        \delta_m\left(\Realet(\Theta^{\mathrm{mot}})\right),
\]
with the usual convention that finite étale coefficients may be written as
\(\mathbb Z/m(p)\) or as \(\mu_m^{\otimes p}\).

This is the formal reason motivic Bockstein trajectories realize to the
Bockstein trajectories used in integral Hodge obstruction theory.  In
Diaz's example, this is the passage from
\[
        \gamma\in H^3(V,\mathbb Z/2(2))
\]
to
\[
        \delta(\gamma)\in H^4(V,\mathbb Z(2)).
\]

\subsection{Localization and support}

Let
\[
        i:Z\hookrightarrow X
\]
be a closed immersion, and let
\[
        j:U=X\setminus Z\hookrightarrow X
\]
be the complementary open immersion.  In a six-functor formalism, closed-open
decompositions give localization triangles, schematically of the form
\[
        i_*i^!\mathbf 1_X
        \longrightarrow
        \mathbf 1_X
        \longrightarrow
        j_*j^*\mathbf 1_X
        \overset{+1}{\longrightarrow},
\]
as well as variants involving \(j_!\).  Applying the same construction to
finite-coefficient objects gives localization triangles for
\[
        \mathbf 1_X(p)/m.
\]
The localization formalism is part of the six-functor setting used in the
motivic and realization-compatible frameworks of Tubach and Ruimy--Tubach
\cite{Tubach2025NoriHodgeRealizations,RuimyTubach24}.

This is the motivic form of support exact sequences and residue phenomena.
It is the formal mechanism behind local-to-global transport.  A local
package on a stratum \(Z\) can be represented by an object such as
\[
        i_*(\mathbf 1_Z(p)/m)[s],
\]
and its global image is controlled by the forget-support maps coming from
localization.

For the Coble boundary package along a stratum \(\Sigma\), the full
order-four package is modeled by
\[
        i_*(\mathbf 1_\Sigma/4)[2],
\]
while the order-two shadow is modeled by
\[
        i_*(\mathbf 1_\Sigma/2)[2]
        \longrightarrow
        i_*(\mathbf 1_\Sigma/4)[2].
\]
Betti realization sends these to
\[
        i_*(\mathbb Z/4)_\Sigma[2]
\]
and
\[
        i_*(\mathbb Z/2)_\Sigma[2]
        \longrightarrow
        i_*(\mathbb Z/4)_\Sigma[2].
\]
Thus the Coble/Benoist--Ottem comparison is an example of a supported motivic torsion package and a selected finite-coefficient shadow.

\subsection{Products and the cup-product hierarchy}

The tensor product is the operation responsible for the Coble--Diaz
hierarchy.  Because realization is compatible with tensor products in the
six-functor realization frameworks cited above
\cite{Tubach2025NoriHodgeRealizations,RuimyTubach24}, motivic products
realize to the cup products used in finite-coefficient cohomology.

\begin{remark}[Toward multiple cup products]
\label{rem:toward-multiple-cup-products}
This observation suggests a possible higher cup-product extension of the
Coble--Diaz hierarchy.  One may ask whether several finite-coefficient torsion
inputs
\[
   \theta_1,\ldots,\theta_n
\]
can be combined into a product
\[
   \Theta=\theta_1\cup\cdots\cup\theta_n
\]
whose Bockstein carries a higher integral obstruction.  The present paper does
not develop this theory.  Such an extension would require a separate analysis
of product survivability, residue conditions, and algebraic-control
constraints.  We leave this higher cup-product Bockstein hierarchy to future
work.
\end{remark}

For Diaz, the relevant finite-coefficient classes are
\[
        \alpha\in H^1(S,\mathbb Z/2(1)),
        \qquad
        \beta\in H^2(S_2,\mathbb Z/2(1)).
\]
Motivically, one records lifts
\[
        \alpha^{\mathrm{mot}}
        \in
        H^1_{\mathrm{mot}}(S,\mathbf 1_S(1)/2),
\]
and
\[
        \beta^{\mathrm{mot}}
        \in
        H^2_{\mathrm{mot}}(S_2,\mathbf 1_{S_2}(1)/2).
\]
The product on \(V=S\times S_2\) is
\[
        \gamma^{\mathrm{mot}}
        =
        \pi_S^*\alpha^{\mathrm{mot}}
        \cup
        \pi_{S_2}^*\beta^{\mathrm{mot}}.
\]
Its Betti realization is
\[
        \RealB(\gamma^{\mathrm{mot}})
        =
        \pi_S^*\alpha
        \cup
        \pi_{S_2}^*\beta
        =
        \gamma.
\]
The Bockstein then gives
\[
        \delta(\gamma)\in H^4(V,\mathbb Z(2)).
\]

For the level-three candidate, one uses
\[
        X=S_1\times S_2\times S_3
\]
and
\[
        \Theta
        =
        \pi_1^*\alpha_1
        \cup
        \pi_2^*\beta_2
        \cup
        \pi_3^*\beta_3
        \in
        H^5(X,\mathbb Z/2(3)).
\]
The same motivic pattern gives
\[
        \Theta^{\mathrm{mot}}
        =
        \pi_1^*\alpha_1^{\mathrm{mot}}
        \cup
        \pi_2^*\beta_2^{\mathrm{mot}}
        \cup
        \pi_3^*\beta_3^{\mathrm{mot}},
\]
followed by
\[
        \delta^{\mathrm{mot}}(\Theta^{\mathrm{mot}}).
\]
Thus \(\otimes\) is the operation that turns isolated torsion packages into
higher-level product mechanisms.

\subsection{Duality, pairings, and internal Hom}

The operation \(\mathcal Hom\), together with duality, becomes important
when torsion packages carry pairings.  Singular packages often come with
discriminant forms
\[
        q:E\times E\to \mathbb Q/\mathbb Z.
\]
For normal surface singularities, this pairing is induced by the exceptional
lattice and is compared with the torsion linking pairing on the link
\cite{GoreskySiegel83,Nikulin80,RahmanIntegralPerverseObstructions}.

Categorically, such pairings should be understood as maps from a package to
its dual or as pairings mediated by internal Hom.  Although the Diaz
construction does not require a detailed duality formalization, later
comparisons of discriminant forms, Brauer pairings, and stacky finite-group
pairings will require this part of the six-operation formalism.

The distinction between the finite group \(E\) and the full discriminant
package \((E,q)\) is important.  Two examples may have the same underlying
finite abelian group but different pairings.  Internal Hom and duality give
the natural categorical language for recording this extra structure.

\subsection{Example dictionary}

We summarize how the principal examples fit the six-operation dictionary.

\subsubsection*{Coble/Benoist--Ottem}

The source package is the Coble local discriminant package
\[
        E_{\frac14(1,1)}\cong\mathbb Z/4.
\]
The visible Enriques direction is
\[
        2E\cong\mathbb Z/2.
\]
Motivically, the full package is represented by
\[
        i_*(\mathbf 1_\Sigma/4)[2],
\]
and the visible order-two shadow by
\[
        i_*(\mathbf 1_\Sigma/2)[2].
\]
The active operations are \(i_*\), \(i^!\), coefficient Bockstein, and
realization.  This is a local boundary/support trajectory.  The local
Coble package and its order-two shadow are computed in the local
discriminant framework of \cite{RahmanIntegralPerverseObstructions}.

\subsubsection*{Diaz}

The source packages are smooth finite-coefficient classes
\[
        \alpha\in H^1(S,\mathbb Z/2(1)),
        \qquad
        \beta\in H^2(S_2,\mathbb Z/2(1)).
\]
The active operations are
\[
        \pi_S^*,\qquad
        \pi_{S_2}^*,\qquad
        \otimes,\qquad
        i^*,\qquad
        \delta_2.
\]
The trajectory is
\[
        \alpha,\beta
        \longrightarrow
        \gamma=\pi_S^*\alpha\cup\pi_{S_2}^*\beta
        \longrightarrow
        \delta(\gamma),
\]
with non-algebraicity controlled by unramified survival.  This is a smooth finite-coefficient cup-product trajectory, and it is precisely the construction used by Diaz \cite[Section 4]{Diaz2023IHCTrivialChow}.

\subsubsection*{Level three}

The candidate source packages are
\[
        \alpha_1\in H^1(S_1,\mathbb Z/2(1)),
        \qquad
        \beta_2\in H^2(S_2,\mathbb Z/2(1)),
        \qquad
        \beta_3\in H^2(S_3,\mathbb Z/2(1)).
\]
The product is
\[
        \Theta
        =
        \pi_1^*\alpha_1
        \cup
        \pi_2^*\beta_2
        \cup
        \pi_3^*\beta_3.
\]
The active operations are the same as in Diaz, but with one additional
tensor factor.  The unresolved station is the restricted nonvanishing
problem on
\[
        E\times S_2\times S_3.
\]
This is the first higher-level test of the Coble--Diaz cup-product
hierarchy.

\subsubsection*{Kummer}

The local singular package is
\[
        E_{A_1}^{\mathrm{sing}}\cong\mathbb Z/2.
\]
The stacky stabilizer package is
\[
        H^2(B\mu_2,\mathbb Z)\cong\mathbb Z/2.
\]
The compatibility
\[
        E_{A_1}^{\mathrm{sing}}
        \cong
        H^2(B\mu_2,\mathbb Z)
\]
is the local singular-stacky calibration.  The active operations are
stack/coarse comparison, realization, and supported transport of fixed-point
packages.

\subsubsection*{Atiyah--Hirzebruch and Godeaux--Serre}

In smooth finite-group quotient examples, ordinary singular \(E\)-packages may vanish because the final variety is smooth.  This means that torsion computations using $E$ constructed in \cite{RahmanIntegralPerverseObstructions} is trivial. Hence, due to the variety's smoothness, the expected torsion source is instead equivariant or stacky:
\[
        E_G,\qquad
        H^*(BG,\mathbb Z),\qquad
        H^*_G(V,\mathbb Z).
\]
Thus, the active operations are quotient-stack realization, finite-group cohomology, pullback to covers, and descent to the coarse space.  These integral hodge counterexamples are expected to test whether stacky \(E_G\)-packages are the correct torsion detection analogs for singular \(E\)-packages in smooth quotient constructions.

\subsection{The pre-rationalization obstruction calculus}

The long-term purpose of this appendix is to explain the operating system of
the torsion-trajectory program.  The integral Hodge obstruction group
\[
        \mathfrak O^p_{\mathbb Z}(X)
        :=
        \operatorname{coker}
        \left(
        CH^p(X)
        \longrightarrow
        \operatorname{Hdg}^p_{\mathbb Z}(X)
        \right)
\]
contains torsion phenomena that disappear after tensoring with
\(\mathbb Q\).  The rational obstruction is
\[
        \mathfrak O^p_{\mathbb Z}(X)\otimes\mathbb Q.
\]
Thus torsion classes are killed by rationalization, but before they die they record how integral algebraicity fails. The motivic finite-coefficient setting is useful because it records the pre-rationalization data: \emph{source package}, \emph{transport operations}, \emph{cup product}, \emph{support or residue survival}, \emph{Bockstein}, and \emph{realization}. This is the obstruction calculus one needs in order to compare integral
Hodge counterexamples.
\\
\\
\noindent
\emph{Guiding principle.}
Torsion mechanisms should be compared before tensoring with \(\mathbb Q\); the six operations provide the transport calculus for that comparison. Hence, each known counterexample should be lifted, or at least represented, by finite-coefficient motivic packages and their realization data.  Only then can one ask whether apparently different counterexamples are genuinely different, or whether they are different realizations of a smaller number of underlying motivic torsion mechanisms.
%
%
\printbibliography

\end{document}